\documentclass[11pt]{amsart}
\usepackage{amssymb}
\usepackage{mathrsfs}
\usepackage{stmaryrd}
\usepackage{bbm}
\usepackage{oldgerm}
\usepackage[frenchb,english]{babel}
\usepackage[utf8]{inputenc}
\usepackage[T1]{fontenc}
\usepackage{lmodern}
\usepackage[all]{xy}
\usepackage{hyperref}
\usepackage{upref}

\def\Zar{\mathbf{Zar}}
\def\aa{\mathfrak{a}}

\def\PP{P}
\def\QQ{Q}
\def\TT{T}
\def\RR{R}
\def\wT{E}
\def\XX{\mathfrak{X}}
\def\DD{\mathfrak{D}}
\def\UU{\mathfrak{U}}
\def\YY{\mathfrak{Y}}
\def\ZZ{\mathfrak{Z}}
\def\SS{\mathscr{S}}

\def\Mod{\mathbf{Mod}}

\def\Ob{\mathbf{Ob}}
\def\FHom{\mathscr{H}om}

\def\MF{\mathbf{MF}}
\def\MFb{\mathbf{MF}_{\textnormal{big}}}
\def\pCRIS{\mathscr{E}}
\def\CRIS{\underline{\mathscr{E}}}

\def\pMIC{p\textnormal{-MIC}}
\def\lMIC{\lambda\textnormal{-MIC}}
\def\FEnd{\mathscr{E}nd}
\def\CC{\mathscr{C}}

\DeclareMathOperator{\length}{lg}
\DeclareMathOperator{\IM}{Im}
\DeclareMathOperator{\cris}{crys}
\DeclareMathOperator{\Cris}{Crys}
\DeclareMathOperator{\tot}{tot}
\DeclareMathOperator{\Tot}{Tot}
\DeclareMathOperator{\gr}{gr}
\DeclareMathOperator{\qn}{qn}

\DeclareMathOperator{\Cov}{Cov}
\DeclareMathOperator{\fppf}{fppf}
\DeclareMathOperator{\Supp}{Supp}
\DeclareMathOperator{\MIC}{MIC}
\DeclareMathOperator{\HIGG}{HIG}
\DeclareMathOperator{\qcoh}{qcoh}
\DeclareMathOperator{\can}{can}

\DeclareMathOperator{\zar}{zar}

\DeclareMathOperator{\Ker}{Ker}
\DeclareMathOperator{\Coker}{Coker}

\DeclareMathOperator{\Spec}{Spec}
\DeclareMathOperator{\Spf}{Spf}
\DeclareMathOperator{\Proj}{Proj}

\DeclareMathOperator{\rH}{H}
\DeclareMathOperator{\rI}{I}

\DeclareMathOperator{\rmC}{C}
\DeclareMathOperator{\rS}{S}
\DeclareMathOperator{\rT}{T}

\DeclareMathOperator{\rD}{D}
\DeclareMathOperator{\rW}{W}
\DeclareMathOperator{\rF}{F}
\DeclareMathOperator{\rE}{E}
\DeclareMathOperator{\rR}{R}
\DeclareMathOperator{\Gr}{Gr}

\DeclareMathOperator{\btor}{-tor}
\DeclareMathOperator{\id}{id}

\newtheorem{theorem}{Theorem}[section]
\newtheorem{prop}[theorem]{Proposition}
\newtheorem{lemma}[theorem]{Lemma}

\newtheorem{coro}[theorem]{Corollary}

\theoremstyle{definition}
\newtheorem{rem}[theorem]{Remark}
\newtheorem{definition}[theorem]{Definition}
\newtheorem{nothing}[theorem]{}

\numberwithin{equation}{section}
\numberwithin{equation}{theorem}

\makeatletter
\newsavebox{\@brx}
\newcommand{\llangle}[1][]{\savebox{\@brx}{\(\m@th{#1\langle}\)}%
  \mathopen{\copy\@brx\kern-0.5\wd\@brx\usebox{\@brx}}}
\newcommand{\rrangle}[1][]{\savebox{\@brx}{\(\m@th{#1\rangle}\)}%
  \mathclose{\copy\@brx\kern-0.5\wd\@brx\usebox{\@brx}}}
\makeatother

\makeatletter
\newcommand*{\relrelbarsep}{.386ex}
\newcommand*{\relrelbar}{%
  \mathrel{%
    \mathpalette\@relrelbar\relrelbarsep
  }%
}
\newcommand*{\@relrelbar}[2]{%
  \raise#2\hbox to 0pt{$\m@th#1\relbar$\hss}%
  \lower#2\hbox{$\m@th#1\relbar$}%
}
\providecommand*{\rightrightarrowsfill@}{%
  \arrowfill@\relrelbar\relrelbar\rightrightarrows
}
\providecommand*{\leftleftarrowsfill@}{%
  \arrowfill@\leftleftarrows\relrelbar\relrelbar
}
\providecommand*{\xrightrightarrows}[2][]{%
  \ext@arrow 0359\rightrightarrowsfill@{#1}{#2}%
}
\providecommand*{\xleftleftarrows}[2][]{%
  \ext@arrow 3095\leftleftarrowsfill@{#1}{#2}%
}
\makeatother

\title{Lifting the Cartier transform of Ogus-Vologodsky modulo $p^{n}$}
\author{Daxin Xu}
\date{\today}

\AtEndDocument{\bigskip{\footnotesize%
  \textsc{Daxin Xu, Department of Mathematics, California Institute of Technology, Pasadena, CA 91125, USA} \par
  \textit{E-mail address}: \texttt{daxinxu@caltech.edu} \par
}}

\begin{document}
\selectlanguage{english}
\maketitle
\begin{abstract}
	Let $\rW$ be the ring of the Witt vectors of a perfect field of characteristic $p$, $\XX$ a smooth formal scheme over $\rW$, $\XX'$ the base change of $\XX$ by the Frobenius morphism of $\rW$, $\XX_{2}'$ the reduction modulo $p^{2}$ of $\XX'$ and $X$ the special fiber of $\XX$. 
	We lift the Cartier transform of Ogus--Vologodsky defined by $\XX_{2}'$ modulo $p^{n}$. More precisely, we construct a functor from the category of $p^{n}$-torsion $\mathscr{O}_{\XX'}$-modules with integrable $p$-connection to the category of $p^{n}$-torsion $\mathscr{O}_{\XX}$-modules with integrable connection, each subject to suitable nilpotence conditions. Our construction is based on Oyama's reformulation of the Cartier transform of Ogus--Vologodsky in characteristic $p$. 
	If there exists a lifting $F:\XX\to \XX'$ of the relative Frobenius morphism of $X$, our functor is compatible with a functor constructed by Shiho from $F$. As an application, we give a new interpretation of Faltings' relative Fontaine modules and of the computation of their cohomology. 
\end{abstract}
\tableofcontents
\section{Introduction}
\begin{nothing} \label{Intro general}
	In his seminal work \cite{Sim92}, Simpson established a deep relation between complex representations of the fundamental group of a projective complex manifold $X$ and Higgs modules on $X$, leading to a theory called \textit{nonabelian Hodge theory}. Recall that a Higgs module on $X$ is a coherent sheaf $M$ together with an $\mathscr{O}_{X}$-linear morphism $\theta:M\to M\otimes_{\mathscr{O}_X}\Omega_{X/\mathbb{C}}^1$ such that $\theta\wedge \theta=0$. (Simpson's result uses, but is much deeper than, the Riemann-Hilbert correspondence relating representations of the fundamental group and modules with integrable connection.)
	In \cite{Fal05}, Faltings developed a partial $p$-adic analogue of Simpson correspondence for $p$-adic local systems on varieties over $p$-adic fields. 

	On the other hand, in \cite{OV07}, Ogus and Vologodsky constructed a version of nonabelian Hodge theory in characteristic $p$. 
	If $X$ is a smooth scheme over a perfect field $k$ of characteristic $p>0$, they established an equivalence, called \textit{Cartier transform}, between certain modules with integrable connection on $X/k$ and certain Higgs modules on $X/k$, depending on a lifting of $X'$ (the base change of $X$ by the Frobenius morphism of $k$) to $\rW_2(k)$. 
	They also constructed a canonical quasi-isomorphism betweeen certain truncations of the de Rham complex of a module with integrable connection and of the Higgs complex of its Cartier transform. 
	This result generalizes the Cartier isomorphism and the decomposition of the de Rham complex given by Deligne--Illusie \cite{DL}; it is also an analogue of a corresponding result in Simpson's theory.

	The relation between Faltings' $p$-adic Simpson correspondence and the Cartier transform is not yet understood. The first difficulty is to lift the Cartier transform modulo $p^{n}$. 
	This is our main goal in the present article. 
	Shiho \cite{Shiho} constructed a ``local'' lifting of the Cartier transform modulo $p^{n}$ under the assumption of a lifting of the relative Frobenius morphism modulo $p^{n+1}$.
	In \cite{Oy}, Oyama gave a new construction of the Cartier transform of Ogus--Vologodsky as the inverse image by a morphism of topoi. His work is inspired by Tsuji's approach to the $p$-adic Simpson correspondence (\cite{AGT} IV). 
	In this article, we use Oyama topoi to ``glue'' Shiho's functor and obtain a lifting of the Cartier transform modulo $p^{n}$ under the (only) assumption that $X$ lifts to a smooth formal scheme over $\rW$.
\end{nothing}


\begin{nothing} \label{lamda connection}
	Shiho's construction applies to modules with \textit{$\lambda$-connection}, a notion of introduced by Deligne. 
	Let $f:X\to S$ be a smooth morphism of schemes, $M$ an $\mathscr{O}_{X}$-module and $\lambda\in \Gamma(S,\mathscr{O}_{S})$. A \textit{$\lambda$-connection on $M$ relative to $S$} is an $f^{-1}(\mathscr{O}_{S})$-linear morphism $\nabla:M\to M\otimes_{\mathscr{O}_{X}}\Omega_{X/S}^{1}$ such that $\nabla(xm)=x\nabla(m)+\lambda m\otimes d(x)$ for every local sections $x$ of $\mathscr{O}_{X}$ and $m$ of $M$.
	$1$-connections correspond to the classical notion of connections, and $0$-connections to Higgs fields.
	The integrability of $\lambda$-connections is defined in the same way as for connections. 
	We denote by $\MIC(X/S)$ (resp. $\lMIC(X/S)$) the category of $\mathscr{O}_{X}$-modules with integrable connection (resp. $\lambda$-connection) relative to $S$.
\end{nothing}

\begin{nothing} \label{Shiho functor}
	In the following, if we use a gothic letter $\mathfrak{T}$ to denote an adic formal $\rW$-scheme, the corresponding roman letter $T$ will denote its special fiber. 
	Let $\XX$ be a smooth formal scheme over $\rW$ and $n$ an integer $\ge 1$. We denote by $\sigma:\rW\to \rW$ the Frobenius automorphism of $\rW$, by $\XX'$ the base change of $\XX$ by $\sigma$ and by $\XX_{n}$ the reduction of $\XX$ modulo $p^{n}$. 
	In \cite{Shiho}, Shiho constructed a ``local'' lifting modulo $p^{n}$ of the Cartier transform of Ogus--Vologodsky defined by $\XX_{2}'$, using a lifting $F_{n+1}:\XX_{n+1}\to \XX'_{n+1}$ of the relative Frobenius morphism $F_{X/k}:X\to X'$ of $X$.  
	
	The image of the differential morphism $dF_{n+1}:F_{n+1}^{*}(\Omega^{1}_{\XX_{n+1}'/\rW_{n+1}})\to \Omega^{1}_{\XX_{n+1}/\rW_{n+1}}$ of $F_{n+1}$ is contained in $p\Omega^{1}_{\XX_{n+1}/\rW_{n+1}}$. Dividing by $p$, it induces an $\mathscr{O}_{\XX_{n}}$-linear morphism
	\begin{displaymath}
		dF_{n+1}/p:F_{n}^{*}(\Omega^{1}_{\XX_{n}'/\rW_{n}})\to \Omega^{1}_{\XX_{n}/\rW_{n}}.
	\end{displaymath}
	Shiho defined a functor (depending on $F_{n+1}$) (\cite{Shiho} 2.5)	
	\begin{eqnarray}\label{Shiho functor eq}
	\Phi_{n}: \pMIC(\XX_{n}'/\rW_{n})&\to& \MIC(\XX_{n}/\rW_{n}) \\
	(M',\nabla') &\mapsto & (F_{n}^{*}(M'),\nabla) \nonumber
\end{eqnarray}
where $\nabla:F_{n}^{*}(M')\to \Omega_{\XX_{n}/\rW_{n}}^{1}\otimes_{\mathscr{O}_{\XX_{n}}}F_{n}^{*}(M')$ is the integrable connection defined for every local section $e$ of $M'$ by
	\begin{equation} \label{shiho formula intro}
		\nabla(F_{n}^{*}(e))= (\id\otimes \frac{dF_{n+1}}{p})(F_{n}^{*}(\nabla'(e))).
	\end{equation}

	Shiho showed that the functor $\Phi_{n}$ induces an equivalence of categories between the full subcategories of $\pMIC(\XX_{n}'/\rW_{n})$ and of $\MIC(\XX_{n}/\rW_{n})$ consisting of quasi-nilpotent objects (\cite{Shiho} Thm. 3.1).
	When $n=1$, Ogus and Vologodsky proved that the functor $\Phi_{1}$ is compatible with the Cartier transform defined by $\XX_{2}'$ (\cite{OV07} Thm. 2.11; \cite{Shiho} 1.12). 
\end{nothing}

\begin{nothing} \label{intro Hopf algebra}
	The categories of connections and their analogues we will be studying can be understood geometrically using the language of groupoids. Our groupoids
will be relatively affine and hence correspond to Hopf algebras.
	If $(\mathscr{T},A)$ is a ringed topos, \textit{a Hopf $A$-algebra} is the data of a ring $B$ of $\mathscr{T}$ together with five homomorphisms 
	\begin{displaymath}
		A\xrightrightarrows[d_{1}]{d_{2}}B, \quad \delta:B\to B\otimes_{A}B~ \textnormal{(comultiplication)},\quad \pi:B\to A~\textnormal{(counit)},\quad \sigma:B\to B~\textnormal{(antipode)},
	\end{displaymath}
	where the tensor product $B\otimes_{A}B$ is taken on the left (resp. right) for the $A$-algebra structure of $B$ defined by $d_{2}$ (resp. $d_{1}$), satisfying the compatibility conditions for coalgebras (cf. \ref{def of Hopf alg}, \cite{Ber} II 1.1.2).
%
%
%
%

	\textit{A $B$-stratification on an $A$-module $M$} is a $B$-linear isomorphism 
	\begin{equation}
		\varepsilon:B\otimes_{A}M\xrightarrow{\sim}M\otimes_{A}B
	\end{equation}
	where the tensor product is taken on the left (resp. right) for the $A$-algebra structure defined by $d_{2}$ (resp. $d_{1}$), satisfying $\pi^{*}(\varepsilon)=\id_{M}$ and a cocycle condition (cf. \ref{def stratification}).
\end{nothing}

\begin{nothing} \label{intro PX}
	A classical example of a Hopf algebra is given by the PD-envelope of the diagonal immersion.
	Let $\XX$ be a smooth formal $\rW$-scheme, $\XX^{2}$ the product of two copies of $\XX$ over $\rW$. For any $n\ge 1$, we denote by $\PP_{\XX_n}$ the PD-envelope of the diagonal immersion $\XX_{n}\to \XX_{n}^{2}$ compatible with the canonical PD-structure on $(\rW_{n},p\rW_{n})$ and by $\PP_{\XX}$ the associated adic formal $\rW$-scheme. 
	The $\mathscr{O}_{\XX}$-bialgebra $\mathscr{O}_{\PP_{\XX}}$ of $\XX_{\zar}$ is naturally equipped with a formal Hopf $\mathscr{O}_{\XX}$-algebra structure (i.e. for every $n\ge 1$, a Hopf $\mathscr{O}_{\XX_{n}}$-algebra structure on $\mathscr{O}_{\PP_{\XX_{n}}}$, which is compatible) (cf. \ref{def of groupoid str}, \ref{groupoid PX}).

	A quasi-nilpotent integrable connection relative to $\rW_{n}$ on an $\mathscr{O}_{\XX_{n}}$-module $M$ (cf. \ref{def quasinilpotent}) is equivalent to an $\mathscr{O}_{\PP_{\XX}}$-stratification on $M$ (\cite{BO} 4.12). 
	Following Shiho \cite{Shiho}, we give below an analogous description of $p$-connections; the relevant Hopf algebra is constructed by dilatation (certain distinguished open subset of admissible blow-up) in formal geometry. 
\end{nothing}

\begin{nothing}\label{intro RX}
	We define by dilatation an adic formal $\XX^{2}$-scheme $\RR_{\XX}$ satisfying the following conditions (\ref{prop univ RQ}). 
	
	(i) The canonical morphism $\RR_{\XX,1}\to X^{2}$ factors through the diagonal immersion $X\to X^{2}$.
	
	(ii) Let $X\to \XX^{2}$ be the morphism induced by the diagonal immersion. For any flat formal $\rW$-scheme $\YY$ and any $\rW$-morphisms $f:\YY\to \XX^{2}$ and $g:Y\to X$ which fit into the following commutative diagram
\begin{displaymath}
	\xymatrix{
		Y \ar[r] \ar[d]_{g} & \YY \ar[d]^{f}\\
		X \ar[r] & \XX^{2},
	}
\end{displaymath}
	there exists a unique $\rW$-morphism $f':\YY\to \RR_{\XX}$ lifting $f$. 

	We denote abusively by $\mathscr{O}_{\RR_{\XX}}$ the direct image of $\mathscr{O}_{\RR_{\XX}}$ via the morphism $\RR_{\XX,\zar}\to \XX_{\zar}$ (i). Using the universal property of $\RR_{\XX}$, we show that $\mathscr{O}_{\RR_{\XX}}$ is equipped with a formal Hopf $\mathscr{O}_{\XX}$-algebra structure (\ref{prop R Q Hopf alg}).
	
	The diagonal immersion $\XX\to \XX^{2}$ induces a closed immersion $\iota: \XX\to \RR_{\XX}$ \eqref{prop univ RQ}. 
	For any $n\ge 1$, we denote by $\TT_{\XX,n}$ the PD-envelope of $\iota_{n}:\XX_{n}\to \RR_{\XX,n}$ compatible with the canonical PD-structure on $(\rW_{n},p\rW_{n})$. 
	The schemes $\{\TT_{\XX,n}\}_{n\ge 1}$ form an adic inductive system and we denote by $\TT_{\XX}$ the associated adic formal $\rW$-scheme. 
	By the universal property of PD-envelope, the formal Hopf algebra structure on $\mathscr{O}_{\RR_{\XX}}$ extends to a formal Hopf $\mathscr{O}_{\XX}$-algebra structure on the $\mathscr{O}_{\XX}$-bialgebra $\mathscr{O}_{\TT_{\XX}}$ of $\XX_{\zar}$ \eqref{groupoid TX}.
	
	In (\cite{Shiho} Prop. 2.9), Shiho showed that for any $n\ge 1$ and any $\mathscr{O}_{\XX_{n}}$-module $M$, an $\mathscr{O}_{\TT_{\XX}}$-stratification on $M$ is equivalent to a quasi-nilpotent integrable $p$-connection on $M$ (cf. \ref{qn pconnection T}).
\end{nothing}

\begin{nothing}\label{intro QX}
	Shiho's local construction deals with modules with $p$-connection and connection, which is different to the (global) Cartier transform of Ogus--Vologodsky. 
	We need a fourth Hopf algebra, introduced by Oyama \cite{Oy}, and we will use it to define a notion of stratification that will enable us to globalise Shiho's construction. 
	
	For any $k$-scheme $Y$, we denote by $\underline{Y}$ the closed subscheme of $Y$ defined by the ideal sheaf of $\mathscr{O}_{Y}$ consisting of the sections of $\mathscr{O}_{Y}$ whose $p$th power is zero.
	In (\ref{prop univ RQ}, \ref{notations R Q}), we construct an adic formal $\XX^{2}$-scheme $\QQ_{\XX}$ satisfying the following conditions.
	
	(i) The canonical morphism $\underline{\QQ_{\XX,1}}\to X^{2}$ factors through the diagonal immersion $X\to X^{2}$.
	
	(ii) For any flat formal $\rW$-scheme $\YY$ and any $\rW$-morphisms $f:\YY\to \XX^{2}$ and $g:\underline{Y}\to X$ which fit into the following commutative diagram
\begin{displaymath}
	\xymatrix{
		\underline{Y} \ar[r] \ar[d]_{g} & \YY \ar[d]^{f}\\
		X \ar[r] & \XX^{2},
	}
\end{displaymath}
there exists a unique $\rW$-morphism $f':\YY\to \QQ_{\XX}$ lifting $f$. 

	We denote abusively by $\mathscr{O}_{\QQ_{\XX}}$ the direct image of $\mathscr{O}_{\QQ_{\XX}}$ via the morphism $\QQ_{\XX,\zar}\to \XX_{\zar}$ (i). It is also equipped with a formal Hopf $\mathscr{O}_{\XX}$-algebra structure (\ref{prop R Q Hopf alg}).

Let $\PP_{\XX}$ be the formal $\XX^{2}$-scheme defined in \ref{intro PX}, $\iota:\XX\to \PP_{\XX}$ the canonical morphism lifting the diagonal immersion $\XX\to \XX^{2}$ and $\mathscr{J}$ the PD-ideal of $\mathscr{O}_{\PP_{X}}$ associated to $\iota_{1}$. 
For any local section of $\mathscr{J}$, we have $x^{p}=p!x^{[p]}=0$. Then we deduce a closed immersion $\underline{\PP_{X}}\to X$ over $\XX^{2}$. By the universal property of $\QQ_{\XX}$, we obtain an $\XX^{2}$-morphism $\lambda:\PP_{\XX}\to \QQ_{\XX}$. 
\end{nothing}

\begin{nothing} \label{intro F induces phi psi}
	The groupoids and Hopf algebras constructed above give a geometric interpretation of Shiho's functor $\Phi$ and of a variation of $\Phi$ which can be globalized. 
	Let $F:\XX\to \XX'$ be a lifting of the relative Frobenius morphism $F_{X/k}$ of $X$. By the universal properties of $\RR_{\XX'}$ and of PD-envelopes, the morphism $F^{2}:\XX^{2}\to \XX'^{2}$ induce morphisms $\psi:\QQ_{\XX}\to \RR_{\XX'}$ \eqref{prop Q to R'} and $\varphi:\PP_{\XX}\to \TT_{\XX'}$ \eqref{prop P to T'} above $F^{2}$ which fit into a commutative diagram \eqref{square Frob}
	\begin{equation}
		\xymatrixcolsep{4pc}\xymatrix{
			\PP_{\XX} \ar[r]^{\varphi} \ar[d]_{\lambda}& \TT_{\XX'} \ar[d]^{\varpi}\\
			\QQ_{\XX} \ar[r]^{\psi} & \RR_{\XX'}},
	\end{equation}
	where $\varpi:\TT_{\XX'}\to \RR_{\XX'}$ \eqref{intro RX} and $\lambda:\PP_{\XX}\to \QQ_{\XX}$ \eqref{intro QX} are independent of $F$. Moreover, $\psi$ and $\varphi$ induce homomorphisms of formal Hopf algebras $\mathscr{O}_{R_{\XX'}}\to F_{*}(\mathscr{O}_{\QQ_{\XX}})$ and $\mathscr{O}_{\TT_{\XX'}}\to F_{*}(\mathscr{O}_{\PP_{\XX}})$. The above diagram induces a commutative diagram \eqref{square stra}
	\begin{equation} \label{intro carre stratification}
		\xymatrixcolsep{4pc}\xymatrix{
		\Big\{ \txt{	\textnormal{category of $\mathscr{O}_{\XX'_{n}}$-modules}\\
		\textnormal{with $\mathscr{O}_{\RR_{\XX'}}$-stratification}} \Big\}
		\ar[r]^{\psi_{n}^{*}} \ar[d]_{\varpi_{n}^{*}}&
		\Big\{	\txt{\textnormal{category of $\mathscr{O}_{\XX_{n}}$-modules}\\
		\textnormal{with $\mathscr{O}_{\QQ_{\XX}}$-stratification}} \Big\}
		\ar[d]^{\lambda_{n}^{*}}\\
		\Big\{ \txt{	\textnormal{category of $\mathscr{O}_{\XX'_{n}}$-modules} \\
		\textnormal{with $\mathscr{O}_{\TT_{\XX'}}$-stratification}} \Big\}
		\ar[r]^{\varphi_{n}^{*}} &
		\Big\{	\txt{\textnormal{category of $\mathscr{O}_{\XX_{n}}$-modules}\\
		\textnormal{with $\mathscr{O}_{\PP_{\XX}}$-stratification}} \Big\}}
\end{equation}
In (\cite{Shiho} 2.17), Shiho showed that the functor $\varphi_{n}^{*}$ is compatible with the functor $\Phi_{n}$ defined by $F$ \eqref{Shiho functor eq}, via the equivalence between the category of modules with quasi-nilpotent integrable connection (resp. $p$-connection) and the category of modules with $\mathscr{O}_{\PP_{\XX}}$-stratification (resp. $\mathscr{O}_{\TT_{\XX}}$-stratification). 
\end{nothing} 

\begin{nothing} \label{Oyama topos intro}
	Let us explain the Oyama sites $\mathscr{E}$ and $\underline{\mathscr{E}}$ whose crystals corresponding to $\mathscr{O}_{\QQ_{\XX}}$ and $\mathscr{O}_{\RR_{\XX}}$ stratification, and a morphism of topoi which will be used to lift the Cartier transform and to globalize the funtor $\psi_n^{*}$. 
	
	Let $X$ be a scheme over $k$. 	
	An object of $\mathscr{E}$ (resp. $\underline{\mathscr{E}}$) is a triple $(U,\mathfrak{T},u)$ consisting of an open subscheme $U$ of $X$, a flat formal $\rW$-scheme $\mathfrak{T}$ and an affine $k$-morphism $u:T\to U$ (resp. $u:\underline{T}\to U$ \eqref{intro QX}). Morphisms are defined in a natural way (cf. \ref{definition Oyama cat}). We denote by $\pCRIS'$ Oyama's category associated to the $k$-scheme $X'$. 
	We denote by $\widetilde{\mathscr{E}}$ (resp. $\widetilde{\underline{\mathscr{E}}}$) the topos of sheaves of sets on $\mathscr{E}$ (resp. $\underline{\mathscr{E}}$) with respect to the Zarisiki topology \eqref{usual covering}. 

	Let $(U,\mathfrak{T},u)$ be an object of $\CRIS$.
	The relative Frobenius morphism $F_{T/k}:T\to T'$ factors through a $k$-morphism $f_{T/k}:T\to \underline{T}'$. 
	We have a commutative diagram 
	\begin{equation}
		\xymatrixcolsep{4pc}\xymatrix{
			U \ar[d]_{F_{U/k}} & \underline{T}~\ar[l]_{u} \ar@{^{(}->}[r] \ar[d]_{F_{\underline{T}/k}} & T \ar[d]^{F_{T/k}} \ar[ld]_{f_{T/k}}\\
			U' & \underline{T}'~\ar[l]_{u'} \ar@{^{(}->}[r] & T' } \label{diagram rho intro}
\end{equation}
where the vertical arrows denote the relative Frobenius morphisms. 
Then $(U',\mathfrak{T},u'\circ f_{T/k})$ is an object of $\pCRIS'$. We obtain a functor \eqref{functor rho}
	\begin{equation} \label{functor rho intro}
		\rho:\underline{\mathscr{E}}\to \mathscr{E}',\qquad (U,\mathfrak{T},u)\mapsto (U',\mathfrak{T},u'\circ f_{T/k}).
	\end{equation}
	The functor $\rho$ is continuous and cocontinuous \eqref{lemma con con} and induces a morphism of topoi \eqref{morphism of topoi Cartier} 
	\begin{equation} \label{morphism C}
		\rmC_{X/\rW}:\widetilde{\underline{\mathscr{E}}}\to \widetilde{\mathscr{E}}'
	\end{equation}
	such that its inverse image functor is induced by the composition with $\rho$. 
\end{nothing}
\begin{nothing} \label{descent data Oyama intro}
	Let $n$ be an integer $\ge 1$. The contravariant functor $(U,\mathfrak{T},u)\mapsto \Gamma(\mathfrak{T},\mathscr{O}_{\mathfrak{T}_{n}})$ defines a sheaf of rings on $\pCRIS$ (resp. $\CRIS$) that we denote by $\mathscr{O}_{\pCRIS,n}$ (resp. $\mathscr{O}_{\CRIS,n}$). By definition, we have $\rmC^{*}_{X/\rW}(\mathscr{O}_{\mathscr{E}',n})=\mathscr{O}_{\underline{\mathscr{E}},n}$.
	To give an $\mathscr{O}_{\mathscr{E},n}$-module (resp. $\mathscr{O}_{\underline{\mathscr{E}},n}$-module) $\mathscr{F}$ amounts to give the following data \eqref{lin descent data}: 
	\begin{itemize}
		\item[(i)] For every object $(U,\mathfrak{T},u)$ of $\pCRIS$ (resp. $\CRIS$), an $u_{*}(\mathscr{O}_{\mathfrak{T}_{n}})$-module $\mathscr{F}_{(U,\mathfrak{T})}$ of $U_{\zar}$. 

		\item[(ii)] For every morphism $f:(U_{1},\mathfrak{T}_{1},u_{1})\to (U_{2},\mathfrak{T}_{2},u_{2})$ of $\pCRIS$ (resp. $\CRIS$), an $u_{1*}(\mathscr{O}_{\mathfrak{T}_{1,n}})$-linear morphism
			\begin{displaymath}
				c_{f}:u_{1*}(\mathscr{O}_{\mathfrak{T}_{1,n}})\otimes_{(u_{2*}(\mathscr{O}_{\mathfrak{T}_{2,n}}))|_{U_{1}}} (\mathscr{F}_{(U_{2},\mathfrak{T}_{2})})|_{U_{1}}\to \mathscr{F}_{(U_{1},\mathfrak{T}_{1})},
			\end{displaymath}
	\end{itemize}
	satisfying a cocycle condition for the composition of morphisms as in (\cite{BO} 5.1). 
	
	Following (\cite{BO} 6.1), we say that $\mathscr{F}$ is a \textit{crystal} if $c_{f}$ is an isomorphism for every morphism $f$ and that $\mathscr{F}$ is \textit{quasi-coherent} if $\mathscr{F}_{(U,\mathfrak{T})}$ is a quasi-coherent $u_{*}(\mathscr{O}_{\mathfrak{T}_{n}})$-module of $U_{\zar}$ for every object $(U,\mathfrak{T},u)$. 
	We denote by $\mathscr{C}^{\qcoh}(\mathscr{O}_{\mathscr{E},n})$ (resp. $\mathscr{C}^{\qcoh}(\mathscr{O}_{\underline{\mathscr{E}},n})$) the category of quasi-coherent crystals of $\mathscr{O}_{\pCRIS,n}$-modules (resp. $\mathscr{O}_{\CRIS,n}$-modules). 
\end{nothing}

The following are the main results of this article. 

\begin{prop}[\ref{equi crystals stratification}] \label{crystal stratification data intro}
	Let $\XX$ be a smooth formal $\SS$-scheme and $X$ its special fiber. 
	There exists a canonical equivalence of categories between the category $\CC^{\qcoh}(\mathscr{O}_{\pCRIS,n})$ (resp. $\CC^{\qcoh}(\mathscr{O}_{\CRIS,n})$) and the category of quasi-coherent $\mathscr{O}_{\XX_{n}}$-modules with $\mathscr{O}_{\RR_{\XX}}$-stratification (resp. $\mathscr{O}_{\QQ_{\XX}}$-stratification) (\ref{intro Hopf algebra}, \ref{intro RX}, \ref{intro QX}). 
\end{prop}
\begin{theorem}[\ref{thm pullback Cartier}] \label{thm 1 intro}
	Let $X$ be a smooth $k$-scheme. Then, for any $n\ge 1$, the inverse image and the direct image functors of the morphism $\rmC_{X/\rW}$ \eqref{morphism C} induce equivalences of categories quasi-inverse to each other
	\begin{equation} \label{Cartier inverse}
		\CC^{\qcoh}(\mathscr{O}_{\pCRIS',n})\rightleftarrows\CC^{\qcoh}(\mathscr{O}_{\CRIS,n}).
	\end{equation}
\end{theorem}

The theorem is proved by fppf descent for quasi-coherent modules.

We call \textit{Cartier equivalence modulo $p^{n}$} the equivalence of categories $\rmC_{X/\rW}^{*}$ \eqref{Cartier inverse}. 
Indeed, given a smooth formal $\rW$-scheme $\XX$ with special fiber $X$, Oyama proved \ref{thm 1 intro} in the case $n=1$ and showed that $\rmC_{X/\rW}^{*}$ is compatible with the Cartier transform of Ogus--Vologodsky defined by the lifting $\XX_{2}'$ of $X'$ (cf. \cite{Oy} section 1.5). 
In section \ref{Comparison of Cartier}, we reprove the later result in a different way \eqref{comparison Cartier}. 

The following result explains the relation between the Cartier equivalence $\rmC_{X/\rW}^{*}$ and Shiho's construction, in the presence of a lifting of Frobenius. 

\begin{prop}[\ref{Cartier global local}]\label{Shiho Cartier compatible}
	Let $\XX$ be a smooth formal $\rW$-scheme, $X$ its special fiber, $F:\XX\to \XX'$ a lifting of the relative Frobenius morphism $F_{X/k}$ of $X$ and $\psi_{n}^{*}$ the functor defined by $F$ in \eqref{intro carre stratification}. Then, the following diagram \eqref{crystal stratification data intro}
	\begin{equation} \label{intro diag Cartier Shiho}
	\xymatrix{
		\mathscr{C}(\mathscr{O}_{\mathscr{E}',n})\ar[r]^{\rmC^{*}_{X/\rW}} \ar[d]_{\mu}^{\wr} &
		\mathscr{C}(\mathscr{O}_{\underline{\mathscr{E}},n}) \ar[d]^{\nu}_{\wr}\\
			\Big\{ \txt{\textnormal{$\mathscr{O}_{\XX'_{n}}$-modules}\\
		\textnormal{with $\mathscr{O}_{\RR_{\XX'}}$-stratification}} \Big\}
		\ar[r]^{\psi_{n}^{*}} &
		\Big\{	\txt{\textnormal{$\mathscr{O}_{\XX_{n}}$-modules}\\
		\textnormal{with $\mathscr{O}_{\QQ_{\XX}}$-stratification}} \Big\} }
	\end{equation}
is commutative up to a functorial isomorphism. That is, for every crystal $\mathscr{M}$ of $\mathscr{O}_{\pCRIS',n}$-modules of $\widetilde{\pCRIS}'$, we have a canonical functorial isomorphism
		\begin{equation}
			\eta_{F}:\psi_{n}^{*}(\mu(\mathscr{M}))\xrightarrow{\sim} \nu(\rmC^{*}_{X/\rW}(\mathscr{M})).
		\end{equation} 
	\end{prop}

	In the diagram \eqref{intro diag Cartier Shiho}, while $\rmC_{X/\rW}^{*}$ does not depend on models of $X$, $\psi_{n}^{*}$ depends on the lifting $F$ of the relative Frobenius morphism and the vertical functors $\mu$, and $\nu$ depend on the formal model $\XX$ of $X$. The isomorphism $\eta_{F}$ depends also on $F$. 
	For different choice of liftings of Frobenius, $\eta_F$ can be related by an explicit formula encorded in Oyama topos \eqref{alpha iso stra}. 
	
	By \ref{intro F induces phi psi}, the equivalence of categories $\rmC_{X/\rW}^{*}$ \eqref{Cartier inverse} is compatible with Shiho's functor $\Phi_{n}$ defined by $F$ \eqref{Shiho functor eq}.
	In the case $n=1$, an analogous relation between the Cartier transform and $\Phi_1$ is shown in (\cite{OV07} 2.11). 

\begin{nothing} \label{fil mod W}
	In \cite{FL}, Fontaine and Laffaille introduced the notion of \textit{Fontaine module} to study $p$-adic Galois representations. It is inspired by the work of Mazur (\cite{Ma72}, \cite{Ma73}) and Ogus (\cite{BO}, \S~8) on the Katz conjecture. 
	Let $\sigma:\rW\to \rW$ be the Frobenius endomorphism and $K_{0}=\rW[\frac{1}{p}]$. A Fontaine module is a triple $(M,M^{\bullet},\varphi^{\bullet})$ made of a $\rW$-module of finite length $M$, a decreasing filtration $\{M^{i}\}_{i\in \mathbb{Z}}$ such that $M^{0}=M$, $M^{p}=0$ and $\rW$-linear morphisms
	\begin{equation} \label{Frob div}
		\varphi^{i}:\rW\otimes_{\sigma,\rW}M^{i} \to M, \quad 0\le i\le p-1,
	\end{equation}
	such that $\varphi^{i}|_{M^{i+1}}=p\varphi^{i+1}$ and $\sum_{i=0}^{p-1}\varphi^{i}(M^{i})=M$. The $\varphi^{i}$'s are called \textit{divided Frobenius morphisms}. 
	
	The main result of Fontaine--Laffaille is the construction of a fully faithful and exact functor from the category of Fontaine modules of length $\le p-2$ to the category of torsion $\mathbb{Z}_{p}$-representations of the Galois group $G_{K_{0}}$ of $K_{0}$ (\cite{Wach} Thm. 2). Its essential image consists of torsion crystalline representations of $G_{K_{0}}$ with weights $\le p-2$ (cf. \cite{BM06} 3.1.3.3). 

	Fontaine and Messing showed that there exists a natural Fontaine module structure on the crystalline cohomology of a smooth proper scheme $\mathcal{X}$ over $\rW$ of relative dimension $\le p-1$ (\cite{FM} II.2.7). Then they deduced the degeneration of the Hodge to de Rham spectral sequence of $\mathcal{X}_{n}/\rW_{n}$. 
\end{nothing}
\begin{nothing} \label{fil mod Faltings}
	A generalisation of Fontaine modules in a relative situation was proposed by Faltings in \cite{Fal89}. Relative Fontaine modules can be viewed as an analogue of variation of Hodge structures on smooth formal schemes over $\rW$. 
	Let $\XX=\Spf(R)$ be an affine smooth formal scheme over $\rW$, $X$ its special fiber and $F:\XX\to \XX$ a $\sigma$-lifting of the absolute Frobenius morphism $F_{X}$ of $X$. 
	\textit{A Fontaine module over $\XX$ with respect to $F$} is a quadruple  $(M,\nabla,M^{\bullet},\varphi_{\rF}^{\bullet})$ made of a coherent, torsion $\mathscr{O}_{\XX}$-module $M$, an integrable connection $\nabla$ on $M$, a decreasing exhaustive filtration $M^{\bullet}$ on $M$ of length at most $p-1$ satisfying Griffiths' transversality, and a family of divided Frobenius morphisms $\{\varphi_{F}^{\bullet}\}$ as in \eqref{Frob div} satisfying a compatibility condition between $\{\varphi_{F}^{\bullet}\}$ and $\nabla$ (cf. \cite{Fal89} 2.c, 2.d).

	Using the connection, Faltings glued the categories of Fontaine modules with respect to different Frobenius liftings by a \textit{Taylor formula} (cf. \cite{Fal89} Thm. 2.3). By gluing local data, he defined Fontaine modules over a general smooth formal $\rW$-scheme $\XX$, even if there is no lifting of $F_{X}$. 

	If $\XX$ is the $p$-adic completion of a smooth, proper $\rW$-scheme $\mathcal{X}$, Faltings associated to each Fontaine module of length $\le p-2$ over $\XX$ a representation of the \'etale fundamental group of $\mathcal{X}_{K_{0}}$ on a torsion $\mathbb{Z}_{p}$-module. 
	Moreover, Faltings generalised Fontaine--Messing's result for the crystalline cohomology of a relative Fontaine module.
\end{nothing}
\begin{nothing} \label{fil mod ptorsion}
	Let $\XX$ be a smooth formal $\rW$-scheme. 
	As an application of their Cartier transform \cite{OV07}, Ogus and Vologodsky proposed an interpretation of $p$-torsion Fontaine modules over $\XX$ (\cite{OV07}, 4.16). \textit{A $p$-torsion Fontaine module over $\XX$} is a triple $(M,\nabla,M^{\bullet})$ as in \ref{fil mod Faltings} such that $pM=0$ and equipped with a horizontal isomorphism
	\begin{equation} \label{morphism phi intro}
	\phi:\rmC^{-1}_{\XX_{2}'}(\pi^{*}(\Gr(M),\kappa))\xrightarrow{\sim} (M,\nabla),
\end{equation}
where $\kappa$ is the Higgs field on $\Gr(M)$ induced by $\nabla$ and Griffiths' transversality, and $\pi:X'\to X$ is the base change of the Frobenius morphism of $k$ to $X$. 
	Given a $\sigma$-lifting $F:\XX\to \XX$ of $F_{X}$, such a morphism $\phi$ is equivalent to a family of divided Frobenius morphisms $\{\varphi_{F}^{\bullet}\}$ with respect to $F$ \eqref{fil mod Faltings} (cf. \ref{phiF familly phi}). 

	By Griffiths' transversality, the de Rham complex $M\otimes_{\mathscr{O}_{X}}\Omega_{X/k}^{\bullet}$ is equipped with a decreasing filtration
	\begin{equation}
		\rF^{i}(M\otimes_{\mathscr{O}_{X}}\Omega_{X/k}^{q})=M^{i-q}\otimes_{\mathscr{O}_{X}}\Omega_{X/k}^{q}.
	\end{equation}
	Let $\ell$ be the length of the filtration $M^{\bullet}$ (i.e. $M^{0}=M, M^{\ell+1}=0$) and $d$ the relative dimension of $\XX$ over $\rW$. 

	(i) For any $i,m$, the canonical morphism $\mathbb{H}^{m}( \rF^{i+1}(M\otimes_{\mathscr{O}_{X}}\Omega_{X/k}^{\bullet}))\to \mathbb{H}^{m}( \rF^{i}(M\otimes_{\mathscr{O}_{X}}\Omega_{X/k}^{\bullet}))$ is injective. 
	The morphism $\phi$ \eqref{morphism phi intro} induces a family of divided Frobenius morphisms on $(\mathbb{H}^{m}(M\otimes_{\mathscr{O}_{X}}\Omega_{X/k}^{\bullet}),\{\mathbb{H}^{m}( \rF^{i})\}_{i\le p-1})$ which make it into a Fontaine module over $\rW$ \eqref{fil mod W}. 
	
	(ii) The hypercohomology spectral sequence of the filtered de Rham complex $(M\otimes_{\mathscr{O}_{X}}\Omega_{X/k}^{\bullet},\rF^{i})$ degenerates at $\rE_{1}$.
\end{nothing}
\begin{nothing} \label{fil mod pn}
	For any $n\ge 1$, using the Cartier transform modulo $p^{n}$, we reformulate Faltings' definition of $p^{n}$-torsion Fontaine modules over $\XX$ following Ogus--Vologodsky (\ref{definition MF}). The Taylor formula used by Faltings to glue the data relative to different liftings of $F_{X}$ is naturally encoded in Oyama topos (\ref{coro MF1 MF2}). Following Faltings' strategy, we prove the analogue of the previous results (\ref{fil mod ptorsion}(i-ii)) on the crystalline cohomology of a $p^{n}$-torsion Fontaine module over $\XX$ \eqref{fil mod coh fil mod}. 
\end{nothing}

\begin{nothing}
	Section \ref{not and con} contains the main notation and general conventions. In section \ref{blowup and dila}, we recall the notion of dilatation in formal geometry. In section \ref{Hopf groupoids}, after recalling the notions of Hopf algebras and groupoids, we present the constructions of groupoids $\RR_{\XX}$ and $\QQ_{\XX}$ (\ref{intro RX}, \ref{intro QX}). In section \ref{conn and stra}, we recall the notions of modules with integrable $\lambda$-connection \eqref{lamda connection} and of modules with stratification \eqref{intro Hopf algebra} and we discuss the relation between them. Following \cite{Shiho}, we present the construction of Shiho's functor $\Phi_{n}$ \eqref{Shiho functor eq} in section \ref{local Shiho}. 
	In section \ref{Oyama topos}, we explain the Oyama topoi $\widetilde{\mathscr{E}}$ and $\widetilde{\underline{\mathscr{E}}}$ \eqref{Oyama topos intro} and their fppf variants. Section \ref{crystals} is devoted to the study of crystals in the Oyama topoi (\ref{descent data Oyama intro}). In section \ref{Cartier equiv}, we study the morphism of topoi $\rmC_{X/\rW}$ \eqref{morphism C} and prove our main results \ref{thm 1 intro} and \ref{Shiho Cartier compatible}. 
	We recall the construction of the Cartier transform of Ogus--Vologodsky \cite{OV07} in section \ref{Cartier OV}. Section \ref{D ring} is devoted to several rings of differential operators after Oyama and serves as a preparation for next section. In section \ref{Comparison of Cartier}, we compare the Cartier equivalence $\rmC_{X/\rW}^{*}$ \eqref{Cartier inverse} and the Cartier transform of Ogus--Vologodsky. 
	In section \ref{Fil modules}, we introduce a notion of relative Fontaine modules using Oyama topoi \eqref{fil mod pn}. We compare it with Faltings' definition \cite{Fal89} and Tsuji's definition \cite{Tsu96}. In section \ref{coh fil mod}, we construct a Fontaine module structure on the crystalline cohomology of a relative Fontaine module.
\end{nothing}

After finishing this article, I learned from Arthur Ogus that Vadim Vologodsky has skecthed a similar approach for lifting the Cartier transform in a short note \footnote{available at \href{http://pages.uoregon.edu/vvologod/papers/p-adiccartier.pdf}{http://pages.uoregon.edu/vvologod/papers/p-adiccartier.pdf}} without providing any detail.

A different approach to the formulation of a Cartier transform modulo $p^n$ and its relationship to Fontaine modules was taken in \cite{LSZ}.\\

\textbf{Acknowledgement.} This article is a part of my thesis prepared at Universit\'e Paris-Sud and IHÉS. I would like to express my great gratitude to my thesis advisor Ahmed Abbes for leading me to this question and for his helpful comments on earlier versions of this work. 
I would like to thank Arthur Ogus, Takeshi Tsuji and an anonymous referee for their careful reading and suggestions. 

\section{Notations and conventions} \label{not and con}
\begin{nothing} \label{notations}
	In this article, $p$ denotes a prime number, $k$ a perfect field of characteristic $p$, $\rW$ the ring of Witt vectors of $k$ and $\sigma:\rW\to \rW$ the Frobenius automorphism of $\rW$. For any integer $n\ge 1$, we set $\rW_{n}=\rW/p^{n}\rW$ and $\SS=\Spf(\rW)$.
\end{nothing}

\begin{nothing} \label{notations Yk}
	Let $X$ be a scheme over $k$. We denote by $F_{X}$ the absolute Frobenius morphism of $X$ and by $F_{X/k}:X\to X'=X\otimes_{k,F_{k}}k$ the relative Frobenius morphism. Then we have a commutative diagram
	\begin{equation}
		\xymatrix{
			X\ar[r]^{F_{X/k}} \ar[rd] & X' \ar[r] \ar@{}[dr]|{\Box} \ar[d]& X \ar[d] \\
			& \Spec{k} \ar[r]^{F_{k}} & \Spec{k}
		}
	\end{equation}
\end{nothing}

\begin{nothing}\label{notation underline}
	Let $X$ be a scheme over $k$. We denote by $\underline{X}$ the scheme theoretic image of $F_{X}:X\to X$ (\cite{EGAInew} 6.10.1 and 6.10.5). By (\cite{EGAInew} 6.10.4), $\underline{X}$ is the closed subscheme of $X$ defined by the ideal sheaf of $\mathscr{O}_{X}$ consisting of the sections of $\mathscr{O}_{X}$ whose $p$th power is zero. It is clear that the correspondence $X\mapsto \underline{X}$ is functorial. Note that the canonical morphism $\underline{X}\to X$ induces an isomorphism of the underlying topological spaces. 
	
	The relative Frobenius morphism $F_{X/k}:X\to X'$ factors through $\underline{X}'$. We denote the induced morphism by $f_{X/k}:X\to \underline{X}'$. By definition, the homomorphism $\mathscr{O}_{\underline{X}'}\to f_{X/k*}(\mathscr{O}_{X})$ is injective, i.e. $f_{X/k}$ is scheme theoretically dominant (\cite{EGAInew} 5.4.2).

	If $Y\to X$ and $Z\to X$ are two morphisms of $k$-schemes, by functoriality, we have a canonical morphism
	\begin{equation}
		\underline{Y\times_{X}Z}\to \underline{Y}\times_{\underline{X}}\underline{Z}. \label{underline product}
	\end{equation}
	Since $X$ is affine if and only if $\underline{X}$ is affine (cf. \cite{EGAInew} 2.3.5), we verify that \eqref{underline product} is an affine morphism.
\end{nothing}
\begin{nothing} \label{conventions adic}
	In this article, we follow the conventions of \cite{Ab10} for adic rings (\cite{Ab10} 1.8.4) and adic formal schemes (\cite{Ab10} 2.1.24). Note that these notions are stronger than those introduced by Grothendieck in (\cite{EGAInew} 0.7.1.9 and 10.4.2).

	If $\XX$ is an adic formal scheme such that $p\mathscr{O}_{\XX}$ is an ideal of definition of $\XX$,  for any integer $n\ge 1$, we denote the usual scheme $(\XX,\mathscr{O}_{\XX}/p^{n}\mathscr{O}_{\XX})$ by $\XX_{n}$.
\end{nothing}
\begin{nothing} \label{notations S flat}
	We say that an adic formal $\SS$-scheme $\XX$ (\cite{Ab10} 2.2.7) is \textit{flat over $\SS$} (or that $\XX$ is a flat formal $\SS$-scheme) if the morphism $\mathscr{O}_{\XX}\to \mathscr{O}_{\XX}$ induced by multiplication by $p$ is injective (i.e. if $\mathscr{O}_{\XX}$ is rig-pur in the sense of (\cite{Ab10} 2.10.1.4.2)). It is clear that the above condition is equivalent to the fact that, for every affine open formal subscheme $U$ of $\XX$, the algebra $\Gamma(U,\mathscr{O}_{\XX})$ is flat over $\rW$.

	Let $A$ be an adic $\rW$-algebra (\cite{Ab10} 1.8.4.5). Then $A$ is flat over $\rW$ if and only if $A_{n}=A/p^{n}A$ is flat over $\rW_{n}$ for all integers $n\ge 1$. Indeed, we only need to prove that the condition is sufficient. Let $a$ be an element of $A$ such that $pa=0$. For any integer $n\ge 1$, by the flatness of $A_{n}$ over $\rW_{n}$, the image of $a$ in $A_{n}$ is contained in $p^{n-1}A/p^{n}A$. Since $A$ is separated, we see that $a=0$, i.e. $A$ is flat over $\rW$. We deduce that an adic formal $\SS$-scheme $\XX$ is flat over $\SS$ if and only if $\XX_{n}$ is flat over $\SS_{n}$ \eqref{conventions adic} for all integers $n\ge 1$.
\end{nothing}

\section{Blow-ups and dilatations} \label{blowup and dila}
\begin{nothing} \label{basic setting top rings}
	Let $A$ be an adic ring \eqref{conventions adic}, $J$ an ideal of definition of finite type of $A$. We put $X=\Spec(A)$, $X'=\Spec(A/J)$ and $\XX=\Spf(A)$. The formal scheme $\XX$ is the completion of $X$ along $X'$. For any $A$-module $M$, we denote by $\widetilde{M}$ the associated $\mathscr{O}_{X}$-module and by $M^{\Delta}$ the completion of $\widetilde{M}$ along $X'$ (\cite{Ab10} 2.7.1), which is an $\mathscr{O}_{\XX}$-module.

	Let $\aa$ be an open ideal of finite type of $A$. We denote by $\aa\mathscr{O}_{\XX}$ the ideal sheaf of $\mathscr{O}_{\XX}$ associated to the presheaf defined by $U\mapsto \aa\Gamma(U,\mathscr{O}_{\XX})$. By (\cite{Ab10} 2.1.13), we have $\aa^{\Delta}=\aa\mathscr{O}_{\XX}$.

	Let $B$ be an adic ring, $u:A\to B$ an adic homomorphism (\cite{Ab10} 1.8.4.5) and $f:\mathfrak{Y}=\Spf(B)\to \XX=\Spf(A)$ the associated morphism. In view of (\cite{Ab10} 2.5.11), we have a canonical functorial isomorphism
	\begin{equation}
		f^{*}(\aa^{\Delta})\xrightarrow{\sim} (\aa\otimes_{A}B)^{\Delta}. \label{isomorphism pullback Delta}
	\end{equation}
	Then, we deduce a canonical isomorphism
	\begin{equation} \label{lemma ideal ring Delta}
		f^{*}(\aa^{\Delta})\mathscr{O}_{\YY}\xrightarrow{\sim}(\aa B)^{\Delta}.		
	\end{equation}
	Indeed, by definition, $f^{*}(\aa^{\Delta})\mathscr{O}_{\YY}$ is the image of the morphism $f^{*}(\aa^{\Delta})\to \mathscr{O}_{\YY}=f^{*}(\mathscr{O}_{\XX})$, which clearly factors through $(\aa B)^{\Delta}$. The isomorphism \eqref{lemma ideal ring Delta} follows from the fact that $\aa^{\Delta}=\aa\mathscr{O}_{\XX}$ and $(\aa B)^{\Delta}=(\aa B)\mathscr{O}_{\YY}$.
\end{nothing}

\begin{nothing} \label{blowup formal sch}
	Let $\XX$ be an adic formal scheme, $\mathscr{J}$ an ideal of definition of finite type of $\XX$ and $\mathscr{A}$ an open ideal of finite type of $\XX$ (\cite{Ab10} 2.1.19). For any $n\ge 1$, we denote by $\XX_{n}$ the usual scheme $(\XX,\mathscr{O}_{\XX}/\mathscr{J}^{n})$ and we set
	\begin{equation} \label{blowup mod In}
		\XX'_{n}=\Proj(\oplus_{m\ge 0}\mathscr{A}^{m}\otimes_{\mathscr{O}_{\XX}}\mathscr{O}_{\XX_{n}}).
	\end{equation}
	The sequence $(\XX_{n}')$ forms an adic inductive $(\XX_{n})$-system (\cite{Ab10} 2.2.13). We call its inductive limit $\XX'$ the \textit{admissible blow-up of $\mathscr{A}$ in $\XX$} (\cite{Ab10} 3.1.2). By (\cite{Ab10} 2.2.14, 2.3.13), $\XX'$ is an adic formal $\XX$-scheme of finite type. 
	Note that $\XX_{n}'$ is different from the blow-up of $\XX_{n}$ along $(\mathscr{A}+\mathscr{J}^n)/\mathscr{J}^n$. 
\end{nothing}

\begin{nothing}\label{dilatation def}
	Let $\XX$ be a flat formal $\SS$-scheme locally of finite type (\ref{notations S flat}, \cite{Ab10} 2.3.13) and $\mathscr{A}$ an open ideal of finite type of $\XX$ containing $p$. Let $\varphi:\XX'\to \XX$ be the admissible blow-up of $\mathscr{A}$ in $\XX$. Then the ideal $\mathscr{A}\mathscr{O}_{\XX'}$ is invertible (\cite{Ab10} 3.1.4(i)), and $\XX'$ is flat over $\SS$ (\cite{Ab10} 3.1.4(ii)). We denote by $\XX_{(\mathscr{A}/p)}$ the maximal open formal subscheme of $\XX'$ on which 
	\begin{equation}
		(\mathscr{A}\mathscr{O}_{\mathfrak{X}'})|\mathfrak{X}_{(\mathscr{A}/p)}=(p\mathscr{O}_{\mathfrak{X}'})|\mathfrak{X}_{(\mathscr{A}/p)}
	\end{equation}
	and we call it \textit{the dilatation of $\mathscr{A}$ with respect to $p$}. 
	
	Note that $\XX_{(\mathscr{A}/p)}$ is the complement of $\Supp(\mathscr{A}\mathscr{O}_{\XX'}/p\mathscr{O}_{\XX'})$ in $\XX'$ (\cite{EGAInew} 0.5.2.2). In view of (\cite{Ab10} 3.1.5 and 3.2.7), the above definition coincides with the notion of dilatation of $\mathscr{A}$ with respect to $p$ introduced in (\cite{Ab10} 3.2.3.4). We denote the restriction of $\varphi:\XX'\to \XX$ to $\XX_{(\mathscr{A}/p)}$ by 
	\begin{equation}
		\psi:\XX_{(\mathscr{A}/p)}\to \XX. \label{morphism dilatation}
	\end{equation}

	We set $\mathscr{A}^{\sharp}=\mathscr{A}^p+p\mathscr{O}_{\XX}$ the open ideal of $\XX$. If $\mathscr{A}$ is locally generated by $\{a_1,\cdots,a_n\}$ then $\mathscr{A}^{\sharp}$ is locally generated by $\{p,a_1^p,\cdots,a_n^p\}$. 
\end{nothing}

\begin{nothing} \label{local description blow up}
		Keep the notation of \ref{dilatation def} and assume moreover that $\XX=\Spf(A)$ is affine. There exists an open ideal of finite type $\aa$ of $A$ containing $p$ such that $\aa^{\Delta}=\mathscr{A}$ (\cite{Ab10} 2.1.10 and 2.1.13). Let $X=\Spec(A)$, $Y=\Spec(A/pA)$, $\phi:X'\to X$ be the blow-up of $\widetilde{\aa}$ in $X$ and $Y'=\phi^{-1}(Y)$; so $\XX$ is the completion of $X$ along $Y$. Then $\XX'$ is canonically isomorphic to the completion of $X'$ along $Y'$ and $\varphi$ is the extension of $\phi$ to the completions (\cite{Ab10} 3.1.3).

		Let $(a_{i})_{0\le i\le n}$ be a finite set of generators of $\aa$ such that $a_{0}=p$. For any $0\le i\le n$, let $U_{i}$ be the maximal open subset of $X'$ where $a_{i}$ generates $\widetilde{\aa}\mathscr{O}_{X'}$. Since $\widetilde{\aa}\mathscr{O}_{X'}$ is invertible, $(U_{i})_{0\le i\le n}$ form an open covering of $X'$. It is well known that $U_{i}$ is the affine scheme associated to the $A$-algebra $A_{i}$ defined as follows:
		\begin{displaymath}
			A_{i}'=A\Big[\frac{a_{0}}{a_{i}},\cdots,\frac{a_{n}}{a_{i}}\Big]=\frac{A[x_{0},\ldots,x_{i-1},x_{i+1},\cdots,x_{n}]}{(a_{j}-a_{i}x_{j})_{j\neq i}},\quad
			A_{i}=\frac{A_{i}'}{(A_{i}')_{a_{i}\btor}},
		\end{displaymath}
		where $(A_{i}')_{a_{i}\btor}$ denotes the ideal of $a_{i}$-torsion elements of $A_{i}'$. Let $\widehat{A}_{i}$ be the separated completion of $A_{i}$ for the $p$-adic topology and $\widehat{U}_{i}=\Spf(\widehat{A}_{i})$. Then $(\widehat{U}_{i})_{0\le i\le n}$ form a covering of $\XX'$; for any $0\le i\le n$, $\widehat{U}_{i}$ is the maximal open of $\XX'$ where $a_{i}$ generates the invertible ideal $\aa^{\Delta}\mathscr{O}_{\XX'}$ (\cite{Ab10} 3.1.7(ii)). The open formal subscheme $\XX_{(\mathscr{A}/p)}$ of $\XX'$ is equal to $\widehat{U}_{0}=\Spf(\widehat{A}_{0})$. In particular, we see that, in the general setting of \ref{dilatation def}, $\psi:\XX_{(\mathscr{A}/p)}\to \XX$ is affine (\cite{Ab10} 2.3.4).

		Let $A\{x_{1},\cdots,x_{n}\}$ denote the $p$-adic completion of the polynomial ring in $n$ variables $A[x_{1},\cdots,x_{n}]$, which is flat over $A[x_{1},\cdots,x_{n}]$ by (\cite{Ab10} 1.12.12). If $A_0'$ is flat over $\rW$, then we have $A_{0}=A_{0}'$ and deduce a canonical isomorphism (\cite{Ab10} 1.12.16(iv))
	\begin{equation} \label{hat A0 blowup}
		A\Big\{\frac{a_{1}}{p},\cdots,\frac{a_{n}}{p}\Big\}=\frac{A\{x_{1},\cdots,x_{n}\}}{(a_{i}-px_{i})_{1\le i\le n}}\xrightarrow{\sim} \widehat{A}_{0}.
	\end{equation}
\end{nothing}

\begin{prop} \label{prop univ RQ}
	Let $\XX$ be a flat formal $\SS$-scheme locally of finite type, $i:T\to \XX_{1}$ an immersion (not necessary closed). There exists a formal $\XX$-scheme $\psi:\XX_{(T/p)}\to \XX$ (resp. $\psi:\XX^{\sharp}_{(T/p)}\to \XX$) unique up to canonical isomorphisms satisfying the following conditions: 

	\textnormal{(i)} The canonical morphism $(\XX_{(T/p)})_{1}\to \XX_{1}$ (resp. $\underline{(\XX_{(T/p)}^{\sharp})_{1}}\to (\XX_{(T/p)}^{\sharp})_{1}\to \XX_{1}$) factors through the immersion $T\to \XX_{1}$. 
	
	\textnormal{(ii)} Let $\YY$ be a flat formal $\SS$-scheme, $Y=\YY_{1}$ and $f:\YY\to \XX$ an $\SS$-morphism. 
	Suppose that there exists a $k$-morphism $g:Y\to T$ (resp. $g:\underline{Y}\to T$) which fits into the following diagram:
	\begin{equation} \label{diag univ RQ}
		\xymatrix{
			Y\ar[r] \ar[d]_{g} & \YY \ar[d]^{f} \\
			T\ar[r] & \XX} 
		\qquad \textnormal{(resp. } 
		\xymatrix{
			\underline{Y}\ar[r] \ar[d]_{g} & \YY \ar[d]^{f} \\
			T\ar[r] & \XX} )
	\end{equation}
	Then there exists a unique $\SS$-morphism $f': \YY\to \XX_{(T/p)}$ (resp. $f': \YY\to \XX_{(T/p)}^{\sharp}$) such that $f=\psi\circ f'$. 
	If $T\to \XX$ and $f$ are moreover closed immersion, then so is $f':\YY\to \XX_{(T/p)}$. 
\end{prop} 
\begin{proof} It suffices to prove the existence. The uniqueness follows from (i) and the universal property (ii). 
	We first prove the case where $T\to \XX$ is a closed immersion and denote the associated ideal sheaf by $\mathscr{A}$.

	In the first situation, we take $\XX_{(T/p)}$ to be the dilatation $\XX_{(\mathscr{A}/p)}$. 
	Since $p\in \mathscr{A}$, the commutativity of the first diagram of \eqref{diag univ RQ} is equivalent to $\mathscr{A}\mathscr{O}_{\YY}=p\mathscr{O}_{\YY}$. 
	To verify condition (ii), we can reduce to the case where $\XX=\Spf(A)$ is affine and then to the case where $\YY=\Spf(B)$ is affine and the morphism $f$ is associated to an adic homomorphism $u:A\to B$. We take again the notation of \ref{local description blow up}. By \eqref{lemma ideal ring Delta}, we have $(pB)^{\Delta}=(\aa B)^{\Delta}$. The open ideals of finite type $pB$ and $\aa B$ are complete by (\cite{Comalg} III \S 2.12 Cor 1 of Prop 16) and separated as submodules of $B$. We deduce that $pB=\aa B$ by taking $\Gamma(\YY,-)$. Since $B$ is flat over $\rW$, the homomorphism $u$ extends uniquely to an $A$-homomorphism $A_{0}\to B$ and then to an adic $A$-homomorphism $w:\widehat{A}_{0}\to B$ by $p$-adic completion.
	In the first situation, we take for $f'$ the morphism induced by $w$ which is uniquely determined by $f$. 
	If $u:A\to B$ is moreover surjective, then so is $w$. 

	In the second situation, the commutativity of the second diagram of \eqref{diag univ RQ} means the $p$-th power of every local section of $\mathscr{A}\mathscr{O}_{\YY}/p\mathscr{O}_{\YY}$ is zero in $\mathscr{O}_{Y}$, which is equivalent to the fact that $\mathscr{A}^{\sharp}\mathscr{O}_{\YY}=p\mathscr{O}_{\YY}$ \eqref{dilatation def}. 
	We take $\XX^{\sharp}_{(T/p)}$ to be the dilatation $\XX_{(\mathscr{A}^{\sharp}/p)}$.
	Then condition (ii) in this situation follows from the first situation.

	In general, let $\UU$ be an open formal subscheme of $\XX$ such that $i(T)\subset \UU$ and that $T\to \UU$ is a closed immersion and, let $\mathscr{A}$ be the open ideal of finite type associated to $T\to \UU$. 
	We take $\XX_{(T/p)}$ (resp. $\XX^{\sharp}_{(T/p)}$) to be the dilatation $\UU_{(\mathscr{A}/p)}$ (resp. $\UU_{(\mathscr{A}^{\sharp}/p)}$). 
	For any morphism $f:\YY\to \XX$ as in (ii), by \eqref{diag univ RQ}, the morphism $f$ factors through the open subscheme $\UU$ of $\XX$. Then the assertion follows from the case where $T\to \XX$ is a closed immersion.
\end{proof}
\begin{rem}
	The formal $\XX$-scheme $\XX_{(T/p)}$ (resp. $\XX^{\sharp}_{(T/p)}$ is the same as the Higgs envelope $\widetilde{R}_{T}(\mathfrak{X})$ (resp. $\widetilde{Q}_{T}(\mathfrak{X})$) introduced in (\cite{Oy} page 6).
\end{rem}

The next result shows that the constructions of \ref{prop univ RQ} are compatible with \'etale localization.

\begin{prop}[\cite{Oy} 1.1.3]\label{prop dilatation etale}
	Let $\XX$, $\YY$ be two flat formal $\SS$-schemes locally of finite type, $f:\XX\to \YY$ an \'etale $\SS$-morphism \textnormal{(\cite{Ab10} 2.4.5)}. Suppose that there exist a $k$-scheme $T$ and two immersions $i:T\to \XX_{1}$, $j:T\to \YY_{1}$ such that $j=f\circ i$. Then $f$ induces canonical isomorphisms of the formal schemes
	\begin{equation}
		\XX_{(T/p)}\xrightarrow{\sim} \YY_{(T/p)}, \qquad \XX_{(T/p)}^{\sharp}\xrightarrow{\sim} \YY_{(T/p)}^{\sharp}.
	\end{equation}
\end{prop}
\begin{proof} By the universal property (\ref{prop univ RQ}), the composition $\XX_{(T/p)}\to \XX\to \YY$ induces a canonical $\YY$-morphism $u:\XX_{(T/p)}\to \YY_{(T/p)}$. 
We consider the commutative diagram
\begin{equation} \label{construction of inverse etale}
	\xymatrix{
		(\YY_{(T/p)})_{1}\ar[r] \ar[d] & T \ar[r] & \XX_{1} \ar[d]\\
		\YY_{(T/p)} \ar[rr] && \YY}
\end{equation}
Since $f:\XX\to \YY$ is \'etale, the composition of the top horizontal morphisms and $\XX_{1}\to \XX$ lifts uniquely to a $\YY$-morphism $g:\YY_{(T/p)}\to \XX$. By \eqref{construction of inverse etale} and the universal property, $g$ induces an $\XX$-morphism $v:\YY_{(T/p)}\to \XX_{(T/p)}$. The following diagram
\begin{equation}
	\xymatrix{
		\YY_{(T/p)} \ar[r]^{v} \ar[rd]^{g}& \XX_{(T/p)} \ar[r]^{u} \ar[d] & \YY_{(T/p)}\ar[d]\\
		& \XX \ar[r]^{f} & \YY }
\end{equation}
is commutative. By the universal property, we deduce that $u\circ v=\id$. 

We consider the diagram
\begin{equation}\label{diag comm to prove etale}
	\xymatrix{
		\XX_{(T/p)} \ar[r]^{u} \ar[d]_{\psi} & \YY_{(T/p)}\ar[d] \ar[ld]_{g}\\
		\XX \ar[r]^{f} & \YY }
\end{equation}
where the lower triangle is commutative. Since $\psi_{1}$ and $g_{1}$ factor through $T$, we have $\psi_{1}=(g\circ u)_{1}$. Since $f$ is \'etale and the square of \eqref{diag comm to prove etale} commutes, there exists one and only one lifting of $(\XX_{(T/p)})_{1}\to \XX_{1}$ to a $\YY$-morphism $\XX_{(T/p)}\to \XX$. We deduce that $\psi=g\circ u$, i.e. the diagram \eqref{diag comm to prove etale} commutes. Then we have $\psi\circ v\circ u=\psi$. By the universal property, we deduce that $v\circ u=\id$. The first isomorphism follows.

The second isomorphism can be verified in the same way.
\end{proof}

The next result shows that the constructions of \ref{prop univ RQ} are compatible with flat base change.

\begin{prop}[\cite{Oy} 1.1.4]\label{prop dilatation flat}
	Let $\XX$, $\YY$ be two flat formal $\SS$-schemes locally of finite type, $f:\XX\to \YY$ a flat $\SS$-morphism, $T\to \YY_{1}$ an immersion and $S=T\times_{\YY}\XX$. Then $f$ induces canonical isomorphisms of formal schemes
	\begin{equation}
		\XX_{(S/p)}\xrightarrow{\sim} \YY_{(T/p)}\times_{\YY}\XX, \qquad \XX_{(S/p)}^{\sharp}\xrightarrow{\sim} \YY_{(T/p)}^{\sharp}\times_{\YY}\XX.
	\end{equation}
\end{prop}
\begin{proof} By \ref{prop univ RQ}, we can reduce to the case where $T\to \YY$ is a closed immersion. Let $\mathscr{A}$ (resp. $\mathscr{B}$) be the open ideal of finite type associated to the closed immersion $T\to \YY$ (resp. $S\to \XX$). 
Put $\ZZ=\YY_{(\mathscr{B}/p)}\times_{\YY}\XX$. By the universal property, the morphism $f$ induces a $\YY$-morphism $\XX_{(\mathscr{A}/p)}\to \YY_{(\mathscr{B}/p)}$ and then an $\XX$-morphism $u: \XX_{(\mathscr{A}/p)}\to \ZZ$.

Since $f$ is flat, $\ZZ$ is flat over $\SS$. We have $\mathscr{A}\mathscr{O}_{\ZZ}=\mathscr{B}\mathscr{O}_{\ZZ}=p\mathscr{O}_{\ZZ}$. By the universal property, we deduce an $\XX$-morphism $v:\ZZ\to \XX_{(\mathscr{A}/p)}$. Since $u,v$ are $\XX$-morphisms, we deduce that $u\circ v=\id$ and $v\circ u=\id$ by the universal property as in the proof of \ref{prop dilatation etale}.

The second isomorphism can be verified in the same way.
\end{proof}

\section{Hopf algebras and groupoids} \label{Hopf groupoids}
\begin{nothing}\label{tensor product OX bimodule}
	In this section, we review the notion of Hopf algebras and groupoids following \cite{Ber} and the construction of certain Hopf algebras used in \cite{Oy}. 

	Let $(\mathscr{T},A)$ be a ringed topos. For any $A$-bimodules (resp. $A$-bialgebras) $M$ and $N$ of $\mathscr{T}$, $M\otimes_{A}N$ denotes the tensor product of the right $A$-module $M$ and the left $A$-module $N$, and we regard $M\otimes_{A}N$ as an $A$-bimodule (resp. $A$-bialgebra) through the left $A$-action on $M$ and the right $A$-action on $N$.
\end{nothing}

\begin{definition}[\cite{Ber} II 1.1.2 and \cite{Oy} 1.2.1] \label{def of Hopf alg}
	Let $(\mathscr{T},A)$ be a ringed topos. \textit{A Hopf $A$-algebra} is the data of an $A$-bialgebra $B$ and three ring homomorphisms
	\begin{equation}
		\delta:B\to B\otimes_{A}B~ \textnormal{(comultiplication)},\quad \pi:B\to A~\textnormal{(counit)},\quad \sigma:B\to B~\textnormal{(antipode)}	
	\end{equation}
	satisfying the following conditions.
	\begin{itemize}	
		\item[(a)] $\delta$ and $\pi$ are $A$-bilinear and the following diagrams are commutative:
		\begin{equation} \label{diag def Hopf}
		\xymatrix{
			B\ar[r]^-{\delta} \ar[d]_{\delta}& B\otimes_{A}B \ar[d]^{\delta\otimes \id_{B}} \\ 
			B\otimes_{A}B \ar[r]^-{\id_{B}\otimes \delta} & B\otimes_{A}B\otimes_{A}B} \quad
		\xymatrix{
			B\ar[r]^-{\delta}\ar[d]_{\delta} \ar[rd]^{\id_{B}}& B\otimes_{A}B \ar[d]^{\pi\cdot\id_{B}}\\
			B\otimes_{A}B\ar[r]^{\id_{B}\cdot \pi}& B
		}
	\end{equation}

		\item[(b)] $\sigma$ is a homomorphism of $A$-algebras for the left (resp. right) $A$-action on the source and the right (resp. left) $A$-action on the target, and satisfies $\sigma^{2}=\id_{B}$, $\pi \circ \sigma=\pi$.

		\item[(c)] The following diagrams are commutative:
			\begin{equation} \label{diag def Hopf 2}
	\xymatrix{
		B\ar[r]^{\pi}\ar[d]_{\delta} & A \ar[d]^{d_{1}} \\
		B\otimes_{A}B \ar[r]^-{\id_{B}\cdot \sigma} & B} \qquad
	\xymatrix{
		B\ar[r]^{\pi}\ar[d]_{\delta} & A \ar[d]^{d_{2}} \\
		B\otimes_{A}B \ar[r]^-{\sigma\cdot \id_{B}} & B}
	\end{equation}
	\end{itemize}
	where $d_{1}$ (resp. $d_{2}$) is the structural homomorphism of the left (resp. right) $A$-algebra $B$.
\end{definition}

Such a data is also called an affine groupoid of $(\mathscr{T},A)$ by Berthelot (\cite{Ber} II 1.1.2).

\begin{definition}[\cite{Ber} II 1.1.6] \label{Hopf algebra homomorphism}
	Let $f:(\mathscr{T}',A')\to (\mathscr{T},A)$ be a morphism of ringed topoi, $B$ a Hopf $A$-algebra and $B'$ a Hopf $A'$-algebra. \textit{A homomorphism of Hopf algebras} is an $A$-bilinear homomorphism $B\to f_{*}(B')$ compatible with comultiplications, counits and antipodes.
\end{definition}

\begin{nothing} \label{Hopf alg dual}
	Let $(\mathscr{T},A)$ be a ringed topos, $M$ and $N$ two $A$-bimodules. We denote by $\FHom_{ll}(M,N)$ (resp. $\FHom_{lr}(M,N)$) the sheaf of $A$-linear homomorphisms from the left $A$-module $M$ to the left (resp. right) $A$-module $N$. We put $M^{\vee}=\FHom_{ll}(M,A)$. The actions of $A$ on $M$ induces a natural $A$-bimodule structure on $M^{\vee}$. 
	There exists a canonical $A$-bilinear morphism
	\begin{equation} \label{dual tensor product dual}
		M^{\vee}\otimes_{A}N^{\vee}\to (M\otimes_{A}N)^{\vee}
	\end{equation}
	which sends $\varphi\otimes \psi$ to $\theta$ defined by $\theta(m\otimes n)=\varphi(m\psi(n))$ for all local sections $m$ of $M$ and $n$ of $N$.

	Let $B$ be a Hopf $A$-algebra. By \eqref{dual tensor product dual}, we obtain a morphism
	\begin{equation}
		B^{\vee}\times B^{\vee} \to (B\otimes_{A}B)^{\vee} \xrightarrow{\delta^{\vee}} B^{\vee}.
	\end{equation}
	Letting $\pi:B\to A$ be the unit element, the above morphism induces a non-commutative ring structure on $B^{\vee}$. 
	The homomorphism $\pi:B\to A$ induces a ring homomorphism $i:A=A^{\vee}\to B^{\vee}$. In this way, we regard $B^{\vee}$ as a non commutative $A$-algebra. 

	A homomorphism of Hopf $A$-algebras $\varphi:B\to C$ induces a homomorphism of $A$-algebras $\varphi^{\vee}:C^{\vee}\to B^{\vee}$. 
\end{nothing}
\begin{nothing} \label{homo Hopf alg dual}
	Let $f:(\mathscr{T}',A')\to (\mathscr{T},A)$ be a morphism of ringed topoi, $(B,\delta_{B},\sigma_{B},\pi_{B})$ a Hopf $A$-algebra. Suppose that the left and the right $A$-algebra structures on $B$ are equal. Then $(f^{*}(B),f^{*}(\delta_{B}),f^{*}(\pi_{B})$ $,f^{*}(\sigma_{B}))$ form a Hopf $A'$-algebra. 
	
	Let $(B',\delta_{B'},\sigma_{B'},\pi_{B'})$ be a Hopf $A'$-algebra and $u:B\to f_{*}(B')$ a homomorphism of Hopf algebras \eqref{Hopf algebra homomorphism}.
	By adjunction, we obtain homomorphisms of $A'$-algebras
	\begin{equation}
		u^{\sharp}:f^{*}(B)\to B', \qquad \widetilde{u}:f^{*}(B)\otimes_{A'} f^{*}(B)\to B'\otimes_{A'}B'
	\end{equation}
	for the left $A'$-algebra structures on the targets. Then the diagrams 
	\begin{equation} \label{diag commu homo Hopf}
		\xymatrixcolsep{4pc}\xymatrix{
			f^{*}(B)\ar[r]^-{f^{*}(\delta_{B})}\ar[d]_{u^{\sharp}} & f^{*}(B)\otimes_{A'}f^{*}(B) \ar[d]^{\widetilde{u}}\\
			B' \ar[r]^-{\delta_{B'}} & B'\otimes_{A'}B'}\qquad 
		\xymatrixcolsep{4pc}\xymatrix{
			f^{*}(B)\ar[r]^{f^{*}(\pi_{B})} \ar[d]_{u^{\sharp}}& A'\ar@{=}[d]\\
			B' \ar[r]^{\pi_{B'}} & A' }
\end{equation}
are commutative. The restriction of $\widetilde{u}$ on $f^{-1}(B\otimes_{A}B)$ is given by $u^{\sharp}|_{f^{-1}(B)}$. In view of \eqref{Hopf alg dual} and \eqref{diag commu homo Hopf}, the homomorphism $u^{\sharp}$ induces a homomorphism of $A'$-algebras $(u^{\sharp})^{\vee}:(B')^{\vee}\to (f^{*}(B))^{\vee}$.
\end{nothing}
%
%

\begin{nothing} \label{notations XX projections}
	In the following, we review the notion of \textit{formal groupoid}. 
	Let $\XX$ be an adic formal $\SS$-scheme. For any integer $r \ge 1$, let $\XX^{r+1}$ be the fiber product of $(r+1)$-copies of $\XX$ over $\SS$. We consider $\XX$ as an adic formal $(\XX^{r+1})$-scheme by the diagonal immersion $\Delta(r):\XX\to \XX^{r+1}$.
	
	We denote by $\tau:\XX^{2}\to \XX^{2}$ the morphism which exchanges the factors of $\XX^{2}$, by $p_{i}:\XX^{2}\to \XX$ ($i=1,2$) the canonical projections and by $p_{ij}:\XX^{3}\to\XX^{2}$ ($1\le i<j\le 3$) be the projection whose composition with $p_{1}$ (resp. $p_{2}$) is the projection $\XX^{3}\to \XX$ on the $i$-th (resp. $j$-th) factor.

	For any formal $\XX^{2}$-schemes $\YY$ and $\ZZ$, $\YY\times_{\XX}\ZZ$ denotes the product of $\YY\to \XX^{2}\xrightarrow{p_{2}} \XX$ and $\ZZ\to \XX^{2}\xrightarrow{p_{1}}\XX$, and we regard $\YY\times_{\XX}\ZZ$ as a formal $\XX^{2}$-scheme by the projection $\YY\times_{\XX}\ZZ\to \XX^{2}\times_{\XX}\XX^{2}=\XX^{3}\xrightarrow{p_{13}}\XX^{2}$.
\end{nothing}

\begin{definition} \label{def of groupoid str}
	Let $\XX$ be an adic formal $\SS$-scheme. \textit{A formal $\XX$-groupoid over $\SS$} is the data of an adic formal $\XX^{2}$-scheme $\mathfrak{G}$ and three adic $\SS$-morphisms \eqref{notations XX projections}
	\begin{equation}
	\alpha:\mathfrak{G}\times_{\XX}\mathfrak{G}\to \mathfrak{G}, \quad \iota:\XX\to \mathfrak{G}, \quad \eta:\mathfrak{G}\to \mathfrak{G}
	\end{equation}
	satisfying the following conditions.

	(i) $\alpha$ and $\iota$ are $\XX^{2}$-morphisms and the following diagrams are commutative:
	\begin{equation}
		\xymatrix{
			\mathfrak{G}\times_{\XX}\mathfrak{G}\times_{\XX} \mathfrak{G}\ar[r]^-{\id\times \alpha} \ar[d]_{\alpha\times \id} & \mathfrak{G}\times_{\XX}\mathfrak{G} \ar[d]^{\alpha}\\
			\mathfrak{G}\times_{\XX}\mathfrak{G} \ar[r]^-{\alpha} & \mathfrak{G}} \qquad
		\xymatrix{
			\mathfrak{G}\ar[r]^-{\iota\times \id} \ar[d]_{\id \times \iota} \ar[rd]^{\id} & \mathfrak{G}\times_{\XX}\mathfrak{G} \ar[d]^{\alpha}\\
			\mathfrak{G}\times_{\XX}\mathfrak{G} \ar[r]^-{\alpha} & \mathfrak{G}}
	\end{equation}

	(ii) The morphism $\eta$ is compatible with $\tau:\XX^{2}\to \XX^{2}$ and we have $\eta^{2}=\id_{\mathfrak{G}}$, $\eta\circ \iota=\iota$. 
	
	(iii) Let $q_{1}$ (resp. $q_{2}$) be the projection $\mathfrak{G}\to \XX$ induced by $p_{1}$ (resp. $p_{2}$). The following diagrams are commutative:
	\begin{equation}
		\xymatrix{
			\mathfrak{G} \ar[r]^-{q_{1}} \ar[d]_{\id\times \eta} & \XX\ar[d]^{\iota} \\
			\mathfrak{G}\times_{\XX}\mathfrak{G} \ar[r]^-{\alpha} & \mathfrak{G}
		} \qquad
		\xymatrix{
			\mathfrak{G} \ar[r]^-{q_{2}} \ar[d]_{\eta\times \id} & \XX\ar[d]^{\iota} \\
			\mathfrak{G}\times_{\XX}\mathfrak{G} \ar[r]^-{\alpha} & \mathfrak{G}
		}		
	\end{equation}

	(iv) The morphism of underlying topological spaces $|\mathfrak{G}|\to |\XX^{2}|$ factors through $\Delta:|\XX|\to |\XX^{2}|$.
\end{definition}

Let $\varpi:\mathfrak{G}_{\zar}\to \XX_{\zar}$ be the morphism of topoi induced by $\mathfrak{G}_{\zar}\to \XX^{2}_{\zar}$ and its factorization through $\Delta$. Then we have $\varpi_{*}(\mathscr{O}_{\mathfrak{G}})=q_{1*}(\mathscr{O}_{\mathfrak{G}})=q_{2*}(\mathscr{O}_{\mathfrak{G}})$. In this way, we regard $\mathscr{O}_{\mathfrak{G}}$ as an $\mathscr{O}_{\XX}$-bialgebra of $\XX_{\zar}$. Then, the formal $\XX$-groupoid structure on $\mathfrak{G}$ induces a \textit{formal Hopf $\mathscr{O}_{\XX}$-algebra structure} on $\mathscr{O}_{\mathfrak{G}}$, that is for every $n\ge 1$, a Hopf $\mathscr{O}_{\XX_{n}}$-algebra structure on $\mathscr{O}_{\mathfrak{G}_{n}}$ which are compatible. 

\begin{definition} \label{morphism of groupoids}
	Let $\XX, \YY$ be two adic formal $\SS$-schemes, $f:\YY\to \XX$ an $\SS$-morphism, $(\mathfrak{G},\alpha_{\mathfrak{G}},\iota_{\mathfrak{G}},\eta_{\mathfrak{G}})$ a formal $\XX$-groupoid and $(\mathfrak{H},\alpha_{\mathfrak{H}},\iota_{\mathfrak{H}},\eta_{\mathfrak{H}})$ a formal $\YY$-groupoid. 
	\textit{A morphism of formal groupoids above $f$} is an $\XX^{2}$-morphism $\varphi:\mathfrak{H}\to \mathfrak{G}$ compatible with $\alpha$'s, $\iota$'s and $\eta$'s. 
\end{definition}

A morphism of formal groupoids $\varphi:\mathfrak{H}\to \mathfrak{G}$ induces a homomorphism of formal Hopf algebras $\mathscr{O}_{\mathfrak{G}}\to f_{*}(\mathscr{O}_{\mathfrak{H}})$, that is for every $n\ge 1$, $\mathscr{O}_{\mathfrak{G}_{n}}\to f_{*}(\mathscr{O}_{\mathfrak{H}_{n}})$ is a homomorphism of Hopf algebras.

%

\begin{nothing} \label{notations R Q}
	In the remainder of this section, $\XX$ denotes a smooth formal $\SS$-scheme. We put $X=\XX_{1}$. We denote by $\RR_{\XX}(r)$ (resp. $\QQ_{\XX}(r)$) the dilatation $(\XX^{r+1})_{(X/p)}$ (resp. $(\XX^{r+1})_{(X/p)}^{\sharp}$) with respect to the diagonal immersion $X\to \XX^{r+1}$ \eqref{prop univ RQ}. 	
	By \ref{prop univ RQ}(i), the canonical morphisms $(\RR_{\XX}(r))_{1}\to X^{r+1}$ and $\underline{(\QQ_{\XX}(r))_{1}}\to (\QQ_{\XX}(r))_{1} \to X^{r+1}$ factor through the diagonal immersion $X\to X^{r+1}$. 

	To simplify the notation, we put $\RR_{\XX}=\RR_{\XX}(1)$, $\QQ_{\XX}=\QQ_{\XX}(1)$, $\mathcal{R}_{\XX}=\mathscr{O}_{\RR_{\XX}}$ and $\mathcal{Q}_{\XX}=\mathscr{O}_{\QQ_{\XX}}$.
	Following Oyama \cite{Oy}, we will present the formal $\XX$-groupoid structure on $\RR_{\XX}$ and $\QQ_{\XX}$. 

	Our notations are different to those of \cite{Oy}. 
	In (\cite{Oy} 1.2), $\widetilde{\RR}_{\XX}$ (resp. $\widetilde{\QQ}_{\XX}$) denotes the formal $\XX^{2}$-scheme constructed by dilatation and $\RR_{\XX}$ (resp. $\QQ_{\XX}$) denotes its reduction modulo $p$.
\end{nothing}

\begin{prop}[\cite{Oy} 1.2.5 and 1.2.6] \label{prop iso R2 R times R}
	Let $r,r'$ be two integers $\ge 1$. There exists canonical isomorphisms 
	\begin{equation}
		\RR_{\XX}(r)\times_{\XX}\RR_{\XX}(r')\xrightarrow{\sim} \RR_{\XX}(r+r'), \qquad \QQ_{\XX}(r)\times_{\XX}\QQ_{\XX}(r')\xrightarrow{\sim} \QQ_{\XX}(r+r'), \label{iso R2 R times R}
	\end{equation}
	where the projections $\RR_{\XX}(r)\to \XX$ and $\QQ_{\XX}(r)\to \XX$ (resp. $\RR_{\XX}(r')\to \XX$ and $\QQ_{\XX}(r')\to \XX$) are induced by the projection $\XX^{r+1}\to \XX$ on the last factor (resp. $\XX^{r'+1}\to \XX$ on the first factor).
\end{prop}
\begin{proof} By \ref{notations R Q}, we have a commutative diagram
	\begin{displaymath}
		\xymatrix{
			(\RR_{\XX}(r))_{1}\times_{X}(\RR_{\XX}(r'))_{1}\ar[r] \ar[d]& \RR_{\XX}(r)\times_{\XX}\RR_{\XX}(r')\ar[d]\\
			X\ar[r] & \XX^{r+1}\times_{\XX}\XX^{r'+1}}
	\end{displaymath}
	By the universality of $\RR_{\XX}(r+r')$ (\ref{prop univ RQ}), we obtain an $(\XX^{r+r'+1})$-morphism
	\begin{equation}
		\varphi: \RR_{\XX}(r)\times_{\XX}\RR_{\XX}(r')\to \RR_{\XX}(r+r').
	\end{equation}

	On the other hand, by the universal property of $\RR_{\XX}(r)$ and $\RR_{\XX}(r')$, the projection $\XX^{r+r'+1}\to \XX^{r+1}$ on the first ($r+1$)-factors (resp. $\XX^{r+r'+1}\to \XX^{r'+1}$ on the last ($r'+1$)-factors) induces a morphism $\RR_{\XX}(r+r')\to \RR_{\XX}(r)$ (resp. $\RR_{\XX}(r+r')\to \RR_{\XX}(r')$) and hence an $(\XX^{r+r'+1})$-morphism
	\begin{equation}
		\psi:\RR_{\XX}(r+r')\to \RR_{\XX}(r)\times_{\XX}\RR_{\XX}(r').
	\end{equation}
	
	The composition $\RR_{\XX}(r+r')\xrightarrow{\psi} \RR_{\XX}(r)\times_{\XX}\RR_{\XX}(r')\xrightarrow{\varphi} \RR_{\XX}(r+r')\to \XX^{r+r'+1}$ is the canonical morphism $\RR_{\XX}(r+r')\to \XX^{r+r'+1}$. By the universal property of $\RR_{\XX}(r+r')$, we have $\varphi\circ \psi=\id$. Let $q_{1}$ (resp. $q_{2}$) denote the projection on the first (resp. second) factor of $\RR_{\XX}(r)\times_{\XX}\RR_{\XX}(r')$. We consider the commutative diagram
	\begin{equation}
		\xymatrix{
			\RR_{\XX}(r)\times_{\XX}\RR_{\XX}(r') \ar[r]^-{\varphi} \ar[rd]&\RR_{\XX}(r+r') \ar[r]^-{\psi} \ar[rd] \ar[d] & \RR_{\XX}(r)\times_{\XX}\RR_{\XX}(r')\ar[d]^{q_{1}}\\
			& \XX^{r+r'+1}\ar[rd] & \RR_{\XX}(r)\ar[d]\\
			& & \XX^{r+1} }
	\end{equation}
	By the universal property of $\RR_{\XX}(r)$, we see that $q_{1}\circ \psi\circ \varphi$ is equal to $q_{1}$. Similarly, we verify that $q_{2}\circ \psi\circ \varphi$ is equal to $q_{2}$. Hence, we have $\psi\circ \varphi=\id$.

	Since $X$ is reduced, we have $\underline{X}=X$. We have a canonical morphism $\underline{(\QQ_{\XX}(r))_{1}\times_{X}(\QQ_{\XX}(r'))_{1}}\to \underline{(\QQ_{\XX}(r))_{1}}\times_{X}\underline{(\QQ_{\XX}(r'))_{1}}$ \eqref{underline product} and the following commutative diagram \eqref{notations R Q}
	\begin{displaymath}
		\xymatrix{
			\underline{(\QQ_{\XX}(r))_{1}\times_{X}(\QQ_{\XX}(r'))_{1}}\ar[r]& \underline{(\QQ_{\XX}(r))_{1}}\times_{X}\underline{(\QQ_{\XX}(r'))_{1}} \ar[r] \ar[d] & \QQ_{\XX}(r)\times_{\XX}\QQ_{\XX}(r')\ar[d]\\
			& X \ar[r]& \XX^{r+1}\times_{\XX}\XX^{r'+1}}
	\end{displaymath}
	By the universal property of $\QQ_{\XX}(r+r')$ \eqref{prop univ RQ}, we obtain an $(\XX^{r+r'+1})$-morphism
	\begin{equation}
		\QQ_{\XX}(r)\times_{\XX}\QQ_{\XX}(r')\to \QQ_{\XX}(r+r'). \label{iso Q2 Q times Q}
	\end{equation}
	By repeating the proof for $\RR_{\XX}(r+r')$, we verify that the above morphism is an isomorphism.
\end{proof}
\begin{prop} \label{prop R Q Hopf alg}
	The formal $\XX^{2}$-scheme $\RR_{\XX}$ (resp. $\QQ_{\XX}$) has a natural formal $\XX$-groupoid structure. 
\end{prop}

\begin{proof} We follow the proof of (\cite{Oy} 1.2.7) where the author proves the analogous results for the $X^{2}$-schemes $\RR_{\XX,1}$ and $\QQ_{\XX,1}$. By \ref{notations R Q}, the morphism of the underlying topological spaces $|\RR_{\XX}|\to |\XX^{2}|$ (resp. $|\QQ_{\XX}|\to |\XX^{2}|$) factors through the diagonal immersion $|\XX|\to |\XX^{2}|$. We consider the following commutative diagram \eqref{notations R Q}:
\begin{equation} \label{diag univ alpha R}
		\xymatrix{
		(\RR_{\XX}(2))_{1} \ar@{^{(}->}[r] \ar[dd]& \RR_{\XX}(2)\ar[d] \\
		& \XX^{3}\ar[d]^{p_{13}}\\
		X\ar[r] & \XX^{2}}
	\end{equation}
	By the universal property of $\RR_{\XX}$ \eqref{prop univ RQ}, we deduce an adic $\XX^{2}$-morphism
	\begin{equation} \label{alpha R}
		\alpha_{R}:\RR_{\XX}(2)\to \RR_{\XX}.
	\end{equation}
	We identify $\RR_{\XX}(2)$ (resp. $\RR_{\XX}(3)$) and $\RR_{\XX}\times_{\XX}\RR_{\XX}$ (resp. $\RR_{\XX}\times_{\XX}\RR_{\XX}\times_{\XX}\RR_{\XX}$) by \ref{prop iso R2 R times R}. The diagrams
	\begin{displaymath}
		\xymatrix{
			\RR_{\XX}\times_{\XX}\RR_{\XX}\times_{\XX}\RR_{\XX}  \ar[r]^-{\id\times \alpha_{R}} \ar[d] & \RR_{\XX}\times_{\XX}\RR_{\XX} \ar[r]^-{\alpha_{R}} \ar[d] & \RR_{\XX} \ar[d] \\
			\XX^{4} \ar[r]^{p_{124}} & \XX^{3} \ar[r]^{p_{13}}& \XX^{2}} \quad
		\xymatrix{
			\RR_{\XX}\times_{\XX}\RR_{\XX}\times_{\XX}\RR_{\XX}  \ar[r]^-{\alpha_{R}\times \id} \ar[d] & \RR_{\XX}\times_{\XX}\RR_{\XX} \ar[r]^-{\alpha_{R}} \ar[d] & \RR_{\XX} \ar[d] \\
			\XX^{4} \ar[r]^{p_{134}} & \XX^{3} \ar[r]^{p_{13}}& \XX^{2}}
	\end{displaymath}
	are commutative and the compositions of the lower horizontal arrows coincide. By the universal property of $\RR_{\XX}$, we deduce that $\alpha_{R}\circ(\id\times \alpha_{R})=\alpha_{R}\circ(\alpha_{R}\times \id)$.

	For any integer $r\ge 1$, we consider the following commutative diagram:
	\begin{equation}
		\xymatrix{
			& \XX \ar[d]^{\Delta(r)} \\
			X \ar@{^{(}->}[ru] \ar[r]& \XX^{r+1}}
	\end{equation}
	By the universal property of $\RR_{\XX}(r)$, we deduce a $\XX^{r+1}$-morphism
	\begin{equation} \label{iota R}
		\iota_{R}(r):\XX\to \RR_{\XX}(r).
	\end{equation}
	In view of \ref{prop univ RQ}, since $\XX\to\XX^{r+1}$ is an immersion, $\iota_{R}(r)$ is a closed immersion. When $r=1$, the diagrams
	\begin{displaymath}
		\xymatrix{
			\RR_{\XX} \ar[r]^-{\id\times \iota_{R}} \ar[d] & \RR_{\XX}\times_{\XX}\RR_{\XX} \ar[r]^-{\alpha_{R}} \ar[d] & \RR_{\XX} \ar[d] \\
			\XX^{2} \ar[r]^{(p_{1},p_{2},p_{2})} & \XX^{3} \ar[r]^{p_{13}}& \XX^{2}}\quad
		\xymatrix{
			\RR_{\XX} \ar[r]^-{\iota_{R}\times \id} \ar[d] & \RR_{\XX}\times_{\XX}\RR_{\XX} \ar[r]^-{\alpha_{R}} \ar[d] & \RR_{\XX} \ar[d] \\
			\XX^{2} \ar[r]^{(p_{1},p_{1},p_{2})} & \XX^{3} \ar[r]^{p_{13}}& \XX^{2}}
\end{displaymath}
	are commutative and the compositions of the lower horizontal arrows are equal to $\id_{\XX^{2}}$. By the universal property of $\RR_{\XX}$, we deduce that $\alpha_{R}\circ(\id\times \iota_{R})=\alpha_{R}\circ(\iota_{R}\times \id)=\id$.

	We consider the following commutative diagram \eqref{notations XX projections}:
	\begin{equation}
		\xymatrix{
		\RR_{\XX,1} \ar@{^{(}->}[r] \ar[dd]& \RR_{\XX}\ar[d] \\
		& \XX^{2}\ar[d]^{\tau}\\
		X\ar[r] & \XX^{2}}
	\end{equation}
	By the universal property of $\RR_{\XX}$, we deduce an adic morphism
	\begin{equation}
		\eta_{R}: \RR_{\XX}\to \RR_{\XX}. \label{eta R}
	\end{equation}
	Since $\tau\circ \Delta=\Delta$, we deduce that $\eta_{R}\circ \iota_{R}=\iota_{R}$ by the universal property of $\RR_{\XX}$. By construction, $\eta_{R}$ satisfies the condition \ref{def of groupoid str}(ii). The diagrams
	\begin{equation}
		\xymatrix{
			\RR_{\XX} \ar[r]^-{\id\times \eta_{R}} \ar[d] & \RR_{\XX}\times_{\XX}\RR_{\XX} \ar[r]^-{\alpha_{R}} \ar[d] & \RR_{\XX} \ar[d] \\
			\XX^{2} \ar[r]^{(p_{1},p_{2},p_{1})}  & \XX^{3} \ar[r]^{p_{13}}& \XX^{2}
		}\qquad
		\xymatrix{
			\RR_{\XX} \ar[r]^-{q_{1}} \ar[d] & \XX \ar[r]^-{\iota_{R}} \ar[d] & \RR_{\XX} \ar[d] \\
			\XX^{2} \ar[r]^{p_{1}}  & \XX \ar[r]^{\Delta}& \XX^{2}
		}
	\end{equation}
	are commutative and the compositions of the lower horizontal arrows coincide. By the universal property of $\RR_{\XX}$, we deduce that $\alpha_{R}\circ(\id\times \eta_{R})=\iota_{R}\circ q_{1}$. We prove the equality $\alpha_{R}\circ(\eta_{R}\times \id)=\iota_{R}\circ q_{2}$ in the same way.

	The proposition for $\mathcal{Q}_{\XX}$ can be verified in exactly the same way using \ref{notations R Q} and \ref{prop iso R2 R times R}.
\end{proof}
\begin{nothing} \label{local description of R}
		We put $\YY=\Spf(\rW\{T_{1},\cdots,T_{d}\})$ and we present local descriptions for $\mathcal{R}_{\YY}$ and $\mathcal{Q}_{\YY}$ \eqref{notations R Q}. Put $\xi_{i}=1\otimes T_{i}-T_{i}\otimes 1\in \mathscr{O}_{\YY^{2}}$. The ideal $\mathscr{A}$, associated to the diagonal closed immersion $\YY_{1}\to \YY^{2}$, is generated by $p,\xi_{1},\xi_{2},\cdots,\xi_{d}$. 
		The algebra $W[1\otimes T_1,\cdots,1\otimes T_d, T_1\otimes 1,\cdots,T_d\otimes 1,x_1,\cdots,x_d]/(\xi_i-px_i)_{1\le i\le d}$ is free over $\rW$. 
		Hence, we have an isomorphism of $\mathscr{O}_{\YY^{2}}$-algebras (\ref{hat A0 blowup}, \ref{notations R Q})
	\begin{equation} \label{RX X2 alg}
		\mathscr{O}_{\YY^{2}}\Big\{\frac{\xi_{1}}{p},\cdots,\frac{\xi_{d}}{p}\Big\}=\frac{\mathscr{O}_{\YY^{2}}\{x_{1},\cdots,x_{d}\}}{(\xi_{i}-px_{i})_{1\le i\le d}}\xrightarrow{\sim} \mathcal{R}_{\YY}.
	\end{equation}

	By \ref{dilatation def}, the ideal $\mathscr{A}^{\sharp}$ is generated by $p,\xi_{1}^{p},\cdots,\xi_{d}^{p}$.
	The algebra $W[1\otimes T_1,\cdots,1\otimes T_d, T_1\otimes 1,\cdots,T_d\otimes 1,x_1,\cdots,x_d]/(\xi_i^p-px_i)_{1\le i\le d}$ is free over $\rW$. 
	Hence, we have an isomorphism of $\mathscr{O}_{\YY^{2}}$-algebras \eqref{hat A0 blowup}
	\begin{equation}\label{QX X2 alg}
		\mathscr{O}_{\YY^{2}}\Big\{\frac{\xi_{1}^{p}}{p},\cdots,\frac{\xi_{d}^{p}}{p}\Big\}=\frac{\mathscr{O}_{\YY^{2}}\{x_{1},\cdots,x_{d}\}}{(\xi_{i}^{p}-px_{i})_{1\le i\le d}}\xrightarrow{\sim} \mathcal{Q}_{\YY}.
	\end{equation}
\end{nothing}

\begin{lemma}\label{local description R Ad}
	Let $d$ be the relative dimension of $\XX$ over $\SS$ and $\widehat{\mathbb{A}}_{\XX}^{d}$ the $d$-dimensional affine space. Assume that there exists an \'etale $\SS$-morphism $f:\XX\to \YY=\Spf(\rW\{T_{1},\cdots,T_{d}\})$. Considering $\RR_{\XX}$ as a formal $\XX$-scheme via the morphism $q_{1}:\RR_{\XX}\to \XX$ (resp. $q_{2}:\RR_{\XX}\to \XX$), then $f$ induces an isomorphism over $\XX$:
		\begin{equation} \label{iso RX d affine}
			\lambda:\RR_{\XX}\xrightarrow{\sim} \widehat{\mathbb{A}}_{\XX}^{d},
		\end{equation}
		such that $\lambda\circ \iota_{R}:\XX\to \widehat{\mathbb{A}}_{\XX}^{d}$ is the closed immersion associated to the zero section of $\widehat{\mathbb{A}}_{\XX}^{d}$.
\end{lemma}
\begin{proof} We follow the proof of (\cite{Oy} 1.1.8) and we first prove the assertions for $\YY$. Observe that $q_{1}$ and $q_{2}$ are affine. For any $1\le i\le d$, we have $1\otimes T_{i}=p(\frac{\xi_{i}}{p})+T_{i}\otimes 1$ in $\mathcal{R}_{\YY}$. By \eqref{RX X2 alg}, we deduce the isomorphisms
	\begin{equation} 
		\mathscr{O}_{\YY}\{x_{1},\cdots,x_{d}\}\xrightarrow{\sim} q_{1*}(\mathcal{R}_{\YY}),\qquad \mathscr{O}_{\YY}\{x_{1},\cdots,x_{d}\}\xrightarrow{\sim} q_{2*}(\mathcal{R}_{\YY})
	\end{equation}
	where $x_{i}$ is sent to $\frac{\xi_{i}}{p}$ in both cases. The isomorphisms \eqref{iso RX d affine} for $\YY$ follows. We put $Y=\YY_{1}$ and we consider the following commutative diagram
	\begin{equation} \label{diag etale f}
	\xymatrix{
		& \XX^{2}\ar[d]^{\id \times f}\\
		X\ar[ru]^{\Delta} \ar[r] \ar[d]_{f_{1}} & \XX\times_{\SS}\YY\ar[d]^{f\times \id}\\
		Y\ar[r]^{\Delta}& \YY^{2}
	}
\end{equation}
where the square is Cartesian. By \ref{prop dilatation etale} and \ref{prop dilatation flat}, we deduce an isomorphism
\begin{equation}
	\RR_{\XX}\xrightarrow{\sim} \RR_{\YY}\times_{\YY^{2}}(\XX\times_{\SS}\YY)=\RR_{\YY}\times_{q_{1},\YY}\XX.
\end{equation}
Considering $\RR_{\XX}$ as a formal $\XX$-scheme via $q_{1}$, the isomorphism \eqref{iso RX d affine} follows from that of $\RR_{\YY}$. The another isomorphism can be verified in the same way. In view of the construction of $\iota_{R}$ \eqref{iota R}, the composition $\lambda\circ \iota_{R}$ corresponds to the zero section.
\end{proof}
\begin{coro}\label{local description R lemma}
		Keep assumptions of \ref{local description R Ad} and consider $\frac{\xi_{i}}{p}$'s as sections of $\mathcal{R}_{\XX}$. We have the following isomorphisms of $\mathscr{O}_{\XX}$-algebras
	\begin{equation} \label{local p adic poly R}
		\mathscr{O}_{\XX}\{x_{1},\cdots,x_{d}\}\xrightarrow{\sim} q_{1*}(\mathcal{R}_{\XX}),\qquad \mathscr{O}_{\XX}\{x_{1},\cdots,x_{d}\}\xrightarrow{\sim} q_{2*}(\mathcal{R}_{\XX})
	\end{equation}
	where $x_{i}$ is sent to $\frac{\xi_{i}}{p}$ in both cases.
\end{coro}

\begin{nothing}\label{local description Q lemma}
		We put $\YY=\Spf(\rW\{T_{1},\cdots,T_{d}\})$. By \eqref{QX X2 alg}, we have following isomorphisms
	\begin{equation} 
		\frac{\mathscr{O}_{\YY}\{x_{1},\cdots,x_{d},y_{1},\cdots,y_{d}\}}{(y_{i}^{p}-px_{i})_{1\le i\le d}}\xrightarrow{\sim} q_{1*}(\mathcal{Q}_{\YY}),\qquad \frac{\mathscr{O}_{\YY}\{x_{1},\cdots,x_{d},y_{1},\cdots,y_{d}\}}{(y_{i}^{p}-px_{i})_{1\le i\le d}}\xrightarrow{\sim} q_{2*}(\mathcal{Q}_{\YY})
	\end{equation}
	where $x_{i}$ is sent to $\frac{\xi_{i}^{p}}{p}$ and $y_{i}$ is sent to $\xi_{i}$ in both cases. 

	Assume that there exists an \'etale $\SS$-morphism $f:\XX\to \YY=\Spf(\rW\{T_{1},\cdots,T_{d}\})$. We consider the $\xi_{i}$'s and $\frac{\xi_{i}^{p}}{p}$'s as sections of $\mathcal{Q}_{\XX}$. By \ref{prop dilatation etale}, \ref{prop dilatation flat} and \eqref{diag etale f}, we deduce the following isomorphisms
	\begin{equation} 
		\frac{\mathscr{O}_{\XX}\{x_{1},\cdots,x_{d},y_{1},\cdots,y_{d}\}}{(y_{i}^{p}-px_{i})_{1\le i\le d}}\xrightarrow{\sim} q_{1*}(\mathcal{Q}_{\XX}),\qquad \frac{\mathscr{O}_{\XX}\{x_{1},\cdots,x_{d},y_{1},\cdots,y_{d}\}}{(y_{i}^{p}-px_{i})_{1\le i\le d}}\xrightarrow{\sim} q_{2*}(\mathcal{Q}_{\XX})
	\end{equation}
	where $x_{i}$ is sent to $\frac{\xi_{i}^{p}}{p}$ and $y_{i}$ is sent to $\xi_{i}$ in both cases. 
\end{nothing}

\begin{nothing} \label{local description Hopf RQ}
	Let $n$ be an integer $\ge 1$. We describe the Hopf algebra structure of $\mathcal{R}_{\XX,n}$ and $\mathcal{Q}_{\XX,n}$ in terms of a system of local coordinates. Keep the assumption and notation of \ref{local description R Ad}. 
	The homomorphism $\mathscr{O}_{\XX^{2}_n}\to \mathscr{O}_{\XX^{2}_{n}}\otimes_{\mathscr{O}_{\XX_{n}}}\mathscr{O}_{\XX^{2}_{n}}$ induced by $p_{13}:\XX^{3}\to \XX^{2}$ sends $\xi_{i}$ to $1\otimes \xi_{i}+\xi_{i}\otimes 1$. The homomorphism $\mathscr{O}_{\XX^{2}_{n}}\to \mathscr{O}_{\XX^{2}_{n}}$ induced by $\tau:\XX^{2}\to \XX^{2}$, sends $\xi_{i}$ to $-\xi_{i}$. In view of the proof of \ref{prop R Q Hopf alg}, we have following descriptions:
	\begin{eqnarray}
		&&\left\{ 
			\begin{array}{ll}
				\delta: \mathcal{R}_{\XX,n}\to \mathcal{R}_{\XX,n}\otimes_{\mathscr{O}_{\XX_{n}}}\mathcal{R}_{\XX,n} & \frac{\xi_{i}}{p}\mapsto 1\otimes\frac{\xi_{i}}{p}+\frac{\xi_{i}}{p}\otimes1 \\
				\sigma:\mathcal{R}_{\XX,n}\to \mathcal{R}_{\XX,n} & \frac{\xi_{i}}{p} \mapsto -\frac{\xi_{i}}{p}  \\
				\pi:\mathcal{R}_{\XX,n}\to \mathscr{O}_{\XX_{n}} & \frac{\xi_{i}}{p}\mapsto 0
			\end{array}
			\right. \\
		&&\left\{ 
			\begin{array}{ll}
				\delta: \mathcal{Q}_{\XX,n}\to \mathcal{Q}_{\XX,n}\otimes_{\mathscr{O}_{\XX_{n}}}\mathcal{Q}_{\XX,n} & \xi_{i}\mapsto 1\otimes \xi_{i} +\xi_{i}\otimes 1 \\
				& \frac{\xi_{i}^{p}}{p}\mapsto 1\otimes \frac{\xi_{i}^{p}}{p} +\sum_{j=1}^{p-1}\frac{(p-1)!}{j!(p-j)!}\xi_{i}^{j}\otimes \xi_{i}^{p-j} +\frac{\xi_{i}^{p}}{p}\otimes 1 \\
				\sigma:\mathcal{Q}_{\XX,n}\to \mathcal{Q}_{\XX,n} & \xi_{i}\mapsto -\xi_{i}, \qquad \frac{\xi_{i}^{p}}{p} \mapsto \frac{(-\xi_{i})^{p}}{p}  \\
				\pi:\mathcal{Q}_{\XX,n}\to \mathscr{O}_{\XX_{n}} & \xi_{i}\mapsto 0, \qquad \frac{\xi_{i}^{p}}{p}\mapsto 0
			\end{array} \right.
	\end{eqnarray}
\end{nothing}

\section{Connections and stratifications} \label{conn and stra}
\begin{nothing}\label{notation MIC}
	Let $S$ be a scheme, $f:X\to S$ a smooth morphism, $M$ an $\mathscr{O}_{X}$-module and $\lambda\in \Gamma(S,\mathscr{O}_{S})$. 
	We say (abusively) that a morphism of $\mathscr{O}_{X}$-modules $u:M\to N$ is \textit{$\mathscr{O}_{S}$-linear} if it is $f^{-1}(\mathscr{O}_{S})$-linear.
	\textit{A $\lambda$-connection on $M$} relative to $S$ is an $\mathscr{O}_{S}$-linear morphism 
	\begin{equation}
		\nabla:M\to M\otimes_{\mathscr{O}_{X}}\Omega_{X/S}^{1}
	\end{equation}
	such that for every local sections $f$ of $\mathscr{O}_{X}$ and $e$ of $M$, we have $\nabla(fe)=\lambda e\otimes d(f)+ f\nabla(e)$. 
	We will simply call $\nabla$ a $\lambda$-connection on $M$ when there is no risk of confusion. For any $q\ge 0$, the morphism $\nabla$ extends to a unique $\mathscr{O}_{S}$-linear morphism 
\begin{equation}
	\nabla_{q}:M\otimes_{\mathscr{O}_{X}}\Omega_{X/S}^{q} \to M\otimes_{\mathscr{O}_{X}}\Omega_{X/S}^{q+1}
\end{equation}
such that for every local sections $\omega$ of $\Omega_{X/S}^{q}$ and $e$ of $M$, we have $\nabla_{q}(e\otimes\omega)=\lambda e\otimes d(\omega) + \nabla(e)\wedge\omega$. 
The composition $\nabla_{1}\circ \nabla$ is $\mathscr{O}_{X}$-linear. We say that $\nabla$ is \textit{integrable} if $\nabla_{1}\circ \nabla =0$.


Let $(M,\nabla)$ and $(M',\nabla')$ be two $\mathscr{O}_{X}$-modules with $\lambda$-connection. A morphism from $(M,\nabla)$ to $(M',\nabla')$ is an $\mathscr{O}_{X}$-linear morphism $u:M\to M'$ such that $(\id\otimes u)\circ \nabla=\nabla'\circ u$.

Classically, $1$-connections are called \textit{connections} and integrable $0$-connections are called \textit{Higgs fields}. \textit{A Higgs module} is an $\mathscr{O}_{X}$-module equipped with a Higgs field. 

We denote by $\MIC(X/S)$ (resp. $\lMIC(X/S)$, resp. $\HIGG(X/S)$) the category of $\mathscr{O}_{X}$-modules with integrable connection (resp. $\lambda$-connection, resp. Higgs field) relative to $S$.

Let $(M,\nabla)$ be an object of $\lMIC(X/S)$. We deduce that $\nabla_{q+1}\circ \nabla_{q}=0$ for all $q\ge 0$. 
Then we can associate to $(M,\nabla)$ a \textit{$\lambda$-de Rham complex}:
\begin{equation} \label{lambla dR complex}
	M\xrightarrow{\nabla} M\otimes_{\mathscr{O}_{X}}\Omega_{X/S}^{1} \xrightarrow{\nabla_{1}} M\otimes_{\mathscr{O}_{X}}\Omega_{X/S}^{2} \xrightarrow{\nabla_{2}}\cdots.
\end{equation}
Classically, $0$-de Rham complexes are called \textit{Dolbeault complexes} in (\cite{Sim92} p.~24, \cite{AGT} I.2.3) or \textit{Higgs complexes} in (\cite{OV07} p.~2).
\end{nothing}

\begin{nothing}\label{quasi-nilpotent coords}
	Let $d,n$ be integers $\ge 1$, $S$ a $(\mathbb{Z}/p^{n}\mathbb{Z})$-scheme, $f:X\to\mathbb{A}_{S}^{d}=\Spec(\mathscr{O}_{S}[T_{1},\cdots,T_{d}])$ an \'{e}tale $S$-morphism. For any $1\le i\le d$, we denote by $t_{i}$ the image of $T_{i}$ in $\mathscr{O}_{X}$. Let $m$ be an integer $\ge 0$ and $(M,\nabla)$ an $\mathscr{O}_{X}$-module with integrable $p^{m}$-connection relative to $S$. There are $\mathscr{O}_{S}$-linear endomorphisms $\nabla_{\partial_{1}},\cdots,\nabla_{\partial_{d}}$ of $M$ such that for every local section $e$ of $M$, we have
	\begin{equation}
		\nabla(e)=\sum_{i=1}^{d} \nabla_{\partial_{i}}(e)\otimes dt_{i}.
	\end{equation}
	Since $\nabla$ is integrable, we have $\nabla_{\partial_{i}}\circ \nabla_{\partial_{j}}=\nabla_{\partial_{j}}\circ \nabla_{\partial_{i}}$ for all $1\le i,j\le d$. Therefore, for every multi-index $I=(i_{1},\cdots,i_{d})\in \mathbb{N}^{d}$, the endomorphism $\nabla_{\partial^{I}}=\prod_{j=1}^{d} (\nabla_{\partial_{j}})^{i_{j}}$ is well-defined.

	Following (\cite{BO} 4.10, \cite{Shiho} Definition 1.5), we say that $(M,\nabla)$ \textit{is quasi-nilpotent with respect to $f$} if, for any open subscheme $U$ of $X$ and any section $e\in M(U)$, there exists a Zariski covering $\{U_{j}\to U\}_{j\in J}$ and a family of integers $\{N_{j}\}_{j\in J}$ such that $\nabla_{\partial^{I}}(e|_{U_{j}})=0$ for all $j\in J$ and $I\in \mathbb{N}^{d}$ with $|I|\ge N_{j}$. 
	
	If $f':X\to \mathbb{A}_{S}^{d}$ is another \'{e}tale $S$-morphism, $(M,\nabla)$ is quasi-nilpotent with respect to $f$ if and only if it is quasi-nilpotent with respect to $f'$ (\cite{BO} 4.13, \cite{Shiho} Lemma 1.6). Note that this result requires that $p^{n}\mathscr{O}_{S}=0$ for some $n>0$.
\end{nothing}

\begin{definition}[\cite{BO} 4.13; \cite{Shiho} Definition 1.8] \label{def quasinilpotent}
	Let $n$ be an integer $\ge 1$, $S$ a $(\mathbb{Z}/p^{n}\mathbb{Z})$-scheme, $X$ a smooth $S$-scheme and $(M,\nabla)$ an $\mathscr{O}_{X}$-module with integrable $p^{m}$-connection relative to $S$. We say that $(M,\nabla)$ is \textit{quasi-nilpotent} if for any point $x$ of $X$, there exists a Zariski neighborhood $U$ of $x$ in $X$ and an \'{e}tale $S$-morphism $f:U\to \mathbb{A}_{S}^{d}$ such that $(M,\nabla)|_{U}$ is quasi-nilpotent with respect to $f$ \eqref{quasi-nilpotent coords}. 
\end{definition}

We denote by $\MIC^{\qn}(X/S)$ (resp. $\lambda\textnormal{-}\MIC^{\qn}(X/S)$) the full subcategory of $\MIC(X/S)$ (resp. $\lMIC(X/S)$) consisting of the quasi-nilpotent objects.



\begin{definition} \label{def stratification}
	Let $(\mathscr{T},A)$ be a ringed topos, $(B,\delta,\pi,\sigma)$ a Hopf $A$-algebra \eqref{def of Hopf alg} and $M$ an $A$-module. \textit{A $B$-stratification on $M$} is a $B$-linear isomorphism $\varepsilon:B\otimes_{A}M\xrightarrow{\sim}M\otimes_{A}B$ \eqref{tensor product OX bimodule} such that:
	
	(i) $\pi^{*}(\varepsilon)=\id_{M}$.

	(ii) (cocycle condition) The following diagram is commutative:
	\begin{equation}
		\xymatrix{
			B\otimes_{A}B\otimes_{A}M\ar[rd]^{\id_{B}\otimes \varepsilon} \ar[rr]^{\delta^{*}(\varepsilon)}& & M\otimes_{A}B\otimes_{A}B\\
			& B\otimes_{A} M\otimes_{A} B \ar[ru]^{\varepsilon\otimes \id_{B}}&}
	\end{equation}
\end{definition}

Given two $A$-modules with $B$-stratification $(M_{1},\varepsilon_{1})$ and $(M_{2},\varepsilon_{2})$, a morphism from $(M_{1},\varepsilon_{1})$ to $(M_{2},\varepsilon_{2})$ is an $A$-linear morphism $f:M_{1}\to M_{2}$ compatible with $\varepsilon_{1}$ and $\varepsilon_{2}$. The $A$-module $M_{1}\otimes_{A}M_{2}$ has a canonical $B$-stratification
\begin{equation}
	B\otimes_{A}M_{1}\otimes_{A}M_{2}\xrightarrow{\varepsilon_{1}\otimes \id_{M_{2}}} M_{1}\otimes_{A}B\otimes_{A}M_{2}\xrightarrow{\id_{M_{1}}\otimes \varepsilon_{2}}M_{1}\otimes_{A}M_{2}\otimes_{A}B. \label{stratification tensor product}
\end{equation}
The above stratification on the tensor product makes the category of $A$-modules with $B$-stratification into a tensor category.

Let $n$ be an integer $\ge 1$, $\XX$ an adic formal $\SS$-scheme, $M$ an $\mathscr{O}_{\XX_{n}}$-module and $\mathfrak{G}$ a formal $\XX$-groupoid over $\SS$ \eqref{def of groupoid str}. We call abusively \textit{$\mathscr{O}_{\mathfrak{G}}$-stratification on $M$} instead of $\mathscr{O}_{\mathfrak{G}_{n}}$-stratification on $M$. 

We have a simpler description of a stratification.

\begin{lemma}[\cite{Oy} 1.2.4] \label{lemma stratification}
	Let $(\mathscr{T},A)$ be a ringed topos, $B$ a Hopf $A$-algebra and $M$ an $A$-module. A $B$-stratification on $M$ is equivalent to an $A$-linear morphism $\theta:M\to M\otimes_{A}B$ for the right $A$-action on the target satisfying the following conditions:
	
	\textnormal{(i)} The composition $(\id_{M}\otimes \pi)\circ \theta:M\to M\otimes_{A}B \to M$ is the identity morphism.

	\textnormal{(ii)} The following diagram is commutative
	\begin{equation}
		\xymatrix{
			M\ar[r]^{\theta} \ar[d]_-{\theta} & M\otimes_{A}B \ar[d]^{\theta\otimes \id_{B}}\\
			M\otimes_{A}B \ar[r]^-{\id_{M}\otimes \delta} & M\otimes_{A}B\otimes_{A}B }
	\end{equation}
\end{lemma}


	Let $\theta:M\to M\otimes_{A} B$ be an $A$-linear morphism satisfying the conditions of \ref{lemma stratification} and 
	\begin{displaymath}
		\alpha(\theta):B^{\vee}\otimes_{A}M\to M,\quad \varphi\otimes m=(\id \otimes \varphi)\bigl(\theta(m)\bigr)
	\end{displaymath}
	the associated $A$-linear morphism. In view of conditions (i-ii) of \ref{lemma stratification}, the morphism $\alpha(\theta)$ makes $M$ into a left $B^{\vee}$-module (cf. \cite{Oy} Proof of 1.2.9 page 18 for details).

\begin{nothing}\label{morphism Hopf functor}
	Let $f:(\mathscr{T}',A')\to (\mathscr{T},A)$ be a morphism of ringed topoi, $B$ a Hopf $A$-algebra and $B'$ a Hopf $A'$-algebra. A homomorphism of Hopf algebras $B\to f_{*}(B')$ \eqref{Hopf algebra homomorphism} induces a functor \eqref{def stratification} (cf. \cite{Ber} II 1.2.5):
	\begin{eqnarray}
		\Big\{\begin{array}{c}
	\textnormal{$A$-modules}\\
	\textnormal{with $B$-stratification}
		\end{array}\Big\}
	&\to& 
	\Big\{\begin{array}{c}
	\textnormal{$A'$-modules}\\
	\textnormal{with $B'$-stratification}
		\end{array}\Big\}\\
		(M,\varepsilon)\qquad &\mapsto& (f^{*}(M),f^{-1}(\varepsilon)\otimes_{f^{-1}(B)}B') \nonumber.
\end{eqnarray}
\end{nothing}

\begin{nothing} \label{PD envelop P}
	In the remainder of this section, $\XX$ denotes a smooth formal $\SS$-scheme. For any integer $n\ge 1$, we equip $(\rW_{n},p\rW_{n})$ with the canonical PD-structure $\gamma_{n}$. 
	We briefly review the formal groupoid structure on the PD-envelope of the diagonal immersion following \cite{Ber}. 

	Let $r,n$ be integers $\ge 1$. We denote by $\XX_{n}^{r+1}$ the product of $(r+1)$-copies of $\XX_{n}$ over $\SS_{n}$ \eqref{conventions adic} and by $\PP_{\XX_n}(r)$ the PD-envelope of the diagonal immersion $\XX_{n}\to \XX_{n}^{r+1}$ compatible with the PD-structure $\gamma_{n}$ (\cite{Ber} I 4.3.1). By extension of scalars (\cite{BO} 3.20.8, \cite{Ber} I 2.8.2), we have a canonical PD-isomorphism $\PP_{\XX_n}(r)\times_{\SS_{n}}\SS_{m}\xrightarrow{\sim} \PP_{\XX_m}(r)$ for all integers $1\le m< n$. The inductive limit $\PP_{\XX}(r)$ of the inductive system $(\PP_{\XX_n}(r))_{n\ge 1}$ is an adic affine formal $(\XX^{r+1})$-scheme (\cite{Ab10} 2.3.10). We drop $(r)$ from the notation when $r=1$. 
\end{nothing}

\begin{nothing} \label{notation pd sym algebra}
	For a commutative ring $A$, we denote by $A\langle x_{1},\cdots,x_{d} \rangle$ the PD polynomial ring in $d$ variables (\cite{Ber} I 1.5). If $A$ is an adic ring such that $pA$ is an ideal of the definition, we denote by $A\llangle x_{1},\cdots,x_{d}\rrangle$ the $p$-adic completion of the PD polynomial algebra $A\langle x_{1},\cdots,x_{d}\rangle$. 
\end{nothing}

\begin{nothing} \label{flatness PX}
	Assume that there exists an \'{e}tale $\SS$-morphism $\XX\to \widehat{\mathbb{A}}_{\SS}^{d}=\Spf(\rW\{T_{1},\cdots$ $,T_{d}\})$ and we set $t_{i}$ the image of $T_{i}$ in $\mathscr{O}_{\XX}$ for all $1\le i\le d$. We note $\xi_{i}$ the section $1\otimes t_{i}-t_{i}\otimes 1$ of $\mathscr{O}_{\XX^{2}}$ and also its image in $\mathscr{O}_{\PP_{\XX}}$. 	
	By (\cite{Ber} I 4.4.1 and 4.5.3), we deduce the following PD-isomorphisms \eqref{notation pd sym algebra}
	\begin{equation}\label{local description PX}
		\mathscr{O}_{\XX}\llangle x_{1},\cdots,x_{d}\rrangle \xrightarrow{\sim} q_{1*}(\mathscr{O}_{\PP_{\XX}}), \qquad 	\mathscr{O}_{\XX}\llangle x_{1},\cdots,x_{d}\rrangle \xrightarrow{\sim} q_{2*}(\mathscr{O}_{\PP_{\XX}})
	\end{equation}
	where $q_{1},q_{2}:\PP_{\XX}\to \XX$ are the canonical morphisms and $x_{i}$ is sent to $\xi_{i}$. 
	In general, we deduce that $\PP_{\XX}$ is flat over $\SS$ \eqref{notations S flat}.
	
	For any integers $r,r'\ge 1$, by (\cite{Ber} II 1.3.4 and 1.3.5), we deduce a canonical isomorphism of formal $(\XX^{r+r'+1})$-schemes
	\begin{equation}
		\PP_{\XX}(r)\times_{\XX}\PP_{\XX}(r') \xrightarrow{\sim} \PP_{\XX}(r+r'). \label{P2 to Ptimes P}
	\end{equation}
\end{nothing}

\begin{prop} \label{groupoid PX}
	The formal $\XX^{2}$-scheme $\PP_{\XX}$ has a natural formal $\XX$-groupoid structure. 
\end{prop}

\begin{proof}
	For any $r\ge 1$, the diagonal immersion $\Delta(r):\XX\to \XX^{r+1}$ induces a canonical $(\XX^{r+1})$-morphism $\iota_{P}(r):\XX\to \PP_{\XX}(r)$. 
	Set $X=\XX_{1}$ and let $J$ be the PD-ideal of $\mathscr{O}_{\PP_{X}}$ associated to the closed immersion $X\to \PP_{X}$. For any local section $x$ of $J$, we have $x^{p}=p!x^{[p]}=0$. Hence, we have a closed immersion $\underline{\PP_{X}}\hookrightarrow X$ \eqref{notation underline}. Since $X$ is reduced, the composition $\underline{X}\to \underline{\PP_{X}}\to X$ is an isomorphism. 
	We deduce an isomorphism:
	\begin{equation}
		\underline{\PP_{X}} \xrightarrow{\sim} X, \label{iso X underline PX}
	\end{equation}	
	Hence the morphism of the underlying topological spaces $|\PP_{\XX}|\to |\XX^{2}|$ factors through $\Delta:|\XX|\to |\XX^{2}|$. The canonical morphism 
	$\PP_{\XX}(2)\to \XX^{3}\xrightarrow{p_{13}} \XX^{2}$ is compatible with $\iota_{P}(2)$ and $\Delta$. By the universal property of $(\PP_{\XX_{n}})_{n\ge 1}$, we deduce an $\XX^{2}$-morphism $\alpha_{P}:\PP_{\XX}(2)\to \PP_{\XX}$.  
	Similary, by the universal property of $(\PP_{\XX_{n}})_{n\ge 1}$, the composition $\PP_{\XX}\to \XX^{2}\xrightarrow{\tau}\XX^{2}$ \eqref{notations XX projections} induces a morphism
		$\eta_{P}:\PP_{\XX}\to \PP_{\XX}$. 
	By the universal property of $(\PP_{\XX_{n}})_{n\ge 1}$, we verify that $(\alpha_{P},\iota_{P},\eta_{P})$ defines a formal $\XX$-groupoid structure on $\PP_{\XX}$ (cf. the proof of \ref{prop R Q Hopf alg}).
\end{proof}

	To simplify the notations, we put $\mathcal{P}_{\XX}=\mathscr{O}_{\PP_{\XX}}$ considered as a formal Hopf $\mathscr{O}_{\XX}$-algebra by \ref{groupoid PX}.

\begin{prop}[\cite{BO} 4.12] \label{B-O stratification MIC}
	Let $n$ be an integer $\ge 1$. There is a canonical equivalence of categories between the category of $\mathscr{O}_{\XX_{n}}$-modules with $\mathcal{P}_{\XX}$-stratification and the category $\MIC^{\qn}(\XX_{n}/\SS_{n})$ \eqref{def quasinilpotent}.
\end{prop}

We now explain the relationships among the various groupoids and stratifications we have constructed. 

\begin{prop} \label{prop P to Q}
	Let $\QQ_{\XX}$ be the formal $\XX$-groupoid defined in \ref{prop R Q Hopf alg}. We have a canonical morphism of formal $\XX$-groupoids \eqref{groupoid PX} 
	\begin{equation} \label{pi P to Q} 
		\lambda:\PP_{\XX}\to \QQ_{\XX}.
	\end{equation}
\end{prop}
\begin{proof} 
	The isomorphism $\underline{\PP_{X}}\simeq X$ \eqref{iso X underline PX} fits into a commutative diagram
\begin{equation}
	\xymatrix{
		\underline{\PP_{X}} \ar[r] \ar[d]& \PP_{\XX} \ar[d]\\
		X \ar[r]^{\Delta}& \XX^{2} }
\end{equation}
	By the universal property of $\QQ_{\XX}$ \eqref{prop univ RQ}, we deduce a canonical $\XX^{2}$-morphism $\lambda:\PP_{\XX}\to \QQ_{\XX}$ . We denote by $(\alpha_{P},\iota_{P},\eta_{P})$ (resp. $(\alpha_{Q},\iota_{Q},\eta_{Q})$) the formal groupoid structure on $\PP_{\XX}$ (resp. $\QQ_{\XX}$). The diagrams 
	\begin{equation}
		\xymatrix{
			\PP_{\XX}\times_{\XX}\PP_{\XX} \ar[r]^-{\lambda^{2}} \ar[d]& \QQ_{\XX}\times_{\XX}\QQ_{\XX} \ar[r]^-{\alpha_{Q}} \ar[d] & \QQ_{\XX} \ar[d]\\
			\XX^{3}\ar@{=}[r]& \XX^{3} \ar[r]^{p_{13}} & \XX^{2}
		}\qquad
		\xymatrix{
			\PP_{\XX}\times_{\XX}\PP_{\XX} \ar[r]^-{\alpha_{P}} \ar[d]& \PP_{\XX} \ar[r]^-{\lambda} \ar[d] & \QQ_{\XX} \ar[d]\\
			\XX^{3}\ar[r]^{p_{13}}& \XX^{2} \ar@{=}[r] & \XX^{2} }
	\end{equation}
	are commutative and the compositions of the lower horizontal arrows coincide. By the universal property of $\QQ_{\XX}$, we deduce that $\alpha_{Q}\circ \lambda^{2}=\lambda\circ \alpha_{P}$. The composition $\XX\xrightarrow{\iota_{P}}\PP_{\XX}\to \XX^{2}$ is the diagonal immersion. By the universal property of $\QQ_{\XX}$, we deduce that $\lambda\circ \iota_{P}=\iota_{Q}$. 
	The following diagrams
	\begin{equation}
		\xymatrix{
			\PP_{\XX}\ar[r]^{\eta_{P}} \ar[d]& \PP_{\XX}\ar[r]^{\lambda} \ar[d]& \QQ_{\XX} \ar[d]\\
			\XX^{2}\ar[r]^{\tau} &\XX^{2} \ar@{=}[r] & \XX^{2}} \qquad
		\xymatrix{
			\PP_{\XX}\ar[r]^{\lambda} \ar[d]& \QQ_{\XX}\ar[r]^{\eta_{Q}} \ar[d]& \QQ_{\XX} \ar[d]\\
			\XX^{2}\ar@{=}[r] &\XX^{2}\ar[r]^{\tau} & \XX^{2}}
	\end{equation}
	are commutative and the compositions of the lower arrows coincide. By the universal property of $\QQ_{\XX}$, we deduce that $\lambda\circ \eta_{P}=\eta_{Q}\circ \lambda$. The proposition follows.
\end{proof}
\begin{nothing} \label{basic notation T X}
	We take again the notation of \ref{notations R Q} for $\XX$. For any integers $n,r \ge 1$, we denote by $\TT_{\XX,n}(r)$ the PD-envelope of the closed immersion $\XX_{n}\hookrightarrow (\RR_{\XX}(r))_{n}$ \eqref{iota R} compatible with the PD-structure $\gamma_{n}$. By extension of scalars (\cite{Ber} I 2.8.2), we have a canonical isomorphism of PD-schemes $\TT_{\XX,n}(r)\times_{\SS_{n}}\SS_{m}\xrightarrow{\sim} \TT_{\XX,m}(r)$ for all integers $1\le m< n$. The inductive limit $\TT_{\XX}(r)$ of the inductive system $(\TT_{\XX,n}(r))_{n\ge 1}$ is an adic affine formal $\RR_{\XX}(r)$-scheme. We denote by 
	\begin{equation}
		\varpi(r):\TT_{\XX}(r)\to \RR_{\XX}(r) \label{pi T to R}
	\end{equation}
	the canonical morphism.	We set $\TT_{\XX}=\TT_{\XX}(1)$ and $\varpi=\varpi(1)$. 
\end{nothing}

\begin{nothing}\label{TX2 times}
	Let $n,r,r'$ be integers $\ge 1$. We denote by $J_{(R_{\XX}(r))_n}$ the ideal of $\mathscr{O}_{(\RR_{\XX}(r))_{n}}$ associated to the closed immersion $\XX_{n}\to (\RR_{\XX}(r))_{n}$, which induces an isomorphism on the underlying topological spaces. Via the $(\XX_{n}^{r+r'+1})$-isomorphism \eqref{iso R2 R times R} 
	\begin{displaymath}
		(\RR_{\XX}(r))_{n}\times_{\XX_{n}}(\RR_{\XX}(r'))_{n}\xrightarrow{\sim} (\RR_{\XX}(r+r'))_{n},
	\end{displaymath}
	the ideal $J_{(R_{\XX}(r))_n}\otimes_{\mathscr{O}_{\XX_{n}}}\mathscr{O}_{(\RR_{\XX}(r'))_{n}}+\mathscr{O}_{(\RR_{\XX}(r))_{n}}\otimes_{\mathscr{O}_{\XX_{n}}}J_{(R_{\XX}(r'))_n}$ corresponds to $J_{(R_{\XX}(r+r'))_n}$. In view of (\cite{Ber} II 1.3.5), we deduce a canonical isomorphism of PD-schemes 
	\begin{equation}
		\lambda_{n}:\TT_{\XX,n}(r)\times_{\XX_{n}}\TT_{\XX,n}(r')\xrightarrow{\sim} \TT_{\XX,n}(r+r'), \label{T 2 T times 1 iso mod n}
	\end{equation}
	where the projections $\TT_{\XX,n}(r)\to \XX_{n}$ (resp. $\TT_{\XX,n}(r')\to \XX_{n}$) is induced by the projection $\XX^{r+1}_{n}\to \XX_{n}$ on the last factor (resp. $\XX^{r'+1}_{n}\to \XX_{n}$ on the first factor). In view of the construction of \eqref{T 2 T times 1 iso mod n}, the isomorphisms $\lambda_{m}$ and $\lambda_{n}$ are compatibles for all integers $1\le m< n$. We deduce a canonical isomorphism of formal $(\XX^{r+r'+1})$-schemes
	\begin{equation} \label{T 2 T times 1 iso}
		\TT_{\XX}(r)\times_{\XX}\TT_{\XX}(r')\xrightarrow{\sim} \TT_{\XX}(r+r').
	\end{equation}
\end{nothing}


\begin{prop} \label{groupoid TX}
	The formal $\XX^{2}$-scheme $\TT_{\XX}$ has a natural formal $\XX$-groupoid structure such that the morphism $\varpi:\TT_{\XX}\to \RR_{\XX}$ \eqref{pi T to R} is a morphism of formal $\XX$-groupoids \textnormal{(\ref{morphism of groupoids}, \ref{prop R Q Hopf alg})}.
\end{prop}
\begin{proof} Since the morphism of underlying topological spaces $|\RR_{\XX}|\to |\XX^{2}|$ factors through $\Delta:|\XX|\to |\XX^{2}|$, the same holds for $|\TT_{\XX}|$. The $\XX^{2}$-morphism $\iota_{R}(r):\XX\to \RR_{\XX}(r)$ \eqref{iota R} induces a canonical $\XX^{2}$-morphism
\begin{equation}
	\iota_{T}(r):\XX\to \TT_{\XX}. \label{iota T}
\end{equation}
	The composition of $\varpi(2)$ and $\alpha_{R}$ \eqref{alpha R}
	\begin{equation}
		\TT_{\XX}(2) \to \RR_{\XX}(2)\to \RR_{\XX} 
	\end{equation}
	is compatible with $\iota_{T}(2)$ and $\iota_{R}$. By the universal property of $(\TT_{\XX,n})_{n\ge 1}$, we deduce an $\XX^{2}$-morphism
	\begin{equation}
		\alpha_{T}: \TT_{\XX}(2)\to \TT_{\XX} \label{alpha T}
	\end{equation}
	compatible with $\alpha_{R}$. We identify $\TT_{\XX}(2)$ and $\TT_{\XX}\times_{\XX}\TT_{\XX}$ (resp. $\TT_{\XX}(3)$ and $\TT_{\XX}\times_{\XX}\TT_{\XX}\times_{\XX}\TT_{\XX}$) by \ref{TX2 times}. The diagrams
	\begin{displaymath}
		\xymatrix{
			\TT_{\XX}\times_{\XX}\TT_{\XX}\times_{\XX}\TT_{\XX}  \ar[r]^-{\id\times \alpha_{T}} \ar[d] & \TT_{\XX}\times_{\XX}\TT_{\XX} \ar[r]^-{\alpha_{T}} \ar[d] & \TT_{\XX} \ar[d] \\
			\RR_{\XX}\times_{\XX}\RR_{\XX}\times_{\XX}\RR_{\XX} \ar[r]^-{\id\times \alpha_{R}}  & \RR_{\XX}\times_{\XX}\RR_{\XX} \ar[r]^-{\alpha_{R}} & \RR_{\XX} } \quad
		\xymatrix{
			\TT_{\XX}\times_{\XX}\TT_{\XX}\times_{\XX}\TT_{\XX}  \ar[r]^-{\alpha_{T}\times \id} \ar[d] & \TT_{\XX}\times_{\XX}\TT_{\XX} \ar[r]^-{\alpha_{T}} \ar[d] & \TT_{\XX} \ar[d] \\
			\RR_{\XX}\times_{\XX}\RR_{\XX}\times_{\XX}\RR_{\XX} \ar[r]^-{\alpha_{R}\times \id}  & \RR_{\XX}\times_{\XX}\RR_{\XX} \ar[r]^-{\alpha_{R}} & \RR_{\XX} }
	\end{displaymath}
	are commutative and the compositions of the lower horizontal arrows coincide. By the universal property of $(\TT_{\XX,n})_{n\ge 1}$, we deduce that $\alpha_{T}\circ(\id\times \alpha_{T})=\alpha_{T}\circ(\alpha_{T}\times \id)$. Since $\eta_{R}\circ \iota_{R}=\iota_{R}$ (\ref{def of groupoid str}(ii)), by the universal property of $(\TT_{\XX,n})_{n\ge 1}$, the composition $\TT_{\XX}\to \RR_{\XX}\xrightarrow{\eta_{R}} \RR_{\XX}$ \eqref{eta R} induces a morphism
	\begin{equation}
		\eta_{T}:\TT_{\XX}\to \TT_{\XX}. \label{eta T}
	\end{equation}

	By the universal property of $(\TT_{\XX,n})_{n\ge 1}$, we verify that $(\alpha_{T},\iota_{T},\eta_{T})$ is a formal $\XX$-groupoid structure on $\TT_{\XX}$ (cf. the proof of \ref{prop R Q Hopf alg}). In view of the proof, the morphism $\varpi$ is clearly a morphism of formal groupoids.
\end{proof}
\begin{nothing} \label{Local description T X}
	In the following, we recall Shiho's intepretation of integrable $p$-connections in terms of $\mathcal{T}_{\XX}$-stratifications \cite{Shiho}.

	To simplify the notation, we set $\mathcal{T}_{\XX}=\mathscr{O}_{\TT_{\XX}}$ considered as a formal Hopf $\mathscr{O}_{\XX}$-algebra \eqref{def of groupoid str} and we present a local description for it. Assume that there exists an \'{e}tale $\SS$-morphism $\XX\to \widehat{\mathbb{A}}_{\SS}^{d}=\Spf(\rW\{T_{1},\cdots,T_{d}\})$. We put $\xi_{i}=1\otimes T_{i}-T_{i}-1$ and we consider $\frac{\xi_{i}}{p}$ as the section of $\mathcal{R}_{\XX}$ \eqref{local description of R} and of $\mathcal{T}_{\XX}$. 

	Let $n$ be an integer $\ge 1$. By \ref{local description R Ad} and \ref{local description R lemma}, the closed immersion $\XX_{n}\to \RR_{\XX,n}$ of smooth $\SS_{n}$-schemes is regular (\cite{EGAIV} 17.12.1) and $(\frac{\xi_{1}}{p},\cdots,\frac{\xi_{d}}{p})$ is a regular sequence which generates $J_{\RR_{\XX},n}$ \eqref{TX2 times}. In view of (\cite{Ber} I 4.5.1 and 4.5.2), we deduce the following isomorphisms:
	\begin{equation} \label{des local TX}
		\mathscr{O}_{\XX_{n}}\langle x_{1},\cdots,x_{d}\rangle\xrightarrow{\sim} q_{1*}(\mathscr{O}_{\TT_{\XX,n}})\qquad \mathscr{O}_{\XX_{n}}\langle x_{1},\cdots,x_{d}\rangle\xrightarrow{\sim} q_{2*}(\mathscr{O}_{\TT_{\XX,n}}),
	\end{equation}
	where $q_{1},q_{2}:\TT_{\XX,n}\to \XX_{n}$ are the canonical projections and $x_{i}$ is sent to $\frac{\xi_{i}}{p}$ in both cases.

	Then we deduce the following isomorphisms \eqref{notation pd sym algebra}
	\begin{equation} \label{local p adic PD poly T}
		\mathscr{O}_{\XX} \llangle x_{1},\cdots,x_{d} \rrangle \xrightarrow{\sim} q_{1*}(\mathcal{T}_{\XX})\qquad \mathscr{O}_{\XX}\llangle x_{1},\cdots,x_{d} \rrangle \xrightarrow{\sim} q_{2*}(\mathcal{T}_{\XX}),
	\end{equation}
	where $q_{1},q_{2}:\TT_{\XX}\to \XX$ are the canonical projections and $x_{i}$ is send to $\frac{\xi_{i}}{p}$ in both cases. For any multi-index $I=(i_{1},i_{2},\cdots,i_{d})\in \mathbb{N}^{d}$, we put $(\frac{\xi}{p})^{[I]}=\prod_{j=1}^{d}(\frac{\xi_{j}}{p})^{[i_{j}]}\in \mathcal{T}_{\XX}$.
\end{nothing}

\begin{prop}[\cite{Shiho} Prop. 2.9]\label{qn pconnection T}
	There is a canonical equivalence of categories between the category of $\mathscr{O}_{\XX_{n}}$-modules with $\mathcal{T}_{\XX}$-stratification and the category $\pMIC^{\qn}(\XX_{n}/\SS_{n})$ \eqref{def quasinilpotent}.
\end{prop}

We recall the description of this equivalence in the local case (cf. \cite{Shiho} Prop. 2.9). Suppose that there exists an \'{e}tale $\SS$-morphism $\XX\to \widehat{\mathbb{A}}_{\SS}^{d}=\Spf(\rW\{T_{1},\cdots,T_{d}\})$. We take again the notation of \ref{quasi-nilpotent coords} and \ref{Local description T X}. Let $(M,\nabla)$ be an $\mathscr{O}_{\XX_{n}}$-modules with quasi-nilpotent integrable $p$-connection. The associated stratification $\varepsilon:\mathcal{T}_{\XX}\otimes_{\mathscr{O}_{\XX}}M\xrightarrow{\sim} M\otimes_{\mathscr{O}_{\XX}}\mathcal{T}_{\XX}$ is defined, for every local section $m$ of $M$ by
\begin{equation}
	\varepsilon(1\otimes m)=\sum_{I\in \mathbb{N}^{d}} \nabla_{\partial^{I}}(m)\otimes \biggl(\frac{\xi}{p}\biggr)^{[I]}, \label{stratification given by pconnection}
\end{equation}
	where the right hand side is a locally finite sum since $\nabla$ is quasi-nilpotent.

\begin{lemma}
	There exists a canonical morphism of formal $\XX$-groupoids 
	\begin{equation} \label{R to P}		
		\varsigma: \RR_{\XX} \to \PP_{\XX}.
	\end{equation}
\end{lemma}
\begin{proof} Recall \eqref{prop univ RQ} that we have a commutative diagram
\begin{displaymath}
	\xymatrix{
		\RR_{\XX,1} \ar[r] \ar[d] & \RR_{\XX} \ar[d] \\
		X \ar[r] & \XX^{2}
	}
\end{displaymath}
Since $\RR_{\XX}$ is flat over $\rW$, the ideal $(p)$ of $\mathscr{O}_{\RR_{\XX}}$ has a canonical PD-structure. By the universal property of PD-envelope, we deduce a canonical $\XX^{2}$-morphism $\varsigma:\RR_{\XX}\to \PP_{\XX}$. 
We denote by $(\alpha_{R},\iota_{R},\eta_{R})$ (resp. $(\alpha_{P},\iota_{P},\eta_{P})$) the formal groupoid structure on $\RR_{\XX}$ (resp. $\PP_{\XX}$). 
By the universal property of $(\PP_{\XX,n})_{n\ge 1}$, we verify that $\varsigma$ is a morphism of formal $\XX$-groupoids (cf. the proof of \ref{prop P to Q}).
\end{proof}
\begin{nothing} \label{lemma alg P to R}
	Let $s:\mathcal{P}_{\XX}\to \mathcal{R}_{\XX}$ be the homomorphism of formal Hopf algebras induced by $\varsigma$. We present a local description of $s$. 
	Suppose that there exists an \'{e}tale $\SS$-morphism $\XX\to \widehat{\mathbb{A}}_{\SS}^{d}=\Spf(\rW\{T_{1},\cdots,T_{d}\})$. We take again the notation of \ref{local description R lemma} and \ref{flatness PX}. 
	In view of the construction of $\varsigma$, for any multi-index $I\in \mathbb{N}^{d}$, we have
	\begin{equation} \label{calcul P to R alg}
		s(\xi^{[I]})=\frac{p^{|I|}}{I!}\Big(\frac{\xi}{p}\Big)^{I}.
	\end{equation}

	We denote by $J_{R}$ (resp. $J_{P}$) the ideal sheaf of $\mathcal{R}_{\XX}$ (resp. $\mathcal{P}_{\XX}$) associated to the closed immersion $\iota_{R}$ (resp. $\iota_{P}$). 
	Note that the $p$-adic valuation of $I!$ is $\le \sum_{k\ge 1}\lfloor \frac{|I|}{p^{k}}\rfloor < |I|$. 
	By \eqref{calcul P to R alg}, we deduce that in general, $s$ is injective and $s(J_{P})\subset pJ_{R}$. Then we have $s(J_{P}^{[i]})\subset p^{i}J_{R}^{i}$ for any $1\le i\le p-1$. By dividing by $p^{i}$, we obtain $\mathscr{O}_{\XX}$-bilinear morphisms
	\begin{equation}
		s^{i}:J^{[i]}_{P} \to J_{R}^{i} \qquad \forall~ 0\le i\le p-1. \label{P to R divided}
	\end{equation}
\end{nothing}

\section{Local constructions of Shiho} \label{local Shiho}
	In the section, we review Shiho's local Cartier transform \cite{Shiho} (which depends on a lifting of Frobenius) and explain how it can be understood in terms of the groupoids $\PP_{\XX},\TT_{\XX},\RR_{\XX}$ and $\QQ_{\XX}$. 

	We denote by $\XX$ a smooth formal $\SS$-scheme, by $X$ the special fiber of $\XX$, by $\XX'=\XX\times_{\SS,\sigma}\SS$ the base change of $\XX$ by $\sigma$ \eqref{notations} and by $\pi:\XX'\to \XX$ the canonical projection.

\begin{nothing} \label{lemma verification shiho}
	Let $n$ be an integer $\ge 1$. We assume that there exists an $\SS_{n+1}$-morphism $F_{n+1}:\XX_{n+1}\to \XX'_{n+1}$ \eqref{conventions adic} whose reduction modulo $p$ is the relative Frobenius morphism $F_{X/k}$ of $X$ \eqref{notations Yk} and we denote by $F_{n}$ the reduction modulo $p^{n}$ of $F_{n+1}$. The morphism $F_{n+1}$ induces an $(\mathscr{O}_{\XX_{n+1}})$-linear morphism $dF_{n+1}:F_{n+1}^{*}(\Omega^{1}_{\XX_{n+1}'/\SS_{n+1}})\to \Omega^{1}_{\XX_{n+1}/\SS_{n+1}}$ whose image is contained in $p\Omega^{1}_{\XX_{n+1}/\SS_{n+1}}$. By dividing by $p$, it induces an $\mathscr{O}_{\XX_{n}}$-linear morphism 
	\begin{equation} \label{dFn+1 Omega}
		\frac{dF_{n+1}}{p}:F_{n}^{*}(\Omega^{1}_{\XX_{n}'/\SS_{n}})\to \Omega^{1}_{\XX_{n}/\SS_{n}}.
	\end{equation}
	
	Let $(M',\nabla')$ be an $\mathscr{O}_{\XX'_{n}}$-module with an integrable $p$-connection relative to $\SS_{n}$ \eqref{notation MIC}.
	We denote by $\zeta_{n}$ the composition  
	\begin{equation} \label{zeta n Fn+1}
		\zeta_{n}:F_{n}^{*}(\Omega^{1}_{\XX'_{n}/\SS_{n}}\otimes_{\mathscr{O}_{\XX_{n}'}}M')\xrightarrow{\sim} F_{n}^{*}(\Omega^{1}_{\XX'_{n}/\SS_{n}})\otimes_{\mathscr{O}_{\XX_{n}}}F_{n}^{*}(M')\to \Omega^{1}_{\XX_{n}/\SS_{n}}\otimes_{\mathscr{O}_{\XX_{n}}}F_{n}^{*}(M'),
	\end{equation}
	i.e., the composition of $\frac{dF_{n+1}}{p}\otimes \id$ and the canonical isomorphism.

	Shiho constructs a $\rW_{n}$-linear morphism $\nabla:F_{n}^{*}(M')\to \Omega_{\XX_{n}/\SS_{n}}^{1}\otimes_{\mathscr{O}_{\XX_{n}}}F_{n}^{*}(M')$ as follows. For any local sections $f$ of $\mathscr{O}_{\XX_{n}}$ and $e$ of $M'$, we put 
	\begin{equation} \label{shiho formula}
		\nabla(f F_{n}^{*}(e))= f\zeta_{n}(F_{n}^{*}(\nabla'(e))) +df\otimes(F_{n}^{*}(e)).
	\end{equation}
	The morphism $\nabla$ is well-defined and is an integrable connection on $F_{n}^{*}(M')$ relative to $\SS_{n}$ (cf. \cite{Shiho} page 805-806).
	Shiho defines a functor (cf. \cite{Shiho} 2.5)
	\begin{eqnarray}
		\Phi_{n}: \pMIC(\XX_{n}'/\SS_{n}) &\to& \MIC(\XX_{n}/\SS_{n}), \label{Shiho Phi}\\
		(M',\nabla')&\mapsto& (F_{n}^{*}(M'),\nabla). \nonumber
	\end{eqnarray}
	The functor $\Phi_{1}$ sends quasi-nilpotent objects to quasi-nilpotent objects (cf. \ref{lemma qn} below). By d\'evissage (\cite{Shiho} 1.13), the same holds for $\Phi_{n}$. It induces an equivalence of the categories (\cite{Shiho} 3.1):
	\begin{equation}
		\Phi_{n}: \pMIC^{\qn}(\XX_{n}'/\SS_{n}) \xrightarrow{\sim} \MIC^{\qn}(\XX_{n}/\SS_{n}). \label{Shiho equi}
	\end{equation}

	Let $(M',\nabla')$ be an object of $\pMIC(\XX_{n}'/\SS_{n})$ and $(M,\nabla)=\Phi_{n}(M',\nabla')$. 
	In view of \eqref{shiho formula},  the adjunction morphism of $\id_{M}\otimes \wedge^{\bullet}(\frac{dF_{n+1}}{p})$ \eqref{dFn+1 Omega} induce a $\rW$-linear morphism of complexes \eqref{lambla dR complex}
	\begin{equation} \label{morphism of complexes Shiho}
	\lambda: M'\otimes_{\mathscr{O}_{\XX_{n}'}}\Omega_{\XX_{n}'/\SS_{n}}^{\bullet} \to F_{X/k*}(M\otimes_{\mathscr{O}_{\XX_{n}}}\Omega_{\XX_{n}/\SS_{n}}^{\bullet}).
	\end{equation}
	Indeed, for local coordinates $t_{1},\cdots,t_{d}$ of $\XX'_{n}$ over $\SS_{n}$, any local section of $M'\otimes_{\mathscr{O}_{\XX_{n}'}}\Omega_{\XX_{n}'/\SS_{n}}^{q}$ can be written as a sum of sections of the form $m\otimes dt_{i_{1}}\wedge\cdots\wedge dt_{i_{q}}$. Using \eqref{shiho formula}, one verifies that \eqref{morphism of complexes Shiho} is a morphism of complexes.
\end{nothing}

\begin{lemma} \label{lemma qn}
	Let $(M',\theta)$ be a Higgs module on $X'/k$ \eqref{notation MIC} and $\nabla$ the integrable connection on $M=F_{X/k}^{*}(M')$ constructed in \ref{lemma verification shiho}. If $(M',\theta)$ is quasi-nilpotent \eqref{def quasinilpotent}, then so is $(M,\nabla)$.
\end{lemma}
\begin{proof} The question being local, we can reduce to the case where there exists an \'{e}tale $\SS_{2}$-morphism $\XX_{2}\to \mathbb{A}_{\SS_{2}}^{d}=\Spec(\rW_{2}[T_{1},\cdots,T_{d}])$. 
	For any $1\le i\le d$, let $t_{i}$ be the image of $T_{i}$ in $\mathscr{O}_{X}$ and $t'_{i}=\pi^{*}(t_{i})\in \mathscr{O}_{X'}$. There exists a section $a_{i}$ of $\mathscr{O}_{X}$ such that $\frac{dF_{2}}{p}(dt'_{i})=t_{i}^{p-1}dt_{i}+da_{i}$ and an $\mathscr{O}_{X'}$-linear morphism $\theta_{i}:M'\to M'$ such that for every local section $e$ of $M'$, we have $\theta(e)=\sum_{i=1}^{d}dt'_{i}\otimes\theta_{i}(e)$. 
	We set $\partial_{i}$ the dual of $dt_{i}$. Then we have
\begin{equation}\label{nabla i mod p}
	\nabla_{\partial_{i}}(F_{X/k}^{*}(e))=t_{i}^{p-1}F^{*}_{X/k}(\theta_{i}(e))+\sum_{j=1}^{d}\frac{\partial a_{j}}{\partial t_{i}}F_{X/k}^{*}(\theta_{j}(e)).
\end{equation}

We denote by $\psi:M\to M\otimes_{\mathscr{O}_{X}} F_{X}^{*}(\Omega_{X/k}^{1})$ the $p$-curvature associated to $\nabla$ (\cite{Ka71} 5.0). There exists $\mathscr{O}_{X}$-linear endomorphisms $\psi_{i}:M\to M$ for $1\le i\le d$ such that
	\begin{equation}
		\psi=\sum_{i=1}^{d}\psi_{i}\otimes F_{X}^{*}(dt_{i}).
	\end{equation}
	Recall (\cite{Ka71} 5.2) that $\psi_{i}$ and $\psi_{j}$ commutes for $1\le i,j\le d$. For any $I=(i_{1},\cdots,i_{d})\in \mathbb{N}^{d}$, we put $\psi_{I}=\prod_{j=1}^{d}\psi_{j}^{i_{j}}$ and $\theta_{I}=\prod_{j=1}^{d}\theta_{j}^{i_{j}}$. 
	The $p$th iterate $\partial_{i}^{(p)}$ of $\partial_{i}$ is zero (\cite{Ka71} 5.0). Then we have
\begin{equation}
	\psi_{i}=(\nabla_{\partial_{i}})^{p}. \label{psi nabla p}
\end{equation}

By \eqref{nabla i mod p} and induction, one verifies that for any integer $l\ge 1$, there exist elements $\{a_{l,I}\in \mathscr{O}_{X}\}_{I\in \mathbb{N}^{d}, 1\le |I|\le l}$ such that for every local section $e$ of $M'$, we have
\begin{equation} \label{formula nabla partial l Shiho}
	(\nabla_{\partial_{i}})^{l}(F_{X/k}^{*}(e))=\sum_{1\le |I|\le l} a_{l,I}F_{X/k}^{*}(\theta_{I}(e)).
\end{equation}
Since the $\psi_{i}$'s are $\mathscr{O}_{X}$-linear morphisms, if there exists an integer $N$ such that $\theta_{I}(e)=0$ for all $|I|\ge N$, then $\psi_{I}(F_{X/k}^{*}(e))=0$ for all $|I|\ge N$ by \eqref{psi nabla p} and \eqref{formula nabla partial l Shiho}. We deduce that $\nabla$ is quasi-nilpotent. 
\end{proof}

\begin{rem} \label{Cartier descente}
	Given an $\mathscr{O}_{X'}$-module $M'$, the \textit{Frobenius descent connection} $\nabla_{\can}$ on $F_{X/k}^{*}(M')$ is defined for local sections $m$ of $M'$ and $f$ of $\mathscr{O}_{X}$, by 
	\begin{equation} \label{can connection}
		\nabla_{\can}(fF_{X/k}^{*}(m))=m\otimes df.
	\end{equation}
	It is integrable and of $p$-curvature zero. Cartier descent states that the functor $M'\mapsto (F_{X/k}^{*}(M'),\nabla_{\can})$ induces an equivalence of categories between the category of quasi-coherent $\mathscr{O}_{X'}$-modules and the full subcategory of $\MIC^{\qn}(X/k)$ consisting of quasi-coherent objects whose $p$-curvature is zero (\cite{Ka71} 5.1). Considering $\mathscr{O}_{X'}$-modules as Higgs modules with the zero Higgs field, by \eqref{shiho formula}, we see that $\Phi_{1}$ is compatible with Cartier descent.
\end{rem}
\begin{nothing}
	Let $(M',\theta)$ be a Higgs module on $X'/k$ and $\ell$ an integer $\ge 0$. We suppose that $(M',\theta)$ is nilpotent of level $\le \ell$, i.e. there exists an increasing filtration of $M'$:
	\begin{equation}
		0=N'_{0}\subset N'_{1} \subset \cdots \subset N_{\ell}' \subset N'_{\ell+1}=M'
	\end{equation}
	such that $\theta(N'_{i})\subset N'_{i-1}\otimes_{\mathscr{O}_{X'}}\Omega_{X'/k}^{1}$ for $1\le i\le \ell+1$. Then the induced Higgs field on $\gr^{i}_{N'}(M')$ is trivial. 

	We set $(M,\nabla)=\Phi_{1}(M',\theta)$ and $N_{i}=\Phi_{1}(N'_{i},\theta|N_{i}')$ for $0\le i\le \ell+1$. By \eqref{shiho formula}, we see that $\nabla$ induces an integrable connection on each graded piece $\gr^{i}_{N}(M)$, with zero $p$-curvature. 

	We have a filtration on the de Rham complex $M\otimes_{\mathscr{O}_{X}}\Omega_{X/k}^{\bullet}$ (resp. the Dolbeault complex $M'\otimes_{\mathscr{O}_{X'}}\Omega_{X'/k}^{\bullet}$) defined by: 
	\begin{displaymath}
		N_{i}\otimes_{\mathscr{O}_{X}}\Omega_{X/k}^{\bullet} \qquad \textnormal{(resp. } N'_{i}\otimes_{\mathscr{O}_{X'}}\Omega_{X'/k}^{\bullet}).
	\end{displaymath}
\end{nothing}

\begin{prop} \label{Cartier iso nilpotent}
	The morphism of complexes \eqref{morphism of complexes Shiho} induces for every $i\in [1,\ell+1]$ a quasi-isomorphism
	\begin{equation}\label{morphism Cartier iso nilpotent} 
		N_{i}'\otimes_{\mathscr{O}_{X'}}\Omega_{X'/k}^{\bullet} \to F_{X/k*}(N_{i}\otimes_{\mathscr{O}_{X}}\Omega_{X/k}^{\bullet}).
	\end{equation}
\end{prop}
\begin{proof} 
	We first consider the case where $\ell=0$, i.e. $\theta$ is the zero Higgs field. 
	We follow a similar argument of (\cite{Ogus04} 1.2) where Ogus shows an analogous result in the level of cohomology of complexes. 
	Then $\nabla$ is the Frobenius descent connection on $M=F_{X/k}^{*}(M')$. When $M'=\mathscr{O}_{X'}$, the morphism on the cohomology induced by $\lambda$ \eqref{morphism of complexes Shiho} is the Cartier isomorphism (\cite{Ka71} 7.2)
\begin{equation} \label{Cartier isomorphism}
	\rmC_{X/k}^{-1}:\Omega_{X'/k}^{i}\xrightarrow{\sim} \mathcal{H}^{i}(F_{X/k*}(\Omega_{X/k}^{\bullet})).
\end{equation}
Since $F_{X/k}$ induces an isomorphism on the underlying topological spaces, we have an isomorphism of complexes
\begin{equation}
	M'\otimes_{\mathscr{O}_{X'}}F_{X/k*}(\Omega_{X/k}^{\bullet})\xrightarrow{\sim} F_{X/k*}(M\otimes_{\mathscr{O}_{X}}\Omega_{X/k}^{\bullet}).
\end{equation}
Since $F_{X/k*}(\Omega_{X/k}^{\bullet})$ is a complex of flat $\mathscr{O}_{X'}$-modules whose cohomology sheaves are also flat \eqref{Cartier isomorphism}, the canonical morphism
\begin{equation}
	M'\otimes_{\mathscr{O}_{X'}}\Omega_{X'/k}^{i}\xrightarrow{\sim} M'\otimes_{\mathscr{O}_{X'}}\mathcal{H}^{i}(F_{X/k*}(\Omega_{X/k}^{\bullet})) \to \mathcal{H}^{i}(M'\otimes_{\mathscr{O}_{X'}}F_{X/k*}(\Omega_{X/k}^{\bullet}))
\end{equation}
is an isomorphism. The assertion in the case $\ell=0$ follows. 

We prove the general case by induction on $i$. The assertion for $i=1$ is already proved. If the assertion is true for $i-1$, then the assertion for $i$ follows by dévissage from the induction hypothesis. 
\end{proof}

In the remainder of this section, we suppose that there exists an $\SS$-morphism $F:\XX\to \XX'$ which lifts the relative Frobenius morphism $F_{X/k}$ of $X$. We take again the notation of $\RR_{\XX}$, $\QQ_{\XX}$, $\TT_{\XX}$ and $\PP_{\XX}$ (\ref{prop R Q Hopf alg}, \ref{groupoid TX}, \ref{groupoid PX}).

\begin{prop} \label{prop Q to R'}
	The morphism $F$ induces a morphism of formal groupoids above $F$ \textnormal{(\ref{morphism of groupoids})}
	\begin{equation}\label{Q to R'}
		\psi:\QQ_{\XX}\to \RR_{\XX'}.	
	\end{equation}
\end{prop}
\begin{proof} First, we show that there exists a unique morphism $g:\QQ_{\XX,1}\to X'$ which fits into a commutative diagram
\begin{equation} \label{diag g existence}
	\xymatrix{
		\QQ_{\XX,1}\ar[r] \ar[dd]_{g}& \QQ_{\XX} \ar[d]\\
		&\XX^{2} \ar[d]^{F^{2}} \\
		X'\ar[r]^{\Delta}& \XX'^{2} }
\end{equation}
	where the bottom map is induced by the diagonal immersion. 
	The problem being local on $\XX$, we can assume that there exists an \'{e}tale $\SS$-morphism $\XX\to \widehat{\mathbb{A}}_{\SS}^{d}=\Spf(\rW\{T_{1},\cdots,T_{d}\})$. 
	For any $1\le i\le d$, we put $t_{i}$ the image of $T_{i}$ in $\mathscr{O}_{\XX}$, $t_{i}'=\pi^{*}( t_{i}) \in \mathscr{O}_{\XX'}$, $\xi_{i}=1\otimes t_{i}-t_{i}\otimes 1 \in \mathscr{O}_{\XX^{2}}$ and $\xi_{i}'=1\otimes t'_{i}-t'_{i}\otimes 1 \in \mathscr{O}_{\XX'^{2}}$. Locally, there is a section $a_{i}$ of $\mathscr{O}_{\XX}$ such that $F^{*}(t_{i}')=t_{i}^{p}+p a_{i}$. Then we have 
	\begin{eqnarray}
		F^{2*}(\xi_{i}') &=& 1\otimes t_{i}^{p}-t_{i}^{p}\otimes 1 +p(1\otimes a_{i}-a_{i}\otimes 1) \label{calcul F2 pullback}\\
		&=& (\xi_{i}+t_{i}\otimes 1)^{p} -t_{i}^{p}\otimes 1 + p(1\otimes a_{i}-a_{i}\otimes 1) \nonumber\\
		&=& \xi_{i}^{p} + \sum_{j=1}^{p-1} \binom{p}{j}\xi_{i}^{j}(t_{i}\otimes 1)^{p-j} + p(1\otimes a_{i}-a_{i}\otimes 1). \nonumber
	\end{eqnarray}
	Since $\xi_{i}^{p}=p\cdot\big(\frac{\xi_{i}^{p}}{p}\big)$ in $\mathcal{Q}_{\XX}$, the image of $F^{2*}(\xi_{i}')$ in $\mathcal{Q}_{\XX}$ is contained in $p\mathcal{Q}_{\XX}$. 
	Then the existence and the uniqueness of $g$ follow. 
	By the universal property of $\RR_{\XX'}$ (\ref{prop univ RQ}), we deduce an $\XX'^{2}$-morphism $\psi:\QQ_{\XX}\to \RR_{\XX'}$.
	Using the universal property of $\RR_{\XX'}$, we verify that $\psi$ is a morphism of formal groupoids above $F$ (cf. the proof of \ref{prop P to Q}).
\end{proof}
\begin{nothing} \label{prop P to R'}
	We denote the composition of $\psi:\QQ_{\XX}\to \RR_{\XX'}$ and $\lambda:\PP_{\XX}\to \QQ_{\XX}$ \eqref{pi P to Q} by
	\begin{equation}
		\phi: \PP_{\XX}\to \RR_{\XX'}. \label{P to R'}
	\end{equation}
	The morphism $\phi$ is a morphism of formal groupoids above $F$ \eqref{morphism of groupoids}.
\end{nothing}

\begin{lemma}[\cite{Shiho} 2.14] \label{prop P to T'}
	The morphism $\phi$ induces a morphism of formal groupoids above $F$ 
	\begin{equation}
		\varphi:\PP_{\XX}\to \TT_{\XX'}. \label{P to T'}
	\end{equation}
\end{lemma}
\begin{proof} For any $n\ge 1$, by $\phi\circ \iota_{P}=F\circ \iota_{R'}$ (\ref{def of groupoid str}(ii)) and the universal property of $\TT_{\XX',n}$ \eqref{basic notation T X}, $\phi_{n}:\PP_{\XX_{n}}\to \RR_{\XX',n}$ induces a PD-$\RR_{\XX',n}$-morphism $\varphi_{n}:\PP_{\XX_n}\to \TT_{\XX',n}$. For any $1\le m< n$, since $\phi_{m}$ and $\phi_{n}$ are compatible, we see that $\varphi_{m}$ and $\varphi_{n}$ are compatible. Hence we obtain a $\XX'^{2}$-morphism $\varphi:\PP_{\XX}\to \TT_{\XX'}$. It follows from \ref{prop P to R'} and the universal property of $(\TT_{\XX',n})_{n\ge 1}$ that $\varphi$ is a morphism of formal groupoids above $F$.
\end{proof}
\begin{nothing}
	We have a commutative diagram of formal groupoids 
	\begin{equation} \label{square Frob}
		\xymatrixcolsep{4pc}\xymatrix{
			\PP_{\XX} \ar[r]^{\varphi} \ar[rd]^{\phi} \ar[d]_{\lambda}& \TT_{\XX'} \ar[d]^{\varpi}\\
			\QQ_{\XX} \ar[r]^{\psi} & \RR_{\XX'}}
	\end{equation}
where $\varphi,\phi,\psi$ are induced by $F$. By \ref{morphism Hopf functor}, we deduce a commutative diagram:
	\begin{equation}\label{square stra}
		\xymatrixcolsep{4pc}\xymatrix{
		\Big\{ \txt{	\textnormal{$\mathscr{O}_{\XX'_{n}}$-modules}\\
				\textnormal{with $\mathcal{R}_{\XX'}$-stratification}} \Big\}
		\ar[r]^{\psi_{n}^{*}} \ar[d]_{\varpi_{n}^{*}}&
		\Big\{	\txt{\textnormal{$\mathscr{O}_{\XX_{n}}$-modules}\\
		\textnormal{with $\mathcal{Q}_{\XX}$-stratification}} \Big\}
		\ar[d]^{\lambda_{n}^{*}}\\
		\Big\{ \txt{	\textnormal{$\mathscr{O}_{\XX'_{n}}$-modules} \\
		\textnormal{with $\mathcal{T}_{\XX'}$-stratification}} \Big\}
		\ar[r]^{\varphi_{n}^{*}} &
		\Big\{	\txt{\textnormal{$\mathscr{O}_{\XX_{n}}$-modules}\\
		\textnormal{with $\mathcal{P}_{\XX}$-stratification}} \Big\}}
	\end{equation}
\end{nothing}

\begin{nothing} \label{Shiho two construction} 
	Let $n$ be an integer $\ge 1$ and $F_{n}$ (resp. $F_{n+1}$) the reduction modulo $p^{n}$ (resp. $p^{n+1}$) of $F$. In (\cite{Shiho} Prop. 2.17), Shiho shows that, through the equivalences of categories \ref{B-O stratification MIC} and \ref{qn pconnection T}, the functor
	\begin{eqnarray}
		\varphi_{n}^{*}: \Big\{\begin{array}{c}
			\textnormal{$\mathscr{O}_{\XX'_{n}}$-modules}\\
			\textnormal{with $\mathcal{T}_{\XX'}$-stratification}
		\end{array}\Big\}
	&\to& 
	\Big\{\begin{array}{c}
		\textnormal{$\mathscr{O}_{\XX_{n}}$-modules}\\
		\textnormal{with $\mathcal{P}_{\XX}$-stratification}
		\end{array}\Big\}\\
		(M,\varepsilon) &\mapsto& (F_{n}^{*}(M),\varphi^{-1}(\varepsilon)\otimes_{\varphi^{-1}(\mathcal{T}_{\XX'})}\mathcal{P}_{\XX}). \nonumber 
	\end{eqnarray}
	coincides with the functor induced by $F_{n+1}$ \eqref{Shiho equi}
	\begin{displaymath}
		\Phi_{n}:\pMIC^{\qn}(\XX_{n}'/\SS_{n})\xrightarrow{\sim} \MIC^{\qn}(\XX_{n}/\SS_{n}).
	\end{displaymath} 
	Indeed, if we write down a $p$-connection and the associated connection in terms of local coordinates (as in \eqref{nabla i mod p}), then we can verify the compatibility between $\varphi_n$ and $\Phi_n$ using formula \eqref{calcul F2 pullback} (cf. \cite{Shiho} 2.17 for more details).
\end{nothing}

\section{Oyama topoi} \label{Oyama topos}
	In this section, $X$ denotes a $k$-scheme. We explain two ``crystalline like'' topoi introduced by Oyama \cite{Oy} associated to $X$. 
	When $X$ admits a smooth lifting $\XX$ to $\rW$, crystals on these topoi (introduced in \S~\ref{crystals}) are equivalent to modules with $\mathcal{R}_{\XX}$-stratification and (resp. $\mathcal{Q}_{\XX}$-stratification) discussed in \ref{conn and stra}, and are independent of the choice of any lifting of $X$. 

	If we use a gothic letter $\mathfrak{T}$ to denote an adic formal $\SS$-scheme, the corresponding roman letter $T$ will denote its special fiber $\mathfrak{T}\otimes_{\rW}k$. 

\begin{definition}[\cite{Oy} 1.3.1] \label{definition Oyama cat}
	We define the category $\pCRIS(X/\SS)$ (resp. $\CRIS(X/\SS)$)\footnote{In \cite{Oy}, the category $\pCRIS(X/\SS)$ (resp. $\CRIS(X/\SS)$) is denoted by $\textnormal{HIG}^{\gamma}(X/\SS)$ (resp. $\textnormal{CRIS}^{\gamma}(X/\SS)$).} as follows.
	\begin{itemize}
		\item[(i)] An object of $\pCRIS(X/\SS)$ is a triple $(U,\mathfrak{T},u)$ consisting of an open subscheme $U$ of $X$, a flat formal $\SS$-scheme $\mathfrak{T}$ \eqref{notations S flat} and an affine $k$-morphism $u:T\to U$.
			
		\item[(ii)] An object of $\CRIS(X/\SS)$ is a triple $(U,\mathfrak{T},u)$ consisting of an open subscheme $U$ of $X$, a flat formal $\SS$-scheme $\mathfrak{T}$ and an affine $k$-morphism $u:\underline{T}\to U$, where $\underline{T}$ is the closed subscheme of $T$ defined in \ref{notation underline}.
	
		\item[(iii)] Let $(U_{1},\mathfrak{T}_{1},u_{1})$ and $(U_{2},\mathfrak{T}_{2},u_{2})$ be two objects of $\pCRIS(X/\SS)$ (resp. $\CRIS(X/\SS)$). A morphism from $(U_{1},\mathfrak{T}_{1},u_{1})$ to $(U_{2},\mathfrak{T}_{2},u_{2})$ consists of an $\SS$-morphism $f:\mathfrak{T}_{1}\to \mathfrak{T}_{2}$ and an $X$-morphism $g:U_{1}\to U_{2}$ such that $g\circ u_{1}=u_{2}\circ f_{s}$ (resp. $g\circ u_{1}=u_{2}\circ \underline{f_{s}}$), where $f_{s}$ is the reduction modulo $p$ of $f$.
	\end{itemize}
\end{definition}

We denote an object $(U,\mathfrak{T},u)$ of $\pCRIS$ (resp. $\CRIS$) simply by $(U,\mathfrak{T})$, if there is no risk of confusion.

To simplify the notation, we drop $(X/\SS)$ in the notation $\mathscr{E}(X/\SS)$ (resp. $\underline{\mathscr{E}}(X/\SS)$) and we write simply $\mathscr{E}$ (resp. $\underline{\mathscr{E}}$), if there is no risk of confusion. We put $\pCRIS'=\pCRIS(X'/\SS)$ \eqref{notations Yk}.

\begin{lemma} \label{lemma flatness}
	Let $f:\mathfrak{T}_{1}\to \mathfrak{T}_{2}$ be an $\SS$-morphism of flat formal $\SS$-schemes and $f_{s}:T_{1}\to T_{2}$ its special fiber. 

	\textnormal{(i)} If $f_{s}$ is an isomorphism, then $f$ is an isomorphism.

	\textnormal{(ii)} If $f_{s}$ is flat, then the morphism $\mathfrak{T}_{1,n}\to \mathfrak{T}_{2,n}$ induced by $f$ \eqref{conventions adic} is flat for all integers $n\ge 1$.
\end{lemma}
\begin{proof} We can reduce to the case where $\mathfrak{T}_{1}=\Spf(B)$, $\mathfrak{T}_{2}=\Spf(A)$ are affine (\cite{Ab10} 2.1.37) and $f$ is induced by an adic $\rW$-homomorphism $u:A\to B$. For any integers $n\ge 1$, we put $A_{n}=A/p^{n}A$, $B_{n}=B/p^{n}B$, $\gr^{n}(A)=p^{n}A/p^{n+1}A$ and $\gr^{n}(B)=p^{n}B/p^{n+1}B$. Since $A$ and $B$ are flat over $\rW$, the canonical morphism of $B_{1}$-modules
\begin{equation} \label{gr Bourbaki}
	B_{1}\otimes_{A_{1}} \gr^{n}(A)\to \gr^{n}(B)
\end{equation}
is an isomorphism.

(i) If $u$ induces an isomorphism $A_{1}\xrightarrow{\sim} B_{1}$, we deduce that $u$ is an isomorphism by \eqref{gr Bourbaki} and (\cite{Comalg} III \S~2.8 Cor. 3 to Théo. 1).

(ii) If $B_{1}$ is flat over $A_{1}$, we deduce that $B_{n}$ is flat over $A_{n}$ for all integers $n\ge 1$ by \eqref{gr Bourbaki} and the local criterion of flatness (\cite{Comalg} III \S~5.2 Théo. 1).
\end{proof}
\begin{nothing} \label{fiber product Oyama}
	We say that a morphism $(U_{1},\mathfrak{T}_{1},u_{1})\to (U_{2},\mathfrak{T}_{2},u_{2})$ of $\pCRIS$ (resp. $\CRIS$) is \textit{flat} if the special fiber $T_{1}\to T_{2}$ of the morphism $\mathfrak{T}_{1}\to \mathfrak{T}_{2}$ is flat.

	Let $(U_{1},\mathfrak{T}_{1},u_{1})\to (U,\mathfrak{T},u)$ be a flat morphism and $(U_{2},\mathfrak{T}_{2},u_{2})\to (U,\mathfrak{T},u)$ a morphism of $\pCRIS$ (resp. $\CRIS$). Their fiber product is represented by the fiber products $\mathfrak{T}_{12}=\mathfrak{T}_{1}\times_{\mathfrak{T}}\mathfrak{T}_{2}$ (which is flat over $\SS$ in view of \ref{lemma flatness}(ii) and \ref{notations S flat}) and $U_{12}=U_{1}\times_{U}U_{2}$ endowed with the affine morphism $T_{1}\times_{T}T_{2}\to U_{1}\times_{U}U_{2}$ (resp. composite morphism $\underline{T_{1}\times_{T}T_{2}}\to \underline{T_{1}}\times_{\underline{T}}\underline{T_{2}}\to U_{1}\times_{U}U_{2}$ \eqref{underline product}) induced by $u_{1}$, $u_{2}$ and $u$.
\end{nothing}

\begin{definition}\label{def morphism Car}
	\textnormal{(i)} We say that a morphism $f:(U_{1},\mathfrak{T}_{1})\to (U_{2},\mathfrak{T}_{2})$ of $\pCRIS$ is \textit{Cartesian} if the canonical morphism $T_1\to T_2\times_{U_2}U_1$ is an isomorphism. 
		
	\textnormal{(ii)} We say that a morphism $f:(U_{1},\mathfrak{T}_{1},u_1)\to (U_{2},\mathfrak{T}_{2},u_2)$ of $\CRIS$ is \textit{Cartesian} if $T_1\to T_2$ is an open immersion and the canonical morphism $\underline{T}_1\to \underline{T}_2\times_{U_2}U_1$ is an isomorphism. \footnote{The above definition of Cartesian morphism in $\underline{\mathscr{E}}$ is equivalent to the original definition in (\cite{Oy} 1.3.1), where Oyama requires the canonical morphism $T_2\to U'_2\times_{U'_1,u'_1\circ f_{T_1/k}} T_1$ is an isomorphism \eqref{notation underline}.} 
\end{definition}

If $f:(U_{1},\mathfrak{T}_{1})\to (U_{2},\mathfrak{T}_{2})$ is a Cartesian morphism of $\pCRIS$ (resp. $\CRIS$), $T_{1}$ is an open subscheme of $T_{2}$ and $f$ identifies $\mathfrak{T}_{1}$ with the open formal subscheme induced by $\mathfrak{T}_{2}$ on $T_{1}$ by \ref{lemma flatness}(i).

A Cartesian morphism is clearly flat \eqref{fiber product Oyama}. By \ref{fiber product Oyama}, the base change of a Cartesian morphism of $\pCRIS$ (resp. $\CRIS$) is Cartesian. The composition of Cartesian morphisms of $\pCRIS$ (resp. $\CRIS$) is clearly Cartesian.

\begin{nothing} \label{CRIS to Zar fib}
	Let $(U,\mathfrak{T},u)$ be an object of $\pCRIS$ (resp. $\CRIS$) and $V$ an open subscheme of $U$. 
	Note that the canonical morphism $\underline{T}\to T$ induces an isomorphism on the underlying topological spaces. We denote by $\mathfrak{T}_{V}$ the restriction of $\mathfrak{T}$ to the open subset $u^{-1}(V)$ of the topological space $|T|$. 
	Then, we obtain an object $(V,\mathfrak{T}_{V})$ of $\pCRIS$ (resp. $\CRIS$) and a Cartesian morphism $(V,\mathfrak{T}_{V})\to (U,\mathfrak{T})$ of $\pCRIS$ (resp. $\CRIS$). 


	Any morphism $(U_{1},\mathfrak{T}_{1})\to (U_{2},\mathfrak{T}_{2})$ of $\pCRIS$ (resp. $\CRIS$) factors uniquely through the Cartesian morphism $(U_{1},(\mathfrak{T}_{2})_{U_{1}})\to (U_{2},\mathfrak{T}_{2})$. 
	The category $\pCRIS$ (resp. $\CRIS$) is fibered over the category $\mathbf{Zar}_{/X}$ of open subschemes of $X$ by the functor
	\begin{eqnarray}
		\pi:\pCRIS\to \mathbf{Zar}_{/X} \quad (\textnormal{resp. }\CRIS\to \mathbf{Zar}_{/X}) \quad (U,\mathfrak{T})\mapsto U, \label{pi Cris to Zar}
	\end{eqnarray}
\end{nothing}

\begin{nothing}\label{presheaf to datas}
	Let $(U,\mathfrak{T})$ be an object of $\pCRIS$ (resp. $\CRIS$). By \ref{CRIS to Zar fib}, we have a functor 
	\begin{equation}
		\alpha_{(U,\mathfrak{T})}:\Zar_{/U}\to \pCRIS \quad \textnormal{(resp. $\CRIS$)} \qquad V\mapsto (V,\mathfrak{T}_{V}). \label{alpha UT}
	\end{equation}

	Let $f:(U_{1},\mathfrak{T}_{1})\to (U_{2},\mathfrak{T}_{2})$ be a morphism of $\pCRIS$ (resp. $\CRIS$), $\jmath_{f}:\Zar_{/U_{1}}\to \Zar_{/U_{2}}$ the functor induced by composing with $U_{1}\to U_{2}$. Then, the morphism $f$ induces a morphism of functors:
	\begin{equation}
		\beta_{f}:\alpha_{(U_{1},\mathfrak{T}_{1})}\to \alpha_{(U_{2},\mathfrak{T}_{2})}\circ \jmath_{f}. \nonumber
	\end{equation}

	Let $\mathscr{F}$ be a presheaf on $\pCRIS$ (resp. $\CRIS$). We denote by $\mathscr{F}_{(U,\mathfrak{T})}$ the presheaf $\mathscr{F}\circ \alpha_{(U,\mathfrak{T})}$ on $\Zar_{/U}$. 
	The morphism $\beta_{f}$ induces a morphism of presheaves:
	\begin{equation} \label{morphism cf}
		\gamma_{f}:\mathscr{F}_{(U_{2},\mathfrak{T}_{2})}|_{U_{1}}\to \mathscr{F}_{(U_{1},\mathfrak{T}_{1})}.
	\end{equation}
	It is clear that $\gamma_{\id}=\id$. By construction, if $f$ is a Cartesian morphism, then $\gamma_{f}$ is an isomorphism. If $g:(U_{2},\mathfrak{T}_{2})\to (U_{3},\mathfrak{T}_{3})$ is another morphism of $\pCRIS$ (resp. $\CRIS$), we verify that $\gamma_{g\circ f}=\gamma_{f}\circ \jmath_{f}^{*}(\gamma_{g})$,
	where $\jmath_{f}^{*}$ denotes the localisation functor from the category of presheaves on $\Zar_{/U_{2}}$ to the category of presheaves on $\Zar_{/U_{1}}$.
\end{nothing}

\begin{prop}\label{prop presheaf to datas}
		A presheaf $\mathscr{F}$ on $\pCRIS$ (resp. $\CRIS$) is equivalent to the following data:
		
		\begin{itemize}
			\item[(i)] For every object $(U,\mathfrak{T})$ of $\pCRIS$ (resp. $\CRIS$) a presheaf $\mathscr{F}_{(U,\mathfrak{T})}$ on $\Zar_{/U}$,

			\item[(ii)] For every morphism $f:(U_{1},\mathfrak{T}_{1})\to (U_{2},\mathfrak{T}_{2})$ of $\pCRIS$ (resp. $\CRIS$) a morphism $\gamma_{f}:\mathscr{F}_{(U_{2},\mathfrak{T}_{2})}|_{U_{1}}\to \mathscr{F}_{(U_{1},\mathfrak{T}_{1})}$,
		\end{itemize}
		subject to the following conditions
	\begin{itemize}
		\item[(a)] If $f$ is the identity morphism of $(U,\mathfrak{T})$, then $\gamma_{f}$ is the identity of $\mathscr{F}_{(U,\mathfrak{T})}$.
		\item[(b)] If $f$ is a Cartesian morphism, then $\gamma_{f}$ is an isomorphism.			
		\item[(c)] If $f$ and $g$ are two composable morphisms of $\pCRIS$ (resp. $\CRIS$), then we have $\gamma_{g\circ f}=\gamma_{f}\circ \jmath_{f}^{*}(\gamma_{g})$.
	\end{itemize}
\end{prop}
\begin{proof} Let $\{\mathscr{F}_{(U,\mathfrak{T})},\gamma_{f}\}$ be a data as in the proposition. We associate to it a presheaf on $\pCRIS$ (resp. $\CRIS$) as follows. Let $(U,\mathfrak{T})$ be an object of $\pCRIS$ (resp. $\CRIS$). We define $\mathscr{F}(U,\mathfrak{T})=\mathscr{F}_{(U,\mathfrak{T})}(U)$. For any morphism $f:(V,\mathfrak{Z})\to (U,\mathfrak{T})$ of $\pCRIS$ (resp. $\CRIS$), we deduce a morphism from $\gamma_{f}$
\begin{equation}
	\mathscr{F}_{(U,\mathfrak{T})}(U)\to \mathscr{F}_{(U,\mathfrak{T})}(V)\to \mathscr{F}_{(V,\mathfrak{Z})}(V).
\end{equation}
In view of conditions (a) and (c), the correspondence 
\begin{equation}
	(U,\mathfrak{T})\mapsto \mathscr{F}_{(U,\mathfrak{T})}(U)
\end{equation}
is a presheaf. In view of condition (b), the above construction is quasi-inverse to the construction of \ref{presheaf to datas}. The assertion follows.
\end{proof}

We call \textit{descent data associated to $\mathscr{F}$} the data $\{\mathscr{F}_{(U,\mathfrak{T})},\gamma_{f}\}$ as in the proposition.
\begin{definition}[\cite{Oy} 1.3.1] \label{usual covering}
	Let $(U,\mathfrak{T},u)$ be an object of $\pCRIS(X/\SS)$ (resp. $\CRIS(X/\SS)$). We denote by $\Cov(U,\mathfrak{T},u)$ the set of families of Cartesian morphisms  $\{(U_{i},\mathfrak{T}_{i},u_{i})\to (U,\mathfrak{T},u)\}_{i\in I}$ such that $\{U_{i}\to U\}_{i\in I}$ is a Zariski covering.
\end{definition}

%

\begin{nothing} \label{usual topology}
	By \ref{def morphism Car}, we see that the sets $\Cov(U,\mathfrak{T})$ for $(U,\mathfrak{T})\in \Ob(\pCRIS(X/\SS))$ (resp. $(U,\mathfrak{T})\in \Ob(\CRIS(X/\SS))$) form a pretopology \textnormal{(\cite{SGAIV} II 1.3)}. 
	We call the topology on $\pCRIS(X/\SS)$ (resp. $\CRIS(X/\SS)$) associated to the pretopology defined by the $\Cov(U,\mathfrak{T})$'s \textit{Zariski topology}. We denote by $\widetilde{\mathscr{E}}(X/\SS)$ (resp. $\widetilde{\underline{\mathscr{E}}}(X/\SS)$) the topos of sheaves of sets on $\pCRIS(X/\SS)$ (resp. $\CRIS(X/\SS)$), equipped with the Zariski topology. 
\end{nothing}

\begin{prop}\label{sheaf to datas}
	Let $\mathscr{F}$ be a presheaf on $\pCRIS$ (resp. $\CRIS$) and $\{\mathscr{F}_{(U,\mathfrak{T})},\gamma_{f}\}$ the associated descent data \eqref{prop presheaf to datas}. Then $\mathscr{F}$ is a sheaf for the Zariski topology, if and only if for each object $(U,\mathfrak{T})$ of $\pCRIS$ (resp. $\CRIS$), the presheaf $\mathscr{F}_{(U,\mathfrak{T})}$ \eqref{presheaf to datas} is a sheaf of the Zariski topos $U_{\zar}$. 
\end{prop}

\begin{proof} Let $(U,\mathfrak{T})$ be an object of $\pCRIS$ (resp. $\CRIS$). The functor $\alpha_{(U,\mathfrak{T})}$ \eqref{alpha UT} sends morphisms of $\Zar_{/U}$ to Cartesian morphisms of $\pCRIS$ (resp. $\CRIS$) and commutes with fibered products. It is clearly continuous for the Zariski topologies. Hence, if $\mathscr{F}$ is a sheaf, then $\mathscr{F}_{(U,\mathfrak{T})}$ is a sheaf of $U_{\zar}$ (\cite{SGAIV} III 1.2). 

Conversely, suppose that each presheaf $\mathscr{F}_{(U,\mathfrak{T})}$ is a sheaf of $U_{\zar}$. Let $\{(U_{i},\mathfrak{T}_{i})\to (U,\mathfrak{T})\}_{i\in I}$ be an element of $\Cov(U,\mathfrak{T})$. In view of condition (b) of \ref{prop presheaf to datas}, we deduce that the sequence
\begin{displaymath}
	\mathscr{F}(U,\mathfrak{T})\to \prod_{i\in I}\mathscr{F}(U_{i},\mathfrak{T}_{i})\rightrightarrows \prod_{i,j\in I} \mathscr{F}(U_{ij},\mathfrak{T}_{ij})
\end{displaymath}
is exact, where $(U_{ij},\mathfrak{T}_{ij})=(U_{i},\mathfrak{T}_{i})\times_{(U,\mathfrak{T})}(U_{j},\mathfrak{T}_{j})$. Hence, $\mathscr{F}$ is a sheaf. 
\end{proof}
\begin{nothing} \label{fppf covering sch}
	Recall that a family of morphisms of schemes $\{f_{i}:T_{i}\to T\}_{i\in I}$ is called an \textit{fppf covering} if each morphism $f_{i}$ is flat and locally of finite presentation and if $|T|=\bigcup_{i\in I}f_{i}(|T_{i}|)$ (cf. \cite{SGAIII} IV 6.3.2).

	We say that a family of $\SS$-morphisms of flat formal $\SS$-schemes $\{f_{i}:\mathfrak{T}_{i}\to \mathfrak{T}\}_{i\in I}$ is an \textit{fppf covering} if each morphism $f_{i}$ is locally of finite presentation (\cite{Ab10} 2.3.15) and if its special fiber $\{T_{i}\to T\}_{i\in I}$ is an fppf covering of schemes. By (\cite{Ab10} 2.3.16) and \ref{lemma flatness}(ii), a family of $\SS$-morphisms of flat formal $\SS$-schemes $\{\mathfrak{T}_{i}\to \mathfrak{T}\}_{i\in I}$ is an fppf covering if and only if the family $\{\mathfrak{T}_{i,n}\to \mathfrak{T}_{n}\}_{i\in I}$ is an fppf covering of schemes for all integers $n\ge 1$.  

	Recall that an adic formal $\SS$-scheme $\mathfrak{T}$ is flat over $\SS$ if and only if $\mathfrak{T}_{n}$ is flat over $\SS_{n}$ for all integers $n\ge 1$ \eqref{notations S flat}. Since fppf coverings of schemes are stable by base change and by composition, the same holds for fppf coverings of flat formal $\SS$-schemes.
\end{nothing}

\begin{definition} \label{fppf covering}
	Let $(U,\mathfrak{T},u)$ be an object of $\pCRIS(X/\SS)$ (resp. $\CRIS(X/\SS)$). We denote by $\Cov_{\fppf}(U,\mathfrak{T},u)$ the set of families of flat morphisms $\{(U_{i},\mathfrak{T}_{i},u_{i})\to (U,\mathfrak{T},u)\}_{i\in I}$ \eqref{fiber product Oyama} such that $\{U_{i}\to U\}_{i\in I}$ is a Zariski covering and that $\{\mathfrak{T}_{i}\to \mathfrak{T}\}_{i\in I}$ is an fppf covering \eqref{fppf covering sch}.
\end{definition}

	It is clear that $\Cov_{\fppf}(U,\mathfrak{T},u)$ contains $\Cov(U,\mathfrak{T},u)$ \eqref{usual covering}.


\begin{nothing} \label{fppf topology}
	By \ref{fiber product Oyama} and \ref{fppf covering sch}, we see that the sets $\Cov_{\fppf}(U,\mathfrak{T})$ for $(U,\mathfrak{T})\in \Ob(\pCRIS(X/\SS))$ (resp. $(U,\mathfrak{T})\in \Ob(\CRIS(X/\SS))$) form a pretopology \textnormal{(\cite{SGAIV} II 1.3)}. 
	We call \textit{fppf topology} the topology on $\pCRIS(X/\SS)$ (resp. $\CRIS(X/\SS)$) associated to the pretopology defined by the sets $\Cov_{\fppf}(U,\mathfrak{T})$. We denote by $\widetilde{\pCRIS}(X/\SS)_{\fppf}$ (resp. $\widetilde{\CRIS}(X/\SS)_{\fppf}$) the topos of sheaves of sets on $\pCRIS(X/\SS)$ (resp. $\CRIS(X/\SS)$), equipped with the fppf topology.
\end{nothing}

\begin{nothing}\label{fppf sheaves to descent datas}
	Let $\mathscr{F}$ be a sheaf of $\widetilde{\pCRIS}_{\fppf}$ (resp. $\widetilde{\CRIS}_{\fppf}$) and $(U,\mathfrak{T})$ an object of $\pCRIS$ (resp. $\CRIS$). 
	Since $\mathscr{F}$ is also a sheaf for the Zariski topology, the presheaf $\mathscr{F}_{(U,\mathfrak{T})}$ \eqref{prop presheaf to datas} is a sheaf of $U_{\zar}$.
	
	Let $\{f:(U,\mathfrak{Z})\to (U,\mathfrak{T})\}$ be an element of $\Cov_{\fppf}(U,\mathfrak{T})$ and put $(U,\mathfrak{Z}\times_{\mathfrak{T}}\mathfrak{Z})=(U,\mathfrak{Z})\times_{(U,\mathfrak{T})}(U,\mathfrak{Z})$. Since $\mathscr{F}$ is a sheaf for the fppf topology, we deduce an exact sequence of $U_{\zar}$ \eqref{morphism cf}
	\begin{equation}
		\mathscr{F}_{(U,\mathfrak{T})} \xrightarrow{\gamma_{f}} \mathscr{F}_{(U,\mathfrak{Z})} \rightrightarrows \mathscr{F}_{(U,\mathfrak{Z}\times_{\mathfrak{T}}\mathfrak{Z})}. \label{descent datas sequence}
	\end{equation}
\end{nothing}

\begin{nothing}\label{generality morphism topos}
	Let $\mathscr{C}$ and $\mathscr{D}$ be two categories, $\widehat{\mathscr{C}}$ (resp. $\widehat{\mathscr{D}}$) the category of presheaves of sets on $\mathscr{C}$ (resp. $\mathscr{D}$) and $u:\mathscr{C}\to \mathscr{D}$ a functor. We have a functor
	\begin{equation} \label{pull back presheaves}
		\widehat{u}^{*}:\widehat{\mathscr{D}}\to \widehat{\mathscr{C}}\qquad \mathscr{G}\mapsto \widehat{u}^{*}(\mathscr{G})=\mathscr{G}\circ u.
	\end{equation}
	It admits a right adjoint (\cite{SGAIV} I 5.1)
	\begin{equation}
		\widehat{u}_{*}:\widehat{\mathscr{C}}\to \widehat{\mathscr{D}}.
	\end{equation}

	Let $\mathscr{C}$ and $\mathscr{D}$ be two sites \footnote{We suppose that the site $\mathscr{C}$ is small.}. If $u:\mathscr{C}\to \mathscr{D}$ is a cocontinuous (resp. continuous) functor and $\mathscr{F}$ (resp. $\mathscr{G}$) is a sheaf on $\mathscr{C}$ (resp. $\mathscr{D}$), then $\widehat{u}_{*}(\mathscr{F})$ (resp. $\widehat{u}^{*}(\mathscr{G})$) is a sheaf on $\mathscr{D}$ (resp. $\mathscr{C}$) (\cite{SGAIV} III 1.2, 2.2).

	Let $\widetilde{\mathscr{C}}$ (resp. $\widetilde{\mathscr{D}}$) be the topos of the sheaves of sets on $\mathscr{C}$ (resp. $\mathscr{D}$) and $u:\mathscr{C}\to \mathscr{D}$ a cocontinuous functor. Then $u$ induces a morphism of topoi $g:\widetilde{\mathscr{C}}\to \widetilde{\mathscr{D}}$ defined by $g_{*}=\widehat{u}_{*}$ and $g^{*}=a\circ \widehat{u}^{*}$, where $a$ is the sheafification functor (cf. \cite{SGAIV} III 2.3).
\end{nothing}

\begin{nothing}
	Note that the fppf topology on $\pCRIS$ (resp. $\CRIS$) is finer than the Zariski topology. Equipped with the fppf topology on the source and the Zariski topology on the target, the identical functors $\pCRIS\to \pCRIS$ and $\CRIS\to \CRIS$ are cocontinuous. By \ref{generality morphism topos}, they induce morphisms of topoi
	\begin{eqnarray} \label{functor restriction}
		\sigma:\widetilde{\pCRIS}_{\fppf}\to \widetilde{\pCRIS},\qquad \sigma:\widetilde{\CRIS}_{\fppf}\to \widetilde{\CRIS}.
	\end{eqnarray}
	If $\mathscr{F}$ is a sheaf of $\widetilde{\pCRIS}_{\fppf}$ (resp. $\widetilde{\CRIS}_{\fppf}$), $\sigma_{*}(\mathscr{F})$ is equal to $\mathscr{F}$. If $\mathscr{G}$ is a sheaf of $\widetilde{\pCRIS}$ (resp. $\widetilde{\CRIS}$), then $\sigma^{*}(\mathscr{G})$ is the sheafification of $\mathscr{G}$ with respect to the fppf topology.
\end{nothing}


%

\section{Crystals in Oyama topoi}\label{crystals}
In this section, we keep the notation in \S~\ref{Oyama topos} and we study crystals in Oyama topoi and explain their interpretations in terms of modules with stratification \eqref{equi crystals stratification}. 
Let $n$ an integer $\ge 1$.

\begin{nothing}\label{ring usual topology}
	We define a presheaf of rings $\mathscr{O}_{\mathscr{E}(X/\SS),n}$ on $\pCRIS(X/\SS)$ (resp. $\mathscr{O}_{\underline{\mathscr{E}}(X/\SS),n}$ on $\CRIS(X/\SS)$) \eqref{definition Oyama cat} by 
	\begin{equation}
		(U,\mathfrak{T})\mapsto \Gamma(\mathfrak{T},\mathscr{O}_{\mathfrak{T}_{n}}).
	\end{equation}

	For any object $(U,\mathfrak{T})$ of $\pCRIS(X/\SS)$ (resp. $\CRIS(X/\SS)$) and any element $\{(U_{i},\mathfrak{T}_{i})\to (U,\mathfrak{T})\}_{i\in I}$ of $\Cov_{\fppf}(U,\mathfrak{T})$ \eqref{fppf covering}, $\{\mathfrak{T}_{i,n}\to \mathfrak{T}_{n}\}_{i\in I}$ is an fppf covering of schemes \eqref{fppf covering sch}. 
	By fppf descent for quasi-coherent modules (\cite{SGAI} VIII 1.2), $\mathscr{O}_{\mathscr{E}(X/\SS),n}$ (resp. $\mathscr{O}_{\underline{\mathscr{E}}(X/\SS),n}$) is a sheaf for the fppf topology \eqref{fppf topology}. Since the fppf topology is finer than the Zariski topology, it is also a sheaf for the Zariski topology \eqref{usual topology}.	
\end{nothing}
\begin{nothing} \label{lin descent data}
	For any object $(U,\mathfrak{T},u)$ of $\pCRIS$ (resp. $\CRIS$), we consider $\mathscr{O}_{\mathfrak{T}_{n}}$ as a sheaf of $T_{\zar}$ (resp. $\underline{T}_{\zar}$). We have \eqref{sheaf to datas}
	\begin{displaymath}
		(\mathscr{O}_{\mathscr{E},n})_{(U,\mathfrak{T},u)}=u_{*}(\mathscr{O}_{\mathfrak{T}_{n}}) \quad \textnormal{(resp. $(\mathscr{O}_{\underline{\mathscr{E}},n})_{(U,\mathfrak{T},u)}=u_{*}(\mathscr{O}_{\mathfrak{T}_{n}})$)}.
	\end{displaymath}
	A morphism $f:(U_{1},\mathfrak{T}_{1},u_{1})\to (U_{2},\mathfrak{T}_{2},u_{2})$ of $\pCRIS$ (resp. $\CRIS$) induces a morphism of ringed topoi 
	\begin{equation} \label{tilde fn morphism}
		\widetilde{f}_{n}:(U_{1,\zar},u_{1*}(\mathscr{O}_{\mathfrak{T}_{1,n}}))\to (U_{2,\zar},u_{2*}(\mathscr{O}_{\mathfrak{T}_{2,n}})).
	\end{equation}

	By \ref{prop presheaf to datas} and \ref{sheaf to datas} and a standard argument, an $\mathscr{O}_{\pCRIS,n}$-module of $\widetilde{\pCRIS}$ (resp. $\mathscr{O}_{\CRIS,n}$-module of $\widetilde{\CRIS}$) is equivalent to the following data:
	\begin{itemize}
		\item[(i)] For every object $(U,\mathfrak{T},u)$ of $\pCRIS$ (resp. $\CRIS$), an $u_{*}(\mathscr{O}_{\mathfrak{T}_{n}})$-module $\mathscr{F}_{(U,\mathfrak{T})}$ of $U_{\zar}$,

		\item[(ii)] For every morphism $f:(U_{1},\mathfrak{T}_{1},u_{1})\to (U_{2},\mathfrak{T}_{2},u_{2})$ of $\pCRIS$ (resp. $\CRIS$), an $u_{1*}(\mathscr{O}_{\mathfrak{T}_{1,n}})$-linear morphism $c_{f}:\widetilde{f}_{n}^{*}(\mathscr{F}_{(U_{2},\mathfrak{T}_{2})})\to \mathscr{F}_{(U_{1},\mathfrak{T}_{1})}$,
	\end{itemize}
	which is subject to the following conditions
	\begin{itemize}
		\item[(a)] If $f$ is the identity morphism, then $c_{f}$ is the identity.
		\item[(b)] If $f$ is a Cartesian morphism, then $c_{f}$ is an isomorphism.
		\item[(c)] If $f$ and $g$ are two composable morphisms of $\pCRIS$ (resp. $\CRIS$), then we have $c_{g\circ f}=c_{f}\circ \widetilde{f}^{*}_{n}(c_{g})$. 
	\end{itemize}

We call the data $\{\mathscr{F}_{(U,\mathfrak{T})},c_{f}\}$ \textit{the linearized descent data associated to the $\mathscr{O}_{\pCRIS,n}$-module (resp. $\mathscr{O}_{\CRIS,n}$-module) $\mathscr{F}$} of $\widetilde{\pCRIS}$ (resp. of $\widetilde{\CRIS}$).

An $\mathscr{O}_{\pCRIS,n}$-module (resp. $\mathscr{O}_{\CRIS,n}$-module) $\mathscr{F}$ of $\widetilde{\pCRIS}_{\fppf}$ (resp. $\widetilde{\CRIS}_{\fppf}$) gives rise to an $\mathscr{O}_{\pCRIS,n}$-module (resp. $\mathscr{O}_{\CRIS,n}$-module) $\sigma_{*}(\mathscr{F})$ of $\widetilde{\pCRIS}$ (resp. $\widetilde{\CRIS}$) \eqref{functor restriction}. 
We can associate to $\mathscr{F}$ a linearized descent data $\{\mathscr{F}_{(U,\mathfrak{T})},c_{f}\}$ by that of $\sigma_{*}(\mathscr{F})$.
\end{nothing}

\begin{definition}[\cite{Oy} 1.3.3] \label{def crystals}
	Let $\mathscr{F}$ be an $\mathscr{O}_{\pCRIS,n}$-module of $\widetilde{\pCRIS}$ (resp. an $\mathscr{O}_{\pCRIS,n}$-module of $\widetilde{\pCRIS}_{\fppf}$, resp. an $\mathscr{O}_{\CRIS,n}$-module of $\widetilde{\CRIS}$, resp. an $\mathscr{O}_{\CRIS,n}$-module of $\widetilde{\CRIS}_{\fppf}$).

	\textnormal{(i)} We say that $\mathscr{F}$ is \textit{quasi-coherent} if for every object $(U,\mathfrak{T},u)$ of $\pCRIS$ (resp. $\CRIS$), the $u_{*}(\mathscr{O}_{\mathfrak{T}_{n}})$-module $\mathscr{F}_{(U,\mathfrak{T})}$ of $U_{\zar}$ \eqref{lin descent data} is quasi-coherent (\cite{EGAInew} 0.5.1.3).

	\textnormal{(ii)} We say that $\mathscr{F}$ is a \textit{crystal} or a \textit{crystal of $\mathscr{O}_{\pCRIS,n}$-modules of $\widetilde{\pCRIS}$} (resp. \textit{$\mathscr{O}_{\pCRIS,n}$-modules of $\widetilde{\pCRIS}_{\fppf}$}, resp. \textit{$\mathscr{O}_{\CRIS,n}$-modules of $\widetilde{\CRIS}$}, resp. \textit{$\mathscr{O}_{\CRIS,n}$-modules of $\widetilde{\CRIS}_{\fppf}$}) if for every morphism $f:(U_{1},\mathfrak{T}_{1})\to (U_{2},\mathfrak{T}_{2})$ of $\pCRIS$ (resp. $\CRIS$), the canonical morphism $c_{f}$ \eqref{lin descent data} is an isomorphism.
\end{definition}

We denote by $\CC(\mathscr{O}_{\pCRIS,n})$ (resp. $\CC_{\fppf}(\mathscr{O}_{\pCRIS,n})$, resp. $\CC(\mathscr{O}_{\CRIS,n})$, resp. $\CC_{\fppf}(\mathscr{O}_{\CRIS,n})$) the category of crystals of $\mathscr{O}_{\pCRIS,n}$-modules of $\widetilde{\pCRIS}$ (resp. $\mathscr{O}_{\pCRIS,n}$-modules of $\widetilde{\pCRIS}_{\fppf}$, resp. $\mathscr{O}_{\CRIS,n}$-modules of $\widetilde{\CRIS}$, resp. $\mathscr{O}_{\CRIS,n}$-modules of $\widetilde{\CRIS}_{\fppf}$) and we use the notation $\CC^{\qcoh}(-)$ or $\CC^{\qcoh}_{\fppf}(-)$ to denote the full subcategory consisting of quasi-coherent crystals.

\begin{prop} \label{crystal lin descent data}
	A quasi-coherent $\mathscr{O}_{\mathscr{E},n}$-module $\mathscr{F}$ of $\widetilde{\pCRIS}$ (resp. $\mathscr{O}_{\underline{\mathscr{E}},n}$-module of $\widetilde{\CRIS}$) is equivalent to
	
	\textnormal{(i)} For every object $(U,\mathfrak{T})$ of $\pCRIS$ (resp. $\CRIS$), a quasi-coherent $\mathscr{O}_{\mathfrak{T}_{n}}$-module $\mathscr{F}_{(U,\mathfrak{T})}$ of $\mathfrak{T}_{n,\zar}$;
	
	\textnormal{(ii)} For every morphism $f:(U_{1},\mathfrak{T}_{1})\to (U_{2},\mathfrak{T}_{2})$ of $\pCRIS$ (resp. $\CRIS$), an $(\mathscr{O}_{\mathfrak{T}_{1,n}})$-linear morphism of $\mathfrak{T}_{1,\zar}$:
	\begin{equation} \label{transition morphism qcoh}
		c_{f}:f_{n}^{*}(\mathscr{F}_{(U_{2},\mathfrak{T}_{2})})\to \mathscr{F}_{(U_{1},\mathfrak{T}_{1})},
	\end{equation}
	where $f_{n}$ denotes the morphism $\mathfrak{T}_{1,n}\to \mathfrak{T}_{2,n}$; 

	which are subject to similar conditions (a-c) in \ref{lin descent data}.
\end{prop}

The assertion follows from the following proposition.

\begin{prop} \label{datas crystals coro}
	Let $u:T\to U$ be an affine morphism of schemes, $i:T\to \mathcal{T}$ a nilpotent closed immersion. We denote by $v:(\mathcal{T},\mathscr{O}_{\mathcal{T}})\to (U,u_{*}(\mathscr{O}_{\mathcal{T}}))$ the morphism of ringed topoi induced by $u$. Then, the inverse image and direct image functors of $v$ induce equivalences of categories quasi-inverse to each other between the category of quasi-coherent $\mathscr{O}_{\mathcal{T}}$-modules of $\mathcal{T}_{\zar}$ and the category of quasi-coherent $u_{*}(\mathscr{O}_{\mathcal{T}})$-modules of $U_{\zar}$.
\end{prop}

The proposition follows from \ref{pushforward affine nil} and \ref{pullback affine nil}.

\begin{lemma} \label{exacteness of direct image of v}
	We keep the assumption of \ref{datas crystals coro}. 
	The restriction of the functor $v_{*}$ to the category of quasi-coherent $\mathscr{O}_{\mathcal{T}}$-modules is exact.
\end{lemma}
\begin{proof} The functor $v_{*}$ is left exact. Let $\mathscr{M}\to \mathscr{N}$ be a surjection of quasi-coherent $\mathscr{O}_{\mathcal{T}}$-modules. To show that $v_{*}(\mathscr{M})\to v_{*}(\mathscr{N})$ is surjective, it suffices to show that for any affine open subscheme $V$ of $U$, the morphism $v_{*}(\mathscr{M})(V)\to v_{*}(\mathscr{N})(V)$ is surjective. Since $u$ is affine, the open subscheme $\mathcal{T}_{V}$ of $\mathcal{T}$ associated to the open subset $u^{-1}(V)$ of $T$ is affine. Since $\mathscr{M}$, $\mathscr{N}$ are quasi-coherent and $\mathcal{T}_{V}$ is affine, the morphism $\mathscr{M}(\mathcal{T}_{V})\to \mathscr{N}(\mathcal{T}_{V})$ is surjective. Then the assertion follows.
\end{proof}
\begin{lemma} \label{pushforward affine nil}
	We keep the assumption of \ref{datas crystals coro}. 	
	For any quasi-coherent $\mathscr{O}_{\mathcal{T}}$-module $\mathscr{M}$, $v_{*}(\mathscr{M})$ is quasi-coherent and the adjunction morphism $v^{*}(v_{*}(\mathscr{M}))\to \mathscr{M}$ is an isomorphism. 
\end{lemma}
\begin{proof} The question being local, we may assume that $U$ is affine and so is $\mathcal{T}$. Then the quasi-coherent $\mathscr{O}_{\mathcal{T}}$ module $\mathscr{M}$ admits a presentation
\begin{equation} \label{presentation of M}
	\mathscr{O}_{\mathcal{T}}^{\oplus J}\to \mathscr{O}_{\mathcal{T}}^{\oplus I}\to \mathscr{M}\to 0.
\end{equation}
The canonical morphism $(v_{*}(\mathscr{O}_{\mathcal{T}}))^{\oplus I}\to v_{*}(\mathscr{O}_{\mathcal{T}}^{\oplus I})$ is clearly an isomorphism. Since $v_{*}$ is exact \eqref{exacteness of direct image of v}, we obtain a presentation of $v_{*}(\mathscr{M})$
\begin{displaymath}
	(v_{*}(\mathscr{O}_{\mathcal{T}}))^{\oplus J}\to (v_{*}(\mathscr{O}_{\mathcal{T}}))^{\oplus I}\to v_{*}(\mathscr{M})\to 0.
\end{displaymath}
The first assertion follows. 

The exact sequence \eqref{presentation of M} induces a commutative diagram
\begin{displaymath}
	\xymatrix{
		(v^{*}(v_{*}(\mathscr{O}_{\mathcal{T}})))^{\oplus J} \ar[r] \ar[d]& (v^{*}(v_{*}(\mathscr{O}_{\mathcal{T}})))^{\oplus I} \ar[r] \ar[d] & v^{*}(v_{*}(\mathscr{M}))\ar[r] \ar[d] & 0\\
		\mathscr{O}_{\mathcal{T}}^{\oplus J}\ar[r] & \mathscr{O}_{\mathcal{T}}^{\oplus I}\ar[r] & \mathscr{M}\ar[r] & 0.
	}
\end{displaymath}
Since $v_{*}$ is exact and $v^{*}$ is right exact, horizontal arrows are exact. 
The first two vertical arrows are isomorphisms. Then the second assertion follows.
\end{proof}
\begin{lemma} \label{pullback affine nil}
	We keep the assumption of \ref{datas crystals coro}. 
	For any quasi-coherent $u_{*}(\mathscr{O}_{\mathcal{T}})$-module $\mathscr{N}$, $v^{*}(\mathscr{N})$ is quasi-coherent and the adjunction morphism $\mathscr{N}\to v_{*}(v^{*}(\mathscr{N}))$ is an isomorphism.
\end{lemma}
\begin{proof} The pull-back functor sends quasi-coherent objects to quasi-coherent objects (\cite{EGAInew} 5.1.4). 

The second assertion being local, we may assume that $U$ is affine and that $\mathscr{N}$ admits a presentation:
\begin{displaymath}
	(u_{*}(\mathscr{O}_{\mathcal{T}}))^{\oplus J}\to (u_{*}(\mathscr{O}_{\mathcal{T}}))^{\oplus I} \to \mathscr{N}\to 0
\end{displaymath}
Then we deduce a commutative diagram
\begin{displaymath}
	\xymatrix{
		(u_{*}(\mathscr{O}_{\mathcal{T}}))^{\oplus J}\ar[r] \ar[d]& (u_{*}(\mathscr{O}_{\mathcal{T}}))^{\oplus I} \ar[r]  \ar[d]& \mathscr{N}\ar[r] \ar[d]& 0 \\
		(v_{*}(v^{*}(u_{*}(\mathscr{O}_{\mathcal{T}}))))^{\oplus J} \ar[r]& (v_{*}(v^{*}(u_{*}(\mathscr{O}_{\mathcal{T}}))))^{\oplus I} \ar[r] & v_{*}(v^{*}(\mathscr{N})) \ar[r] & 0 \\		
	}
\end{displaymath}
Since $v_{*}$ is exact \eqref{exacteness of direct image of v} and $v^{*}$ is right exact, the horizontal arrows are exact.
The first two vertical arrows are isomorphisms. Then the second assertion follows.
\end{proof}

\begin{prop} \label{qcoh fppf descent}
	Let $X$ be a $k$-scheme. 
	The direct image functors $\sigma_{*}:\widetilde{\pCRIS}_{\fppf}\to \widetilde{\pCRIS}$ and $\sigma_{*}:\widetilde{\CRIS}_{\fppf}\to \widetilde{\CRIS}$ \eqref{functor restriction} induce equivalences of categories
	\begin{equation} \label{functor qcoh Cris}
		\CC^{\qcoh}_{\fppf}(\mathscr{O}_{\pCRIS,n})\xrightarrow{\sim} \CC^{\qcoh}(\mathscr{O}_{\pCRIS,n}), \qquad \CC^{\qcoh}_{\fppf}(\mathscr{O}_{\CRIS,n})\xrightarrow{\sim} \CC^{\qcoh}(\mathscr{O}_{\CRIS,n}).
	\end{equation}
\end{prop}
\begin{proof} The functor $\sigma_{*}$ sends quasi-coherent crystals of $\mathscr{O}_{\pCRIS,n}$-modules of $\widetilde{\pCRIS}_{\fppf}$ (resp. $\mathscr{O}_{\CRIS,n}$-modules of $\widetilde{\CRIS}_{\fppf}$) to quasi-coherent crystals of $\mathscr{O}_{\pCRIS,n}$-modules of $\widetilde{\pCRIS}$ (resp. $\mathscr{O}_{\CRIS,n}$-modules of $\widetilde{\CRIS}$). The functors \eqref{functor qcoh Cris} are clearly fully faithful. It suffices to show the essential surjectivity. 

Let $\mathscr{F}$ be a quasi-coherent crystal of $\mathscr{O}_{\pCRIS,n}$-modules of $\widetilde{\pCRIS}$ and $\{\mathscr{F}_{(U,\mathfrak{T})},c_{f}\}$ the associated linearized descent data. We consider $\mathscr{F}_{(U,\mathfrak{T})}$ as an $\mathscr{O}_{\mathfrak{T}_n}$-module of $T_{\zar}$ \eqref{crystal lin descent data}. To show $\mathscr{F}$ is a sheaf for fppf topology, we prove that, for any element $\{(U_{i},\mathfrak{T}_{i})\to (U,\mathfrak{T})\}_{i\in I}$ of $\Cov_{\fppf}(U,\mathfrak{T})$ \eqref{fppf covering}, the sequence 
\begin{equation} \label{exact fppf}
	0\to \Gamma(T,\mathscr{F}_{(U,\mathfrak{T})})\to \prod_{i\in I} \Gamma(T_{i},\mathscr{F}_{(U_{i},\mathfrak{T}_{i})})\to \prod_{i,j\in I} \Gamma(T_{ij},\mathscr{F}_{(U_{ij},\mathfrak{T}_{ij})})
\end{equation}
is exact, where $(U_{ij},\mathfrak{T}_{ij})=(U_{i},\mathfrak{T}_{i})\times_{(U,\mathfrak{T})}(U_{j},\mathfrak{T}_{j})$ \eqref{fiber product Oyama}. Then, we have
\begin{equation}
	\mathscr{F}_{(U_{i},\mathfrak{T}_{i})}\simeq f_{i}^{*}(\mathscr{F}_{(U,\mathfrak{T})}),\quad \mathscr{F}_{(U_{ij},\mathfrak{T}_{ij})}\simeq f_{ij}^{*}(\mathscr{F}_{(U,\mathfrak{T})}),
\end{equation}
where $f_{i}$ (resp. $f_{ij}$) denotes the morphism $\mathfrak{T}_{i}\to \mathfrak{T}$ (resp. $\mathfrak{T}_{ij}\to \mathfrak{T}$). It suffices to show that the sequence 
\begin{equation} \label{exact fppf 2}
	0\to \mathscr{F}_{(U,\mathfrak{T})}\to \prod_{i\in I} f_{i*}(\mathscr{F}_{(U_{i},\mathfrak{T}_{i})})\to \prod_{i,j\in I}f_{ij*}(\mathscr{F}_{(U_{ij},\mathfrak{T}_{ij})})
\end{equation}
is exact. The question being local, we can suppose that $U$ is quasi-compact and quasi-separated and so is $\mathfrak{T}$. Hence we can reduce to the case where the set $I$ is finite. Note that the schemes $\mathfrak{T}_{i,n}$ and $\mathfrak{T}_{n}$ are quasi-compact and  quasi-separated and that the morphisms $f_{i}$ and $f_{ij}$ are quasi-compact. Then the exactness of \eqref{exact fppf 2} follows from the fppf descent for quasi-coherent modules.  

The assertion for quasi-coherent crystals of $\mathscr{O}_{\CRIS,n}$-modules can be verified in the same way.
\end{proof}

\begin{prop}\label{equi crystals stratification}
	Let $\XX$ be a smooth formal $\SS$-scheme and $X$ its special fiber. 
	There exists a canonical equivalence of tensor categories between the category $\CC(\mathscr{O}_{\pCRIS,n})$ (resp. $\CC(\mathscr{O}_{\CRIS,n})$) and the category of $\mathscr{O}_{\XX_{n}}$-modules with $\mathcal{R}_{\XX}$-stratification (resp. $\mathcal{Q}_{\XX}$-stratification) \textnormal{(\ref{prop R Q Hopf alg}, \ref{def stratification})}.
\end{prop}
\begin{proof} The proof is standard. We briefly explain the construction following (\cite{Oy} 1.3.4) where the author deals with the case $n=1$. 

	The triple $(X,\XX,\id)$ (resp. $(X,\XX,\underline{X}\to X)$) is an object of $\pCRIS$ (resp. $\CRIS$). By \ref{notations R Q}, $(X,\RR_{\XX}(r))$ (resp. $(X,\QQ_{\XX}(r))$) is an object of $\pCRIS$ (resp. $\CRIS$) for all integers $r\ge 1$.	
	We take again the notation of the proof of \ref{prop R Q Hopf alg}. We denote by $q_{1}, q_{2}:(X,\RR_{\XX})\to (X,\XX)$ the canonical morphisms.
	Let $\mathscr{F}$ be a crystal of $\mathscr{O}_{\pCRIS,n}$-modules of $\widetilde{\pCRIS}$. We set $\mathcal{F}$ to be the $\mathscr{O}_{\XX_{n}}$-module $\mathscr{F}_{(X,\XX)}$. Then we deduce isomorphisms of $\mathcal{R}_{\XX}$-modules of $X_{\zar}$ \eqref{lin descent data}
\begin{equation}
	c_{q_{1}}:\widetilde{q}_{1,n}^{*}(\mathcal{F})\xrightarrow{\sim} \mathscr{F}_{(X,\RR_{\XX})} \qquad c_{q_{2}}:\widetilde{q}_{2,n}^{*}(\mathcal{F})\xrightarrow{\sim} \mathscr{F}_{(X,\RR_{\XX})}. \label{stratification M}
\end{equation}
We denote by $\varepsilon$ the composition of $c_{q_{2}}$ and the inverse of $c_{q_{1}}$. By a standard argument, we can show that $\varepsilon$ defines a $\mathcal{R}_{\XX}$-stratification on $\mathcal{F}$.
It is clear that the correspondence $\mathscr{F}\to (\mathcal{F},\varepsilon)$ is functorial and is compatible with tensor products \eqref{def stratification}. 

	Conversely, let $\mathcal{F}$ be an $\mathscr{O}_{\XX_{n}}$-module with an $\mathcal{R}_{\XX}$-stratification $\varepsilon:\widetilde{q}_{2,n}^{*}(\mathcal{F})\xrightarrow{\sim} \widetilde{q}_{1,n}^{*}(\mathcal{F})$. 
	Let $(U,\mathfrak{T},u)$ be an object of $\pCRIS$ such that $U$ is affine. Since $\XX$ is smooth over $\SS$ and $\mathfrak{T}$ is affine, the $k$-morphism $u:T\to U$ extends to a morphism $\varphi:(U,\mathfrak{T})\to (X,\XX)$ of $\pCRIS$.
	We define $\mathscr{F}_{(U,\mathfrak{T})}$ to be the $\mathscr{O}_{\mathfrak{T}_{n}}$-module $\widetilde{\varphi}^{*}_{n}(\mathcal{F})$ of $U_{\zar}$. 
	By a standard arguement, we can show that this definition of $\mathscr{F}_{(U,\mathfrak{T})}$ is independent of the choice of the deformation $\mathfrak{T}\to \XX$ of $u:T\to U$ up to a canonical isomorphism which comes from the stratification.

Let $g:(U_{1},\mathfrak{T}_{1})\to (U_{2},\mathfrak{T}_{2})$ be a morphism of $\pCRIS$ such that $U_{1}$ and $U_{2}$ are affine. There exists morphisms $\varphi_{1}:(U_{1},\mathfrak{T}_{1})\to (X,\XX)$ and $\varphi_{2}:(U_{2},\mathfrak{T}_{2})\to (X,\XX)$ of $\pCRIS$. Then there exists a unique morphism $h:(U_{1},\mathfrak{T}_{1})\to (X,\RR_{\XX})$ such that $q_{1}\circ h=\varphi_{1}$ and $q_{2}\circ h=\varphi_{2}\circ g$. We deduce a canonical $\mathscr{O}_{\mathfrak{T}_{1,n}}$-linear isomorphism of $U_{1,\zar}$
\begin{equation}
	c_{g}:\widetilde{g}_{n}^{*}(\mathscr{F}_{(U_{2},\mathfrak{T}_{2})})=\widetilde{h}_{n}^{*}(\widetilde{q}_{2,n}^{*}(\mathcal{F}))\xrightarrow[\sim]{\widetilde{h}_{n}^{*}(\varepsilon)} \widetilde{h}_{n}^{*}(\widetilde{q}_{1,n}^{*}(\mathcal{F}))=\mathscr{F}_{(U_{1},\mathfrak{T}_{1})}.
\end{equation}
By gluing the constructions for affine objects, we obtain an isomorphism $c_{g}$ for a general morphism $g$ of $\pCRIS$. 
We deduce the cocycle properties for $c_{g}$ by the cocycle condtion of $\varepsilon$. 

Hence we obtain a linearized descent data $\{\mathscr{F}_{(U,\mathfrak{T},u)},c_{f}\}$ such that each morphism $c_{f}$ is an isomorphism. Then, we get a crystal of $\mathscr{O}_{\pCRIS,n}$-modules $\mathscr{F}$ of $\widetilde{\pCRIS}$ by \ref{lin descent data}. The correspondence $(\mathcal{F},\varepsilon)\mapsto \mathscr{F}$ is clearly functorial and quasi-inverse to the previous construction.

The assertion for crystals of $\mathscr{O}_{\CRIS,n}$-modules can be verified in the same way.
\end{proof}

\section{Cartier equivalence} \label{Cartier equiv}

In this section, we show that the relative Frobenius morphism of $X$ induces an equivalence of topoi between $\widetilde{\underline{\mathscr{E}}}_{\fppf}$ and $\widetilde{\mathscr{E}}'_{\fppf}$. Then we prove that this equivalence globalizes Shiho's local Cartier transform modulo $p^n$ explained in \S~\ref{local Shiho}. 

\begin{nothing}\label{morphism topos C}
	Let $(U,\mathfrak{T},u)$ be an object of $\CRIS$. We have a commutative diagram 
	\begin{equation}
		\xymatrixcolsep{4pc}\xymatrix{
			U \ar[d]_{F_{U/k}} & \underline{T}~ \ar[l]_{u} \ar@{^{(}->}[r] \ar[d]_{F_{\underline{T}/k}} & T \ar[d]^{F_{T/k}} \ar[ld]_{f_{T/k}}\\
			U' & \underline{T}' ~\ar[l]_{u'} \ar@{^{(}->}[r] & T' } \label{diagram rho}
\end{equation}
where the vertical arrows denote the relative Frobenius morphisms (\ref{notations Yk}, \ref{notation underline}). It is clear that $(U',\mathfrak{T},u'\circ f_{T/k})$ is an object of $\pCRIS'$ \eqref{definition Oyama cat}. Moreover, the correspondence $(U,\mathfrak{T},u)\mapsto (U',\mathfrak{T},u'\circ f_{T/k})$ is functorial. We denote by $\rho$ the functor defined as above:
	\begin{equation}\label{functor rho}
		\rho:\CRIS\to \pCRIS',\qquad (U,\mathfrak{T},u)\mapsto (U',\mathfrak{T},u'\circ f_{T/k}).
	\end{equation}

	We will show in \ref{lemma con con}(ii) that $\rho$ is continuous and cocontinuous with respect to either Zariski or fppf topology. By \ref{generality morphism topos}, the functor $\rho$ \eqref{functor rho} induces morphisms of topoi
	\begin{eqnarray}
		\label{morphism of topoi Cartier} &\rmC_{X/\SS}&:\widetilde{\CRIS}\to \widetilde{\pCRIS}',\\
		\label{morphism Cartier fppf} &\rmC_{X/\SS,\fppf}&:\widetilde{\CRIS}_{\fppf} \to \widetilde{\pCRIS}'_{\fppf}.
	\end{eqnarray}
	such that the pullback functor is induced by the composition with $\rho$. They fit into a commutative diagram \eqref{functor restriction}
	\begin{equation} \label{Cartier sigma fppf Zar}
		\xymatrixcolsep{5pc}\xymatrix{
			\widetilde{\CRIS}_{\fppf} \ar[r]^{\rmC_{X/\SS,\fppf}} \ar[d]_{\sigma}& \widetilde{\pCRIS}'_{\fppf} \ar[d]^{\sigma'} \\
			\widetilde{\CRIS} \ar[r]^{\rmC_{X/\SS}} & \widetilde{\pCRIS}'}
	\end{equation}
\end{nothing}

One of the main results in this section is the following.

\begin{theorem}[cf. \cite{Oy} 1.4.6] \label{eq topos fppf}
	For any smooth $k$-scheme $X$, the morphism $\rmC_{X/\SS,\fppf}:\widetilde{\CRIS}_{\fppf}\to \widetilde{\pCRIS}'_{\fppf}$ is an equivalence of topoi.
\end{theorem}

The theorem follows from \ref{lemma con con}, \ref{lemma rho 4} and \ref{lemma adjunction iso}.

\begin{prop}[cf. \cite{Oy} 1.4.1] \label{lemma con con}
	\textnormal{(i)} The functor $\rho$ \eqref{functor rho} is fully faithful.

	\textnormal{(ii)} Equipped with the Zariski topology \eqref{usual topology} (resp. fppf topology \eqref{fppf topology}) on both sides, the functor $\rho$ is continuous and cocontinuous \textnormal{(\cite{SGAIV} III 1.1, 2.1)}.
\end{prop}

\begin{lemma}\label{full faithful of rho}
	Let $Y$, $Z$ be two $k$-schemes and $g_{1},g_{2}:\underline{Y}\to Z$ two $k$-morphisms. We put $h_{i}=g'_{i}\circ f_{Y/k}:Y\to \underline{Y}'\to Z'$ \eqref{diagram rho} for $i=1,2$. If $h_{1}=h_{2}$, then $g_{1}=g_{2}$.
\end{lemma}
\begin{proof} Let $U$ be an affine open subscheme of $Z$. Since $f_{Y/k}$ is a homeomorphism and $h_{1}=h_{2}$, we have $g_{1}^{-1}(U)=g_{2}^{-1}(U)$. Hence we can reduce to case where $Z$ is affine.

Since the morphism $f_{Y/k}$ is scheme theoretically dominant \eqref{notation underline}, we deduce that $g_{1}'=g_{2}'$ by (\cite{EGAInew} 5.4.1). 
The functor $X\mapsto X'$ from the category of $k$-schemes to itself is clearly faithful. Then the lemma follows.
\end{proof}

\begin{lemma} \label{lemma cocon}
	Let $(U,\mathfrak{T},u)$ be an object of $\CRIS$ and $g:(V',\mathfrak{Z},w)\to \rho(U,\mathfrak{T},u)$ a morphism of $\pCRIS'$. Then there exist an object $(V,\mathfrak{Z},v)$ of $\CRIS$ and a morphism $f:(V,\mathfrak{Z},v)\to (U,\mathfrak{T},u)$ of $\CRIS$ such that $g=\rho(f)$.
	If $g$ is Cartesian (resp. flat), so is $f$. 
\end{lemma}
\begin{proof} Put $V=F_{U/k}^{-1}(V')$. Since the composition $Z\to T\xrightarrow{u'\circ f_{T/k}} U'$ factors through $V'\subset U'$, the composition $\underline{Z}\to \underline{T}\xrightarrow{u} U$ factors through $V$. We obtain a $k$-morphism $v:\underline{Z}\to V$ such that $w=v'\circ f_{Z/k}$. Since $f_{Z/k}$ is separated, $v'$ is affine (\cite{EGAII} 1.6.2(v)) and so is $v$. Hence, we get an object $(V,\mathfrak{Z},v)$ of $\CRIS$ and a morphism $f:(V,\mathfrak{Z},v)\to (U,\mathfrak{T},u)$ of $\CRIS$ such that $g=\rho(f)$.
\end{proof}

\begin{nothing}
	\textit{Proof of \ref{lemma con con}}. 
	(i) The functor $\rho$ is clearly faithful. We prove its fullness. Let $(U_{1},\mathfrak{T}_{1},u_{1})$, $(U_{2},\mathfrak{T}_{2},u_{2})$ be two objects of $\CRIS$ and $g:\rho(U_{1},\mathfrak{T}_{1},u_{1})\to \rho(U_{2},\mathfrak{T}_{2},u_{2})$ a morphism of $\pCRIS'$. Since $U_{1}'\subset U_{2}'\subset U'$, we have $U_{1}\subset U_{2}\subset U$. It suffices to show that the diagram 
\begin{displaymath}
	\xymatrix{
		\underline{T}_{1} \ar[r]^{\underline{g}_{s}} \ar[d]_{u_{1}} & \underline{T}_{2} \ar[d]^{u_{2}} \\
		U_{1} \ar[r] & U_{2}
	}
\end{displaymath}
is commutative. It follows from \ref{full faithful of rho} applied to the compositions $\underline{T_{1}}\xrightarrow{u_{1}}U_{1}\to U_{2}$ and $\underline{T_{1}}\to \underline{T_{2}}\xrightarrow{u_{2}}U_{2}$.

(ii) A family of morphisms $\{(U_{i},\mathfrak{T}_{i})\to (U,\mathfrak{T})\}_{i\in I}$ of $\CRIS$ belongs to $\Cov(U,\mathfrak{T})$ (resp. $\Cov_{\fppf}(U,\mathfrak{T})$) if and only if, its image by $\rho$ belongs to $\Cov(\rho(U,\mathfrak{T}))$ (resp. $\Cov_{\fppf}(\rho(U,\mathfrak{T}))$). 
	The functor $\rho$ sends flat morphisms to flat morphisms and it commutes with the fibered product of a flat morphism and a morphism of $\CRIS$. Indeed, by functoriality, if $T_{1}\to T$ and $T_{2}\to T$ are two morphisms of $k$-schemes, the canonical morphism $f_{T_{1}/k}\times f_{T_{2}/k}: T_{1}\times_{T}T_{2}\to \underline{T_{1}}'\times_{\underline{T}'}\underline{T_{2}}'$ is equal to the composition $T_{1}\times_{T}T_{2}\to (\underline{T_{1}\times_{T}T_{2}})'\to (\underline{T_{1}}\times_{\underline{T}}\underline{T}_{2})'=\underline{T_{1}}'\times_{\underline{T}'}\underline{T_{2}}'$.
Then the continuity of $\rho$ follows from (\cite{SGAIV} III 1.6). 

Let $\{(U_{i}',\mathfrak{T}_{i})\to \rho(U,\mathfrak{T})\}_{i\in I}$ be an element of $\Cov(\rho(U,\mathfrak{T}))$ (resp. $\Cov_{\fppf}(\rho(U,\mathfrak{T}))$). By \ref{lemma cocon}, there exists an element $\{(U_{i},\mathfrak{T}_{i})\to (U,\mathfrak{T})\}_{i\in I}$ of $\Cov(U,\mathfrak{T})$ (resp. $\Cov_{\fppf}(U,\mathfrak{T})$) mapping by $\rho$ to the given element. Then, $\rho$ is cocontinuous by \textnormal{(\cite{SGAIV} III 2.1)}. \hfill$\qed$
\end{nothing}

\begin{nothing} \label{lifting of Frobenius}
	To prove \ref{eq topos fppf}, we use (local) liftings of the Frobenius morphism to construct fppf converings. 
	Let $X$ be an affine scheme and $h:X\to Y=\Spec(k[T_{1},\cdots,T_{d}])$ an \'etale morphism. By (\cite{Hodge} 3.2), the following diagram is Cartesian \eqref{notations Yk}
	\begin{equation}\label{Car etale Frobenius}
		\xymatrix{
			X\ar[r]^{F_{X}}\ar[d]_{h} \ar@{}[dr]|{\Box} & X\ar[d]^{h}\\
			Y\ar[r]^{F_{Y}} & Y.}
	\end{equation}

	We put $\YY=\Spf(\rW\{\rT_{1},\cdots,T_{d}\})$ and we denote by $F_{\YY}:\YY\to \YY$ the affine morphism defined by $\sigma:\rW\to \rW$ and $T_{i}\mapsto T_{i}^{p}$. 
	Since $X$ is \'etale over $Y$, there exists a unique deformation $\XX$ of $X$ over $\YY$ up to a unique isomorphism. 
	The formal scheme $\XX\times_{\YY,F_{\YY}}\YY$ is also a deformation of $X=X\times_{Y,F_{Y}}Y$ over $\YY$. 
	Then we deduce a Cartesian diagram  
	\begin{equation}
		\xymatrix{
			\XX\ar[r] \ar[d] \ar@{}[dr]|{\Box}& \XX\ar[d]\\
			\YY\ar[r]^{F_{\YY}}& \YY}
	\end{equation}
	In particular, we obtain a morphism of finite type $F_{\XX}:\XX\to \XX$ above $\sigma$, which lifts the absolute Frobenius morphism $F_{X}$ of $X$. We put $\XX'=\XX\times_{\SS,\sigma}\SS$. Then we obtain an $\SS$-morphism of finite type $F_{\XX/\SS}:\XX\to \XX'$ which lifts the relative Frobenius morphism $F_{X/k}$. 
	Since $X$ is smooth, the morphism $F_{X/k}:X\to X'$ is faithfully flat (cf. \cite{Hodge} 3.2). Hence $\{F_{\XX/\SS}:\XX\to \XX'\}$ is an fppf covering in the sense of \ref{fppf covering sch}. 
\end{nothing}

\begin{lemma} \label{lemma rho 4}
	Let $(U',\mathfrak{T},u)$ be an object of $\pCRIS'$, $U=F_{X/k}^{-1}(U')$. Suppose that $U$ is affine and that there exists an \'etale $k$-morphism $U\to \Spec(k[T_{1},\cdots,T_{d}])$. 
	
	\textnormal{(i)} There exists an object $(U,\mathfrak{Z})$ of $\CRIS$ and an element $\{f:\rho(U,\mathfrak{Z})\to (U',\mathfrak{T})\}$ of $\Cov_{\fppf}(U',\mathfrak{T})$.

	\textnormal{(ii)} Let $g:(U'_{1},\mathfrak{T}_{1})\to (U',\mathfrak{T})$ be a morphism of $\pCRIS'$. Then there exists a morphism $h:(U_{1},\mathfrak{Z}_{1})\to (U,\mathfrak{Z})$ of $\CRIS$ and an element $\{\varphi:\rho(U_{1},\mathfrak{Z}_{1})\to (U'_{1},\mathfrak{T}_{1})\}$ of $\Cov_{\fppf}(U'_{1},\mathfrak{T}_{1})$ such that the following diagram is Cartesian:
	\begin{equation}\label{diag rho U1U2 pre}
	\xymatrix{
		\rho(U_{1},\mathfrak{Z}_{1})\ar[r]^{\varphi} \ar[d]_{\rho(h)} \ar@{}[dr]|{\Box}& (U_{1}',\mathfrak{T}_{1}) \ar[d]^{g}\\
		\rho(U,\mathfrak{Z})\ar[r]^{f} & (U',\mathfrak{T}) }
	\end{equation}
	If $g$ is a Cartesian morphism so is $h$. 
\end{lemma}

\begin{proof} (i) We follow the proof of (\cite{Oy} 1.4.5). 
	Let $\UU$ be a smooth lifting of $U$ over $\SS$ and $F:\UU\to \UU'$ a lifting of $F_{U/k}$ as in \ref{lifting of Frobenius}. Note that $U'$ is affine. Since the morphism $u$ is affine, $T$ and $\mathfrak{T}$ are affine. Since $\UU'$ is smooth over $\SS$, there exists an $\SS$-morphism $\tau:\mathfrak{T}\to \UU'$ which lifts $u$. We consider the following commutative diagram:
\begin{equation} \label{diag lifting Frob}
\xymatrix{
		& T\times_{U'}U  \ar[rr]\ar'[d][dd] \ar[ld] && \mathfrak{T}\times_{\mathfrak{U}'}\mathfrak{U} \ar[dd] \ar[ld]\\
		T \ar[rr]\ar[dd]_{u}
		&& \mathfrak{T}\ar[dd]\\
		& U \ar[ld]^{F_{U/k}} \ar'[r][rr]
		&& \mathfrak{U}\ar[ld]^{F}\\
		U' \ar[rr]
		&& \mathfrak{U}'}
\end{equation}
We set $\mathfrak{Z}=\mathfrak{T}\times_{\UU'}\UU$, $Z=T\times_{U'}U$ and we denote the composition $\underline{Z}\to Z \to U$ by $v$. Then we obtain an object $(U,\mathfrak{Z},v)$ of $\CRIS$. By \eqref{diag lifting Frob}, one verifies that the diagram
\begin{equation}
	\xymatrix{
		Z\ar[r] \ar[d]_{f_{Z/k}} & T \ar[d]^{u}\\
		\underline{Z}' \ar[r]^{v'}& U'}
\end{equation}
is commutative. Then, we obtain a morphism $\rho(U,\mathfrak{Z},v)\to (U',\mathfrak{T},u)$ of $\pCRIS'$. Since $\{F:\UU\to \UU'\}$ is an fppf covering, $\{\rho(U,\mathfrak{Z},v)\to (U',\mathfrak{T},u)\}$ is an element of $\Cov_{\fppf}(U',\mathfrak{T})$.

(ii) The morphism $f$ is flat. We denote by $(U'_{1},\mathfrak{Z}_{1})$ the fibered product $\rho(U,\mathfrak{Z})\times_{(U',\mathfrak{T})}(U'_{1},\mathfrak{T}_{1})$ in $\pCRIS$. By applying \ref{lemma cocon} to the projection $(U'_{1},\mathfrak{Z}_{1})\to \rho(U,\mathfrak{Z})$, we obtain the Cartesian diagram \eqref{diag rho U1U2 pre}. Since $\varphi$ is the base change of $f$, $\varphi$ is an element of $\Cov_{\fppf}(U'_{1},\mathfrak{T}_{1})$.
\end{proof}

The following lemma is a complement of \ref{lemma rho 4} and will be used in the proof of \ref{q-c crystal} below. 

\begin{lemma} \label{lemma rho 4 comp}
	Let $X$ be a $k$-scheme and $(U',\mathfrak{T})$ an object of $\pCRIS'$, $(U,\mathfrak{Z})$ an object of $\CRIS$ and $\{\rho(U,\mathfrak{Z})\to (U',\mathfrak{T})\}$ an element of $\Cov_{\fppf}(U',\mathfrak{T})$. Then there exists an object $(U,\mathfrak{Z}\times_{\mathfrak{T}}\mathfrak{Z})$ of $\CRIS$ and two morphisms $p_{1},p_{2}:(U,\mathfrak{Z}\times_{\mathfrak{T}}\mathfrak{Z})\to (U,\mathfrak{Z})$ such that $\rho(U,\mathfrak{Z}\times_{\mathfrak{T}}\mathfrak{Z})=\rho(U,\mathfrak{Z})\times_{(U',\mathfrak{T})}\rho(U,\mathfrak{Z})$ and that $\rho(p_{1})$ (resp. $\rho(p_{2})$) is the projection $\rho(U,\mathfrak{Z})\times_{(U',\mathfrak{T})}\rho(U,\mathfrak{Z})\to \rho(U,\mathfrak{Z})$ on the first (resp. second) component. 
\end{lemma}

\begin{proof} By applying \ref{lemma cocon} to the projection $\rho(U,\mathfrak{Z})\times_{(U',\mathfrak{T})}\rho(U,\mathfrak{Z})\to \rho(U,\mathfrak{Z})$ on the first component, we obtain an object $(U,\mathfrak{Z}\times_{\mathfrak{T}}\mathfrak{Z})$ of $\CRIS$ and a morphism $p_{1}:(U,\mathfrak{Z}\times_{\mathfrak{T}}\mathfrak{Z})\to (U,\mathfrak{Z})$ as in the proposition. The existence of $p_{2}$ follows from the fullness of $\rho$ (\ref{lemma con con}(i)).
\end{proof}

We conclude theorem \ref{eq topos fppf} by a general result on topoi due to Oyama \cite{Oy} which we do not repeat the proof. 

\begin{prop}[\cite{Oy} 4.2.1] \label{lemma adjunction iso}
	Let $\mathscr{C}$ be a site, $\mathscr{D}$ a site whose topology is defined by a pretopology and $u:\mathscr{C}\to \mathscr{D}$ a functor. Assume that:
	\begin{itemize}
		\item[(i)] $u$ is fully faithful,
		\item[(ii)] $u$ is continuous and cocontinuous,
		\item[(iii)] For every object $V$ of $\mathscr{D}$, there exists a covering of $V$ in $\mathscr{D}$ of the form $\{u(U_{i})\to V\}_{i\in I}$ where $U_{i}$ is an object of $\mathscr{C}$. 
	\end{itemize}
	Then the morphism of topoi $g:\widetilde{\mathscr{C}}\to \widetilde{\mathscr{D}}$ defined by $g^{*}=\widehat{u}^{*}$ and $g_{*}=\widehat{u}_{*}$ \eqref{generality morphism topos} is an equivalence of topoi.
\end{prop}

\begin{nothing}\label{morphism ringed topos}
	
	Let $\mathscr{F}$ be a sheaf of $\widetilde{\pCRIS}'$ (resp. $\widetilde{\pCRIS}'_{\fppf}$), $\{\mathscr{F}_{(U,\mathfrak{T})},\gamma_{f,\mathscr{F}}\}$ the descent data associated to $\mathscr{F}$ \eqref{prop presheaf to datas} and $\{\rmC^{*}(\mathscr{F})_{(U,\mathfrak{T})},\gamma_{f,\rmC^{*}(\mathscr{F})}\}$ the descent data associated to $\rmC^{*}_{X/\SS}(\mathscr{F})$ (resp. $\rmC^{*}_{X/\SS,\fppf}(\mathscr{F})$).  
	Since $\rho$ takes Cartesian morphisms to Cartesian morphisms \eqref{def morphism Car}, for any object $(U,\mathfrak{T})$ of $\CRIS$, we have 
	\begin{equation}
		\rmC^{*}(\mathscr{F})_{(U,\mathfrak{T})}=\pi_{U*}(\mathscr{F}_{\rho(U,\mathfrak{T})}) \label{C-1 calcul}
	\end{equation}
	where $\pi_{U}$ denotes the equivalence of topoi $U'_{\zar}\xrightarrow{\sim} U_{\zar}$. 
	For any morphism $f:(U_{1},\mathfrak{T}_{1})\to (U_{2},\mathfrak{T}_{2})$ of $\CRIS$, we verify that 
	\begin{equation}
		\gamma_{f,\rmC^{*}(\mathscr{F})}=\pi_{U_{1}*}(\gamma_{\rho(f),\mathscr{F}}). \label{C-1 calcul cf}
	\end{equation}

	We verify that, by definition 
	\begin{equation}
		\rmC^{*}_{X/\SS}(\mathscr{O}_{\pCRIS',n})=\mathscr{O}_{\CRIS,n} \qquad \textnormal{(resp. $\rmC^{*}_{X/\SS,\fppf}(\mathscr{O}_{\pCRIS',n})=\mathscr{O}_{\CRIS,n}$)}. 
	\end{equation}

	The morphism $\rmC_{X/\SS}$ (resp. $\rmC_{X/\SS,\fppf}$) is therefore underlying a morphism of ringed topoi, which we denote also by
	\begin{eqnarray}
		\label{morphism of ringed topos Cartier} &\rmC_{X/\SS}&:( \widetilde{\CRIS}, \mathscr{O}_{\CRIS,n}) \to ( \widetilde{\pCRIS}',\mathscr{O}_{\pCRIS',n}),\\
		\label{morphism of ringed topos Cartier fppf} \textnormal{(resp.} &\rmC_{X/\SS,\fppf}&:( \widetilde{\CRIS}_{\fppf}, \mathscr{O}_{\CRIS,n}) \to ( \widetilde{\pCRIS}'_{\fppf},\mathscr{O}_{\pCRIS',n})).
	\end{eqnarray}
\end{nothing}

	Let $\mathscr{F}'$ be a quasi-coherent $\mathscr{O}_{\mathscr{E}',n}$-module of $\widetilde{\pCRIS}'$. By \eqref{C-1 calcul}, $\rmC_{X/\SS}^{*}(\mathscr{F}')$ is also quasi-coherent. For any object $(U,\mathfrak{T})$ of $\CRIS$, we have an equality of $\mathscr{O}_{\mathfrak{T}_{n}}$-modules of $\mathfrak{T}_{n,\zar}$ \eqref{C-1 calcul}
	\begin{equation} \label{inverse image qc}
		(\rmC^{*}_{X/\SS}(\mathscr{F}'))_{(U,\mathfrak{T})}=\mathscr{F}'_{\rho(U,\mathfrak{T})}.
	\end{equation}

\begin{theorem}[cf. \cite{Oy} 1.4.3 for $n=1$] \label{thm pullback Cartier}
	Let $X$ be a smooth scheme over $k$. The inverse image and the direct image functors of $\rmC_{X/\SS}$ induce equivalences of categories quasi-inverse to each other \eqref{def crystals}
	\begin{equation} \label{equi crystal qcoh}
  		\CC^{\qcoh}(\mathscr{O}_{\pCRIS',n})\rightleftarrows\CC^{\qcoh}(\mathscr{O}_{\CRIS,n}).
	\end{equation}
\end{theorem}

The theorem follows from \ref{qcoh fppf descent} and \ref{q-c crystal} below. 

\begin{prop} \label{q-c crystal}
	Let $X$ be a smooth $k$-scheme.
	The inverse image and the direct image functors of the morphism $\rmC_{X/\SS,\fppf}$ \eqref{morphism of ringed topos Cartier fppf} induce equivalences of categories quasi-inverse to each other:
	\begin{equation}\label{Cartier modpn}
		\CC^{\qcoh}_{\fppf}(\mathscr{O}_{\pCRIS',n})\rightleftarrows\CC^{\qcoh}_{\fppf}(\mathscr{O}_{\CRIS,n}).
	\end{equation}
\end{prop}
\begin{proof} We write simply $\rmC$ for $\rmC_{X/\SS,\fppf}$. By \ref{eq topos fppf}, it suffices to show that the functors $\rmC^{*}$ and $\rmC_{*}$ preserve quasi-coherent crystals. The assertion for $\rmC^{*}$ follows from \eqref{C-1 calcul} and \eqref{C-1 calcul cf}. 

Let $\mathscr{F}$ be a quasi-coherent crystal of $\mathscr{O}_{\CRIS,n}$-modules of $\widetilde{\CRIS}_{\fppf}$ and $(U',\mathfrak{T},u)$ an object of $\pCRIS'$. We first show that $(\rmC_{*}(\mathscr{F}))_{(U',\mathfrak{T})}$ is quasi-coherent. 
The statement is local, therefore we may assume that $U=F_{X/k}^{-1}(U')$ satisifies the condition of \ref{lemma rho 4}, i.e. $U$ is affine and there exists an \'etale $k$-morphism $U\to \mathbb{A}_{k}^{d}$.
Then, by \ref{lemma rho 4}(i) and \ref{lemma rho 4 comp}, there exist objects $(U,\mathfrak{Z})$ and $(U,\mathfrak{Z}\times_{\mathfrak{T}}\mathfrak{Z})$ of $\CRIS$, an element $\{f:\rho(U,\mathfrak{Z})\to (U',\mathfrak{T})\}$ of $\Cov_{\fppf}(U',\mathfrak{T})$ and two morphisms $p_{1},p_{2}:(U,\mathfrak{Z}\times_{\mathfrak{T}}\mathfrak{Z})\to (U,\mathfrak{Z})$ such that $\rho(U,\mathfrak{Z}\times_{\mathfrak{T}}\mathfrak{Z})=\rho(U,\mathfrak{Z})\times_{(U',\mathfrak{T})}\rho(U,\mathfrak{Z})$ and that $\rho(p_{1})$ and $\rho(p_{2})$ are the canonical projections of $\rho(U,\mathfrak{Z})\times_{(U',\mathfrak{T})}\rho(U,\mathfrak{Z})$ \eqref{lemma rho 4 comp}. In particular, the morphism $\mathfrak{Z}\times_{\mathfrak{T}}\mathfrak{Z}\to \mathfrak{Z}$ attached to $p_{1}$ (resp. $p_{2}$) is the projection on the first (resp. second) component.

Since the adjunction morphism $\rmC^{*}\rmC_{*}\to \id$ is an isomorphism \eqref{eq topos fppf}, we have \eqref{C-1 calcul}
\begin{equation}\label{C-1 calcul app}
	\pi_{U*}((\rmC_{*}(\mathscr{F}))_{\rho(U,\mathfrak{Z})})=\mathscr{F}_{(U,\mathfrak{Z})}, \qquad \pi_{U*}((\rmC_{*}(\mathscr{F}))_{\rho(U,\mathfrak{Z}\times_{\mathfrak{T}}\mathfrak{Z})})=\mathscr{F}_{(U,\mathfrak{Z}\times_{\mathfrak{T}}\mathfrak{Z})}. 
\end{equation}
By \ref{crystal lin descent data}, we consider $\mathscr{F}_{(U,\mathfrak{Z})}$ (resp. $\mathscr{F}_{(U,\mathfrak{Z}\times_{\mathfrak{T}}\mathfrak{Z})}$) as a quasi-coherent $\mathscr{O}_{\mathfrak{Z}_{n}}$-module of $Z_{\zar}$ (resp. $(\mathscr{O}_{(\mathfrak{Z}\times_{\mathfrak{T}}\mathfrak{Z})_{n}})$-module of $(Z\times_{T}Z)_{\zar}$). Since $\mathscr{F}$ is a crystal, we have $(\mathscr{O}_{(\mathfrak{Z}\times_{\mathfrak{T}}\mathfrak{Z})_{n}})$-linear isomorphisms
\begin{equation}
	p_{2,n}^{*}(\mathscr{F}_{(U,\mathfrak{Z})})\xrightarrow[\sim]{c_{p_{2}}} \mathscr{F}_{(U,\mathfrak{Z}\times_{\mathfrak{T}}\mathfrak{Z})}\xleftarrow[\sim]{c_{p_{1}}} p_{1,n}^{*}(\mathscr{F}_{(U,\mathfrak{Z})}).
\end{equation}
We define $\varphi:p_{2}^{*}(\mathscr{F}_{(U,\mathfrak{Z})}) \to p_{1}^{*}(\mathscr{F}_{(U,\mathfrak{Z})})$ to be the composition of $c_{p_{2}}$ and the inverse of $c_{p_{1}}$. Thus, we obtain a effective descent datum $(\mathscr{F}_{(U,\mathfrak{Z})},\varphi)$ for the fppf covering $\{f_{n}:\mathfrak{Z}_{n}\to \mathfrak{T}_{n}\}$. Therefore there exists a quasi-coherent $\mathscr{O}_{\mathfrak{T}_{n}}$-module $\mathscr{M}$ of $T_{\zar}$, a canonical $(\mathscr{O}_{\mathfrak{Z}_{n}})$-linear isomorphism 
\begin{equation}
	f_{n}^{*}(\mathscr{M})\xrightarrow{\sim} \mathscr{F}_{(U,\mathfrak{Z})} \label{descent fppf}
\end{equation}
and an exact sequence of $U_{\zar}$
\begin{equation}
	0\to \pi_{U*}(u_{*}(\mathscr{M}))\to \mathscr{F}_{(U,\mathfrak{Z})}\to \mathscr{F}_{(U,\mathfrak{Z}\times_{\mathfrak{T}}\mathfrak{Z})}
\end{equation}
where $\mathscr{F}_{(U,\mathfrak{Z})}$ and $\mathscr{F}_{(U,\mathfrak{Z}\times_{\mathfrak{T}}\mathfrak{Z})}$ are now considered as sheaves of $U_{\zar}$.

On the other hand, since $\rmC_{*}(\mathscr{F})$ is a sheaf, there exists an exact sequence of $U'_{\zar}$ \eqref{descent datas sequence}
\begin{equation}
	0\to (\rmC_{*}(\mathscr{F}))_{(U',\mathfrak{T})}\to (\rmC_{*}(\mathscr{F}))_{\rho(U,\mathfrak{Z})}\to (\rmC_{*}(\mathscr{F}))_{\rho(U,\mathfrak{Z}\times_{\mathfrak{T}}\mathfrak{Z})}.
\end{equation}
By \eqref{C-1 calcul app}, we obtain an $u_{*}(\mathscr{O}_{\mathfrak{T}_{n}})$-linear isomorphism $u_{*}(\mathscr{M})\xrightarrow{\sim}(\rmC_{*}(\mathscr{F}))_{(U',\mathfrak{T})}$ of $U'_{\zar}$. In particular, $(\rmC_{*}(\mathscr{F}))_{(U',\mathfrak{T})}$ is quasi-coherent. Hence $\rmC_{*}(\mathscr{F})$ is quasi-coherent.

Let $g:(U_{1}',\mathfrak{T}_{1})\to (U'_{2},\mathfrak{T}_{2})$ be a morphism of $\pCRIS'$. We prove that the morphism \eqref{transition morphism qcoh}
\begin{equation}
	c_{g}:g_{n}^{*}(\rmC_{*}(\mathscr{F})_{(U'_{2},\mathfrak{T}_{2})})\to \rmC_{*}(\mathscr{F})_{(U_{1}',\mathfrak{T}_{1})}
\end{equation}
associated to $\rmC_{*}(\mathscr{F})$ is an isomorphism. 
Since the problem is Zariski local, we may assume that that $U'_{2}$ satisfies the condition of \ref{lemma rho 4}. 
By \ref{lemma rho 4}(ii), there exists a morphism $h:(U_{1},\mathfrak{Z}_{1})\to (U_{2},\mathfrak{Z}_{2})$ of $\CRIS$, an element $\{f_{1}:\rho(U_{1},\mathfrak{Z}_{1})\to (U'_{1},\mathfrak{T}_{1})\}$ of $\Cov_{\fppf}(U'_{1},\mathfrak{T}_{1})$ and an element $\{f_{2}:\rho(U_{2},\mathfrak{Z}_{2})\to (U'_{2},\mathfrak{T}_{2})\}$ of $\Cov_{\fppf}(U'_{2},\mathfrak{T}_{2})$ such that the following diagram is Cartesian
\begin{equation}\label{diag rho U1U2}
	\xymatrix{
		\rho(U_{1},\mathfrak{Z}_{1})\ar[r]^{f_{1}} \ar[d]_{\rho(h)}& (U_{1}',\mathfrak{T}_{1}) \ar[d]^{g}\\
		\rho(U_{2},\mathfrak{Z}_{2})\ar[r]^{f_{2}} & (U_{2}',\mathfrak{T}_{2})
	}
\end{equation}

By \ref{lemma rho 4 comp} and repeating the previous fppf descent argument, we have a canonical isomorphism \eqref{descent fppf}
\begin{equation}\label{iso descent datas}
	f_{i,n}^{*}(\rmC_{*}(\mathscr{F})_{(U_{i}',\mathfrak{T}_{i})}) \xrightarrow{\sim} \mathscr{F}_{(U_{i},\mathfrak{Z}_{i})}
\end{equation}
for $i=1,2$. Furthermore, since $\mathscr{F}$ is a crystal, we have an isomorphism
\begin{equation}\label{iso crystal h pullback}
	c_{h}:h_{n}^{*}(\mathscr{F}_{(U_{2},\mathfrak{Z}_{2})})\xrightarrow{\sim} \mathscr{F}_{(U_{1},\mathfrak{Z}_{1})}.
\end{equation}
In view of \eqref{C-1 calcul cf}, \eqref{C-1 calcul app}, \eqref{diag rho U1U2} and \eqref{iso descent datas}, the morphism
\begin{equation}
	f_{1,n}^{*}(c_{g}):f_{1,n}^{*}(g_{n}^{*}(\rmC_{*}(\mathscr{F})_{(U'_{2},\mathfrak{T}_{2})}))\to f_{1,n}^{*}(\rmC_{*}(\mathscr{F})_{(U_{1}',\mathfrak{T}_{1})})
\end{equation}
is identical to $c_{h}$ \eqref{iso crystal h pullback} and hence is an isomorphism. Since $f_{1,n}:\mathfrak{Z}_{1,n}\to \mathfrak{T}_{1,n}$ is faithfully flat, we deduce that $c_{g}$ is an isomorphism. The proposition follows.
\end{proof}

\begin{definition} \label{Cartier transform def}
	Let $X$ be a smooth $k$-scheme. 
	We call \textit{Cartier equivalence} (\textit{modulo $p^{n}$}) the equivalence of categories \eqref{equi crystal qcoh}
	\begin{equation} \label{Cartier transform}
		\rmC_{X/\SS}^{*}:\CC^{\qcoh}(\mathscr{O}_{\pCRIS',n})\xrightarrow{\sim} \CC^{\qcoh}(\mathscr{O}_{\CRIS,n}).
	\end{equation}
\end{definition}

The above equivalence depends only on $X$ and is different from the Cartier transform of Ogus-Vologodsky \cite{OV07} which depends on a lifting of $X'$ to $\rW_2$. 
We will compare two constructions under the assumption that there exists a smooth lifting of $X$ to $\rW$ in \S~ \ref{Comparison of Cartier}. 

To simplify the notation, we write $\rmC^{*}$ for $\rmC_{X/\SS}^{*}$ if there is no confusion.

\begin{nothing} \label{Cartier localization}
	The Cartier equivalence is compatible with localization with respect to an open subscheme. 
	More precisely, let $U$ be an open subscheme of $X$. The category $\pCRIS(U/\SS)$ (resp. $\CRIS(U/\SS)$) forms naturally a full subcategory of $\pCRIS(X/\SS)$ (resp. $\CRIS(X/\SS)$). Equipped with the Zariski topologies on both sides, the canonical functor $\pCRIS(U/\SS)\to \pCRIS(X/\SS)$ (resp. $\CRIS(U/\SS)\to \CRIS(X/\SS)$) is continuous and cocontinuous. It induces a morphism of topoi 
	\begin{equation} \label{morphism resctriction on U}
		j_{U}: \widetilde{\pCRIS}(U/\SS) \to \widetilde{\pCRIS}(X/\SS) \qquad \textnormal{(resp. } \underline{j}_{U}: \widetilde{\CRIS}(U/\SS) \to \widetilde{\CRIS}(X/\SS))
	\end{equation}
	such that the inverse image functor is given by resctricting a sheaf of $\widetilde{\pCRIS}(X/\SS)$ to $\pCRIS(U/\SS)$ (resp. $\widetilde{\CRIS}(X/\SS)$ to $\CRIS(U/\SS)$). The above morphisms fit into a commutative diagram
	\begin{equation} \label{localisation of Cartier morphism}
		\xymatrixcolsep{4pc}\xymatrix{
			\widetilde{\CRIS}(U/\SS) \ar[r]^{\rmC_{U/\SS}} \ar[d]_{\underline{j}_{U}}& \widetilde{\pCRIS}(U'/\SS) \ar[d]^{j_{U'}} \\
			\widetilde{\CRIS}(X/\SS) \ar[r]^{\rmC_{X/\SS}} & \widetilde{\pCRIS}(X'/\SS) }
	\end{equation}
Then we have a commutative diagram
\begin{equation} \label{Cartier localisation}
	\xymatrix{
		\CC(\mathscr{O}_{\pCRIS'(X/\SS),n})\ar[r]^{\rmC_{X/\SS}^{*}} \ar[d]_{j_{U}^{*}}& \CC(\mathscr{O}_{\CRIS(X/\SS),n}) \ar[d]^{\underline{j}_{U}^{*}} \\
		\CC(\mathscr{O}_{\pCRIS'(U/\SS),n})\ar[r]^{\rmC_{U/\SS}^{*}} & \CC(\mathscr{O}_{\CRIS(U/\SS),n}) }
\end{equation}
\end{nothing}

\begin{nothing} \label{prep model X lifting F}
	In the remainder of this section, $\XX$ denotes a smooth formal $\SS$-scheme with special fiber $X$. 
	We set $\XX'=\XX\times_{\SS,\sigma}\SS$ \eqref{notations} and we suppose that there exists an $\SS$-morphism $F:\XX\to \XX'$ which lifts the relative Frobenius morphism $F_{X/k}:X\to X'$ of $X$. We show that the Cartier equivalence $\rmC^{*}$ \eqref{Cartier transform} globalises Shiho's local construction in \S~\ref{local Shiho} defined by $F$. 
	
	The morphism $F$ induces a morphism of $\pCRIS'$ that we denote also by
	\begin{equation}\label{F in E'}
		F:\rho(X,\XX)=(X',\XX,F_{X/k})\to (X',\XX',\id)
	\end{equation}
	Recall that \eqref{prop Q to R'} the morphism $F$ induces a morphism of formal groupoids $\psi:\QQ_{\XX}\to \RR_{\XX'}$ above $F$. 

	The following result explains the relation between the Cartier equivalence $\rmC^{*}$ and the functor $\psi_n^{*}$ induced by $\psi$ \eqref{morphism Hopf functor}. 
\end{nothing}

\begin{prop} \label{Cartier global local}
	Keep the assumption of \ref{prep model X lifting F}. The diagram \eqref{q-c crystal}:
	\begin{equation} \label{diag Cartier F}
		\xymatrix{
			\CC(\mathscr{O}_{\pCRIS',n}) \ar[r]^{\rmC^{*}} \ar[d]_{\wr} & \CC(\mathscr{O}_{\CRIS,n}) \ar[d]^{\wr}\\
		\Big\{ \txt{	\textnormal{$\mathscr{O}_{\XX'_{n}}$-modules}\\
				\textnormal{with $\mathcal{R}_{\XX'}$-stratification}} \Big\}
		\ar[r]^{\psi_{n}^{*}} &
		\Big\{	\txt{\textnormal{$\mathscr{O}_{\XX_{n}}$-modules}\\
		\textnormal{with $\mathcal{Q}_{\XX}$-stratification}} \Big\} }
	\end{equation}
	where the vertical arrows are the equivalences of categories defined in \ref{equi crystals stratification}, is 2-commutative.
	That is, there exists a functorial isomorphism of $\mathscr{O}_{\XX_{n}}$-modules, depending on $F$
	\begin{equation} \label{eta F}
		\eta_{F}:\psi^{*}_{n}(\mathscr{M}_{(X',\XX')})\xrightarrow{\sim} \rmC^{*}(\mathscr{M})_{(X,\XX)}
	\end{equation}
	compatible with the $\mathcal{Q}_{\XX}$-stratifications, for every crystal $\mathscr{M}$ of $\mathscr{O}_{\mathscr{E}',n}$-modules of $\widetilde{\pCRIS}'$.
\end{prop}
\begin{proof} Let $\mathscr{M}$ be a crystal of $\mathscr{O}_{\mathscr{E}',n}$-modules of $\widetilde{\pCRIS}'$. By \eqref{C-1 calcul}, we have 
\begin{eqnarray}
	\rmC^{*}(\mathscr{M})_{(X,\XX)}=\pi_{X*}(\mathscr{M}_{\rho(X,\XX)})
\end{eqnarray}
Since $\mathscr{M}$ is a crystal, $F$ \eqref{F in E'} induces a functorial isomorphism of $\mathscr{O}_{\XX_{n}}$-modules
\begin{equation} \label{eta F modules}
	\eta_{F}:F^{*}_{n}(\mathscr{M}_{(X',\XX')})=\pi_{X*}(\widetilde{F}^{*}_{n}(\mathscr{M}_{(X',\XX')}))\xrightarrow{\sim} \rmC^{*}(\mathscr{M})_{(X,\XX)}.
\end{equation}

	The composition $\QQ_{\XX,1}\to \underline{\QQ_{\XX,1}}'\to X'$ \eqref{diagram rho} identifies with the morphism $g:\QQ_{\XX,1}\to X'$ induced by $F^{2}:\XX^{2}\to \XX'^{2}$ \eqref{diag g existence}. 
	Then, it follows from the proof of \ref{prop Q to R'} that $\psi$ induces a morphism of $\pCRIS'$ that we denote also by
	\begin{equation}\label{psi in E'}
		\psi:\rho(X,\QQ_{\XX})\to (X',\RR_{\XX'})
	\end{equation}
	which fits into the following commutative diagrams
\begin{equation} \label{double diag proj}
	\xymatrix{
		\rho(X,\QQ_{\XX}) \ar[r]^{\psi} \ar[d]_{\rho(q_{1})} & (X',\RR_{\XX'}) \ar[d]^{q_{1}'}\\
		\rho(X,\XX) \ar[r]^{F} & (X',\XX') }\qquad
	\xymatrix{
		\rho(X,\QQ_{\XX}) \ar[r]^{\psi} \ar[d]_{\rho(q_{2})} & (X',\RR_{\XX'}) \ar[d]^{q_{2}'}\\
		\rho(X,\XX) \ar[r]^{F} & (X',\XX') }
\end{equation}
where $q_{1},q_{2}$ (resp. $q_{1}',q_{2}'$) are the canonical projections of $(X,\QQ_{\XX})$ to $(X,\XX)$ (resp. $(X',\RR_{\XX'})$ to $(X',\XX')$). 
Hence, $\psi$ \eqref{psi in E'} induces an isomorphism
\begin{equation} \label{psi pullback cal}
	\pi_{X*}(\widetilde{\psi}^{*}_{n}(\mathscr{M}_{(X',\RR_{\XX'})}))\xrightarrow{\sim} \rmC^{*}(\mathscr{M})_{(X,\QQ_{\XX})}.
\end{equation}

Recall that the left vertical functor of \eqref{diag Cartier F} is given by $\mathscr{M}\mapsto (\mathscr{M}_{(X',\XX')},\varepsilon')$ where $\varepsilon'$ is induced by isomorphisms $\widetilde{q}'^{*}_{2,n} (\mathscr{M}_{(X',\XX')})\xrightarrow{\sim} \mathscr{M}_{(X',\RR_{\XX'})} \xleftarrow{\sim} \widetilde{q}'^{*}_{1,n}(\mathscr{M}_{(X',\XX')})$. By \ref{morphism Hopf functor}, we have
\begin{equation}
	\psi_{n}^{*}( \mathscr{M}_{(X',\XX')},\varepsilon')=(F^{*}_{n}(\mathscr{M}_{(X',\XX')}),\pi_{X*}(\widetilde{\psi}^{*}_{n}(\varepsilon'))).
\end{equation}

On the other hand, the $\mathscr{O}_{\XX_{n}}$-module associated to $\rmC^{*}(\mathscr{M})$ is $\rmC^{*}(\mathscr{M})_{(X,\XX)}$ and the associated $\mathcal{Q}_{\XX}$-stratification is induced by the isomorphisms $\widetilde{q}^{*}_{2,n} (\rmC^{*}(\mathscr{M})_{(X,\XX)})\xrightarrow{\sim} \rmC^{*}(\mathscr{M})_{(X,\QQ_{\XX})} \xleftarrow{\sim} \widetilde{q}^{*}_{1,n}(\rmC^{*}(\mathscr{M})_{(X,\XX)})$. The proposition follows in view of \eqref{double diag proj}, \eqref{eta F} and \eqref{psi pullback cal}.
\end{proof}
\begin{nothing} \label{Cartier transform shiho}
	Keep the assumption of \ref{prep model X lifting F}. 
	Recall that the morphisms of formal $\XX$-groupoids $\varpi':\TT_{\XX'}\to \RR_{\XX'}$ \eqref{pi T to R} and $\lambda:\PP_{\XX}\to \QQ_{\XX}$ \eqref{pi P to Q} induce functors
	\begin{eqnarray}
		\label{varpi pullback}	\varpi'^{*}_{n}:\Big\{ \txt{ \textnormal{$\mathscr{O}_{\XX'_{n}}$-modules}\\ \textnormal{with $\mathcal{R}_{\XX'}$-stratification}} \Big\} 
		\to \Big\{ \txt{ \textnormal{$\mathscr{O}_{\XX'_{n}}$-modules}\\	\textnormal{with $\mathcal{T}_{\XX'}$-stratification}} \Big\} \\
		\label{pi pullback} \lambda_{n}^{*}:\Big\{ \txt{ \textnormal{$\mathscr{O}_{\XX_{n}}$-modules}\\	\textnormal{with $\mathcal{Q}_{\XX}$-stratification}} \Big\} \to \Big\{ \txt{ \textnormal{$\mathscr{O}_{\XX_{n}}$-modules}\\	\textnormal{with $\mathcal{P}_{\XX}$-stratification}} \Big\}
	\end{eqnarray}
	By \ref{B-O stratification MIC}, \ref{qn pconnection T} and \ref{equi crystals stratification}, they further induce functors
	\begin{eqnarray}
		\mu: \CC(\mathscr{O}_{\pCRIS',n})&\to& \pMIC^{\qn}(\XX'_{n}/\SS_{n}),\label{functor mu}\\
		\nu: \CC(\mathscr{O}_{\CRIS,n})&\to& \MIC^{\qn}(\XX_{n}/\SS_{n}).\label{functor nu}	
	\end{eqnarray}
	By \eqref{square stra}, \ref{Shiho two construction} and \ref{Cartier global local}, the diagrams
		\begin{equation} \label{big diag}
		\xymatrix{
			\CC(\mathscr{O}_{\pCRIS',n}) \ar[r]^{\rmC^{*}} \ar[d]_{\wr} & \CC(\mathscr{O}_{\CRIS,n}) \ar[d]^{\wr}\\
		\Big\{ \txt{	\textnormal{$\mathscr{O}_{\XX'_{n}}$-modules}\\
				\textnormal{with $\mathcal{R}_{\XX'}$-stratification}} \Big\}
				\ar[r]^{\psi_{n}^{*}} \ar[d]_{\varpi'^{*}_{n}}&
		\Big\{	\txt{\textnormal{$\mathscr{O}_{\XX_{n}}$-modules}\\
		\textnormal{with $\mathcal{Q}_{\XX}$-stratification}} \Big\}  \ar[d]^{\lambda^{*}_{n}}\\
			\Big\{ \txt{	\textnormal{$\mathscr{O}_{\XX'_{n}}$-modules}\\
				\textnormal{with $\mathcal{T}_{\XX'}$-stratification}} \Big\}
				\ar[r]^{\varphi_{n}^{*}}  \ar[d]_{\wr}&
		\Big\{	\txt{\textnormal{$\mathscr{O}_{\XX_{n}}$-modules}\\
		\textnormal{with $\mathcal{P}_{\XX}$-stratification}} \Big\} \ar[d]^{\wr} \\
		\pMIC^{\qn}(\XX_{n}'/\SS_{n}) \ar[r]^{\Phi_{n}} & \MIC^{\qn}(\XX_{n}/\SS_{n}) }
	\end{equation}
	where the functors $\psi_{n}^{*}$, $\varphi_{n}^{*}$ and $\Phi_{n}$ are induced by $F$, are commutative up to a functorial isomorphism of $\MIC^{\qn}(\XX_{n}/\SS_{n})$. For every object $\mathscr{M}$ of $\mathscr{C}(\mathscr{O}_{\pCRIS',n})$, we have a functorial isomorphism 
	\begin{equation} \label{eta F MIC}
		\eta_{F}: \Phi_{n}(\mu(\mathscr{M})) \xrightarrow{\sim} \nu(\rmC^{*}(\mathscr{M})).
	\end{equation}
We see that the Cartier equivalence $\rmC^{*}$ \eqref{Cartier transform} is compatible with Shiho's functor $\Phi_{n}$ \eqref{Shiho equi}.
\end{nothing}

In the remainder of this section, we will explain how to relate Shiho's local constructions with respect to different liftings of Frobenius morphism using Cartier equivalence. 

Let $F_{1},F_{2}:\XX\to \XX'$ denote two liftings of the relative Frobenius morphism $F_{X/k}$ of $X$. 

\begin{lemma} \label{construction psi12 alpha}
	The morphisms $F_{1},F_{2}$ induce a morphism of $\pCRIS'$
	\begin{equation}
		\psi_{12}:\rho(X,\QQ_{\XX})\to (X',\RR_{\XX'}).
	\end{equation}
\end{lemma}
\begin{proof} The proof is similar to that of \ref{prop Q to R'}. By the universal property of $\RR_{\XX'}$, it suffices to show that there exists a unique $k$-morphism $g:\QQ_{\XX,1}\to X'$ which fits into a commutative diagram
\begin{equation} 
	\xymatrix{
		\QQ_{\XX,1}\ar[r] \ar[dd]_{g}& \QQ_{\XX} \ar[d]\\
		&\XX^{2} \ar[d]^{(F_{1},F_{2})} \\
		X'\ar[r]^{\Delta}& \XX'^{2}}
\end{equation}
The problem being local on $\XX$, we can assume that there exists an \'{e}tale $\SS$-morphism $\XX\to \widehat{\mathbb{A}}_{\SS}^{d}=\Spf(\rW\{T_{1},\cdots,T_{d}\})$. We put $t_{i}$ the image of $T_{i}$ in $\mathscr{O}_{\XX}$, $\xi_{i}=1\otimes t_{i}-t_{i}\otimes 1$, $t_{i}'=\pi^{*}(t_{i})\in \mathscr{O}_{\XX'}$ and $\xi_{i}'=1\otimes t'_{i}-t'_{i}\otimes 1$ for all $1\le i\le d$. Locally, there exists sections $a_{i},b_{i}$ of $\mathscr{O}_{\XX}$ such that $F_{1}^{*}(t_{i}')=t_{i}^{p}+p a_{i}, F_{2}^{*}(t_{i}')=t_{i}^{p}+p b_{i}$. 
By a similar calculation of \eqref{calcul F2 pullback}, we have
\begin{equation} \label{calcul F1F2 pullback}
	(F_{1},F_{2})^{*}(\xi_{i}')= \xi_{i}^{p}+\sum_{j=1}^{p-1}\binom{p}{j}\xi_{k}^{j}(t_{k}\otimes1)^{p-j}+p(1\otimes b_{k}-a_{k}\otimes 1).
\end{equation}
	Since $\xi_{i}^{p}=p\cdot\big(\frac{\xi_{i}^{p}}{p}\big)$ in $\mathcal{Q}_{\XX}$, the assertion follows.
\end{proof}
\begin{nothing} \label{alpha morphism local}
	Keep the assumption of \ref{construction psi12 alpha}, we denote by $\alpha$ the composition
	\begin{equation} \label{morphism alpha}
		\alpha:(X',\XX,F_{X/k})=\rho(X,\XX) \xrightarrow{\rho(\iota_{Q})} \rho(X,\QQ_{\XX}) \xrightarrow{\psi_{12}} (X',\RR_{\XX'}).
	\end{equation}
	We set $q'_{1}, q'_{2}:(X',\RR_{\XX'})\to (X',\XX')$ the canonical morphisms of $\pCRIS'$ and we have $q'_{i}\circ \alpha=F_{i}$ \eqref{F in E'}.

	Considering $\mathcal{R}_{\XX'}$ as a formal Hopf $\mathscr{O}_{\XX'}$-algebra of $\XX_{\zar}$, $\alpha$ induces a $\rW$-homomorphism of $\XX_{\zar}$:
	\begin{equation} \label{homomorphism a}
		a:\mathcal{R}_{\XX'}\to \mathscr{O}_{\XX}.
	\end{equation}
	Equipped with the left (resp. right) $\mathscr{O}_{\XX'}$-linear action on the source and the $\mathscr{O}_{\XX'}$-linear action induced by $F_{1}$ (resp. $F_{2}$) on the target, $a$ is $\mathscr{O}_{\XX'}$-linear.

	Suppose that there exists an \'etale $\SS$-morphism $\XX\to \widehat{\mathbb{A}}_{\SS}^{d}=\Spf(\rW\{T_{1},\cdots,$$T_{d}\})$ and we take again the notation of the proof of \ref{construction psi12 alpha}. 
	In view of \eqref{calcul F1F2 pullback}, the homomorphism $a$ is determined by
	\begin{equation}
		a\biggl(\frac{\xi'_{i}}{p}\biggr)=\frac{F_{2}^{*}(t'_{i})-F_{1}^{*}(t'_{i})}{p}\qquad \forall ~ 1\le i\le d. \label{homomorphism a formula}
	\end{equation}
\end{nothing}

\begin{nothing} \label{diff MIC to liftings}
	Let $\mathscr{M}$ be a crystal of $\mathscr{O}_{\pCRIS',n}$-modules, $\mathcal{M}=\mathscr{M}_{(X',\XX')}$ and $\varepsilon':\widetilde{q}'^{*}_{2,n}(\mathcal{M})\xrightarrow{\sim} \widetilde{q}'^{*}_{1,n}(\mathcal{M})$ the associated $\mathcal{R}_{\XX'}$-stratification \eqref{equi crystals stratification}. Recall \eqref{Cartier global local} that the morphisms $F_{1}$ and $F_{2}$ induce morphisms $\psi_{1},\psi_{2}:\rho(X,\QQ_{\XX})\to (X',\RR_{\XX'})$ of $\pCRIS'$ respectively. We associate to $(\mathcal{M},\varepsilon')$ two different $\mathscr{O}_{\XX_{n}}$-modules with $\mathcal{Q}_{\XX}$-stratification 
	\begin{equation}
		(F_{1,n}^{*}(\mathcal{M}),\pi_{X*}(\widetilde{\psi}_{1,n}^{*}(\varepsilon')))\quad \textnormal{and} \quad(F_{2,n}^{*}(\mathcal{M}),\pi_{X*}(\widetilde{\psi}_{2,n}^{*}(\varepsilon'))).
	\end{equation}
	
	Let $\Phi_{1,n}$ (resp. $\Phi_{2,n}$) be Shiho's functor induced by $F_{1}$ (resp. $F_{2}$) \eqref{Shiho equi}. We associate to $\mathscr{M}$ two different objects of $\MIC^{\qn}(\XX_{n}/\SS_{n})$ \eqref{big diag}
	\begin{equation} \label{two MIC shiho}
		\Phi_{1,n}(\mu(\mathscr{M})) \quad \textnormal{and} \quad \Phi_{2,n}(\mu(\mathscr{M})),
	\end{equation}
	whose underlying $\mathscr{O}_{\XX_{n}}$-modules are $F_{1,n}^{*}(\mathcal{M})$ and $F_{2,n}^{*}(\mathcal{M})$ respectively.

	The morphism $\alpha$ gives a natural way to glue $\Phi_{1,n}(\mu(\mathscr{M})),\Phi_{2,n}(\mu(\mathscr{M}))$. 
\end{nothing}

\begin{prop} \label{alpha iso stra}
	Keep the assumption of \ref{diff MIC to liftings}. 
	The morphism $\alpha$ and the $\mathcal{R}_{\XX'}$-stratification $\varepsilon'$ induce an isomorphism of $\mathscr{O}_{\XX_{n}}$-modules with $\mathcal{Q}_{\XX}$-stratification:
	\begin{equation}
		\alpha^{*}(\varepsilon'):(F_{2,n}^{*}(\mathcal{M}),\pi_{X*}(\widetilde{\psi}_{2,n}^{*}(\varepsilon')))\xrightarrow{\sim} (F_{1,n}^{*}(\mathcal{M}),\pi_{X*}(\widetilde{\psi}_{1,n}^{*}(\varepsilon'))) \label{alpha connect Fi}
	\end{equation}
such that $\eta_{F_{2}}=\eta_{F_{1}}\circ \alpha^{*}(\varepsilon')$ \eqref{eta F}. 
\end{prop}
\begin{proof} We denote by $q_{1},q_{2}:(X,\QQ_{\XX})\to (X,\XX)$ the canonical morphisms of $\CRIS$, by $q'_{1}, q'_{2}:(X',\RR_{\XX'})\to (X',\XX')$ the canonical morphisms of $\pCRIS'$ and by $q'_{ij}: (X',\RR_{\XX'}(2))\to (X',\RR_{\XX'})$ the morphism induced by $\RR_{\XX'}(2)\to \XX'^{3}\xrightarrow{p'_{ij}} \XX'^{2}$ and the universal property of $\RR_{\XX}$ for all $1\le i<j\le 3$ (cf. \eqref{diag univ alpha R}).

Since $F_{i}=q'_{i}\circ \alpha$ for $i=1,2$, the isomorphism $\varepsilon':\widetilde{q}'^{*}_{2,n}(\mathcal{M})\xrightarrow{\sim} \widetilde{q}'^{*}_{1,n}(\mathcal{M})$ induces an isomorphism
\begin{equation}
	\pi_{X*}(\widetilde{\alpha}_{n}^{*}(\varepsilon')):\widetilde{F}_{2,n}^{*}(\mathcal{M})\xrightarrow{\sim} \widetilde{F}_{1,n}^{*}(\mathcal{M}).
\end{equation}
We write simply $\alpha^{*}(\varepsilon')$ for $\pi_{X*}(\widetilde{\alpha}_{n}^{*}(\varepsilon'))$. Then we have $\eta_{F_{2}}=\eta_{F_{1}}\circ \alpha^{*}(\varepsilon')$. It remains to show that the $\alpha^{*}(\varepsilon')$ is compatible with the $\mathcal{Q}_{\XX}$-stratifications on both sides. 

	The compositions of morphisms
	\begin{eqnarray*}
		&\xymatrix{\XX^{2} \ar[r]^{p_{1}}& \XX \ar[r]^{\Delta} & \XX^{2} \ar[r]^{(F_{1},F_{2})}& \XX'^{2}\ar[r]^{p_{2}'}&\XX'}&\\
		&\xymatrix{\XX^{2} \ar[r]^{(F_{2},F_{2})} &\XX'^{2} \ar[r]^{p_{1}'} &\XX'} &
	\end{eqnarray*}
	are equal. We deduce that the compositions of morphisms of $\pCRIS'$
	\begin{eqnarray*}
		&&\xymatrix{ \rho(X,\QQ_{\XX})\ar[r]^{\rho(q_{1})}& \rho(X,\XX) \ar[r]^{\alpha}& (X',\RR_{\XX'}) \ar[r]^{q_{2}'}& (X',\XX')}\\
		&&\xymatrix{\rho(X,\QQ_{\XX}) \ar[r]^{\psi_{2}}& (X',\RR_{\XX'}) \ar[r]^{q_{1}'}& (X',\XX')}
	\end{eqnarray*}
	are equal. In view of the isomorphism $(X',\RR_{\XX'}(2))\simeq (X',\RR_{\XX'})\times_{(X',\XX')}(X',\RR_{\XX'})$ of $\pCRIS'$ \eqref{iso R2 R times R}, we obtain a morphism of $\pCRIS'$
\begin{equation}
	u: \rho(X,\QQ_{\XX}) \to (X',\RR_{\XX'}(2))
\end{equation}
such that $\alpha\circ \rho(q_{1})=q_{12}'\circ u$ and $\psi_{2}=q_{23}'\circ u$. Symmetrically, we construct a morphism of $\pCRIS'$
\begin{equation}
	v: \rho(X,\QQ_{\XX}) \to (X',\RR_{\XX'}(2))
\end{equation}
such that $\psi_{1}=q_{12}'\circ v$ and $\alpha\circ \rho(q_{2})=q_{23}'\circ v$. The compositions 
\begin{eqnarray*}
	\XX^{2}\xrightarrow{(\iota_{1},\iota_{1},\iota_{2})} \XX^{3} \xrightarrow{(F_{1},F_{2},F_{2})} \XX'^{3} \xrightarrow{p_{13}'} \XX'^{2} \\
	\XX^{2}\xrightarrow{(\iota_{1},\iota_{2},\iota_{2})} \XX^{3} \xrightarrow{(F_{1},F_{1},F_{2})} \XX'^{3} \xrightarrow{p_{13}'} \XX'^{2}
\end{eqnarray*}
are equal to $(F_{1},F_{2}):\XX^{2}\to \XX'^{2}$. By the universal property of $\RR_{\XX'}$ (\ref{prop univ RQ}), we deduce that $q_{13}'\circ u= q_{13}'\circ v=\psi_{12}$ \eqref{construction psi12 alpha}, i.e. the compositions 
\begin{equation}
	\xymatrix{\rho(X,\QQ_{\XX}) \ar@<0.5ex>[r]^-{u} \ar@<-0.5ex>[r]_-{v}& (X',\RR_{\XX'}(2)) \ar[r]^-{q_{13}'}& (X',\RR_{\XX'})}
\end{equation}
are equal to $\psi_{12}$. By \eqref{double diag proj}, we have $\widetilde{\psi}_{i,n}^{*}(\widetilde{q}'^{*}_{j,n}(\mathcal{M}))\simeq \widetilde{q}_{j,n}^{*}(\widetilde{F}_{i,n}^{*}(\mathcal{M}))$ for all $i,j=1,2$. By the cocycle condition $\widetilde{q}'^{*}_{12,n}(\varepsilon')\circ \widetilde{q}'^{*}_{23,n}(\varepsilon')=\widetilde{q}'^{*}_{13,n}(\varepsilon')$, we deduce a commutative diagram:
\begin{equation}
	\xymatrixcolsep{5pc}\xymatrix{
		\widetilde{\psi}_{2,n}^{*}(\widetilde{q}'^{*}_{2,n}(\mathcal{M})) \ar[r]^{\widetilde{\psi}_{2,n}^{*}(\varepsilon')} \ar[d]_{\widetilde{q}_{2,n}^{*}(\widetilde{\alpha}^{*}_{n}(\varepsilon'))} & \widetilde{\psi}_{2,n}^{*}(\widetilde{q}'^{*}_{1,n}(\mathcal{M})) \ar[d]^{\widetilde{q}_{1,n}^{*}(\widetilde{\alpha}_{n}^{*}(\varepsilon'))} \\
		\widetilde{\psi}_{1,n}^{*}(\widetilde{q}'^{*}_{2,n}(\mathcal{M})) \ar[r]^{\widetilde{\psi}_{1,n}^{*}(\varepsilon')} & \widetilde{\psi}_{1,n}^{*}(\widetilde{q}'^{*}_{1,n}(\mathcal{M}))}
\end{equation}
That is the isomorphism $\alpha^{*}(\varepsilon')$ is compatible with the $\mathcal{Q}_{\XX}$-stratifications $\pi_{X*}(\widetilde{\psi}_{2,n}^{*}(\varepsilon'))$ and $\pi_{X*}(\widetilde{\psi}_{1,n}^{*}(\varepsilon'))$.
\end{proof}
\begin{coro} \label{alpha iso conn}
	Keep the assumption of \ref{diff MIC to liftings}. The isomorphism \eqref{alpha connect Fi} induces an isomorphism of $\MIC^{\qn}(\XX_{n}/\SS_{n})$ \eqref{two MIC shiho}
	\begin{equation}\label{alpha iso MIC}
		\alpha^{*}(\varepsilon'): \Phi_{2,n}(\mu(\mathscr{M})) \xrightarrow{\sim} \Phi_{1,n}(\mu(\mathscr{M}))
	\end{equation}
	such that $\eta_{F_{2}}=\eta_{F_{1}}\circ \alpha^{*}(\varepsilon')$ \eqref{eta F MIC}.
\end{coro}

It follows from \eqref{big diag} and \ref{alpha iso stra}.

\begin{rem}
	Given a lifting of the Frobenius morphism, Shiho's functor \eqref{Shiho Phi} applies to $\mathscr{O}_{\XX'_{n}}$-modules with integrable $p$-connection. However, the isomorphisme \eqref{alpha iso MIC}, which glues local constructions of Shiho, depends on the $\mathcal{R}_{\XX'}$-stratification on the $\mathscr{O}_{\XX'_{n}}$-module $\mathscr{M}_{(X',\XX')}$.
\end{rem}

\section{Cartier transform of Ogus--Vologodsky} \label{Cartier OV}
For the convenience of the reader, we shall review the original construction of the Cartier transform of Ogus--Vologodsky \cite{OV07}. 
In particular, we shall clarify some details, especially in regard to sheaf of affine functions on a torsor. 

In this section, $X$ denotes a scheme. Starting from \ref{crystal of affine functions}, we will suppose that $X$ is smooth over $k$.
\begin{nothing} \label{sym alg general not}
	Let $\wT$ a locally free $\mathscr{O}_{X}$-module of finite type. We denote by $\rS(\wT)$ (resp. $\Gamma(\wT)$) the symmetric algebra (resp. PD-algebra) of $\wT$ over $\mathscr{O}_{X}$ (\cite{Il71} I 4.2.2.6) and for any integer $n\ge 0$, by $\rS^{n}(\wT)$ (resp. $\Gamma_{n}(\wT)$) its homogeneous part of degree $n$. 
	There exists a unique homomorphism of $\mathscr{O}_{X}$-algebras
	\begin{equation} \label{delta ST}
		\delta:\rS(\wT)\to \rS(\wT)\otimes_{\mathscr{O}_{X}}\rS(\wT) 
	\end{equation}
	such that for every local section $e$ of $\wT$, we have $\delta(e)=1\otimes e+e\otimes 1$. This homomorphisms makes $\rS(\wT)$ into a Hopf commutative $\mathscr{O}_{X}$-algebra.

	Let $I$ (resp. $J$) be the ideal $\oplus_{n\ge 1}\rS^{n}(\wT)$ of $\rS(\wT)$ (resp. PD-ideal $\oplus_{n\ge 1}\Gamma_{n}(\wT)$ of $\Gamma(\wT)$). 	We denote by $\widehat{\rS}(\wT)$ (resp. $\widehat{\Gamma}(\wT)$) the completion of $\rS(\wT)$ (resp. $\Gamma(\wT)$) with respect to the filtration $\{I^{n}\}_{n\ge 1}$ (resp. PD-filtration $\{J^{[n]}\}_{n\ge 1}$).

	Let $M$ be a $\widehat{\Gamma}(\wT)$-module. 
	We say that $M$ is \textit{quasi-nilpotent} if for any open subscheme $U$ of $X$ and any $e\in M(U)$, there exists a Zariski covering $\{U_{i}\to U\}_{i\in I}$ and a family of integers $\{N_{i}\}_{i\in I}$ such that for each $i\in I$, $e|_{U_{i}}$ is annihilated by the ideal $J^{[N_{i}]}$. For any integer $n\ge 0$, we say that $M$ is \textit{nilpotent of level $\le n$} if $E$ is annihilated by $J^{[n+1]}$. 
\end{nothing}

\begin{nothing} \label{action Gamma T on SO}
	We set $\Omega=\FHom_{\mathscr{O}_{X}}(\wT,\mathscr{O}_{X})$. The pairing $\wT\otimes_{\mathscr{O}_{X}}\Omega\to \mathscr{O}_{X}$ induces a canonical morphism 
	\begin{equation} \label{pairing T Om}
		\Gamma_{n}(\wT)\otimes_{\mathscr{O}_{X}}\rS^{m+n}(\Omega)\to \rS^{m}(\Omega)
	\end{equation}
	which is perfect if $m=0$ (\cite{BO} A.10). If we equip $\Gamma(\wT)$ with the topology defined by the PD-filtration $\{J^{[n]}\}$ and $\rS(\Omega)$ with the discrete topology, the $\mathscr{O}_{X}$-linear morphism $\Gamma(\wT)\otimes_{\mathscr{O}_{X}} \rS(\Omega)\to \rS(\Omega)$ is continuous. It extends by continuity to an action of $\widehat{\Gamma}(\wT)$ on $\rS(\Omega)$:
	\begin{equation}
		\widehat{\Gamma}(\wT)\otimes_{\mathscr{O}_{X}} \rS(\Omega)\to \rS(\Omega).
	\end{equation}
	We have an increasing exhaustive filtration of $\widehat{\Gamma}(\wT)$-submodules $\{\oplus_{m\le n} \rS^{m}(\Omega)\}_{n\ge 0}$ of $\rS(\Omega)$ such that $\oplus_{m\le n} \rS^{m}(\Omega)$ is nilpotent of level $\le n$. 
	Then $\rS(\Omega)$ is quasi-nilpotent. The above morphism induces an $\mathscr{O}_{X}$-linear isomorphism
	\begin{eqnarray} \label{Gamma S dual}
		\widehat{\Gamma}(\wT)&=&\varprojlim_{n} (\oplus_{m\le n}\Gamma_{m}(\wT))\\
		&\simeq& \varprojlim_{n} \FHom_{\mathscr{O}_{X}}(\oplus_{m\le n} \rS^{m}(\Omega),\mathscr{O}_{X}) \nonumber \\
		&=& \FHom_{\mathscr{O}_{X}}(\rS(\Omega),\mathscr{O}_{X}). \nonumber
	\end{eqnarray}
	We equip $\FHom_{\mathscr{O}_{X}}(\rS(\Omega),\mathscr{O}_{X})$ with the $\mathscr{O}_{X}$-algebra structure induced by the Hopf algebra $\rS(\Omega)$ \eqref{Hopf alg dual}. The above isomorphism is an isomorphism of $\mathscr{O}_{X}$-algebras.
\end{nothing}

\begin{nothing}
	Let $f:Y\to X$ be a morphism of schemes. We put $\wT_{Y}=f^{*}(\wT)$ and $\Omega_{Y}=f^{*}(\Omega)$. By the universal property of the symmetric algebra, we have a canonical isomorphism of $\mathscr{O}_{Y}$-algebras
	\begin{equation}
		f^{*}(\rS_{\mathscr{O}_{X}}(\Omega))\xrightarrow{\sim} \rS_{\mathscr{O}_{Y}}(\Omega_{Y}).
	\end{equation}
	Since $\rS_{\mathscr{O}_{X}}(\Omega)$ is a direct sum of locally free $\mathscr{O}_{X}$-modules of finite type, by duality \eqref{Gamma S dual}, we deduce a canonical isomorphism of $\mathscr{O}_{Y}$-algebras
	\begin{eqnarray} \label{base change Gamma T}	
		\widehat{\Gamma}_{\mathscr{O}_{Y}}(\wT_{Y})&\simeq & \FHom_{\mathscr{O}_{Y}}(\rS_{\mathscr{O}_{Y}}(\Omega_{Y}),\mathscr{O}_{Y})\\
		&\simeq& f^{*}(\FHom_{\mathscr{O}_{X}}(\rS_{\mathscr{O}_{X}}(\Omega),\mathscr{O}_{X})) \nonumber\\
		&\simeq& f^{*}(\widehat{\Gamma}_{\mathscr{O}_{X}}(\wT)). \nonumber
	\end{eqnarray}
\end{nothing}

\begin{nothing} \label{affine functions}
	Let $\mathscr{L}$ be an $\wT$-torsor of $X_{\zar}$. \textit{An affine function on $\mathscr{L}$} is a morphism $f:\mathscr{L}\to \mathscr{O}_{X}$ of $X_{\zar}$ satisfying the following equivalent conditions (\cite{AGT} II.4.7):
	\begin{itemize}
		\item[(i)] For every open subscheme $U$ of $X$ and every $s\in \mathscr{L}(U)$, the morphism:
			\begin{equation}
				\wT(U)\to \mathscr{O}(U), \qquad t\mapsto f(s+t)-f(s)
			\end{equation}
			is $\mathscr{O}_{X}(U)$-linear.

		\item[(ii)] There exists a section $\omega_{f}\in\Omega(X)$, called the \textit{linear term of $f$}, such that for every open subscheme $U$ of $X$ and all $s\in \mathscr{L}(U)$ and $t\in \wT(U)$, we have
		\begin{equation} \label{affine functions linear term}
			f(s+t)=f(s)+\omega_{f}(t).
		\end{equation}
\end{itemize}

The condition (i) is clearly local for the Zariski topology on $X$. We denote by $\mathscr{F}$ the subsheaf of $\FHom_{X_{\zar}}(\mathscr{L},\mathscr{O}_{X})$ consisting of affine functions on $\mathscr{L}$; in other words, for any open subscheme $U$ of $X$, $\mathscr{F}(U)$ is the set of affine functions on $\mathscr{L}|_{U}$. It is naturally endowed with an $\mathscr{O}_{X}$-module structure. We call $\mathscr{F}$ \textit{the sheaf of affine functions on $\mathscr{L}$}.

	We have a canonical $\mathscr{O}_{X}$-linear morphism $c:\mathscr{O}_{X}\to \mathscr{F}$ whose image consists of constant functions. The ``linear term'' defines an $\mathscr{O}_{X}$-linear morphism $\omega:\mathscr{F}\to \Omega$. One verifies that the sequence
	\begin{equation} \label{exact seq affine functions}
		0\to \mathscr{O}_{X}\xrightarrow{c} \mathscr{F} \xrightarrow{\omega} \Omega\to 0
	\end{equation}
	is exact. By (\cite{Il71} I 4.3.1.7), the above sequence induces, for any integer $n\ge 1$, an exact sequence:
	\begin{equation}
		0\to \rS^{n-1}(\mathscr{F})\to \rS^{n}(\mathscr{F})\to \rS^{n}(\Omega)\to 0.
	\end{equation}
	The $\mathscr{O}_{X}$-modules $(\rS^{n}(\mathscr{F}))_{n\ge 0}$ form a inductive system. We denote its inductive limit by
	\begin{equation}
		\mathscr{A}=\varinjlim_{n\ge 0} \rS^{n}(\mathscr{F})
	\end{equation}
	which is naturally endowed with a structure of an $\mathscr{O}_{X}$-algebra. For any integer $n\ge 0$, the canonical morphism $\rS^{n}(\mathscr{F})\to \mathscr{A}$ is injective. By letting $N_{n}(\mathscr{A})=\rS^{n}(\mathscr{F})$ for all $n\ge 0$, we obtain an increasing exhaustive filtration of $\mathscr{A}$.

	There exists a unique homomorphism of $\mathscr{O}_{X}$-algebras
	\begin{equation} \label{mu A to SOA}
		\mu:\mathscr{A}\to \rS(\Omega)\otimes_{\mathscr{O}_{X}}\mathscr{A}
	\end{equation}
	such that for every local section $m$ of $\mathscr{F}$, we have	$\mu(m)=1\otimes m+ \omega(m)\otimes 1$. 
	For $n\ge 0$, we have
	\begin{equation} \label{mu filtration}
		\mu(N_{n}(\mathscr{A})) \subset \oplus_{i+j=n}\rS^{i}(\Omega)\otimes N_{j}(\mathscr{A}).
	\end{equation}
	By construction, the following diagram is commutative
	\begin{equation} \label{diag delta mu}
		\xymatrix{
			\mathscr{A}\ar[r]^-{\mu} \ar[d]_{\mu} & \rS(\Omega)\otimes_{\mathscr{O}_{X}}\mathscr{A} \ar[d]^{\delta\otimes \id}\\
			\rS(\Omega)\otimes_{\mathscr{O}_{X}}\mathscr{A} \ar[r]^-{\id\otimes\mu}& \rS(\Omega)\otimes_{\mathscr{O}_{X}}\rS(\Omega)\otimes_{\mathscr{O}_{X}}\mathscr{A} }
	\end{equation}
\end{nothing}

\begin{nothing} \label{Gamma action on A}
	By \eqref{Gamma S dual} and \eqref{mu A to SOA}, we have an $\mathscr{O}_{X}$-linear morphism:
	\begin{eqnarray} \label{Gamma T action on A}
		\widehat{\Gamma}(\wT)\otimes_{\mathscr{O}_{X}}\mathscr{A}&\to& \mathscr{A}\\
		u\otimes a &\mapsto& (u\otimes\id) (\mu(a)). \nonumber
	\end{eqnarray}
	By \eqref{diag delta mu}, the above morphism makes $\mathscr{A}$ into a $\widehat{\Gamma}(\wT)$-module. The action of $\wT$ on $\mathscr{F}$ is given by $\omega$ \eqref{exact seq affine functions} and duality.  
	By \eqref{mu filtration}, we see that $\mathscr{A}$ is quasi-nilpotent and that for any $n\ge 0$, $N_{n}(\mathscr{A})$ is a $\widehat{\Gamma}(\wT)$-submodule of $\mathscr{A}$ and is nilpotent of level $\le n$. 
	
	The canonical $\rS(\Omega)$-linear isomorphism $\Omega_{\rS(\Omega)/\mathscr{O}_{X}}^{1}\xrightarrow{\sim} \Omega\otimes_{\mathscr{O}_{X}} \rS(\Omega)$ induces an isomorphism
	\begin{equation}
		\Omega^{1}_{\mathscr{A}/\mathscr{O}_{X}}\xrightarrow{\sim} \Omega\otimes_{\mathscr{O}_{X}}\mathscr{A}.
	\end{equation}
	We denote the universal $\mathscr{O}_{X}$-derivation by 
	\begin{equation}\label{univer derivation}
		d_{\mathscr{A}}:\mathscr{A}\to \Omega\otimes_{\mathscr{O}_{X}}\mathscr{A}.
	\end{equation}
	For any local section $m$ of $\mathscr{F}$, we have $d_{\mathscr{A}}(m)=\omega(m)\otimes 1$. 
\end{nothing}	
\begin{nothing} \label{action on A splitting}
	Let $s\in \mathscr{L}(X)$ and let $\rho_{s}:\mathscr{F}\to \mathscr{O}_{X}$ be the associated splitting of the exact sequence \eqref{exact seq affine functions}. The morphism $\Omega\to \mathscr{F}$ deduced from $\id-c\circ \rho_{s}$ extends to an isomorphism of $\mathscr{O}_{X}$-algebras
	\begin{equation}
		\psi:\rS(\Omega)\xrightarrow{\sim} \mathscr{A}
	\end{equation}
	which is compatible with the filtrations $(\oplus_{0\le i\le n}\rS^{i}(\Omega))_{n}$ and $(N_{n}(\mathscr{A}))_{n}$. 
	The diagram \eqref{delta ST}
	\begin{equation}
		\xymatrix{
			\rS(\Omega) \ar[r]^{\psi} \ar[d]_{\delta}& \mathscr{A} \ar[d]^{\mu}\\
			\rS(\Omega)\otimes_{\mathscr{O}_{X}}\rS(\Omega) \ar[r]^{\id \otimes \psi} & \rS(\Omega)\otimes_{\mathscr{O}_{X}}\mathscr{A}
		}
	\end{equation}
	is commutative. Hence the isomorphism $\psi$ is compatible with the $\widehat{\Gamma}(\wT)$-module structures \eqref{action Gamma T on SO}.
\end{nothing}

\begin{nothing} \label{affine inverse image num}
	Let $f:Y\to X$ be a morphism of schemes and $\mathscr{L}$ an $\wT$-torsor of $X_{\zar}$. For $\mathscr{O}_{X}$-modules, we will use the notation $f^{-1}$ to denote the inverse image in the sense of abelian sheaves and will keep the notation $f^{*}$ for the inverse image in the sense of modules. 
	The \textit{affine inverse image of $\mathscr{L}$ under $f$}, denoted by $f^{+}(\mathscr{L})$, is the $f^{*}(\wT)$-torsor of $Y_{\zar}$ deduced from the $f^{-1}(\wT)$-torsor $f^{*}(\mathscr{L})$ by extending its structural group by the canonical homomorphism $f^{-1}(\wT)\to f^{*}(\wT)$,
	\begin{equation} \label{affine inverse image}
		f^{+}(\mathscr{L})=f^{*}(\mathscr{L})\wedge^{f^{-1}(\wT)}f^{*}(\wT);
	\end{equation}
in other words, the quotient of $f^{*}(\mathscr{L})\times f^{*}(\wT)$ by the diagonal action of $f^{-1}(\wT)$ (\cite{Gi} III 1.4.6). 

We denote by $\mathscr{F}$ the sheaf of affine functions on $\mathscr{L}$ \eqref{affine functions} and by $\mathscr{F}^{+}$ the sheaf of affine functions on $f^{+}(\mathscr{L})$.
Let $l:\mathscr{L}\to \mathscr{O}_{X}$ be an affine morphism, $\omega\in \Omega(X)$ its linear term and $\omega'=f^{*}(\omega)\in f^{*}(\Omega)(Y)$. Endowing $\mathscr{O}_{Y}$ with the structure of $f^{*}(\wT)$-object defined by $\omega'$ (\cite{AGT} II.4.8), there exists a unique $f^{*}(\wT)$-equivariant morphism $l':f^{+}(\mathscr{L})\to \mathscr{O}_{Y}$ that fits into the commutative diagram
\begin{equation} \label{diag l inverse affine}
	\xymatrix{
		f^{*}(\mathscr{L}) \ar[r]^{l} \ar[d]& f^{-1}(\mathscr{O}_{X})\ar[d] \\
		f^{+}(\mathscr{L}) \ar[r]^{l'} & \mathscr{O}_{Y} }
\end{equation}
where the vertical arrows are the canonical morphisms (\cite{Gi} III 1.3.6). 
The morphism $l'$ is therefore affine, with linear term $\omega'$. 
The resulting correspondence $l\mapsto l'$ induces an $\mathscr{O}_{X}$-linear morphism
\begin{equation}
	\lambda_{\sharp}:\mathscr{F}\to f_{*}(\mathscr{F}^{+}).
\end{equation}
Its adjoint morphism is an $\mathscr{O}_{Y}$-linear isomorphism (\cite{AGT} II.4.13.4)
	\begin{equation} \label{iso pullback tosor sheaf of affine funs}
		\lambda:f^{*}(\mathscr{F})\xrightarrow{\sim} \mathscr{F}^{+}
	\end{equation}
	which fits into a commutative diagram
	\begin{equation} \label{diag extension compatible}
		\xymatrix{
			0\ar[r] & \mathscr{O}_{Y} \ar[d] \ar[r]& f^{*}(\mathscr{F}) \ar[r] \ar[d]_{\lambda} & f^{*}(\Omega) \ar[d] \ar[r] & 0\\
			0\ar[r] & \mathscr{O}_{Y} \ar[r] & \mathscr{F}^{+} \ar[r] & f^{*}(\Omega) \ar[r] & 0 }
	\end{equation}
	In particular, the isomorphism $\lambda$ is compatible with actions of $f^{*}(\widehat{\Gamma}_{\mathscr{O}_{X}}(E))\simeq \widehat{\Gamma}_{\mathscr{O}_{Y}}(f^{*}(E))$ \eqref{base change Gamma T}.
\end{nothing}


\begin{nothing} \label{crystal of affine functions}
	In the remainder of this section, $X$ denotes a \textit{smooth scheme over $k$}. 
	We denote by $\Cris(X/k)$ the crystalline site of $X$ over $k$ equipped with the PD-ideal $0$, by $(X/k)_{\cris}$ the crystalline topos of $X$ over $k$ and by $\mathscr{O}_{X/k}$ the structure ring of $(X/k)_{\cris}$ defined for every object $(U,T)$ of $\Cris(X/k)$, by $(U,T)\mapsto \Gamma(T,\mathscr{O}_{T})$.

	Let $E$ be a crystal of locally free $\mathscr{O}_{X/k}$-modules of finite type on $\Cris(X/k)$ and $\mathscr{L}$ an $E$-torsor of $(X/k)_{\cris}$. For any object $(U,T)$ of $\Cris(X/k)$, $\mathscr{L}_{(U,T)}$ is an $E_{(U,T)}$-torsor of $T_{\zar}$ (\cite{Ber} III 3.5.1). We define $\mathscr{F}_{(U,T)}$ to be the sheaf of affine functions on $\mathscr{L}_{(U,T)}$ of $T_{\zar}$ \eqref{affine functions}. 
	
	Let $g:(U_{1},T_{1})\to (U_{2},T_{2})$ be a morphism of $\Cris(X/k)$ and $\jmath_{g}:|U_{1}|(=|T_{1}|)\to |U_2|(=|T_2|)$ the morphism of underlying topological spaces. The transition morphism of $\mathscr{L}$ associated to $g$ (\cite{BO} 5.1)
	\begin{equation} \label{cg torsor}
		c_{g}:\jmath^{-1}_{g}(\mathscr{L}_{(U_{2},T_{2})})\to \mathscr{L}_{(U_{1},T_{1})}
	\end{equation}
	is $\jmath^{-1}_{g}(E_{(U_{2},T_{2})})$-equivariant. By (\cite{Gi} III 1.4.6(iii)), we obtain an $E_{(U_{1},T_{1})}$-equivariant isomorphism
	\begin{equation}
		\jmath_{g}^{+}(\mathscr{L}_{(U_{2},T_{2})})\xrightarrow{\sim} \mathscr{L}_{(U_{1},T_{1})}.
	\end{equation}
	By \eqref{iso pullback tosor sheaf of affine funs}, we deduce an $\mathscr{O}_{T_{1}}$-linear isomorphism (\cite{AGT} II 4.14.2)
	\begin{equation}
		\gamma_{g}:g^{*}(\mathscr{F}_{(U_{2},T_{2})})\xrightarrow{\sim} \mathscr{F}_{(U_{1},T_{1})}.
	\end{equation}
	In view of the compatibility conditions of $c_{g}$ \eqref{cg torsor} and (\cite{AGT} II 4.15), the data $\{\mathscr{F}_{(U,T)},\gamma_{g}\}$ satisfy the compatibility conditions of (\cite{BO} 5.1). Hence, they define a crystal of $\mathscr{O}_{X/k}$-modules that we denote by $\mathscr{F}$ and call \textit{the crystal of affine functions on $\mathscr{L}$}.

	We denote by $\mathscr{A}$ the crystal of $\mathscr{O}_{X/k}$-algebras $\varinjlim_{n\ge 0} \rS^n(\mathscr{F})$. 
	It admits a increasing filtration of crystals of $\mathscr{O}_{X/k}$-modules $N_{n}(\mathscr{A})=\rS^{n}(\mathscr{F})$. 
\end{nothing}

\begin{nothing} \label{lemma can connection crystal}
	Let $(U,T,\delta)$ be an object of $\Cris(X/k)$ and $J_{T}$ the PD-ideal associated to the closed immersion $i:U\to T$. For any local section $a$ of $J_{T}$, we have $a^{p}=0$ and hence an isomorphism $U\xrightarrow{\sim} \underline{T}$. We denote by $\varphi_{T/k}$ the composition \eqref{notations Yk}
	\begin{equation} \label{varphi Tk}
		\varphi_{T/k}:\xymatrix{ T\ar[r]^{f_{T/k}}& U' \ar[r]& X'. }
	\end{equation}
	The morphism $\varphi_{X/k}$ is equal to the relative Frobenius morphism $F_{X/k}$. If $g:(U_{1},T_{1})\to (U_{2},T_{2})$ is a morphism of $\Cris(X/k)$, then $\varphi_{T_{2}/k}\circ g=\varphi_{T_{1}/k}$.

	Let $M'$ be an $\mathscr{O}_{X'}$-module. There exists a canonical isomorphism
	\begin{equation}
		\widetilde{c}_{g}:g^{*}(\varphi_{T_{2}/k}^{*}(M'))\xrightarrow{\sim} \varphi_{T_{1}/k}^{*}(M').
	\end{equation}
	The data $\{\varphi_{T/k}^{*}(M'),\widetilde{c}_{g}\}$ defines a crystal of $\mathscr{O}_{X/k}$-modules that we denote by $\mathscr{M}$. 
	Then the $\mathscr{O}_{X}$-module with integrable connection associated to $\mathscr{M}$ \textnormal{(\cite{BO} 6.8)} is $F_{X/k}^{*}(M')$ and the Frobenius descent connection $\nabla_{\can}$ \eqref{can connection} (\cite{OV07} 1.1). 
\end{nothing}

\begin{nothing} \label{torsor of liftings}
	We denote by $\Cris(X/\rW_{2})$ the crystalline site of $X$ over $\rW_{2}$ equipped with the PD-structure $\gamma_{2}$ \eqref{PD envelop P}. 
	Let $(U,\widetilde{T},\delta)$ be an object of $\Cris(X/\rW_{2})$ and $T$ the reduction modulo $p$ of $\widetilde{T}$. Since $\delta$ and $\gamma_{2}$ are compatible (\cite{BO} 3.16), $(U,\widetilde{T},\delta)$ induces an object $(U,T,\bar{\delta})$ of $\Cris(X/k)$.

	Let $\rT_{X/k}$ (resp. $\rT_{X'/k}$) be the $\mathscr{O}_{X}$-dual of $\Omega_{X/k}^{1}$ (resp. $\Omega_{X'/k}^{1}$). 
	Suppose we are given a smooth lifting $\widetilde{X}'$ of $X'/k$ over $\rW_{2}$. Let $(U,\widetilde{T},\delta)$ be an object of $\Cris(X/\rW_{2})$ such that $\widetilde{T}$ is flat over $\rW_{2}$, $T$ the reduction modulo $p$ of $\widetilde{T}$ and $V$ an open subscheme of $T$. Then $(U,T)$ is an object of $\Cris(X/k)$. We denote by $\widetilde{V}$ the open subscheme of $\widetilde{T}$ associated to $V$ and by $\mathscr{L}_{\widetilde{X}',\varphi_{T/k}}(V)$ the set of $\rW_{2}$-morphisms $\widetilde{V}\to \widetilde{X}'$ which make the following diagram commute \eqref{lemma can connection crystal}
	\begin{equation}
		\xymatrix{
			V\ar[r] \ar[d]_{\varphi_{T/k}|_{V}} & \widetilde{V} \ar[d] \\
			X' \ar[r] & \widetilde{X}'
		}
	\end{equation}
	The functor $V\mapsto \mathscr{L}_{\widetilde{X}',\varphi_{T/k}}(V)$ defines a sheaf for the Zariski topology on $T$. The sheaf $\mathscr{L}_{\widetilde{X}',\varphi_{T/k}}$ is a torsor under the $\mathscr{O}_{T}$-module $\FHom_{\mathscr{O}_{T}}(f^{*}_{T/k}(\Omega_{X'/k}^{1}),p\mathscr{O}_{\widetilde{T}})\xrightarrow{\sim} \varphi_{T/k}^{*}(\rT_{X'/k})$.
\end{nothing}

\begin{prop}[\cite{OV07} Thm. 1.1] \label{torsor Frob liftings}
	Suppose we are given a smooth lifting $\widetilde{X}'$ of $X'$ over $\rW_{2}$. 
	Let $\mathscr{T}_{X'/k}$ be the crystal of $\mathscr{O}_{X/k}$-modules associated to the $\mathscr{O}_{X'}$-module $\rT_{X'/k}$ \eqref{lemma can connection crystal}. Then, there exists a unique $\mathscr{T}_{X'/k}$-torsor $\mathscr{L}_{\widetilde{X}'}$ of $(X/k)_{\cris}$ satisfying following conditions:

	\textnormal{(i)} For every object $(U,T)$ of $\Cris(X/k)$ admitting a flat lifting $(U,\widetilde{T})$ in $\Cris(X/\rW_{2})$, the abelian sheaf $\mathscr{L}_{\widetilde{X}',(U,T)}$ of $T_{\zar}$ is the sheaf $\mathscr{L}_{\widetilde{X}',\varphi_{T/k}}$ \eqref{torsor of liftings}.

	\textnormal{(ii)} For every morphism $\widetilde{g}:(U_{1},\widetilde{T}_{1})\to (U_{2},\widetilde{T}_{2})$ of flat objects in $\Cris(X/\rW_{2})$ and any lifting $\widetilde{F}:\widetilde{T}_{2}\to \widetilde{X}'\in \mathscr{L}_{\widetilde{X}',\varphi_{T_{2}/k}}(T_{2})$, the transition morphism $c_{g}:\jmath_{g}^{-1}(\mathscr{L}_{\widetilde{X}',(U_{2},T_{2})})\to \mathscr{L}_{\widetilde{X}',(U_{1},T_{1})}$ satisfies
	\begin{equation}
		c_{g}(\jmath^{-1}_{g}(\widetilde{F}))=\widetilde{F}\circ \widetilde{g}:\widetilde{T}_{1}\to \widetilde{X}'.
	\end{equation}
\end{prop}

\begin{nothing} \label{alg splitting mod}
	We denote by $\mathscr{F}_{\widetilde{X}'}$ the crystal of affine functions on $\mathscr{L}_{\widetilde{X}'}$ and by $\mathscr{A}_{\widetilde{X}'}$ the quasi-coherent crystal of $\mathscr{O}_{X/k}$-algebras associated to $\mathscr{F}_{\widetilde{X}'}$ \eqref{crystal of affine functions}. We put $\mathcal{A}_{\widetilde{X}'}=\mathscr{A}_{\widetilde{X}',(X,X)}$. There exists an integrable connection $\nabla_{A}$ on $\mathcal{A}_{\widetilde{X}'}$. By \eqref{base change Gamma T} and \ref{Gamma action on A}, $\mathcal{A}_{\widetilde{X}'}$ is a quasi-nilpotent $F_{X/k}^{*}(\widehat{\Gamma}(\rT_{X'/k}))$-module. 
	The $p$-curvature
	\begin{displaymath}
		\psi:\mathcal{A}_{\widetilde{X}'}\to \mathcal{A}_{\widetilde{X}'}\otimes_{\mathscr{O}_{X}}F_{X/k}^{*}(\Omega_{X'/k}^{1})
	\end{displaymath}
	of $\nabla_{A}$ is equal to the universal $\mathscr{O}_{X}$-derivation \eqref{univer derivation} (cf. \cite{OV07} Prop. 1.5)
	\begin{displaymath}
		d:\mathcal{A}_{\widetilde{X}'} \to \mathcal{A}_{\widetilde{X}'}\otimes_{\mathscr{O}_{X}}F_{X/k}^{*}(\Omega_{X'/k}^{1}).
	\end{displaymath}
\end{nothing}

\begin{nothing}\label{PD Higgs}
	A Higgs field \eqref{notation MIC} on an $\mathscr{O}_X$-module $E$ relative to $k$ is equivalent to an $\rS(\rT_{X/k})$-module structure on $E$, which extends its $\mathscr{O}_{X}$-module structure (cf. \cite{OV07} 5.1).
	Let $I_{X}$ be the ideal $\oplus_{m>0}\rS^{m}(\rT_{X/k})$ of $\rS(\rT_{X/k})$. For any integer $n\ge 0$, we say that a Higgs module $E$ over $X$ relative to $k$ is \textit{nilpotent of level $\le n$}, if $E$ is annihilated by $I^{n+1}_{X}$ as an $\rS(\rT_{X/k})$-module.
	
	We call \textit{PD-Higgs module on $X$ relative to $k$} a $\widehat{\Gamma}(\rT_{X/k})$-module $E$, and we say that the structure morphism $\psi:\widehat{\Gamma}(\rT_{X/k})\to \FHom_{\mathscr{O}_{X}}(E,E)$ is the \textit{PD-Higgs field} on $E$. For any local section $\xi$ of $\widehat{\Gamma}(\rT_{X/k})$, we set $\psi_{\xi}=\psi(\xi)$.
	
	Let $n$ be an integer $\ge 1$. We denote by $\HIGG_{\gamma}(X/k)$ the category of $\widehat{\Gamma}(\rT_{X/k})$-modules and by $\HIGG_{\gamma}^{\qn}(X/k)$ (resp. $\HIGG_{\gamma}^{n}(X/k)$) the full subcategory of $\HIGG_{\gamma}(X/k)$ consisting of quasi-nilpotent objects (resp. nilpotent objects of level $\le n$) \eqref{sym alg general not}.
	
	Since $\rS^{n}(\rT_{X/k})\simeq \Gamma^{n}(\rT_{X/k})$ for all $0\le n\le p-1$, a nilpotent Higgs module of level $\le p-1$ induces naturally a nilpotent PD-Higgs module of level $\le p-1$. 
\end{nothing}

\begin{nothing}\label{Hom PD Higgs}
	Let $(E_{1},\psi_{1})$ and $(E_{2},\psi_{2})$ be two objects of $\HIGG_{\gamma}^{\qn}(X/k)$ and $\partial$ a local section of $\rT_{X/k}$. We define a PD-Higgs field $\psi$ on $E_{1}\otimes_{\mathscr{O}_{X}}E_{2}$ by (\cite{OV07} 2.7.1)
	\begin{equation} \label{tensor Higgs}
		\psi_{\partial^{[n]}}=\sum_{i+j=n} \psi_{1,\partial^{[i]}}\otimes \psi_{2,\partial^{[j]}}.
	\end{equation}

	Let $m,n$ be two integers, $(E_{1},\psi_{1})$ an object of $\HIGG_{\gamma}^{n}(X/k)$ and $(E_{2},\psi_{2})$ an object of $\HIGG_{\gamma}^{m}(X/k)$. There exists a unique PD-Higgs field $\psi$ on $\FHom_{\mathscr{O}_{X}}(E_{1},E_{2})$ defined, for every local sections $h$ of $\FHom_{\mathscr{O}_{X}}(E_{1},E_{2})$ and $\partial$ of $\rT_{X/k}$ by (cf. \cite{OV07} page 31)
	\begin{equation} \label{Hom psi}
		\psi_{\partial^{[l]}}(h)=\sum_{i+j=l}(-1)^{i}\psi_{2,\partial^{[i]}}\circ h\circ \psi_{1,\partial^{[j]}} \qquad \forall~l \ge 1.
	\end{equation}
%
\end{nothing}

\begin{nothing} \label{num p curvature}
	We denote by $\rD_{X/k}$ the ring of PD-differential operators on $X$ relative to $k$ (\cite{BO}, \S 4).
	Let $\partial$ be a local section of $\rT_{X/k}$ considered as a derivation of $\mathscr{O}_{X}$ over $k$ and hence as a PD-differential operator of order $\le 1$. 
	The $p$th iterate $\partial^{(p)}$ of $\partial$ is again a derivation of $\mathscr{O}_{X}$ over $k$ (\cite{Ka71} 5.0.2, \cite{BO} 4.5). We denote by $\partial^{p}$ the $p$th power of $\partial$ in $\rD_{X/k}$, which is an operator of order $\le p$. 
	The $p$-curvature morphism $c:F_{X}^{*}(\rT_{X/k})\to \rD_{X/k}$ defined by $F_{X}^{*}(\partial) \mapsto \partial^{p}-\partial^{(p)}$, induces an isomorphism of $\mathscr{O}_{X'}$-algebras (\cite{BMR} 2.2.3; \cite{OV07} Thm. 2.1)
	\begin{equation}
		\rS(\rT_{X'/k})\xrightarrow{\sim} F_{X/k*}(Z_{X/k}). \label{isomorphism ST ZD}
	\end{equation}
	The above morphism makes $F_{X/k*}(\rD_{X/k})$ into an Azumaya algebra over $\rS(\rT_{X'/k})$ of rank $p^{2d}$, where $d$ is the dimension of $X$ over $k$.
\end{nothing}


\begin{nothing} \label{D gamma}
	We denote by $\rD^{\gamma}_{X/k}$ the tensor product
	\begin{equation}
		\rD_{X/k}\otimes_{\rS(\rT_{X'/k})}\widehat{\Gamma}(\rT_{X'/k}). \label{def D gamma algebra}
	\end{equation}
	via the morphism $\rS(\rT_{X'/k})\to F_{X/k*}(\rD_{X/k})$ induced by the $p$-curvature morphism \eqref{isomorphism ST ZD}. To give a left $\rD_{X/k}^{\gamma}$-module is equivalent to give an $\mathscr{O}_{X}$-module $M$ with integrable connection $\nabla$ and a homomorphism
	\begin{equation} \label{PD p curvature}
		\psi:\widehat{\Gamma}(\rT_{X'/k})\to F_{X/k*}(\FEnd_{\mathscr{O}_{X}}(M,\nabla))
	\end{equation}
	which extends the Higgs field $\rS(\rT_{X'/k})\to F_{X/k*}(\FEnd_{\mathscr{O}_{X}}(M,\nabla))$ given by the $p$-curvature of $\nabla$ \eqref{isomorphism ST ZD} (cf. \cite{OV07} page 32). 
\end{nothing}
\begin{nothing} \label{D gamma qn}
	There exists an isomorphism $F_{X/k}^{*}(\widehat{\Gamma}(\rT_{X'/k}))\xrightarrow{\sim} \widehat{\Gamma}(F_{X/k}^{*}(\rT_{X'/k}))$ \eqref{base change Gamma T}. 
	Let $M$ be a left $\rD_{X/k}^{\gamma}$-module and $n$ an integer $\ge 0$. 
	We say that $M$ is \textit{quasi-nilpotent} (resp. \textit{nilpotent of level $\le n$}) if $M$ is quasi-nilpotent (resp. nilpotent of level $\le n$) as a $\widehat{\Gamma}(F_{X/k}^{*}(\rT_{X'/k}))$-module \eqref{sym alg general not}.

	We denote by $\MIC_{\gamma}(X/k)$ the category of left $\rD^{\gamma}_{X/k}$-modules and by $\MIC_{\gamma}^{\qn}(X/k)$ (resp. $\MIC_{\gamma}^{n}(X/k)$) the full subcategory of $\MIC_{\gamma}(X/k)$ consisting of quasi-nilpotent objects (resp. nilpotent objects of level $\le n$). 
	
	Let $(M,\nabla)$ be an $\mathscr{O}_{X}$-module with integrable connection whose $p$-curvature is nilpotent of level $\le p-1$ (\cite{Ka71} 5.6). Since $\rS^{n}(\rT_{X'/k})\simeq \Gamma_{n}(\rT_{X'/k})$ for all $0\le n\le p-1$, $(M,\nabla)$ induces naturally an object of $\MIC^{p-1}_{\gamma}(X/k)$. 
\end{nothing}

\begin{nothing} \label{FHom Dgamma}
	Let $(M_{1},\nabla_{1},\psi_{1})$ and $(M_{2},\nabla_{2},\psi_{2})$ be two objects of $\MIC_{\gamma}^{\qn}(X/k)$. There exists a canonical integrable connection $\nabla$ on $M_{1}\otimes_{\mathscr{O}_{X}} M_{2}$ (\cite{Ka71} 1.1.1). The morphisms $\psi_{1}$ and $\psi_{2}$ induce an action $\psi$ of $F_{X/k}^{*}(\widehat{\Gamma}(\rT_{X'/k}))$ on $M_{1}\otimes_{\mathscr{O}_{X}}M_{2}$ as in \eqref{tensor Higgs}. Then we obtain an object $(M_{1}\otimes_{\mathscr{O}_{X}}M_{2},\nabla,\psi)$ of $\MIC^{\qn}_{\gamma}(X/k)$.

	Let $m,n$ be two integers and $(M_{1},\nabla_{1},\psi_{1})$ an object of $\MIC_{\gamma}^{m}(X/k)$ and $(M_{2},\nabla_{2},\psi_{2})$ an object of $\MIC_{\gamma}^{n}(X/k)$. There exists an integrable connection $\nabla$ on the $\mathscr{O}_{X}$-module $\FHom_{\mathscr{O}_{X}}(M_{1},M_{2})$ defined, for every local sections $\partial $ of $\rT_{X/k}$ and $h$ of $\FHom_{\mathscr{O}_{X}}(M_{1},M_{2})$ by (\cite{Ka71} 1.1.2)
	\begin{equation} \label{formula hom connection}
		\nabla_{\partial}(h)=\nabla_{2,\partial}\circ h-h\circ \nabla_{1,\partial}.
	\end{equation}
	The morphisms $\psi_{1}$ and $\psi_{2}$ induce an action of $\widehat{\Gamma}(\rT_{X'/k})$ on $F_{X/k*}(\FHom_{\mathscr{O}_{X}}(M_{1},M_{2}))$ defined by the same formula as \eqref{Hom psi}. These data make $\FHom_{\mathscr{O}_{X}}(M_{1},M_{2})$ into an object of $\MIC_{\gamma}^{m+n}(X/k)$ (\cite{OV07} 2.1).
\end{nothing}

\begin{nothing} \label{AX' dual}
	By \ref{alg splitting mod}, the $\mathscr{O}_{X}$-algebra $\mathcal{A}_{\widetilde{X}'}$ \eqref{alg splitting mod} is equipped with a quasi-nilpotent left $\rD_{X/k}^{\gamma}$-module structure. Moreover, we have an exhaustive filtration $\{N_{n}(\mathcal{A}_{\widetilde{X}'})\}_{n\ge 0}$ of left $\rD^{\gamma}_{X/k}$-submodules of $\mathcal{A}_{\widetilde{X}'}$ such that $N_{n}(\mathcal{A}_{\widetilde{X}'})$ is nilpotent of level $\le n$ \eqref{crystal of affine functions}. 
	We define $(\mathcal{A}_{\widetilde{X}'})^{\vee}$ to be
	\begin{equation}
		(\mathcal{A}_{\widetilde{X}'})^{\vee}=\FHom_{\mathscr{O}_{X}}(\mathcal{A}_{\widetilde{X}'},\mathscr{O}_{X})\simeq \varprojlim_{n\ge 0}\FHom_{\mathscr{O}_{X}}(N_{n}(\mathcal{A}_{\widetilde{X}'}),\mathscr{O}_{X}).
	\end{equation}
	By \ref{D gamma qn} and \ref{FHom Dgamma}, we see that $(\mathcal{A}_{\widetilde{X}'})^{\vee}$ is an object of $\MIC_{\gamma}(X/k)$.

	The involution morphism $\rT_{X/k}\to \rT_{X/k}$ defined by $x\mapsto -x$, induces a homomorphism
	\begin{equation} \label{sigma Gamma T}
		\iota:\widehat{\Gamma}(\rT_{X'/k})\to \widehat{\Gamma}(\rT_{X'/k}).
	\end{equation}
\end{nothing}

\begin{theorem}[\cite{OV07} 2.8] \label{thm OV}
	Suppose we are given a smooth lifting $\widetilde{X}'$ of $X'$ over $\rW_{2}$. 

	\textnormal{(i)} The left $\rD_{X/k}^{\gamma}$-module $(\mathcal{A}_{\widetilde{X}'})^{\vee}$ is a splitting module for the Azumaya algebra $F_{X/k*}(\rD_{X/k}^{\gamma})$ over $\widehat{\Gamma}(\rT_{X'/k})$.
	
	\textnormal{(ii)} The functors \eqref{sigma Gamma T}
	\begin{eqnarray}
		\rmC_{\widetilde{X}'}&:\MIC_{\gamma}(X/k) \xrightarrow{\sim} \HIGG_{\gamma}(X'/k) \quad &E\mapsto \iota^{*}(\FHom_{\rD_{X/k}^{\gamma}}((\mathcal{A}_{\widetilde{X}'})^{\vee},E))\\
		\rmC_{\widetilde{X}'}^{-1}&:\HIGG_{\gamma}(X'/k)\xrightarrow{\sim} \MIC_{\gamma}(X/k) \quad & E'\mapsto (\mathcal{A}_{\widetilde{X}'})^{\vee}\otimes_{\widehat{\Gamma}(\rT_{X'/k})}\iota^{*}(E'). \label{C-1 OV}
	\end{eqnarray}
	are equivalences of categories quasi-inverse to each other. 
	Furthermore, they induce equivalences of tensor categories between $\MIC_{\gamma}^{\qn}(X/k)$ and $\HIGG_{\gamma}^{\qn}(X'/k)$ \textnormal{(\ref{Hom PD Higgs}, \ref{FHom Dgamma})}.

	\textnormal{(iii)} Let $(E,\nabla,\psi)$ be an object of $\MIC_{\gamma}(X/k)$ and $(E',\theta')=\rmC_{\widetilde{X'}}(E,\nabla,\psi)$. A lifting $\widetilde{F}$ of the relative Frobenius morphism $F_{X/k}$ induces a natural isomorphism of $F_{X/k}^{*}(\widehat{\Gamma}(\rT_{X'/k}))$-modules
	\begin{equation}
		\eta_{\widetilde{F}}:(E,\psi)\xrightarrow{\sim} F_{X/k}^{*}(E',-\theta').
	\end{equation}
\end{theorem}

We call $\rmC_{\widetilde{X}'}$ (resp. $\rmC_{\widetilde{X}'}^{-1}$) \textit{Cartier transform (resp. inverse Cartier transform)}.

\begin{theorem}[\cite{OV07} 2.17] \label{iso complexes OV}
	Suppose we are given a smooth lifting $\widetilde{X}'$ of $X'$ over $\rW_{2}$. 
	Let $(M',\theta')$ be a nilpotent Higgs module on $X'/k$ of level $\ell<p$ \eqref{PD Higgs} and $(M,\nabla)=\rmC_{\widetilde{X}'}^{-1}(M',\theta')$. Then, the lifting $\widetilde{X}'$ induces an isomorphism in the derived category $\rD(\mathscr{O}_{X'})$
	\begin{equation} \label{iso complexes OV eq}
		\tau_{<p-\ell}(M'\otimes_{\mathscr{O}_{X'}}\Omega_{X'/k}^{\bullet})\xrightarrow{\sim} F_{X/k*}(\tau_{<p-\ell}(M\otimes_{\mathscr{O}_{X}}\Omega_{X/k}^{\bullet})),
	\end{equation}
	where $M'\otimes_{\mathscr{O}_{X'}}\Omega_{X'/k}^{\bullet}$ is the Dolbeault complex of $(M',\theta')$, $M\otimes_{\mathscr{O}_{X}}\Omega_{X/k}^{\bullet}$ is the de Rham complex of $(M,\nabla)$ and $\tau_{<\bullet}$ denotes the truncation of a complex.	
\end{theorem}

We will give a partial generalization of this result for certain $p^{n}$-torsion crystals (cf. \ref{fil mod coh fil mod}, \ref{final remark on coh}).

\section{Prelude on rings of differential operators} \label{D ring}

The purpose of this section is to review the description of crystals of Oyama site in term of modules of rings of differential rings introduced in \S~\ref{Cartier OV} following Oyama. 
It serves as a preparation for section \ref{Comparison of Cartier}. 

In this section, $X$ denotes a smooth scheme over $k$. From \ref{local des RQ 1} on, suppose we are given a smooth formal $\SS$-scheme $\XX$ with special fiber $X$. 

\begin{nothing} \label{notation PX}
	Let $n\ge 1$ be an integer. We denote by $\PP_{X}$ (resp. $\PP_{X}^{n}$) the PD-envelope (resp. the $n^{\textnormal{th}}$ PD-neighborhood) of the diagonal immersion $\Delta: X\to X^{2}$ with respect to the zero PD-ideal of $k$ (\cite{BO} 3.31). We put $\mathcal{P}_{X}=\mathscr{O}_{\PP_{X}}$ and $\mathcal{P}_{X}^{n}=\mathscr{O}_{\PP_{X}^{n}}$ and we consider them as sheaves of $X_{\zar}$. By \ref{groupoid PX}, $\mathcal{P}_{X}$ is equipped with a Hopf $\mathscr{O}_{X}$-algebra structure $(\delta,\pi,\sigma)$ \eqref{def of Hopf alg}. 
	
	In the first part of this section, we study the $\mathscr{O}_{X}$-algebra $(\mathcal{P}_{X})^{\vee}$ \eqref{Hopf alg dual} of hyper PD-differential operators of $X$ relative to $k$. 

	Assume that there exists an \'etale $k$-morphism $X\to \mathbb{A}_{k}^{d}=\Spec(k[T_{1},\cdots,T_{d}])$. We set $t_{i}$ the image of $T_{i}$ in $\mathscr{O}_{X}$ and $\xi_{i}=1\otimes t_{i}-t_{i}\otimes 1$. We consider the $\xi_{i}$'s as sections of $\mathcal{P}_{X}$. Regarding $\mathcal{P}_{X}$ as a left (resp. right) $\mathscr{O}_{X}$-algebra, we have an isomorphism of PD-$\mathscr{O}_{X}$-algebras (\cite{Ber} I 4.5.3)
	\begin{equation} \label{local des PX}
		\mathscr{O}_{X}\langle x_{1},\cdots,x_{d}\rangle \xrightarrow{\sim} \mathcal{P}_{X},
	\end{equation}
	where $x_{i}$ is sent to $\xi_{i}$. 
	The homomorphism of PD-algebras $\delta:\mathcal{P}_{X}\to \mathcal{P}_{X}\otimes_{\mathscr{O}_{X}} \mathcal{P}_{X}$ (\cite{Ber} I 1.7.1) sends $\xi_{i}$ to $\xi_{i}\otimes 1+1\otimes \xi_{i}$. 
	For any $\alpha\in \mathbb{N}^{d}$, we set $\xi^{[\alpha]}=\prod \xi_{i}^{[\alpha_{i}]}$. Then we deduce that
	\begin{equation} \label{delta for PX}
		\delta(\xi^{[\alpha]})=\sum_{\beta\in \mathbb{N}^{d},\beta\le \alpha} \xi^{[\beta]}\otimes \xi^{[\alpha-\beta]}.
	\end{equation}
	The left (resp. right) $\mathscr{O}_{X}$-module $\mathcal{P}_{X}^{n}$ is free with a basis $\{\xi^{[\alpha]}, |\alpha|\le n\}$ (\cite{Ber} I 4.5.3).
\end{nothing}

\begin{nothing} \label{notation Z}
	Let $U$ be a open subscheme of $X^{2}$ such that $\Delta(X)\subset U$ and that $X\to U$ is a closed immersion. 
	The canonical morphism $\PP_{X}\to X^{2}$ factors through an affine morphism $\PP_{X}\to U$. We denote by $Z$ the scheme theoretic image of $\PP_{X}\to U$ (\cite{EGAInew} 6.10.1 and 6.10.5).
	Note that the morphisms $X\to \PP_{X}$ and $\PP_{X}\to Z$ induce isomorphisms between the underlying topological spaces. Hence we regard $\mathscr{O}_{Z}$ as an $\mathscr{O}_{X}$-bialgebra of $X_{\zar}$. We obtain an injective homomorphism of $\mathscr{O}$-bialgbras $\mathscr{O}_{Z}\to \mathcal{P}_{X}$ and we consider $\mathscr{O}_{Z}$ as a subalgebra of $\mathcal{P}_{X}$. 
%
%
\end{nothing}

\begin{lemma}\label{welldef Z}
	The scheme $Z$ defined in \ref{notation Z} is independent of the choice of $U$ up to canonical isomorphisms.
\end{lemma}
\begin{proof} Let $U_{1}$ and $U_{2}$ be two open subscheme of $X^{2}$ such that $\Delta(X)\subset U_{i}$ and that $X\to U_{i}$ is a closed immersion $i=1,2$ and $Z_{1}$ (resp. $Z_{2}$) the scheme theoretic image of $\PP_{X}$ in $U_{1}$ (resp. $U_{2}$). We can suppose that $U_{1}\subset U_{2}$. The image of $|Z_{1}|$ and $|\PP_{\XX}|$ in $|U_{1}|$ (resp. $|U_{2}|$) are equal. Then the composition $Z_{2}\to U_{1}\to U_{2}$ is a closed immersion. By (\cite{EGAInew} 6.10.3), we deduce a canonical isomorphism $Z_{1}\xrightarrow{\sim}Z_{2}$.
\end{proof}
\begin{lemma} \label{OZ description}
	Assume that $X$ is separated and that there exists an \'etale $k$-morphism $X\to \mathbb{A}_{k}^{d}$. 

	\textnormal{(i)} The left (resp. right) $\mathscr{O}_{X}$-module $\mathscr{O}_{Z}$ is free with a basis $\{\xi^{[\alpha]},\alpha\in \{0,1,\cdots,p-1\}^{d}\}$ \eqref{local des PX}.

	\textnormal{(ii)} The $\mathscr{O}_{Z}$-module $\mathcal{P}_{X}$ is free with a basis $\{\xi^{[pI]},I\in \mathbb{N}^{d}\}$.
\end{lemma}
\begin{proof} (i) It is clear that $\mathscr{O}_{Z}$ contains $\xi^{[\alpha]}$ for all $\alpha\in \{0,\cdots,p-1\}^{d}$. It suffices to show that $\mathscr{O}_{Z}$ is contained in the $\mathscr{O}_{X}$-submodule of $\mathcal{P}_{X}$ generated by $\{\xi^{[\alpha]},\alpha\in \{0,1,\cdots,p-1\}^{d}\}$.

The question being local, we suppose that $X$ is quasi-compact. Let $J$ be the ideal sheaf associated to the diagonal closed immersion $X\to X^{2}$ and $\varpi:\PP_{X}\to X^{2}$ the canonical morphism. By \ref{welldef Z}, we consider $Z$ as the scheme theoretic image of $\varpi$. The ideal $J$ being of finite type, we suppose that $J$ is generated by $m$ elements $x_{1},\cdots,x_{m}$ of $J(X^{2})$ for some integer $m\ge d$. Since $\mathscr{O}_{\PP_{X}}$ is a PD-algebra, the image of $x_{1}^{p},\cdots,x_{m}^{p}$ in $\varpi_{*}(\mathscr{O}_{\PP_{X}})$ are zero. Put $N=(p-1)m$. Then the image of the ideal $J^{N+1}$ in $\varpi_{*}(\mathscr{O}_{\PP_{X}})$ is zero, i.e. the morphism $\varpi$ factors through the $N^{\textnormal{th}}$ order infinitesimal neighborhood $Y_{N}=\Spec(\mathscr{O}_{X^{2}}/J^{N+1})$ of the diagonal immersion $X\to X^{2}$. Then we obtain a homomorphism
\begin{equation}
	\mathscr{O}_{Y_{N}} \to \mathcal{P}_{X}
\end{equation}
whose image is $\mathscr{O}_{Z}$. 
Recall (\cite{BO} 2.2) that, the left (resp. right) $\mathscr{O}_{X}$-module $\mathscr{O}_{Y_{N}}$ is free with a basis $\{\xi^{I},|I|\le N\}$. For any element $I=(i_{1},\cdots,i_{d})\in \mathbb{N}^{d}$, if one of the components $i_{j}$ is $\ge p$, then the image of $\xi^{I}$ in $\mathcal{P}_{X}$ is zero. Then the assertion follows. 

(ii) The assertion follows from (i) and the local description of $\mathcal{P}_{X}$ \eqref{local des PX}.
\end{proof}

\begin{nothing}
	The $p$-curvature morphism $c':\rS(\rT_{X'/k})\to F_{X/k*}(\rD_{X/k})$ induces an isomorphism between $\rS(\rT_{X'/k})$ and the center $Z_{X/k}$ of $\rD_{X/k}$ \eqref{isomorphism ST ZD}. It makes $\rD_{X/k}$ into an $\rS(\rT_{X'/k})$-module of finite type. Let $I_{X'}$ be the ideal $\oplus_{m\ge 1}\rS^{m}(\rT_{X'/k})$ of $\rS(\rT_{X'/k})$. We denote by $\mathcal{K}$ the two-side ideal of $\rD_{X/k}$ generated by $c'(I_{X'})$ and by $\widehat{\rD}_{X/k}$ the completion of $\rD_{X/k}$ with respect to the filtration $\{\mathcal{K}^{m}\}_{m\ge 1}$:
	\begin{equation}
		\widehat{\rD}_{X/k}=\varprojlim_{m\ge 1} \rD_{X/k}/\mathcal{K}^{m}.
	\end{equation}
	
	The completion $\widehat{\rD}_{X/k}$ is equal to the $(I_{X'})$-adic completion for the $\rS(\rT_{X'/k})$-module $\rD_{X/k}$. Then we obtain an isomorphism of $\widehat{\rS}(\rT_{X'/k})$-algebras
	\begin{equation} \label{DX hatDX}
		\widehat{\rS}(\rT_{X'/k})\otimes_{\rS(\rT_{X'/k})}\rD_{X/k}\xrightarrow{\sim} \widehat{\rD}_{X/k}.
	\end{equation}
\end{nothing}

\begin{nothing}
	Recall that $\rD_{X/k}$ is defined as $\varinjlim_{n\ge 1} (\mathcal{P}_{X}^{n})^{\vee}$, where $(\mathcal{P}_{X}^{n})^{\vee}=\FHom_{\mathscr{O}_{X}}(\mathcal{P}_{X}^{n},\mathscr{O}_{X})$.  
	Then we have a canonical homomorphism of $\mathscr{O}_{X}$-algebras:
	\begin{equation} \label{DX to hatDX}
		\rD_{X/k}\to (\mathcal{P}_{X})^{\vee}.
	\end{equation}

	For any integer $n\ge 0$ and any open subscheme $U$ of $X$, we define
	\begin{equation}
		\rF^{n}(\mathcal{P}_{X})(U)=\{a\in \mathcal{P}_{X}(U)|u(a)=0\quad \forall u\in \mathcal{K}^{n+1}(U)\}.
	\end{equation}
	Then $\rF^{n}(\mathcal{P}_{X})$ is a left $\mathscr{O}_{X}$-submodule of $\mathcal{P}_{X}$. We set $(\rF^{n}(\mathcal{P}_{X}))^{\vee}=\FHom_{\mathscr{O}_{X}}(\rF^{n}(\mathcal{P}_{X}),\mathscr{O}_{X})$. 
\end{nothing}

\begin{nothing} \label{local des D Dhat}
	Suppose that there exists an \'etale $k$-morphism $X\to \mathbb{A}_{k}^{d}=\Spec(k[t_{1},\cdots,t_{d}])$. Let $t_{i}$ be the image of $T_{i}$ in $\mathscr{O}_{X}$ and $\partial_{i}\in \rT_{X/k}(X)$ the dual of $dt_{i}$ that we consider as a PD-differential operator. 
	For $i,j$, $\partial_{i}$ and $\partial_{j}$ commute. We set $\partial^{I}=\prod_{j=1}^{d}\partial_{j}^{i_{j}}$ for all $I=(i_{1},\cdots,i_{d})\in \mathbb{N}^{d}$. 
	Any local section of $\rD_{X/k}$ can be written as a finite sum $\sum_{I} a_{I}\partial^{I}$. 
	The $p$-curvature morphism $c'$ sends $\partial_{i}'$ to $\partial_{i}^{p}$ and the ideal $\mathcal{K}$ is generated by $\{\partial_{1}^{p},\cdots,\partial_{d}^{p}\}$ over $\rD_{X/k}$. Then, any local section of $\widehat{\rD}_{X/k}$ can be written as an infinite sum:
	\begin{equation} \label{local des D Dhat eq}
		\sum_{I\in \mathbb{N}^{d}} a_{I} \partial^{I}\qquad \textnormal{with } a_{I}\in \mathscr{O}_{X}.
	\end{equation}
	Since $c'(\partial'^{\beta})=\partial^{p\beta}$ for $\beta\in \mathbb{N}^{d}$, the above section can be rewritten as
	\begin{equation} \label{local des D Dhat eq p cur}
		\sum_{\alpha\in \{0,\cdots,p-1\}^{d},\beta\in \mathbb{N}^{d}} b_{\alpha,\beta} \partial^{\alpha}\cdot c'(\partial'^{\beta}) \qquad \textnormal{with } b_{\alpha,\beta}\in \mathscr{O}_{X}.
	\end{equation}

	For any $I,J\in \mathbb{N}^{d}$, the image of $\partial^{I}$ in $(\mathcal{P}_{X})^{\vee}$ \eqref{DX to hatDX} satisfies $\partial^{I}(\xi^{[J]})=\delta_{I,J}$.
\end{nothing}

\begin{lemma}\label{lemma Ber for dual P}
	\textnormal{(i)} For any $n\ge 0$, the sheaf $\rF^{n}(\mathcal{P}_{X})$ is an $\mathscr{O}_{Z}$-submodule of $\mathcal{P}_{X}$ \eqref{notation Z}. 

	\textnormal{(ii)} With the assumption and notation of \textnormal{\ref{OZ description}}, $\rF^{n}(\mathcal{P}_{X})$ is a free $\mathscr{O}_{Z}$-module with basis $\{\xi^{[pI]}, |I|\le n\}$. 
\end{lemma}

\begin{proof} (i) The question being local, we take again the assumption of \ref{OZ description}.
Since the ideal $\mathcal{K}$ is generated by $\{\partial_{1}^{p},\cdots,\partial_{d}^{p}\}$, a local section $a=\sum_{I}b_{I}\xi^{[I]}$ of $\mathcal{P}_{X}$ is annihilated by $\mathcal{K}$ if and only if $b_{I}=0$ for all $I\in \mathbb{N}^{d}-\{0,\cdots,p-1\}^{d}$. By \ref{OZ description}(i), $\rF^{0}(\mathcal{P}_{X})$ is equal to the subsheaf $\mathscr{O}_{Z}$ of $\mathcal{P}_{X}$. Then the assertion follows.

(ii) The ideal $\mathcal{K}^{n+1}$ is generated by the set of PD-differential operators $\{\partial^{pI}, |I|=n+1\}$ over $\rD_{X/k}$. Then the assertion follows from \ref{OZ description}(ii) and the duality between $\partial^{I}$ and $\xi^{[I]}$.
\end{proof}
\begin{nothing}
	For any $n\ge 0$, since $\rF^{n}(\mathcal{P}_{X})$ is locally a direct summand of $\mathcal{P}_{X}$, the canonical morphism $(\mathcal{P}_{X})^{\vee}\to (\rF^{n}(\mathcal{P}_{X}))^{\vee}$ is surjective. By \eqref{delta for PX}, \ref{OZ description}(i) and \ref{lemma Ber for dual P}(ii), the homomorphism $\delta:\mathcal{P}_{X}\to \mathcal{P}_{X}\otimes_{\mathscr{O}_{X}}\mathcal{P}_{X}$ \eqref{notation PX} sends $\rF^{n}(\mathcal{P}_{X})$ to $\rF^{n}(\mathcal{P}_{X})\otimes_{\mathscr{O}_{X}}\rF^{n}(\mathcal{P}_{X})$. In view of \eqref{Hopf alg dual}, the $\mathscr{O}_{X}$-algebra $(\mathcal{P}_{X})^{\vee}$ induces an $\mathscr{O}_{X}$-algebra structure on $(\rF^{n}(\mathcal{P}_{X}))^{\vee}$.
\end{nothing}

\begin{prop}[Berthelot\footnote{We learn the proof of a local version of (ii) from a talk note given by Berthelot on Cartier transform.}] \label{dual of P Ber}
	\textnormal{(i)} For any integer $n\ge 1$, the homomorphism \eqref{DX to hatDX} induces a canonical isomorphism of $\mathscr{O}_{X}$-algebras $\rD_{X/k}/\mathcal{K}^{n+1}\xrightarrow{\sim} (\rF^{n}(\mathcal{P}_{X}))^{\vee}$.

	\textnormal{(ii)} The homomorphism \eqref{DX to hatDX} induces a canonical isomorphism of $\mathscr{O}_{X}$-algebras $\widehat{\rD}_{X/k}\xrightarrow{\sim} (\mathcal{P}_{X})^{\vee}$.	
\end{prop}
\begin{proof} (i) Since $\mathcal{K}^{n+1}$ acts trivially on $\rF^{n}(\mathcal{P}_{X})$, we obtain a homomorphism $\rD_{X/k}/\mathcal{K}^{n+1}\to (\rF^{n}(\mathcal{P}_{X}))^{\vee}$. In view of the local description of $\mathcal{K}^{n+1}$ and of $\rF^{n}(\mathcal{P}_{X})$ (\ref{lemma Ber for dual P}(ii)), this homomorphism is an isomorphism.

(ii) We have a canonical isomorphism $(\mathcal{P}_{X})^{\vee}=\FHom_{\mathscr{O}_{X}}(\varinjlim_{n\ge 0} \rF^{n}(\mathcal{P}_{X}),\mathscr{O}_{X})\simeq \varprojlim (\rF^{n}(\mathcal{P}_{X}))^{\vee}$. Then the assertion follows from (i).
\end{proof}
\begin{nothing}\label{local des RQ 1}
	In the remainder of this section, suppose that there exists a smooth formal $\SS$-scheme $\XX$ with special fiber $X$. 
	We put $\mathcal{R}_{\XX,1}=\mathcal{R}_{\XX}/p\mathcal{R}_{\XX}$ and $\mathcal{Q}_{\XX,1}=\mathcal{Q}_{\XX}/p\mathcal{Q}_{\XX}$ \eqref{notations R Q}. 
	We review the Oyama's description of $\widehat{\Gamma}(\rT_{X/k}),\rD_{X/k}^{\gamma}$ in terms of rings of differential operators associated to Hopf algebras $\mathcal{R}_{\XX,1},\mathcal{Q}_{\XX,1}$. 
	
	Recall that the left and the right $\mathscr{O}_{X}$-algebra structures of $\mathcal{R}_{\XX,1}$ are equal \eqref{notations R Q}. We denote by $q_{1},q_{2}:\QQ_{\XX,1}\to X$ the canonical morphisms. 
	Suppose that there exists an \'etale $\SS$-morphism $\XX\to \widehat{\mathbb{A}}_{\SS}^{d}=\Spf(\rW\{T_{1},\cdots,T_{d}\})$. We set $t_{i}$ the image of $T_{i}$ in $\mathscr{O}_{\XX}$ and $\xi_{i}=1\otimes t_{i}-t_{i}\otimes 1$.
	We take again the notation of \ref{local description R lemma} and of \ref{local description Q lemma} and we denote by $\zeta_{i}$ the element $\frac{\xi_{i}}{p}$ of $\mathcal{R}_{\XX,1}$ and by $\eta_{i}$ the element $\frac{\xi_{i}^{p}}{p}$ of $\mathcal{Q}_{\XX,1}$. We have isomorphisms of $\mathscr{O}_{X}$-algebras (\ref{local description R lemma}, \ref{local description Q lemma})
\begin{eqnarray} \label{R1 des local}
	&\mathscr{O}_{X}[x_1,\cdots,x_{d}]\xrightarrow{\sim} \mathcal{R}_{\XX,1}, \qquad &x_{i}\mapsto \zeta_{i};\\
	\label{Q1 des local} &\mathscr{O}_{X}[x_{1},\cdots,x_{d},y_{1},\cdots,y_{d}]/(x_{1}^{p},\cdots,x_{d}^{p})\xrightarrow{\sim} q_{j*}(\mathcal{Q}_{\XX,1}), &x_{i}\mapsto \xi_{i},~ y_{i}\mapsto \eta_{i}\quad j=1,2.
\end{eqnarray}
By \ref{local description Hopf RQ}, we have the following description of the Hopf algebra structures on $\mathcal{R}_{\XX,1}$ and $\mathcal{Q}_{\XX,1}$ (\cite{Oy} 1.2.9)
	\begin{eqnarray}\label{local Hopf RQ1}	
		&&\left\{ 
			\begin{array}{ll}
				\delta: \mathcal{R}_{\XX,1}\to \mathcal{R}_{\XX,1}\otimes_{\mathscr{O}_{X}}\mathcal{R}_{\XX,1} & \zeta_i\mapsto 1\otimes\zeta_i+\zeta_i\otimes1 \\
				\sigma:\mathcal{R}_{\XX,1}\to \mathcal{R}_{\XX,1} & \zeta_i \mapsto -\zeta_i  \\
				\pi:\mathcal{R}_{\XX,1}\to \mathscr{O}_{X} & \zeta_i\mapsto 0
			\end{array}
			\right. \\
		&&\left\{ 
			\begin{array}{ll}
				\delta: \mathcal{Q}_{\XX,1}\to \mathcal{Q}_{\XX,1}\otimes_{\mathscr{O}_{X}}\mathcal{Q}_{\XX,1} & \xi_{i}\mapsto 1\otimes \xi_{i} +\xi_{i}\otimes 1 \\
				& \eta_i\mapsto 1\otimes \eta_{i} +\sum_{j=1}^{p-1}\frac{(p-1)!}{j!(p-j)!}\xi_{i}^{j}\otimes \xi_{i}^{p-j} +\eta_i\otimes 1 \\
				\sigma:\mathcal{Q}_{\XX,1}\to \mathcal{Q}_{\XX,1} & \xi_{i}\mapsto -\xi_{i}, \eta_i \mapsto -\eta_i  \\
				\pi:\mathcal{Q}_{\XX,1}\to \mathscr{O}_{X} & \xi_{i}\mapsto 0, \eta_i\mapsto 0
			\end{array}
			\right. \nonumber
	\end{eqnarray}
\end{nothing}

\begin{prop}[\cite{Oy} 1.2.12] \label{dual of R}
	The Hopf $\mathscr{O}_{X}$-algebra $\mathcal{R}_{\XX,1}$ is canonically isomorphism to $\rS(\Omega_{X/k}^{1})$ \eqref{sym alg general not} and the $\mathscr{O}_{X}$-algebra $(\mathcal{R}_{\XX,1})^{\vee}$ \eqref{Hopf alg dual} is canonically isomorphic to $\widehat{\Gamma}(\rT_{X/k})$ \eqref{sym alg general not}. 
\end{prop}

\begin{nothing} \label{p curvature geometry}
	In the following, we study the $\mathscr{O}_{X}$-algebra $(\mathcal{Q}_{\XX,1})^{\vee}$. 
	Suppose first that there exists an $\SS$-morphism $F:\XX\to \XX'$ which lifts the relative Frobenius morphism $F_{X/k}$ of $X$. We denote by $\YY$ (resp. $Z$) the fiber product of $F^{2}:\XX^{2}\to \XX'^{2}$ and $\RR_{\XX'}\to \XX'^{2}$ (resp. the diagonal immersion $X'\to \XX'^{2}$) and by $S$ the fiber product of $F_{X/k}:X\to X'$ and $\RR_{\XX',1}\to X'$ \eqref{notations R Q}. The morphism $F_{X/k}:X\to X'$ and the diagonal immersion $\Delta:X\to \XX^{2}$ induce a closed immersion $X\hookrightarrow Z$. We have a commutative diagram:
	\begin{equation} \label{diag S Y}
		\xymatrix{
			S \ar[dd] \ar[r] \ar@{}[ddr]|{\Box}&Y \ar[rd] \ar[rr] \ar[dd] && \RR_{\XX',1} \ar'[d][dd] \ar[rd] &\\
			&& \YY \ar[rr] \ar[dd] && \RR_{\XX'} \ar[dd]\\
			X\ar[r]&Z \ar[rd] \ar'[r][rr]&& X' \ar[rd] &\\
			&& \XX^{2} \ar[rr]^{F^{2}} && \XX'^{2} }
	\end{equation}
	where $Y=\YY_{1}$ is equal to $\RR_{\XX',1}\times_{X'}Z$ and the left square is Cartesian. 

	Let $\UU$ be an open formal subscheme of $\XX^{2}$ such that diagonal immersion $X\to \XX^{2}$ factors through a closed immersion $X\to \UU$ and put $\UU'=\UU\times_{\SS,\sigma}\SS$.  
	Let $\mathscr{I}$ (resp. $\mathscr{I}'$) be the ideal associated to the diagonal immersion $X\to \UU$ (resp. $X'\to \UU'$). For any local section $x$ of $\mathscr{I}$, we have $x^{p}\in \mathscr{I}'\mathscr{O}_{\UU}$.  
	Since $Z$ is defined by $\mathscr{I}'\mathscr{O}_{\UU}$ and $X$ is reduced, the closed immersion $X\hookrightarrow Z$ induces an isomorphism $X\xrightarrow{\sim} \underline{Z}$. Then, the closed immersion $\underline{Y}\to Y$ factors through $S\to Y$. By the universal property of $\QQ_{\XX}$ \eqref{prop univ RQ}, we deduce an $\XX^{2}$-morphism 
	\begin{equation} \label{nu Y to QX}
		\nu: \YY\to \QQ_{\XX}.
	\end{equation}
	The composition $\QQ_{\XX}\to \XX^{2}\xrightarrow{F^{2}}\XX'^{2}$ induces a morphism of formal groupoids above $F$ \eqref{Q to R'}:
	\begin{equation}
		\psi:\QQ_{\XX}\to \RR_{\XX'}. 
	\end{equation}
	By \eqref{diag S Y} and the universal property of $\RR_{\XX'}$ \eqref{prop univ RQ}, we deduce that the composition $\psi\circ \nu:\YY\to \QQ_{\XX}\to \RR_{\XX'}$ is the canonical morphism 
	\begin{equation}
		\YY=\RR_{\XX'}\times_{\XX'^{2}}\XX^{2}\to \RR_{\XX'}.
	\end{equation}
\end{nothing}

\begin{nothing} \label{homomorphism v sF}
	Keep the assumption and notation of \ref{p curvature geometry}. The morphism $\nu$ \eqref{nu Y to QX} induces a morphism $S\to \QQ_{\XX,1}$ which makes the following diagram commutes
	\begin{equation}
		\xymatrix{
			S\ar[r] \ar[d] & \QQ_{\XX,1} \ar[r]^{\psi_{1}} \ar'[d][dd] & \RR_{\XX',1} \ar[d]\\
			X \ar[rd] \ar[rr] && X' \ar[d] \\
			& X^{2} \ar[r] & X'^{2}.
		}
	\end{equation}
	Let $g$ be the canonical morphism $S\to X$. Since $\RR_{\XX',1}\to X'$ is affine, we have an isomorphism $F_{X/k}^{*}(\mathcal{R}_{\XX',1})\xrightarrow{\sim} g_{*}(\mathscr{O}_{S})$ (\cite{EGAInew} 9.3.2). Then the morphism $S\to \QQ_{\XX,1}$ induces an $\mathscr{O}_{X}$-bilinear homomorphism
	\begin{equation} \label{homo v}
		v:\mathcal{Q}_{\XX,1}\to F_{X/k}^{*}(\mathcal{R}_{\XX',1}).
	\end{equation}
With the assumption and notation of \ref{local des RQ 1}, $v$ sends $\xi_{i}$ to $0$ and $\eta_{i}$ to $F_{X/k}^{*}(\zeta_{i}')$ (\cite{Oy} 1.2.14). 
Hence, $v$ is independent of the choice of the lifting $F:\XX\to \XX'$ of $F_{X/k}$, and can be defined for a general smooth formal $\SS$-scheme $\XX$ even if the relative Frobenius morphism cannot be lifted over $\SS$. 

By \eqref{local Hopf RQ1}, $v$ is compatible with Hopf $\mathscr{O}_{X}$-algebras structures. By taking $\mathscr{O}_{X}$-duals, we obtain a homomorphism of $\mathscr{O}_{X}$-algebras (\ref{base change Gamma T}, \ref{dual of R})
	\begin{equation} \label{v dual homo}
		v^{\vee}:F_{X/k}^{*}(\widehat{\Gamma}(\rT_{X'/k}))\to (\mathcal{Q}_{\XX,1})^{\vee}.
	\end{equation}

	The morphism $\psi_{1}$ induces a homomorphism of Hopf algebras $s_{F}:\mathcal{R}_{\XX',1}\to F_{X/k*}(\mathcal{Q}_{\XX,1})$ \eqref{Hopf algebra homomorphism}. By adjunction, we obtain a homomorphism of $\mathscr{O}_{X}$-algebras 
	\begin{equation} \label{homo sF}
		s_{F}^{\sharp}:F_{X/k}^{*}(\mathcal{R}_{\XX',1})\to \mathcal{Q}_{\XX,1}
	\end{equation}
	for the left $\mathscr{O}_{X}$-algebra structure on $\mathcal{Q}_{\XX,1}$. 
	The composition $v\circ s_{F}$ is the identical homomorphism.
	
	Recall that we have a canonical morphism of formal $\XX$-groupoids $\PP_{\XX}\to \QQ_{\XX}$ \eqref{pi P to Q}. Then we obtain a homomorphism of Hopf $\mathscr{O}_{X}$-algebras \eqref{Hopf algebra homomorphism}:
	\begin{equation}
		u:\mathcal{Q}_{\XX,1}\to \mathcal{P}_{X}.
	\end{equation}
	With the assumption and notation of \ref{local des RQ 1}, $u$ sends $\xi_{i}$ to $\xi_{i}$ and $\eta_{i}$ to $-\xi_{i}^{[p]}$ \eqref{local des PX} (\cite{Oy} 1.2.14).
	
	By taking duals (\ref{Hopf alg dual} and \ref{dual of P Ber}(ii)), $u$ induces a homomorphism of $\mathscr{O}_{X}$-algebras
	\begin{equation}
		u^{\vee}:\widehat{\rD}_{X/k}\to (\mathcal{Q}_{\XX,1})^{\vee}.
	\end{equation}
\end{nothing}

\begin{lemma}[\cite{Oy} 1.2.15] \label{lemma concides uv}
	Let $\sigma:F_{X/k}^{*}(\mathcal{R}_{\XX,1})\to F_{X/k}^{*}(\mathcal{R}_{\XX,1})$ be the involution homomorphism defined in \eqref{local Hopf RQ1}. 

	\textnormal{(i)} The restriction of $u^{\vee}$ and $(\sigma\circ v)^{\vee}$ to $F_{X/k}^{-1}(\widehat{\rS}(\rT_{X'/k}))$ coincide.

	\textnormal{(ii)} The images of any local sections of $u^{\vee}(\widehat{\rD}_{X/k})$ and of $(\sigma\circ v)^{\vee} (F_{X/k}^{-1}(\widehat{\Gamma}(\rT_{X'/k})))$ commute in $(\mathcal{Q}_{\XX,1})^{\vee}$.
\end{lemma}

\begin{prop}[\cite{Oy} 1.2.13] \label{dual of Q}
	The homomorphisms $(\sigma\circ v)^{\vee}$ and $u^{\vee}$ induce an isomorphism of $F_{X/k}^{*}(\widehat{\Gamma}(\rT_{X'/k}))$-algebras \eqref{D gamma}
	\begin{equation}
		\rD_{X/k}^{\gamma}\xrightarrow{\sim} (\mathcal{Q}_{\XX,1})^{\vee}. 
	\end{equation}
\end{prop}

\begin{prop}[\cite{Oy} 1.2.10] \label{equi Dmod strat}
	Let $\XX$ be a smooth formal $\SS$-scheme, $X$ its special fiber. 
	There exists a canonical equivalence of tensor categories between $\HIGG_{\gamma}^{\qn}(X/k)$ \eqref{PD Higgs} (resp. $\MIC^{\qn}_{\gamma}(X/k)$ \eqref{D gamma qn}) and the category of $\mathscr{O}_{X}$-modules with $\mathcal{R}_{\XX}$-stratification (resp. $\mathcal{Q}_{\XX}$-stratification) \textnormal{(\ref{prop R Q Hopf alg}, \ref{def stratification})}. 
\end{prop}
\begin{proof} We briefly recall here the main construction of the equivalence which will be used in the following and refer to \cite{Oy} for details. 

	Let $(M,\varepsilon)$ be an $\mathscr{O}_{X}$-module with $\mathcal{R}_{\XX}$-stratification and $\theta:M\to M\otimes_{\mathscr{O}_{X}}\mathcal{R}_{\XX,1}$ the $\mathscr{O}_{X}$-linear morphism defined by $\theta(m)=\varepsilon(1\otimes m)$ \eqref{lemma stratification}. By \ref{lemma stratification} and \ref{dual of R}, we deduce a $\widehat{\Gamma}(\rT_{X/k})$-module $\psi$ structure on $M$.
	To show that $\psi$ is quasi-nilpotent, we suppose that there exists an \'etale $\SS$-morphism $\XX\to \widehat{\mathbb{A}}_{\SS}^{d}$. We take again notation of \ref{local des RQ 1}. For any $I\in \mathbb{N}^{d}$, we set $\partial^{[I]}=\prod_{j=1}^{d}\partial_{j}^{[i_{j}]}\in \widehat{\Gamma}(\rT_{X/k})$. The action of $\partial^{[I]}$ on $M$ is given by the composition
\begin{equation} \label{action psi partialI}
	\psi_{\partial^{[I]}}: \xymatrix{
		M\ar[r]^-{\theta} & M\otimes_{\mathscr{O}_{X}} \mathcal{R}_{\XX,1}\ar[r]^-{\id \otimes \partial^{[I]}} & M}.
\end{equation}
Since $\mathcal{R}_{\XX,1}$ is isomorphic to a polynomial algebra over $\mathscr{O}_{X}$ \eqref{local des RQ 1}, for any local section $m$ of $M$ and any point $x$ of $X$, there exists a neighborhood $U$ of $x$ such that $\theta(m)|_{U}$ is a section of $M(U)\otimes_{\mathscr{O}_{X}(U)}\mathcal{R}_{\XX,1}(U)$ and can be written as a finite sum:
\begin{equation} \label{R stratification def}
	\varepsilon(1\otimes m)=\theta(m)=\sum_{I\in \mathbb{N}^{d}} \psi_{\partial^{[I]}}(m)\otimes \zeta^{I}.
\end{equation}
Hence, $(M,\psi)$ is quasi-nilpotent \eqref{sym alg general not}. 

Let $(M,\varepsilon)$ be an $\mathscr{O}_{X}$-module with $\mathcal{Q}_{\XX}$-stratification and $\theta:M\to M\otimes_{\mathscr{O}_{X}}\mathcal{Q}_{\XX,1}$ the morphism induced by $\varepsilon$. By \ref{lemma stratification} and \ref{dual of Q}, we associate to it a $\rD_{X/k}^{\gamma}$-module $(M,\nabla,\psi)$ \eqref{D gamma}. We take again the assumption and notation of the proof of \ref{dual of Q}. For $I\in \{0,\cdots,p-1\}^{d}, J\in \mathbb{N}^{d}$, the action of the local section $\partial^{I}\otimes F_{X/k}^{*}(\partial'^{[J]})$ of $\rD_{X/k}^{\gamma}$ on $M$ is given by the composition
\begin{equation}\label{action nabla partialI}
	\nabla_{\partial^{I}}\circ \psi_{\partial'^{[J]}}:
\xymatrixcolsep{5pc}\xymatrix{
	M\ar[r]^-{\theta} &M\otimes_{\mathscr{O}_{X}}\mathcal{Q}_{\XX,1} \ar[r]^-{\id \otimes w(\partial^{I}\otimes F_{X/k}^{*}(\partial'^{[J]}))} & M}
\end{equation}
Since $\mathcal{Q}_{\XX,1}$ is of finite type over a polynomial algebra over $\mathscr{O}_{X}$ \eqref{local des RQ 1}, for any local section $m$ of $M$ and any point $x$ of $X$, there exists a neighborhood $U$ of $x$ such that $\theta(m)|_{U}$ is a section of $M(U)\otimes_{\mathscr{O}_{X}(U)}\mathcal{Q}_{\XX,1}(U)$ and can be written as a finite sum (cf. \cite{Oy} page 25)
\begin{equation} \label{Q stratification def}
	\varepsilon(1\otimes m)=\theta(m)=\sum_{I\in \{0,\cdots,p-1\}^{d},J\in \mathbb{N}^{d}} (-1)^{|J|} \frac{1}{I!} (\nabla_{\partial^{I}}\circ \psi_{\partial'^{[J]}})(m) \otimes \xi^{I}\eta^{J}.
\end{equation}
Hence $(M,\nabla,\psi)$ is quasi-nilpotent \eqref{D gamma qn}. 
\end{proof}

We deduce from \ref{equi crystals stratification} and \ref{equi Dmod strat} the following equivalences. 

\begin{coro} \label{equi of Dmod crystals}
	Let $\XX$ be a smooth formal $\SS$-scheme, $X$ its special fiber, $\widetilde{\mathscr{E}}$ and $\widetilde{\underline{\mathscr{E}}}$ the Oyama topoi of $X$ \eqref{usual topology} and $\mathscr{O}_{\pCRIS,1}$, $\mathscr{O}_{\CRIS,1}$ the $p$-torsion structure rings of $\widetilde{\pCRIS}$ and $\widetilde{\CRIS}$ respectively \eqref{ring usual topology}. 
	We denote by $\mathscr{C}(\mathscr{O}_{\pCRIS,1})$ (resp. $\mathscr{C}(\mathscr{O}_{\CRIS,1})$) the category of crystals of $\mathscr{O}_{\pCRIS,1}$-modules of $\widetilde{\pCRIS}$ (resp. $\mathscr{O}_{\CRIS,1}$-modules of $\widetilde{\CRIS}$) \eqref{def crystals}. 
	Then, we have equivalences of tensor categories
	\begin{equation}\label{equivalence Dmod crystals} 
		\HIGG_{\gamma}^{\qn}(X/k)\simeq \mathscr{C}(\mathscr{O}_{\pCRIS,1}), \qquad \MIC_{\gamma}^{\qn}(X/k) \simeq \mathscr{C}(\mathscr{O}_{\CRIS,1}).
	\end{equation}
\end{coro}

The following result is an analogue of \ref{thm OV}(iii) for Cariter equivalence.

\begin{coro} \label{lifting F Cartier}
	Let $\rmC_{X/\SS}:\widetilde{\CRIS}\to \widetilde{\pCRIS}'$ be the morphism of topoi associated to $X$ \eqref{morphism of ringed topos Cartier}, $\mathscr{M}'$ a crystal of $\mathscr{O}_{\pCRIS',1}$-modules of $\widetilde{\pCRIS}'$ \eqref{def crystals}, $(M',\theta)$ the associated $\widehat{\Gamma}(\rT_{X'/k})$-module \eqref{equi of Dmod crystals} and $(\nabla,\psi)$ the $\rD_{X/k}^{\gamma}$-module structure on $(\rmC_{X/\SS}^{*}(\mathscr{M}'))_{(X,\XX)}$. Then a lifting $F:\XX\to \XX'$ of the relative Frobenius morphism $F_{X/k}$ of $X$ induces a functorial isomorphism of $F_{X/k}^{*}(\widehat{\Gamma}(\rT_{X'/k}))$-modules:
	\begin{equation} \label{Frob-lift induces iso}
		\eta_{F}: \iota^{*}(F_{X/k}^{*}(M',\theta)) \xrightarrow{\sim} ((\rmC_{X/\SS}^{*}(\mathscr{M}'))_{(X,\XX)},\psi)
	\end{equation}
	where $\iota$ denotes the involution homomorphism $F_{X/k}^{*}(\widehat{\Gamma}(\rT_{X'/k}))\to F_{X/k}^{*}(\widehat{\Gamma}(\rT_{X'/k}))$ \eqref{sigma Gamma T}. 
\end{coro}
\begin{proof}	We take again the notation of \ref{homomorphism v sF}. 
By \ref{homo Hopf alg dual}, \ref{dual of R} and \ref{dual of Q}, the homomorphism $s_{F}^{\sharp}$ \eqref{homo sF} induces a homomorphism of $\mathscr{O}_{X}$-algebras:
	\begin{equation}\label{sF dual}
		(s_{F}^{\sharp})^{\vee}:\rD_{X/k}^{\gamma}\to F_{X/k}^{*}(\widehat{\Gamma}(\rT_{X'/k})).
	\end{equation}
	Since the composition $v\circ s_{F}^{\sharp}$ is the identity \eqref{homomorphism v sF}, the composition \eqref{v dual homo}
	\begin{equation} \label{composition equal to iota}
		F_{X/k}^{*}(\widehat{\Gamma}(\rT_{X'/k}))\xrightarrow{(\sigma\circ v)^{\vee}} \rD_{X/k}^{\gamma}\xrightarrow{ (s_{F}^{\sharp})^{\vee}} F_{X/k}^{*}(\widehat{\Gamma}(\rT_{X'/k})).
	\end{equation}
	is the involution homomorphism $\iota:F_{X/k}^{*}(\widehat{\Gamma}(\rT_{X'/k}))\to F_{X/k}^{*}(\widehat{\Gamma}(\rT_{X'/k}))$ \eqref{sigma Gamma T}.

Let $\varepsilon$ be the associated $\mathcal{R}_{\XX'}$-stratification on $M'$. Recall \ref{Cartier global local} that the morphism $F$ induces a functorial isomorphism of $\mathscr{O}_{X}$-modules compatible with $\mathcal{Q}_{\XX}$-stratifications: 
	\begin{equation}
		\eta_{F}: (F_{X/k}^{*}(M'),s_{F}^{*}(\varepsilon)) \xrightarrow{\sim} (\rmC_{X/\SS}^{*}(\mathscr{M}'))_{(X,\XX)}.
	\end{equation}
	In view of \ref{equi Dmod strat}, the associated $\rD_{X/k}^{\gamma}$-module structure on $F_{X/k}^{*}(M')$ in the left hand side is give by $F_{X/k}^{*}(\theta)$ via $(s_{F}^{\sharp})^{\vee}$ that we denote by $(\nabla_{F},\psi_{F})$.
	Then we obtain a functorial isomorphism of $\rD_{X/k}^{\gamma}$-modules
	\begin{equation} 
		\eta_{F}: (F_{X/k}^{*}(M'),\nabla_{F},\psi_{F})\xrightarrow{\sim} (\rmC_{X/\SS}^{*}(\mathscr{M}'))_{(X,\XX)}.
	\end{equation}
	In view of \eqref{composition equal to iota}, we have an equality of $F_{X/k}^{*}(\widehat{\Gamma}(\rT_{X'/k}))$-modules
	\begin{equation}\label{psi F theta}
		(F_{X/k}^{*}(M'),\psi_{F})= \iota^{*}(F_{X/k}^{*}(M',\theta)).
	\end{equation}
This concludes the proof.
\end{proof}

\section{Comparison with the Cartier transform of Ogus--Vologodsky} \label{Comparison of Cartier}
	In this section, we compare the Cartier equivalence modulo $p$ and the Cartier transform of Ogus--Vologodsky.
	Our approach is different to that of Oyama. 
	We interpret the Cartier equivalence as an admissibility condition à la Fontaine for a pair of crystals \eqref{thm iso ad} and use it to compare with Cartier transform \eqref{comparison Cartier}. 	

	Let $\XX$ be a smooth formal $\SS$-scheme, $X$ its special fiber and $\rT_{X/k}$ the $\mathscr{O}_{X}$-dual of the $\mathscr{O}_{X}$-module of differential forms $\Omega_{X/k}^{1}$. We set $\XX'=\XX\otimes_{\SS,\sigma}\SS$. 

\begin{nothing}  \label{morphism topos v}
	We first interpret the Cartier descent in the context of Cartier equivalence \eqref{Cartier transform descent}. 

	Let $\widetilde{\pCRIS}$ and $\widetilde{\CRIS}$ be the Oyama topoi of $X$ \eqref{usual topology} and $\mathscr{O}_{\pCRIS,1}$, $\mathscr{O}_{\CRIS,1}$ the $p$-torsion structure rings of $\widetilde{\pCRIS}$ and $\widetilde{\CRIS}$ respectively \eqref{ring usual topology}. We keep the conventions and notation of \S~\ref{Oyama topos}-\ref{crystals}. 
	A morphism $g:(U_{1},\mathfrak{T}_{1},u_{1})\to (U_{2},\mathfrak{T}_{2},u_{2})$ of $\pCRIS$ (resp. $\CRIS$) induces a morphism of ringed topoi \eqref{tilde fn morphism}
	\begin{equation}\label{tilde gs}
		\widetilde{g}_{s}:(U_{1},u_{1*}(\mathscr{O}_{T_{1}}))\to (U_{2},u_{2*}(\mathscr{O}_{T_{2}})).
	\end{equation}
	
	If we equip $\pCRIS$ (resp. $\CRIS$) with the Zariski topology, the functor \eqref{pi Cris to Zar}
	\begin{equation}
		\pi:\pCRIS\to \Zar_{/X} \quad \textnormal{(resp. } \pi:\CRIS\to \Zar_{/X}) \quad (U,\mathfrak{T})\to U
	\end{equation}
	is cocontinuous. Since $\pi$ commutes with the fibered product of a flat morphism and a morphism \eqref{fiber product Oyama}, one verifies that $\pi$ is also continuous (\cite{SGAIV} III 1.6). 
	By \ref{generality morphism topos}, it induces a morphism of topoi that we denote by
	\begin{equation} \label{morphism v}
		v:\widetilde{\pCRIS}\to X_{\zar} \qquad \textnormal{(resp. } v:\widetilde{\CRIS}\to X_{\zar} )
	\end{equation}
	such that the inverse image functor is induced by the composition with $\pi$. Moreover, one verifies that the above morphisms fit into a commutative diagram
	\begin{equation}\label{morphism v Cartier}
		\xymatrix{
			\widetilde{\CRIS} \ar[r]^{\rmC_{X/\SS}} \ar[d]_{v} & \widetilde{\pCRIS}' \ar[d]^{v'} \\
			X_{\zar} \ar[r]^{F_{X/k}} & X'_{\zar}
		}
	\end{equation}

	Let $\mathscr{F}$ be a sheaf of $X_{\zar}$. For any object $(U,\mathfrak{T})$ of $\pCRIS$ (resp. $\CRIS$), we have $(v^{*}(\mathscr{F}))_{(U,\mathfrak{T})}=\mathscr{F}|_{U}$. For any morphism $f:(U_{1},\mathfrak{T}_{1})\to (U_{2},\mathfrak{T}_{2})$ of $\pCRIS$ (resp. $\CRIS$), the transition morphism $\gamma_{f}$ of $v^{*}(\mathscr{F})$ \eqref{morphism cf} is the canonical isomorphism
	\begin{equation} \label{transition morphism pullback v}
		\gamma_{f}:(\mathscr{F}|_{U_{2}})|_{U_{1}}\xrightarrow{\sim} \mathscr{F}|_{U_{1}}.
	\end{equation}
\end{nothing}
\begin{nothing}
	For any object $(U,\mathfrak{T},u)$ of $\pCRIS$, the morphism $u:T\to U$ induces a canonical, functorial homomorphism
	\begin{equation} \label{ring homo nu}
		v^{*}(\mathscr{O}_{X})(U,\mathfrak{T},u)=\mathscr{O}_{X}(U)\to \mathscr{O}_{\pCRIS,1}(U,\mathfrak{T},u)=u_{*}(\mathscr{O}_{T})(U). 
	\end{equation}
Then the morphism of topoi $v:\widetilde{\pCRIS}\to X_{\zar}$ underlies a morphism of ringed topoi 
\begin{equation} \label{morphism topos nu}
	\nu:(\widetilde{\pCRIS},\mathscr{O}_{\pCRIS,1})\to (X_{\zar},\mathscr{O}_{X}).
\end{equation}

For any object $(U,\mathfrak{T},u)$ of $\CRIS$, we have a morphism $u'\circ f_{T/k}:T\to U'$ \eqref{diagram rho} and hence a canonical, functorial homomorphism
	\begin{equation}
		v^{*}(F_{X/k}^{*}(\mathscr{O}_{X'}))(U,\mathfrak{T},u)=\mathscr{O}_{X'}(U')\to \mathscr{O}_{\CRIS,1}(U,\mathfrak{T},u)=(u'\circ f_{T/k})_{*}(\mathscr{O}_{T})(U'). 
	\end{equation}
	The composition of morphisms $F_{X/k}\circ v:\widetilde{\CRIS}\to X_{\zar} \to X'_{\zar}$ underlies a morphism of ringed topoi 
	\begin{equation} \label{morphism topos mu}
	\mu:(\widetilde{\CRIS},\mathscr{O}_{\CRIS,1})\to (X'_{\zar},\mathscr{O}_{X'}).
	\end{equation}
	
	By \eqref{morphism v Cartier}, the morphisms $\mu$ and $\nu'$ fit into a commutative diagram
	\begin{equation} \label{Cartier descent ringed topos}
		\xymatrix{
			(\widetilde{\CRIS},\mathscr{O}_{\CRIS,1}) \ar[rr]^{\rmC_{X/\SS}} \ar[rd]_{\mu} && (\widetilde{\pCRIS}',\mathscr{O}_{\pCRIS',1}) \ar[ld]^{\nu'} \\
			&(X'_{\zar},\mathscr{O}_{X'})& }
	\end{equation}
\end{nothing}

\begin{prop} \label{Cartier transform descent}
	Let $\Mod(\mathscr{O}_{X'})$ be the category of $\mathscr{O}_{X'}$-modules and let $\lambda$ be the functor
	\begin{equation}
		\lambda:\Mod(\mathscr{O}_{X'})\to \MIC_{\gamma}^{0}(X/k) \qquad M'\mapsto (F_{X/k}^{*}(M'),\nabla_{\can},0)
	\end{equation}
	where $\MIC_{\gamma}^{0}(X/k)$ denotes the category of nilpotent $\rD_{X/k}^{\gamma}$-modules of level $\le 0$ \eqref{D gamma qn}.
	Then, the following diagram is commutative up to a canonical isomorphism
	\begin{equation}
		\xymatrix{
			\Mod(\mathscr{O}_{X'}) \ar[r]^{\lambda} \ar[d]_{\nu'^{*}} & \MIC_{\gamma}^{0}(X/k) \ar[d]\\
			\mathscr{C}(\mathscr{O}_{\pCRIS',1}) \ar[r]^{\rmC^{*}_{X/\SS}}& \mathscr{C}(\mathscr{O}_{\CRIS,1}) }
	\end{equation}
	where the right vertical arrow is given by \eqref{equi of Dmod crystals}.
\end{prop}

The proposition follows from lemmas \ref{stratification id Higgs nul}, \ref{pullback Cartier descent} below.

\begin{nothing} \label{module 0 to crystal}
	Let $M$ be an $\mathscr{O}_{X}$-module and $\mathscr{M}=\nu^{*}(M)$ \eqref{morphism topos nu}. For any object $(U,\mathfrak{T},u)$ of $\pCRIS$, we set $\widetilde{u}^{*}(M)=M|_{U}\otimes_{\mathscr{O}_{U}}u_{*}(\mathscr{O}_{T})$. Recall \eqref{morphism topos v} that we have $(v^{*}(M))_{(U,\mathfrak{T})}=M|_{U}$. In view of \ref{sheaf to datas} and \eqref{ring homo nu}, we deduce that $\mathscr{M}_{(U,\mathfrak{T},u)}=\widetilde{u}^{*}(M)$.
	For any morphism $g:(U_{1},\mathfrak{T}_{1},u_{1})\to (U_{2},\mathfrak{T}_{2},u_{2})$ of $\pCRIS$, in view of \eqref{transition morphism pullback v}, the transition morphism $c_{g}:\widetilde{g}_{s}^{*}(\mathscr{M}_{(U_{2},\mathfrak{T}_{2})})\to \mathscr{M}_{(U_{1},\mathfrak{T}_{1})}$ is the canonical isomorphism \eqref{tilde gs}
	\begin{equation} \label{def cg M 0}
		\widetilde{g}_{s}^{*}(\widetilde{u}_{2}^{*}(M)) \xrightarrow{\sim} \widetilde{u}_{1}^{*}(M).
	\end{equation}
	Hence $\mathscr{M}$ is a crystal of $\mathscr{O}_{\pCRIS,1}$-modules of $\widetilde{\pCRIS}$. 
	If $M$ is moreover quasi-coherent, then so is $\mathscr{M}$ \eqref{def crystals}. In this case, for any object $(U,\mathfrak{T},u)$ of $\pCRIS$, the $\mathscr{O}_{T}$-module $\mathscr{M}_{(U,\mathfrak{T})}$ of $T_{\zar}$ \eqref{crystal lin descent data} is $u^{*}(M|_{U})$. 
\end{nothing}

\begin{lemma} \label{stratification id Higgs nul}
	Under the assumption of \ref{module 0 to crystal}, the $\widehat{\Gamma}(\rT_{X/k})$-module associated to $\mathscr{M}$ \eqref{equi of Dmod crystals} is the $\mathscr{O}_{X}$-module $M$ equipped with the zero PD-Higgs field \eqref{PD Higgs}.
\end{lemma}
\begin{proof} The underlying $\mathscr{O}_{X}$-module is $\mathscr{M}_{(X,\XX)}=M$. The reduction modulo $p$ of the two canonical morphisms $q_{1},q_{2}:\RR_{\XX}\to \XX$ are equal. By \eqref{def cg M 0}, the $\mathcal{R}_{\XX}$-stratification on $M$ associated to $\mathscr{M}$ \eqref{equi crystals stratification}
\begin{displaymath}
	\varepsilon:\widetilde{q}_{2,s}^{*}(M)\xrightarrow{\sim} \widetilde{q}_{1,s}^{*}(M)
\end{displaymath}
is the identity morphism. In view of \eqref{action psi partialI} and \eqref{R stratification def}, the PD-Higgs field associated to $\varepsilon$ is zero. 
\end{proof}
\begin{lemma} \label{pullback Cartier descent}
	Let $M'$ be an $\mathscr{O}_{X'}$-module and $\mathscr{M}=\mu^{*}(M')$. Then $\mathscr{M}$ is a crystal of $\mathscr{O}_{\CRIS,1}$-modules of $\widetilde{\CRIS}$ and the $\rD_{X/k}^{\gamma}$-module associated to $\mathscr{M}$ is $(F_{X/k}^{*}(M'),\nabla_{\can},0)$ \eqref{equi of Dmod crystals}, where $\nabla_{\can}$ denotes the Frobenius descent connection on $F_{X/k}^{*}(M')$ \eqref{can connection}. 
\end{lemma}
\begin{proof} We set $\mathscr{M}'=\nu'^{*}(M')$ \eqref{Cartier descent ringed topos}. Then $\mathscr{M}=\rmC_{X/\SS}^{*}(\mathscr{M}')$ is a crystal by \ref{module 0 to crystal}. For any object $(U,\mathfrak{T},u)$ of $\CRIS$, we put $\phi_{T/k}=u'\circ f_{T/k}:T\to \underline{T}'\to U'$. Then we have $\rho(U,\mathfrak{T},u)=(U',\mathfrak{T},\phi_{T/k})$ \eqref{functor rho} and \eqref{C-1 calcul} 
\begin{equation} \label{C-1 calcul descent Cartier}
	\rmC^{*}_{X/\SS}(\mathscr{M}')_{(U,\mathfrak{T},u)}=\pi_{U*}(\mathscr{M}'_{(U',\mathfrak{T},\phi_{T/k})})=\pi_{U*}(\widetilde{\phi}_{T/k}^{*}(M')).
\end{equation}
The morphism $\phi_{X/k}$ associated to the object $(X,\XX)$ of $\CRIS$ is $F_{X/k}$. 
Then we have $\mathscr{M}_{(X,\XX)}=F_{X/k}^{*}(M')$. There exists a commutative diagram
\begin{equation} \label{diag com des Cartier Q}
		\xymatrixcolsep{5pc}\xymatrix{
			\QQ_{\XX,1} \ar[r]^{q_{2}} \ar[d]_{q_{1}} \ar[rd]^{\phi_{\QQ_{\XX,1}/k}}& X \ar[d]^{F_{X/k}}\\
			X \ar[r]^{F_{X/k}} & X'}
\end{equation}
The morphisms $q_{1},q_{2}:(X,\QQ_{\XX})\to (X,\XX)$ of $\CRIS$ induce isomorphisms (\ref{def cg M 0}, \ref{C-1 calcul descent Cartier})
\begin{equation}
	c_{q_{1}}:\widetilde{q}_{1,s}^{*}(F_{X/k}^{*}(M'))\xrightarrow{\sim}\pi_{X*}(\widetilde{\phi}_{\QQ_{\XX,1}/k}^{*}(M')), \qquad c_{q_{2}}:\widetilde{q}_{2,s}^{*}(F_{X/k}^{*}(M'))\xrightarrow{\sim}\pi_{X*}(\widetilde{\phi}_{\QQ_{\XX,1}/k}^{*}(M')).
\end{equation}
The $\mathcal{Q}_{\XX}$-stratification $\varepsilon$ on $F_{X/k}^{*}(M')$ associated to the crystal $\mathscr{M}$ \eqref{equi crystals stratification} is given by the composition of $c_{q_{2}}$ and the inverse of $c_{q_{1}}$. In view of \eqref{diag com des Cartier Q}, for any local section $m'$ of $M'$, we have
\begin{equation} \label{connection nabla can}
	\varepsilon(1\otimes F_{X/k}^{*}(m'))=F_{X/k}^{*}(m')\otimes 1.
\end{equation}
Let $(F_{X/k}^{*}(M'),\nabla,\psi)$ be the $\rD_{X/k}^{\gamma}$-module associated to $(F_{X/k}^{*}(M'),\varepsilon)$ \eqref{equi Dmod strat}. In view of \eqref{connection nabla can}, \eqref{action nabla partialI} and \eqref{Q stratification def}, we deduce that $\nabla$ and $\psi$ annihilate the subsheaf $F_{X/k}^{-1}(M')$ of $F_{X/k}^{*}(M')$. Hence $\nabla$ is the Frobenius descent connection and $\psi=0$.
\end{proof}

\begin{nothing} \label{def sheaf L Zar}
	To compare the Cartier equivalence and the Cartier transform, we reconstruct the crystal $\mathcal{A}_{\XX_2'}$ \eqref{alg splitting mod} as a crystal in Oyama topos $\widetilde{\underline{\mathscr{E}}}$ (cf. \ref{C inverse underlying} and \ref{coro two splitting}) which allows us to compare Cartier transform and Cartier equivalence. 
	Compared to Oyama's approach for this construction (\cite{Oy} 1.5.2), our approach use certain torsor of liftings in Oyama topoi and is close the original construction of Ogus--Vologodsky. 
	
	Let $(U,\mathfrak{T},u)$ be an object of $\pCRIS$, $W$ an open subscheme of $T$ and $\mathfrak{W}_{2}$ the open subscheme of $\mathfrak{T}_{2}$ associated to $W$. We define $\mathscr{R}_{\XX,(U,\mathfrak{T},u)}(W)$ to be the set of $\SS_{2}$-morphisms $\mathfrak{W}_{2}\to \mathfrak{X}_{2}$ which make the following diagram commutes
	\begin{equation} \label{normal lifting}
		\xymatrix{
			W\ar[rr] \ar[d]_{u|_{W}}& & \mathfrak{W}_{2} \ar[d]\\ 
			U \ar[r]& X \ar[r]& \XX_{2} }
	\end{equation}
	The functor $W\mapsto \mathscr{R}_{\XX,(U,\mathfrak{T},u)}(W)$ defines a sheaf of $T_{\zar}$. Since $\XX$ is smooth over $\SS$, such morphisms exist locally. The sheaf $\mathscr{R}_{\XX,(U,\mathfrak{T},u)}$ is a torsor under the $\mathscr{O}_{T}$-module $\FHom_{\mathscr{O}_{\mathfrak{T}_{2}}}(u^{*}(\Omega_{U/k}^{1}),$ $p\mathscr{O}_{\mathfrak{T}_{2}})\xrightarrow{\sim}u^{*}(\rT_{U/k})$ of $T_{\zar}$.
\end{nothing}

\begin{nothing} \label{Relevement X}
	Let $g:(V,\mathfrak{Z},v)\to (U,\mathfrak{T},u)$ be a morphism of $\pCRIS$ and $g_{s}:Z\to T$ (resp. $g_{2}:\mathfrak{Z}_{2}\to \mathfrak{T}_{2}$) the reduction of the morphism $\mathfrak{Z}\to \mathfrak{T}$. 
	For $\mathscr{O}_{T}$-modules, we will use the notation $g_{s}^{-1}$ to denote the inverse image in the sense of abelian sheaves and will keep the notation $g_{s}^{*}$ for the inverse image in the sense of modules. By adjunction, the isomorphism $g_{s}^{*}(u^{*}(\rT_{U/k}))\xrightarrow{\sim} v^{*}(\rT_{V/k})$ induces an $\mathscr{O}_{T}$-linear morphism:
	\begin{displaymath}
		\tau:u^{*}(\rT_{U/k})\to g_{s*}(v^{*}(\rT_{V/k})).
	\end{displaymath}
	We have a canonical $\tau$-equivariant morphism $\mathscr{R}_{\XX,(U,\mathfrak{T},u)} \to g_{s*}(\mathscr{R}_{\XX,(V,\mathfrak{Z},v)})$ of $T_{\zar}$ defined for every local section $h:\mathfrak{W}_{2}\to \XX_{2}$ of $\mathscr{R}_{\XX,(U,\mathfrak{T},u)}$ by
	\begin{equation} \label{pre gamma gs cocycle Lifting}
		\mathscr{R}_{\XX,(U,\mathfrak{T},u)} \to g_{s*}(\mathscr{R}_{\XX,(V,\mathfrak{Z},v)})\qquad h\mapsto h\circ g_{2}|_{g_{2}^{-1}(\mathfrak{W}_{2})}.
	\end{equation}
	By adjunction, we obtain a $g_{s}^{-1}(u^{*}(\rT_{U/k}))$-equivariant morphism of $Z_{\zar}$:
	\begin{equation} \label{gamma gs cocycle Lifting}
		\gamma_{g}:g_{s}^{*}(\mathscr{R}_{\XX,(U,\mathfrak{T},u)})\to \mathscr{R}_{\XX,(V,\mathfrak{Z},v)}. 
	\end{equation}
	By (\cite{Gi} III 1.4.6(iii)), we deduce a $v^{*}(\rT_{V/k})$-equivariant isomorphism of $Z_{\zar}$ \eqref{affine inverse image num}:
	\begin{equation} \label{gs + inverse image}
		g_{s}^{+}(\mathscr{R}_{\XX,(U,\mathfrak{T},u)})\xrightarrow{\sim} \mathscr{R}_{\XX,(V,\mathfrak{Z},v)}.
	\end{equation}

	One verifies that the data $\{u_{*}(\mathscr{R}_{\XX,(U,\mathfrak{T},u)}),\gamma_{g}\}$ satisfy the compatibility conditions of \ref{prop presheaf to datas}. Then it defines a sheaf of $\widetilde{\pCRIS}$ that we denote by $\mathscr{R}_{\XX}$ \eqref{sheaf to datas}. 
	We set $\mathscr{T}_{X/k}=\nu^{*}(\rT_{X/k})$ (\ref{morphism topos nu}, \ref{stratification id Higgs nul}). In view of \ref{def sheaf L Zar} and \eqref{gamma gs cocycle Lifting}, $\mathscr{R}_{\XX}$ is a $\mathscr{T}_{X/k}$-torsor of $\widetilde{\pCRIS}$. 
\end{nothing}

\begin{prop}\label{FX Oyama topos}
	\textnormal{(i)} There exists a quasi-coherent crystal $\mathscr{F}_{\XX}$ of $\mathscr{O}_{\mathscr{E},1}$-modules of $\widetilde{\pCRIS}$ such that:
	\begin{itemize}
		\item[(a)] For every object $(U,\mathfrak{T},u)$ of $\pCRIS$, $\mathscr{F}_{\XX,(U,\mathfrak{T})}$ \eqref{crystal lin descent data} is the sheaf of affine functions on the $u^{*}(\rT_{U/k})$-torsor $\mathscr{R}_{\XX,(U,\mathfrak{T},u)}$ of $T_{\zar}$ \eqref{affine functions}.

		\item[(b)] For every morphism $g:(V,\mathfrak{Z})\to (U,\mathfrak{T})$ of $\pCRIS$, any affine function $l: \mathscr{R}_{\XX,(U,\mathfrak{T},u)}\to \mathscr{O}_{T}$ and any section $h\in \mathscr{R}_{\XX,(U,\mathfrak{T},u)}(T)$, the transition morphism $c_{g}:g_{s}^{*}(\mathscr{F}_{\XX,(U,\mathfrak{T})})\xrightarrow{\sim} \mathscr{F}_{\XX,(V,\mathfrak{Z})}$ \eqref{crystal lin descent data} sends $g_{s}^{*}(l)$ to an affine function $l':\mathscr{R}_{\XX,(V,\mathfrak{Z},v)}\to \mathscr{O}_{Z}$ such that 
	\begin{equation} \label{calcul cg affine functions}
		l'(h\circ g_{2})=g_{s}^{*}(l(h)) \in \mathscr{O}_{Z}.
	\end{equation}
	\end{itemize}
	
	\textnormal{(ii)} We have an exact sequence of crystals \eqref{morphism topos nu}:
	\begin{equation} \label{suite exact crystal of affine}
		0\to \mathscr{O}_{\mathscr{E},1}\to \mathscr{F}_{\XX}\to \nu^{*}(\Omega_{X/k}^{1})\to 0. 
	\end{equation}
\end{prop}
\begin{proof} (i) For any object $(U,\mathfrak{T},u)$ of $\pCRIS$, we define $\mathscr{F}_{\XX,(U,\mathfrak{T})}$ as in (i). Recall \eqref{exact seq affine functions} that we have an exact sequence of $\mathscr{O}_{T}$-modules of $T_{\zar}$
	\begin{equation} \label{exact seq affine functions Oyama}
		0\to \mathscr{O}_{T}\xrightarrow{c} \mathscr{F}_{\XX,(U,\mathfrak{T})}\xrightarrow{\omega} u^{*}(\Omega^{1}_{U/k})\to 0.
	\end{equation}
	For any morphism $g:(V,\mathfrak{Z})\to (U,\mathfrak{T})$ of $\pCRIS$, by \eqref{iso pullback tosor sheaf of affine funs} and \eqref{gs + inverse image}, we obtain an $\mathscr{O}_{Z}$-linear isomorphism 
	\begin{equation} \label{cg for F pCRIS}
		c_{g}:g_{s}^{*}(\mathscr{F}_{\XX,(U,\mathfrak{T})})\xrightarrow{\sim} \mathscr{F}_{\XX,(V,\mathfrak{Z})}
	\end{equation}
	which fits into a commutative diagram \eqref{diag extension compatible} 
	\begin{equation} \label{diag commut suite exact of aff fun}
		\xymatrix{
			0\ar[r] & g_{s}^{*}(\mathscr{O}_{T}) \ar[r] \ar[d]^{\wr}& g_{s}^{*}(\mathscr{F}_{\XX,(U,\mathfrak{T})}) \ar[r]^{\omega} \ar[d]_{c_{g}} & g_{s}^{*}(u^{*}(\Omega_{U/k}^{1})) \ar[r] \ar[d]^{\wr}& 0 \\
			0\ar[r] & \mathscr{O}_{Z} \ar[r] & \mathscr{F}_{\XX,(V,\mathfrak{Z})} \ar[r]^{\omega} & v^{*}(\Omega_{V/k}^{1}) \ar[r]& 0 \\
		}
	\end{equation}

	In view of the compatibility conditions of $\gamma_{g}$ \eqref{gamma gs cocycle Lifting} and (\cite{AGT} II 4.15), the data $\{\mathscr{F}_{\XX,(U,\mathfrak{T})},c_{g}\}$ satisfy the compatibility conditions of \ref{lin descent data}. Hence, they define a quasi-coherent crystal of $\mathscr{O}_{\pCRIS,1}$-modules of $\widetilde{\pCRIS}$ that we denote by $\mathscr{F}_{\XX}$. The equality \eqref{calcul cg affine functions} follows from \eqref{diag l inverse affine} and \eqref{pre gamma gs cocycle Lifting}.
	
	The assertion (ii) follows from \ref{module 0 to crystal} and the diagram \eqref{diag commut suite exact of aff fun}.
\end{proof}
\begin{nothing} \label{splitting transition}
	We denote by $\mathscr{B}_{\XX}$ the quasi-coherent crystal of $\mathscr{O}_{\pCRIS,1}$-algebras $\varinjlim_{n\ge 1} \rS_{\mathscr{O}_{\mathscr{E},1}}^n(\mathscr{F})$ of $\widetilde{\pCRIS}$. 
	We set $\mathcal{F}_{\XX}=\mathscr{F}_{\XX,(X,\XX)}$ and $\mathcal{B}_{\XX}=\mathscr{B}_{\XX,(X,\XX)}$. Then we obtain an $\mathcal{R}_{\XX}$-stratification $\varepsilon_{\mathcal{F}}$ on $\mathcal{F}_{\XX}$ (resp. $\varepsilon_{\mathcal{B}}$ on $\mathcal{B}_{\XX}$) and a $\widehat{\Gamma}(\rT_{X/k})$-module structure $\psi_{\mathcal{F}}$ on $\mathcal{F}_{\XX}$ (resp. $\psi_{\mathcal{B}}$ on $\mathcal{B}_{\XX}$). 	
	On the other hand, $\mathcal{F}_{\XX}$ being the sheaf of affine functions on the $\rT_{X/k}$-torsor $\mathscr{R}_{\XX,(X,\XX)}$, $\mathcal{F}_{\XX}$ (resp. $\mathcal{B}_{\XX}$) is equipped by \ref{Gamma action on A} with a $\widehat{\Gamma}(\rT_{X/k})$-module structure that we denote by $\kappa_{\mathcal{F}}$ (resp. $\kappa_{\mathcal{B}}$). 
	
	In the following, we show that $\psi_{\mathcal{F}}$ and $\kappa_{\mathcal{F}}$ are different by a sign \eqref{GammaT FBXX}. 

	Recall that we have an exact sequence \eqref{exact seq affine functions Oyama}
	\begin{equation}
		0\to \mathscr{O}_{X}\to \mathcal{F}_{\XX} \to \Omega_{X/k}^{1}\to 0.
	\end{equation}	
	The element $\id_{\XX_{2}}:\XX_{2}\to \XX_{2}$ of $\mathscr{R}_{\XX,(X,\XX)}(X)$ induces a canonical splitting
	\begin{equation} \label{canonical splitting of FX}
		s_{\id}: \mathcal{F}_{\XX}\xrightarrow{\sim} \mathscr{O}_{X}\oplus \Omega_{X/k}^{1}, \qquad l\mapsto (l(\id), \omega(l)). 
	\end{equation}
	Then it induces an isomorphism of $\mathscr{O}_{X}$-algebras
	\begin{equation} \label{RXX Somega}
		\mathcal{B}_{\XX}\xrightarrow{\sim} \rS(\Omega_{X/k}^{1}).
	\end{equation}
	By \ref{action on A splitting}, the isomorphism $s_{\id}$ \eqref{canonical splitting of FX} (resp. \eqref{RXX Somega}) is compatible with the action $\kappa_{\mathcal{F}}$ (resp. $\kappa_{\mathcal{B}}$) and the canonical action of $\widehat{\Gamma}(\rT_{X/k})$ on $\mathscr{O}_{X}\oplus \Omega_{X/k}^{1}$ (resp. $\rS(\Omega_{X/k}^{1})$) \eqref{action Gamma T on SO}. 

\end{nothing}

\begin{prop}\label{calcul stratification algebra B}
	Assume that there exists an \'etale $\SS$-morphism $\XX\to \widehat{\mathbb{A}}_{\SS}^{d}=\Spf(\rW\{T_{1},\cdots,T_{d}\})$. We take again the notation of \ref{local des RQ 1}. 
	
	\textnormal{(i)} If $l$ is a section of $\mathcal{F}_{\XX}$ such that $\omega(l)=0$ \eqref{canonical splitting of FX}, then $\varepsilon_{\mathcal{F}}(1\otimes l)=l\otimes 1$.

	\textnormal{(ii)} For $1\le i\le d$, let $l_{i}$ be the section of $\mathcal{F}_{\XX}$ such that $l_{i}(\id)=0$ and $\omega(l_{i})=dt_{i}$. Then $\varepsilon_{\mathcal{F}}(1\otimes l_{i})=l_{i}\otimes 1- 1\otimes \zeta_{i}$.
\end{prop}
\begin{proof} Assertion (i) follows from \ref{stratification id Higgs nul} and \ref{FX Oyama topos}(ii). 

(ii) We denote by $q_{s}:\RR_{\XX,1}\to X$ the canonical morphism. The canonical morphisms $q_{1,2},q_{2,2}:\RR_{\XX,2}\to \XX_{2}\in \mathscr{R}_{\XX,(X,\RR_{\XX},q_{s})}(\RR_{\XX,1})$ induce two splittings of $\mathcal{F}_{\XX,(X,\RR_{\XX})}$ \eqref{exact seq affine functions Oyama} 
	\begin{eqnarray}
		&s_{q_{1}}:q_{s*}(\mathscr{F}_{\XX,(X,\RR_{\XX})})\xrightarrow{\sim} \mathcal{R}_{\XX,1}\oplus (\mathcal{R}_{\XX,1}\otimes_{\mathscr{O}_{X}}\Omega_{X/k}^{1}) \quad & f\mapsto (f(q_{1}),\omega(f))\\
		&s_{q_{2}}:q_{s*}(\mathscr{F}_{\XX,(X,\RR_{\XX})})\xrightarrow{\sim} \mathcal{R}_{\XX,1}\oplus (\mathcal{R}_{\XX,1}\otimes_{\mathscr{O}_{X}}\Omega_{X/k}^{1}) \quad & f\mapsto (f(q_{2}),\omega(f)).
	\end{eqnarray}
	We identify $\mathcal{R}_{\XX,1}\oplus (\mathcal{R}_{\XX,1}\otimes_{\mathscr{O}_{X}}\Omega_{X/k}^{1})$ with $\mathcal{R}_{\XX,1}\otimes_{\mathscr{O}_{X}}\mathcal{F}_{\XX}$ by $\id_{\mathcal{R}_{\XX,1}}\otimes s_{\id}$ \eqref{canonical splitting of FX}. In view of \eqref{calcul cg affine functions} and \eqref{diag commut suite exact of aff fun}, the morphism $s_{q_{i}}$ is inverse to $q_{s*}(c_{q_{i}})$ for $i=1,2$. 
	Then, the $\mathcal{R}_{\XX}$-stratification $\varepsilon$ on $\mathcal{F}_{\XX}$ is given by the composition of the inverse of $s_{q_{2}}$ and $s_{q_{1}}$.

We denote by $f$ the local section $s_{q_{2}}^{-1}(1\otimes l_{i})$ of $q_{s*}(\mathscr{F}_{\XX,(X,\RR_{\XX})})$. Then we have $f(q_{2})=0$ and $\omega(f)=1\otimes dt_{i}$. To show the assertion, it suffices to prove that $f(q_{1})=-\zeta_{i}\otimes 1$. 
The morphisms $q_{1},q_{2}$ are induced by two homomorphisms
\begin{equation}
	\iota_{1}:\mathscr{O}_{\XX_{2}}\to \mathcal{R}_{\XX,2},\quad \iota_{2}:\mathscr{O}_{\XX_{2}}\to \mathcal{R}_{\XX,2}
\end{equation}
such that $\iota_{1}$ is equal to $\iota_{2}$ modulo $p$. Then they define a $\rW_{2}$-derivation
\begin{equation}
	D=\iota_{2}-\iota_{1}:\mathscr{O}_{\XX_{2}}\to p\mathcal{R}_{\XX,2}.
\end{equation}
We denote by 
\begin{equation}
	\phi:\mathcal{R}_{\XX,1}\otimes_{\mathscr{O}_{X}}\Omega_{X/k}^{1}\to p\mathcal{R}_{\XX,2}
\end{equation}
the $\mathcal{R}_{\XX,1}$-linear morphism associated to $D$. Then we have 
\begin{displaymath}
	\phi(1\otimes dt_{i})=1\otimes t_{i}-t_{i}\otimes 1= p\biggl(\frac{\xi_{i}}{p}\biggr) \in p\mathcal{R}_{\XX,2}.
\end{displaymath}
Identifying $\FHom_{\mathscr{O}_{X}}(\mathcal{R}_{\XX,1}\otimes_{\mathscr{O}_{X}}\Omega_{X/k}^{1}, p\mathcal{R}_{\XX,2})\xrightarrow{\sim} \mathcal{R}_{\XX,1}\otimes_{\mathscr{O}_{X}} \rT_{X/k}$, we consider $\phi$ as a section of $q_{s}^{*}(\rT_{X/k})$. Then we have $q_{2}=q_{1}+\phi \in \mathscr{R}_{\XX,q_{s}}(\RR_{\XX,1})$ and we deduce that (\ref{affine functions}(ii))
\begin{eqnarray*}
	f(q_{1})= f(q_{2}-\phi)
		= -\omega(f)(\phi)
		= -\zeta_{i}\otimes 1.
\end{eqnarray*}
Then the proposition follows.
\end{proof}
\begin{coro} \label{GammaT FBXX}
	The $\widehat{\Gamma}(\rT_{X/k})$-actions $\psi_{\mathcal{F}}$ and $\kappa_{\mathcal{F}}$ on $\mathcal{F}_{\XX}$ (resp. $\psi_{\mathcal{B}}$ and $\kappa_{\mathcal{B}}$ on $\mathcal{B}_{\XX}$) satisfy $\psi_{\mathcal{F}}=\iota^{*}(\kappa_{\mathcal{F}})$ (resp. $\psi_{\mathcal{B}}=\iota^{*}(\kappa_{\mathcal{B}})$), where $\iota:\widehat{\Gamma}(\rT_{X/k})\to \widehat{\Gamma}(\rT_{X/k})$ denotes the involution homomorphism \eqref{sigma Gamma T}. 
\end{coro}
\begin{proof} The question being local, we take again the assumptions and notation of \ref{calcul stratification algebra B}. 
Let $\partial_{i} \in \rT_{X/k}(X)$ be the dual of $dt_{i}$. In view of \eqref{R stratification def} and \ref{calcul stratification algebra B}(ii), the action of $\psi_{\mathcal{F},\partial_{i}}$ on $\mathcal{F}_{\XX}$ sends $l_{i}$ to $-1$. We equip $\mathscr{O}_{X}\oplus \Omega_{X/k}^{1}$ (resp. $\rS(\Omega_{X/k}^{1})$) with the canonical action of $\widehat{\Gamma}(\rT_{X/k})$ \eqref{action Gamma T on SO}. The isomorphism $s_{\id}$ \eqref{canonical splitting of FX} (resp. \eqref{RXX Somega}) sends $l_{i}$ to $dt_{i}$ and induces an isomorphism of $\widehat{\Gamma}(\rT_{X/k})$-modules
	\begin{equation} \label{BX iso iota SO}
		(\mathcal{F}_{\XX},\psi_{\mathcal{F}}) \xrightarrow{\sim} \iota^{*}(\mathscr{O}_{X}\oplus \Omega_{X/k}^{1}) \qquad \textnormal{(resp. }((\mathcal{B}_{\XX},\psi_{\mathcal{B}}) \xrightarrow{\sim} \iota^{*}(\rS(\Omega_{X/k}^{1})) ).
	\end{equation}
	Then the assertion follows from \ref{splitting transition}.
\end{proof}

\begin{coro} \label{calculs of varepsilonR}
	We denote by $\varepsilon_{R}^{+}$ (resp. $\varepsilon_{R}^{-}$) the $\mathcal{R}_{\XX}$-stratification on $\mathcal{R}_{\XX,1}$ associated to the $\widehat{\Gamma}(\rT_{X/k})$-action $\kappa_{\mathcal{B}}$ (resp. $\psi_{\mathcal{B}}$) on $\mathcal{B}_{\XX}$ via the isomorphisms of $\mathscr{O}_{X}$-algebras $\mathcal{B}_{\XX}\xrightarrow{\sim} \rS(\Omega_{X/k}^{1})\xrightarrow{\sim} \mathcal{R}_{\XX,1}$ \textnormal{(\ref{dual of R}, \ref{RXX Somega})}

	\textnormal{(i)} For any local section $r$ of $\mathcal{R}_{\XX,1}$, we have $\varepsilon_{R}^{+}(1\otimes r)=\delta(r)$ \eqref{local des RQ 1}.

	\textnormal{(ii)} For any local section $r$ of $\mathcal{R}_{\XX,1}$, we have $\varepsilon_{R}^{-}(\delta(r))=r\otimes 1$. 
\end{coro}

\begin{proof} The question being local, we take again the assumption and the notation of \ref{calcul stratification algebra B}. Then for $1\le i\le d$, $l_{i}$ is sent to $\zeta_{i}$ by the isomorphism $\mathcal{B}_{\XX}\xrightarrow{\sim} \mathcal{R}_{\XX,1}$. Since $\varepsilon_{R}^{+}$, $\varepsilon_{R}^{-}$ and $\delta$ are homomorphisms, it suffices to verify the assertion for the local sections $\zeta_{i}$. By \ref{calcul stratification algebra B}(ii) and \ref{GammaT FBXX}(i), we have
\begin{equation}
	\varepsilon_{R}^{-}(1\otimes\zeta_{i})=\zeta_{i}\otimes 1-1\otimes \zeta_{i},\qquad \varepsilon_{R}^{+}(1\otimes \zeta_{i})=\zeta_{i}\otimes 1+1\otimes\zeta_{i}
\end{equation}
Then, the assertion (i) follows from \ref{local des RQ 1}. The assertion (ii) follows from the relations
	\begin{displaymath}
		\varepsilon_{R}^{-}(\delta(\zeta_{i}))=\varepsilon_{R}^{-}(1\otimes \zeta_{i}+\zeta_{i}\otimes 1)=(\zeta_{i}\otimes 1-1\otimes \zeta_{i}) +1\otimes \zeta_{i}=\zeta_{i}\otimes 1.
	\end{displaymath}
\end{proof}

With the above preparation, we can interpret a key calculation in Oyama's paper in the following forms. 

\begin{prop}[\cite{Oy} 1.5.3]\label{admissibility R}
	Let $M$ be a quasi-nilpotent $\widehat{\Gamma}(\rT_{X/k})$-module, $\varepsilon$ the associated $\mathcal{R}_{\XX}$-stratification on $M$ \eqref{equi Dmod strat}. We denote by $M_{0}$ the underlying $\mathscr{O}_{X}$-module of $M$ equipped with the $\widehat{\Gamma}(\rT_{X/k})$-module structure defined by the zero PD-Higgs field. Then the stratification $\varepsilon:\mathcal{R}_{\XX,1}\otimes_{\mathscr{O}_{X}}M\to M\otimes_{\mathscr{O}_{X}}\mathcal{R}_{\XX,1}$ induces two isomorphisms of $\widehat{\Gamma}(\rT_{X/k})$-modules \eqref{Hom psi}
	\begin{eqnarray*}
		\textnormal{(i)}~ \qquad \qquad (\mathcal{B}_{\XX},\psi_{\mathcal{B}})\otimes_{\mathscr{O}_{X}}M_{0}&\xrightarrow{\sim}& M\otimes_{\mathscr{O}_{X}}(\mathcal{B}_{\XX},\psi_{\mathcal{B}}), \\
		\textnormal{(ii)} \qquad \quad (\mathcal{B}_{\XX},\psi_{\mathcal{B}})\otimes_{\mathscr{O}_{X}}\iota^{*}(M)&\xrightarrow{\sim}& M_{0}\otimes_{\mathscr{O}_{X}}(\mathcal{B}_{\XX},\psi_{\mathcal{B}}).
	\end{eqnarray*} 
\end{prop}

\begin{proof} We take again the notation of \ref{calculs of varepsilonR}. To simplify the notation, we write $\mathcal{R}$ for $\mathcal{R}_{\XX,1}$. The $\mathcal{R}_{\XX}$-stratification on $M_{0}$ is the identity morphism $\id_{\mathcal{R}\otimes M}:\mathcal{R}\otimes_{\mathscr{O}_{X}}M\to \mathcal{R}\otimes_{\mathscr{O}_{X}} M$ (cf. the proof of \ref{stratification id Higgs nul}). We denote by $\theta_{0}$ (resp. $\theta$) the morphism $M\to M\otimes_{\mathscr{O}_{X}}\mathcal{R}$ defined by $m\mapsto m\otimes 1$ (resp. $m\mapsto \varepsilon(1\otimes m)$) for every local section $m$ of $M$. 

(i) Since the action $\psi_{\mathcal{B}}$ on $\mathcal{B}_{\XX}$ is compatible with the ring structure of $\mathcal{B}_{\XX}$, it suffices to show that $\theta:M_{0}\to M\otimes(\mathcal{B}_{\XX},\psi_{\mathcal{B}})$ is $\widehat{\Gamma}(\rT_{X/k})$-equivariant. 
In view of \eqref{stratification tensor product}, it suffices to prove that the following diagram is commutative
\begin{equation}
	\xymatrix{
		M \ar[d]_{\theta} \ar[rr]^{\theta_{0}} && M\otimes_{\mathscr{O}_{X}}\mathcal{R} \ar[d]^{\theta\otimes \id_{\mathcal{R}}}\\
		M\otimes_{\mathscr{O}_{X}}\mathcal{R} \ar[r]^-{\theta\otimes \id_{\mathcal{R}}} & M\otimes_{\mathscr{O}_{X}}\mathcal{R}\otimes_{\mathscr{O}_{X}}\mathcal{R} \ar[r]^-{\id_{M}\otimes \varepsilon_{R}^{-}}& M\otimes_{\mathscr{O}_{X}}\mathcal{R}\otimes_{\mathscr{O}_{X}}\mathcal{R} }
\end{equation}
By condition (ii) of \ref{lemma stratification}, the composition
\begin{equation}
	\xymatrix{
		M\ar[r]^-{\theta} & M\otimes_{\mathscr{O}_{X}}\mathcal{R} \ar[r]^-{\theta\otimes \id_{\mathcal{R}}}& M\otimes_{\mathscr{O}_{X}}\mathcal{R}\otimes_{\mathscr{O}_{X}}\mathcal{R} \ar[r]^-{\id_{M}\otimes \varepsilon_{R}^{-}}& M\otimes_{\mathscr{O}_{X}}\mathcal{R}\otimes_{\mathscr{O}_{X}}\mathcal{R}. }	
\end{equation}
in the above diagram is equal to the composition
\begin{equation}\label{composition of theta theta0}
	\xymatrix{
		M\ar[r]^-{\theta} & M\otimes_{\mathscr{O}_{X}}\mathcal{R} \ar[r]^-{\id_{M}\otimes\delta}& M\otimes_{\mathscr{O}_{X}}\mathcal{R}\otimes_{\mathscr{O}_{X}}\mathcal{R} \ar[r]^-{\id_{M}\otimes \varepsilon_{R}^{-}}& M\otimes_{\mathscr{O}_{X}}\mathcal{R}\otimes_{\mathscr{O}_{X}}\mathcal{R}.
	}
\end{equation}
By \ref{calculs of varepsilonR}(ii), the above composition is equal to the composition
\begin{equation}
	M\xrightarrow{\theta_{0}} M\otimes_{\mathscr{O}_{X}}\mathcal{R}\xrightarrow{\theta\otimes \id_{\mathcal{R}}} M\otimes_{\mathscr{O}_{X}}\mathcal{R}\otimes_{\mathscr{O}_{X}}\mathcal{R}.
\end{equation}
The assertion follows. 

(ii) By \ref{GammaT FBXX}(i), it suffices to show that the morphism $(\mathcal{B}_{\XX},\kappa)\otimes_{\mathscr{O}_{X}}M\to M_{0}\otimes_{\mathscr{O}_{X}}(\mathcal{B}_{\XX},\kappa_{\mathcal{B}})$ is $\widehat{\Gamma}(\rT_{X/k})$-equivariant. Similarly to (i), by \eqref{stratification tensor product}, it suffices to prove that the following diagram is commutative
\begin{equation}
	\xymatrix{
		M \ar[d]_{\theta} \ar[rr]^{\theta} && M\otimes_{\mathscr{O}_{X}}\mathcal{R} \ar[d]^{\theta\otimes \id_{\mathcal{R}}}\\
		M\otimes_{\mathscr{O}_{X}}\mathcal{R} \ar[r]^-{\theta_{0}\otimes \id_{\mathcal{R}}} & M\otimes_{\mathscr{O}_{X}}\mathcal{R}\otimes_{\mathscr{O}_{X}}\mathcal{R} \ar[r]^-{\id_{M}\otimes \varepsilon_{R}^{+}}& M\otimes_{\mathscr{O}_{X}}\mathcal{R}\otimes_{\mathscr{O}_{X}}\mathcal{R} }
\end{equation}
By condition (ii) of \ref{lemma stratification}, the composition
\begin{equation}
	\xymatrix{
		M\ar[r]^-{\theta} & M\otimes_{\mathscr{O}_{X}}\mathcal{R} \ar[r]^-{\theta\otimes \id_{\mathcal{R}}}& M\otimes_{\mathscr{O}_{X}}\mathcal{R}\otimes_{\mathscr{O}_{X}}\mathcal{R} }	
\end{equation}
in the above diagram is equal to the composition
\begin{equation}
	\xymatrix{
		M\ar[r]^-{\theta} & M\otimes_{\mathscr{O}_{X}}\mathcal{R} \ar[r]^-{\id_{M}\otimes\delta}& M\otimes_{\mathscr{O}_{X}}\mathcal{R}\otimes_{\mathscr{O}_{X}}\mathcal{R}.}
\end{equation}
The assertion follows from \ref{calculs of varepsilonR}(i).
\end{proof}

\begin{nothing} \label{C inverse underlying}
	The counterpart of the crystal $\mathscr{A}_{\XX_{2}'}$ \eqref{alg splitting mod} in Oyama topos $\widetilde{\underline{\mathscr{E}}}$ is defined by applying Cartier equivalence $\rmC^{*}_{X/\SS}$ to $\mathscr{B}_{\XX'}$. 

	Let $\rmC_{X/\SS}:\widetilde{\CRIS}\to \widetilde{\pCRIS}'$ be the morphism of topoi \eqref{morphism of topoi Cartier}. We put $\underline{\mathscr{R}}_{\XX}=\rmC^{*}_{X/\SS}(\mathscr{R}_{\XX'})$ \eqref{Relevement X}. For any object $(U,\mathfrak{T},u)$ of $\CRIS$, we set $\phi_{T/k}=u'\circ f_{T/k}: T\to \underline{T}'\to U'$ so we have $\rho(U,\mathfrak{T},u)=(U',\mathfrak{T},\phi_{T/k})$ \eqref{functor rho}. By \ref{morphism topos C}, the descent data of the sheaf $\underline{\mathscr{R}}_{\XX}$ of $\widetilde{\CRIS}$ is $\{u_{*}(\mathscr{R}_{\XX',(U',\mathfrak{T},\phi_{T/k})}),\gamma_{\rho(g)}\}$ (\ref{prop presheaf to datas}, \ref{def sheaf L Zar}, \ref{Relevement X}). 
	
	We put $\underline{\mathscr{F}}_{\XX}=\rmC^{*}_{X/\SS}(\mathscr{F}_{\XX'})$ and $\underline{\mathscr{B}}_{\XX}=\rmC_{X/\SS}^{*}(\mathscr{B}_{\XX'})$. For any object $(U,\mathfrak{T},u)$ of $\CRIS$, we have \eqref{inverse image qc}
	\begin{equation} \label{underlineF F}
		\underline{\mathscr{F}}_{\XX,(U,\mathfrak{T},u)}=\mathscr{F}_{\XX',(U',\mathfrak{T},\phi_{T/k})},\qquad \underline{\mathscr{B}}_{\XX,(U,\mathfrak{T},u)}=\mathscr{B}_{\XX',(U',\mathfrak{T},\phi_{T/k})}.
	\end{equation}
	The linearised descent data of the quasi-coherent crystal of $\mathscr{O}_{\CRIS,1}$-modules $\underline{\mathscr{F}}_{\XX}$ (resp. $\mathscr{O}_{\CRIS,1}$-algebras $\underline{\mathscr{B}}_{\XX}$) of $\widetilde{\CRIS}$ is $\{\mathscr{F}_{\XX',(U',\mathfrak{T},\phi_{T/k})},c_{\rho(g)}\}$ (resp. $\{\mathscr{B}_{\XX',(U',\mathfrak{T},\phi_{T/k})},c_{\rho(g)}\}$) (\ref{morphism topos C}, \ref{lin descent data}, \ref{FX Oyama topos}). 

	We set $\underline{\mathcal{F}}_{\XX}=\underline{\mathscr{F}}_{\XX,(X,\XX)}$ and $\underline{\mathcal{B}}_{\XX}=\underline{\mathscr{B}}_{\XX,(X,\XX)}$. By \ref{equi of Dmod crystals}, these $\mathscr{O}_{X}$-modules are equipped with $\rD_{X/k}^{\gamma}$-module structures that we denote by $(\nabla_{\underline{\mathcal{F}}},\psi_{\underline{\mathcal{F}}})$ (resp. $(\nabla_{\underline{\mathcal{B}}},\psi_{\underline{\mathcal{B}}})$). 
	We will show that $\underline{\mathcal{B}}_{\XX}$ is isomorphic to the algebra $\mathcal{A}_{\XX_{2}'}$ \eqref{alg splitting mod} introduced by Ogus--Vologodsky \eqref{crystal of affine functions concides}. 
\end{nothing}

\begin{nothing} \label{prep iso ad}
	We interpret the Cartier equivalence as an admissible isomorphism à la Fontaine for a pair of crystals with respect to the period ring $\underline{\mathcal{B}}_{\XX}$. 

	Let $\mathscr{M}'$ be a crystal of $\mathscr{O}_{\pCRIS',1}$-modules of $\widetilde{\pCRIS}'$, $(M',\theta'_{M})$ the associated $\widehat{\Gamma}(\rT_{X'/k})$-module, $\mathscr{M}'_{0}=\nu^{*}(M')$ \eqref{module 0 to crystal}, $\mathscr{N}=\rmC^{*}_{X/\SS}(\mathscr{M}')$ and $(N,\nabla_{N},\psi_{N})$ the associated $\rD_{X/k}^{\gamma}$-module. By \ref{equi of Dmod crystals}, \ref{stratification id Higgs nul} and \ref{admissibility R}(i), we have an isomorphism of crystals of $\mathscr{O}_{\mathscr{E}',1}$-modules of $\widetilde{\pCRIS}'$
	\begin{equation} \label{iso compatibility crystal 1}
		\mathscr{B}_{\XX'}\otimes_{\mathscr{O}_{\pCRIS',1}}\mathscr{M}'_{0} \xrightarrow{\sim} \mathscr{M}'\otimes_{\mathscr{O}_{\pCRIS',1}}\mathscr{B}_{\XX'}.
	\end{equation}
	Applying $\rmC^{*}_{X/\SS}$, we obtain an isomorphism of crystals of $\mathscr{O}_{\CRIS,1}$-modules of $\widetilde{\CRIS}$:
	\begin{equation} \label{iso compatibility crystal 2}
		\underline{\mathscr{B}}_{\XX}\otimes_{\mathscr{O}_{\CRIS,1}}\rmC^{*}_{X/\SS}(\mathscr{M}'_{0})\xrightarrow{\sim} \mathscr{N}\otimes_{\mathscr{O}_{\CRIS,1}}\underline{\mathscr{B}}_{\XX}.
	\end{equation}
	By \ref{Cartier transform descent}, we deduce an isomorphism of $\rD_{X/k}^{\gamma}$-modules
	\begin{equation} \label{admissibility M N}
		\lambda: \underline{\mathcal{B}}_{\XX}\otimes_{\mathscr{O}_{X}}F_{X/k}^{*}(M')\xrightarrow{\sim} N\otimes_{\mathscr{O}_{X}}\underline{\mathcal{B}}_{\XX}
	\end{equation}
	where the $\rD_{X/k}^{\gamma}$-action on the left hand side is induced by $(\nabla_{\underline{\mathcal{B}}},\psi_{\underline{\mathcal{B}}})$ on $\underline{\mathcal{B}}_{\XX}$ and $(\nabla_{\can},0)$ on $F_{X/k}^{*}(M')$ \eqref{FHom Dgamma}, that we still denote by $(\nabla_{\underline{\mathcal{B}}},\psi_{\underline{\mathcal{B}}})$; and the $\rD_{X/k}^{\gamma}$-action on the right hand side is induced by $(\nabla_{N},\psi_{N})$ on $N$ and $(\nabla_{\underline{\mathcal{B}}},\psi_{\underline{\mathcal{B}}})$, that we denote by $(\nabla_{\tot},\psi_{\tot})$. 

	The $\widehat{\Gamma}(\rT_{X'/k})$-module structures $\psi_{\underline{\mathcal{B}}}$ on $\underline{\mathcal{B}}_{\XX}$ and $\theta_{M'}$ on $M'$ define a $\widehat{\Gamma}(\rT_{X'/k})$-module structure on $\underline{\mathcal{B}}_{\XX}\otimes_{\mathscr{O}_{X'}}M'$ \eqref{Hom PD Higgs}, that we denote by $\theta_{\tot}$. On the other hand, the zero PD-Higgs field on $N$ and the action of $\psi_{\underline{\mathcal{B}}}$ on $\underline{\mathcal{B}}_{\XX}$ define a $\widehat{\Gamma}(\rT_{X'/k})$-module structure on $N\otimes_{\mathscr{O}_{X}}\underline{\mathcal{B}}_{\XX}$, that we still denote by $\psi_{\underline{\mathcal{B}}}$.
\end{nothing}

\begin{theorem} \label{thm iso ad}
	Let $\mathscr{M}'$ be a crystal of $\mathscr{O}_{\pCRIS',1}$-modules of $\widetilde{\pCRIS}'$, $(M',\theta'_{M})$ the associated $\widehat{\Gamma}(\rT_{X'/k})$-module, $\mathscr{N}=\rmC^{*}_{X/\SS}(\mathscr{M}')$ and $(N,\nabla_{N},\psi_{N})$ the associated $\rD_{X/k}^{\gamma}$-module. 
	Then, the isomorphism \eqref{admissibility M N}
	\begin{equation} \label{iso ad M N}
		\lambda: (\underline{\mathcal{B}}_{\XX}\otimes_{\mathscr{O}_{X'}}M',(\nabla_{\underline{\mathcal{B}}},\psi_{\underline{\mathcal{B}}}),\theta_{\tot})\xrightarrow{\sim} (N\otimes_{\mathscr{O}_{X}}\underline{\mathcal{B}}_{\XX},(\nabla_{\tot},\psi_{\tot}), \psi_{\underline{\mathcal{B}}})
	\end{equation}
	is compatible with the $\rD_{X/k}^{\gamma}$-actions and the $\widehat{\Gamma}(\rT_{X'/k})$-actions defined on both sides in \ref{prep iso ad}. 
\end{theorem}
\begin{proof} We only need to prove the compatibility of the $\widehat{\Gamma}(\rT_{X'/k})$-actions. Let $\varepsilon':\mathcal{B}_{\XX'}\otimes_{\mathscr{O}_{X'}}M'\xrightarrow{\sim} M'\otimes_{\mathscr{O}_{X'}}\mathcal{B}_{\XX'}$ be the $\mathcal{R}_{\XX'}$-stratification on $M'$ \eqref{admissibility R}, and $\psi_{\mathcal{B}'}$ the $\widehat{\Gamma}(\rT_{X'/k})$-action on $\mathcal{B}_{\XX'}$ defined in \ref{splitting transition}. 
The question is local. Since the Cartier equivalence is compatible with localisation \eqref{Cartier localisation}, we can suppose that there exists a lifting $F:\XX\to \XX'$ of the relative Frobenius morphism $F_{X/k}$ of $X$. Then it induces a morphism $F:(X',\XX,F_{X/k})\to (X',\XX')$ of $\pCRIS'$. The isomorphisms \eqref{iso compatibility crystal 1}, \eqref{iso compatibility crystal 2} and the transition morphisms $\eta_{F}$ associated to $F$ \eqref{eta F modules} induce a commutative diagram
\begin{displaymath}
	\xymatrixcolsep{4pc}\xymatrix{
		F_{X/k}^{*}((\mathscr{B}_{\XX'}\otimes_{\mathscr{O}_{\pCRIS',1}}\mathscr{M}'_{0})_{(X',\XX')}) \ar[r]^-{\sim} \ar[d]_{\eta_{F}} & F_{X/k}^{*}((\mathscr{M}'\otimes_{\mathscr{O}_{\pCRIS',1}}\mathscr{B}_{\XX'})_{(X',\XX')}) \ar[d]^{\eta_{F}} \\
		(\underline{\mathscr{B}}_{\XX}\otimes_{\mathscr{O}_{\CRIS,1}}\rmC^{*}_{X/\SS}(\mathscr{M}'_{0}))_{(X,\XX)} \ar[r]^-{\sim}&  (\mathscr{N}\otimes_{\mathscr{O}_{\CRIS,1}}\underline{\mathscr{B}}_{\XX})_{(X,\XX)}		
	}
\end{displaymath}
Then we deduce a commutative diagram
\begin{equation} \label{eq thm ad 1}
	\xymatrixcolsep{4pc}\xymatrix{
		F_{X/k}^{*}(\mathcal{B}_{\XX'}\otimes_{\mathscr{O}_{X'}}M') \ar[r]^-{F_{X/k}^{*}(\varepsilon')} \ar[d]_{\eta_{F}} & F_{X/k}^{*}(M'\otimes_{\mathscr{O}_{X'}}\mathcal{B}_{\XX'}) \ar[d]^{\eta_{F}}  \\
		\underline{\mathcal{B}}_{\XX}\otimes_{\mathscr{O}_{X'}}M' \ar[r]^-{\lambda} & N\otimes_{\mathscr{O}_{X}}\underline{\mathcal{B}}_{\XX}		
	}
\end{equation}
By \ref{lifting F Cartier}, the isomorphism $\eta_{F}: F_{X/k}^{*}(\mathcal{B}_{\XX'})\xrightarrow{\sim} \underline{\mathcal{B}}_{\XX}$ underlies an isomorphism of $F_{X/k}^{*}(\widehat{\Gamma}(\rT_{X'/k}))$-modules
\begin{equation} \label{eq thm ad 2}
	F_{X/k}^{*}(\mathcal{B}_{\XX'},\iota^{*}(\psi_{\mathcal{B}'}))\xrightarrow{\sim} (\underline{\mathcal{B}}_{\XX},\psi_{\underline{\mathcal{B}}}).
\end{equation}
By \ref{admissibility R}(ii), the isomorphism $F_{X/k}^{*}(\varepsilon')$ is compatible with actions of $F_{X/k}^{*}(\widehat{\Gamma}(\rT_{X'/k}))$:
\begin{equation} \label{eq thm ad 3}
	F_{X/k}^{*}(\varepsilon'):F_{X/k}^{*}(\mathcal{B}_{\XX'}\otimes_{\mathscr{O}_{X'}}M',\psi_{\mathcal{B}'}\otimes \id +\id \otimes \iota^{*}( \theta_{M'}))\xrightarrow{\sim} F_{X/k}^{*}(M'\otimes_{\mathscr{O}_{X'}}\mathcal{B}_{\XX'}, \id\otimes\psi_{\mathcal{B}'}).
\end{equation}
Then the assertion follows from \eqref{eq thm ad 1}, \eqref{eq thm ad 2} and \eqref{eq thm ad 3}. 
\end{proof}
\begin{rem}
	In (\cite{OV07}, 2.23), Ogus and Vologodsky showed a similar result for certain filtered PD-Higgs modules. 
\end{rem}

	Let $(H,\theta)$ (resp. $(H,\nabla,\psi)$) be a $\widehat{\Gamma}(\rT_{X/k})$-module (resp. a $\rD_{X/k}^{\gamma}$-module). We define its $\widehat{\Gamma}(\rT_{X/k})$ invariants (resp. $\rD_{X/k}^{\gamma}$ invariants) by 
	\begin{equation}
		H^{\theta}=\FHom_{\widehat{\Gamma}(\rT_{X/k})}( (\mathscr{O}_{X},0), H) \qquad \textnormal{(resp. } H^{(\nabla,\psi)}=\FHom_{\rD_{X/k}^{\gamma}}( (\mathscr{O}_{X},d,0), H) ).
	\end{equation}
	
	We can recover the Cartier equivalence by taking $\widehat{\Gamma}(\rT_{X'/k})$ invariants (resp. $\rD_{X/k}^{\gamma}$ invariants) for the isomorphism \eqref{iso ad M N}.

\begin{prop}\label{D-inv theta-inv}
	Keep the notation of \ref{prep iso ad} and suppose that $\mathscr{M}'$ is quasi-coherent. The isomorphism \eqref{iso ad M N} induces
	\begin{itemize}
		\item[(i)] a canonical isomorphism of $\rD_{X/k}^{\gamma}$-modules 
			\begin{displaymath}
				(\underline{\mathcal{B}}_{\XX}\otimes_{\mathscr{O}_{X'}}M',\nabla_{\underline{\mathcal{B}}},\psi_{\underline{\mathcal{B}}})^{\theta_{\tot}}\xrightarrow{\sim} (N,\nabla_{N},\psi_{N});
			\end{displaymath}
		
		\item[(ii)] a canonical isomorphism of $\widehat{\Gamma}(\rT_{X'/k})$-modules 
			\begin{displaymath}
				(M',\theta_{M'})\xrightarrow{\sim} (N\otimes_{\mathscr{O}_{X}}\underline{\mathcal{B}}_{\XX},\psi_{\underline{\mathcal{B}}})^{(\nabla_{\tot},\psi_{\tot})}.
			\end{displaymath}
			
	\end{itemize}
\end{prop}
The assertion follows from \ref{lemma commute} and \ref{calcul nabla inv}.

\begin{lemma} \label{lemma commute}
	Keep the notation of \ref{prep iso ad}.
	
	\textnormal{(i)} The actions $(\nabla_{\underline{\mathcal{B}}},\psi_{\underline{\mathcal{B}}})$ of $\rD_{X/k}^{\gamma}$ and $\theta_{\tot}$ of $\widehat{\Gamma}(\rT_{X'/k})$ on $\underline{\mathcal{B}}_{\XX}\otimes_{\mathscr{O}_{X'}}M'$ commute with each other.
	
	\textnormal{(ii)} The actions $(\nabla_{\tot},\psi_{\tot})$ of $\rD_{X/k}^{\gamma}$ and $\psi_{\underline{\mathcal{B}}}$ of $\widehat{\Gamma}(\rT_{X'/k})$ on $N\otimes_{\mathscr{O}_{X}}\underline{\mathcal{B}}_{\XX}$ commute with each other.
\end{lemma}
\begin{proof} (i) By the formula \eqref{tensor Higgs}, one verifies that the action $\psi_{\underline{\mathcal{B}}}$ of $\widehat{\Gamma}(\rT_{X'/k})\subset \rD_{X/k}^{\gamma}$ and the action $\theta_{\tot}$ of $\widehat{\Gamma}(\rT_{X'/k})$ on $\underline{\mathcal{B}}_{\XX}\otimes_{\mathscr{O}_{X'}}M'$ commute with each other. For any local sections $D$ of $\rT_{X/k}$, $\xi'$ of $\widehat{\Gamma}(\rT_{X'/k})$, $b$ of $\underline{\mathcal{B}}_{\XX}$ and $m$ of $M'$, by \ref{tensor Higgs}, we have
\begin{eqnarray*}
	\nabla_{\underline{\mathcal{B}},D}(\theta_{\tot,\xi'}(b\otimes m)) &=& \nabla_{\underline{\mathcal{B}},D}(\psi_{\underline{\mathcal{B}},\xi'}(b))\otimes m + \nabla_{\underline{\mathcal{B}},D}(b)\otimes \theta_{M',\xi'}(m) \\
	&=& \psi_{\underline{\mathcal{B}},\xi'}(\nabla_{\underline{\mathcal{B}},D}(b))\otimes m + \nabla_{\underline{\mathcal{B}},D}(b)\otimes \theta_{M',\xi'}(m) \\
	&=& \theta_{\tot,\xi'}(\nabla_{\underline{\mathcal{B}},D}(b\otimes m))
\end{eqnarray*}
Since $\rD_{X/k}^{\gamma}$ is generated by $\rT_{X/k}$ and $\widehat{\Gamma}(\rT_{X'/k})$, the assertion follows.

Assertion (ii) follows from (i) and \ref{thm iso ad}.
\end{proof}

\begin{lemma} \label{calcul nabla inv}
	Keep the notation of \ref{prep iso ad} and suppose that $\mathscr{M}'$ is quasi-coherent. 
	The canonical homomorphism of crystals $\mathscr{O}_{\underline{\mathscr{E}},1}\to \underline{\mathscr{B}}_{\XX}$ induces an isomorphism of $\mathscr{O}_{X}$-modules (resp. $\mathscr{O}_{X'}$-modules)
	\begin{displaymath}
		N\xrightarrow{\sim} (N\otimes_{\mathscr{O}_{X}}\underline{\mathcal{B}}_{\XX})^{\psi_{\underline{\mathcal{B}}}},\qquad \textnormal{(resp. } M'\xrightarrow{\sim} (\underline{\mathcal{B}}_{\XX}\otimes_{\mathscr{O}_{X'}}M')^{(\nabla_{\underline{\mathcal{B}}},\psi_{\underline{\mathcal{B}}})}).
	\end{displaymath}
\end{lemma}
\begin{proof} 
The assertion in lemma being local, we may assume that there exists a lifting $F:\XX\to \XX'$ of the relative Frobenius morphism $F_{X/k}$ of $X$. By \ref{lifting F Cartier} and \eqref{BX iso iota SO}, we have isomorphisms of $F_{X/k}^{*}(\widehat{\Gamma}(\rT_{X'/k}))$-modules
\begin{eqnarray*}
	\eta_{F}: F_{X/k}^{*}(\mathcal{B}_{\XX'}\otimes_{\mathscr{O}_{X'}}M',\iota^{*}(\psi_{\mathcal{B}'})\otimes \id)&\xrightarrow{\sim}& (N\otimes_{\mathscr{O}_{X}}\underline{\mathcal{B}}_{\XX},\psi_{\underline{\mathcal{B}}})\\
	\rS(\Omega_{X/k}^{1})\otimes_{\mathscr{O}_{X}}F_{X/k}^{*}(M') &\xrightarrow{\sim}& F_{X/k}^{*}(\mathcal{B}_{\XX'}\otimes_{\mathscr{O}_{X'}}M',\iota^{*}(\psi_{\mathcal{B}'})\otimes\id)
\end{eqnarray*}
where $\rS(\Omega_{X/k}^{1})$ is equipped with the canonical action of $\widehat{\Gamma}(\rT_{X/k})$. 
A local section $u$ of $\rS(\Omega_{X/k}^{1})\otimes_{\mathscr{O}_{X}}F_{X/k}^{*}(M')$ can be written as a finite sum
\begin{equation}
	u=\sum_{i=0}^{m}\omega_{i}\otimes u_{i}
\end{equation}
with $u_{i}\in F_{X/k}^{*}(M')$ and $\omega_{i}\in \rS^{i}(\Omega)$. In view of the perfect pairing $\Gamma_{i}(\rT_{X/k})\otimes_{\mathscr{O}_{X}}\rS^{i}(\Omega_{X/k}^{1})\to \mathscr{O}_{X}$ \eqref{pairing T Om} for $i\ge 1$, the action of $\widehat{\Gamma}(\rT_{X/k})$ on $u$ is trivial if and only if $u_{i}=0$ for $i\ge 1$, i.e. $u$ belongs to the submodule $F_{X/k}^{*}(M')$ of $F_{X/k}^{*}(M')\otimes_{\mathscr{O}_{X}}\rS(\Omega)$. The first isomorphism follows. 

Equipped with the Frobenius descent connection $\nabla_{\can}$ on $F_{X/k}^{*}(M')$, we have an injection of $\rD_{X/k}^{\gamma}$-modules $(F_{X/k}^{*}(M'),\nabla_{\can},0)\to (\underline{\mathcal{B}}_{\XX}\otimes_{\mathscr{O}_{X'}}M',\nabla_{\underline{\mathcal{B}}},\psi_{\underline{\mathcal{B}}})$ and a canonical $\mathscr{O}_{X'}$-linear isomorphism (\cite{Ka71} 5.1)
\begin{equation}
	M'\xrightarrow{\sim} F_{X/k}^{*}(M')^{(\nabla_{\can},0)}. \label{iso DL deg 0}
\end{equation}
In view of the first isomorphism, we deduce that $(\underline{\mathcal{B}}_{\XX}\otimes_{\mathscr{O}_{X'}}M')^{(\nabla_{\underline{\mathcal{B}}},\psi_{\underline{\mathcal{B}}})}$ is contained in the image of $F_{X/k}^{*}(M')$ in $\underline{\mathcal{B}}_{\XX}\otimes_{\mathscr{O}_{X'}}M'$. Then the second isomorphism follows from \eqref{iso DL deg 0}.
\end{proof}

	Now we use the above results to compare the Cartier equivalence $\rmC_{X/\SS}^{*}$ and the Cartier transform of Ogus--Vologodsky. 

\begin{theorem}\label{comparison Cartier}
	Let $\XX$ be a smooth formal $\SS$-scheme and $X$ its special fiber. The following diagram is commutative up to a canonical isomorphism
	\begin{equation}
		\xymatrix{
			\CC^{\qcoh}(\mathscr{O}_{\underline{\mathscr{E}},1})\ar[d]^{\wr} \ar[r]^{\rmC_{X/\SS*}} & \CC^{\qcoh}(\mathscr{O}_{\mathscr{E}',1}) \ar[d]^{\wr}\\
			\MIC_{\gamma}^{\qn}(X/k) \ar[r]^{\rmC_{\XX_{2}'}} & \HIGG_{\gamma}^{\qn}(X/k)}
	\end{equation}
	where $\rmC_{X/\SS*}$ is the direct image functor of the morphism of topoi $\rmC_{X/\SS}$ \eqref{thm pullback Cartier}, which depends only on $X$, and the vertical functors are equivalences of categories \eqref{equi of Dmod crystals} and depend on $\XX$. 
\end{theorem}

\begin{rem}
	The top (resp. lower) horizontal functor in the above diagram depends only on $X$ (resp. $\XX_{2}'$). 
	However, the vertical functors \eqref{equivalence Dmod crystals} are constructed by a formal model $\XX$ of $X$. 
	A natural question is: do vertical equivalences depend only on a lifting of $X$ over the quotient of $\rW$ by a power of $p$? 
\end{rem}

\begin{nothing}
	Recall that Ogus and Vologodsky constructed a torsor $\mathscr{L}_{\XX_{2}'}$ of $(X/k)_{\cris}$ \eqref{torsor Frob liftings}, the crystal of affine functions $\mathscr{F}_{\XX_{2}'}$ on $\mathscr{L}_{\XX_{2}'}$ of $(X/k)_{\cris}$ and the quasi-coherent crystal of $\mathscr{O}_{X/k}$-algebras $\mathscr{A}_{\XX'_{2}}$ associated to $\mathscr{F}_{\XX_{2}'}$ \eqref{alg splitting mod}.
	We put $\mathcal{F}_{\XX'_{2}}=\mathscr{F}_{\XX_{2}',(X,X)}$ and $\mathcal{A}_{\XX'_{2}}=\mathscr{A}_{\XX'_{2},(X,X)}$, which are equipped with $\rD^{\gamma}_{X/k}$-module structures \eqref{AX' dual}. 
	The Cartier transform is defined by (\ref{thm OV}(ii))
	\begin{equation} \label{Cartier OV compare}
		\rmC_{\XX_{2}'}: \MIC_{\gamma}^{\qn}(X/k)\xrightarrow{\sim} \HIGG_{\gamma}^{\qn}(X'/k), \qquad N\mapsto \iota^{*}(\FHom_{\rD_{X/k}^{\gamma}}((\mathcal{A}_{\XX_{2}'})^{\vee}, N)). 
	\end{equation}
\end{nothing}

\begin{nothing} \label{lifting object CRIS of cris}
	Let $(U,T,\delta)$ be an object of $\Cris(X/k)$ such that there exists a flat formal $\SS$-scheme $\mathfrak{T}$ with special fiber $T$. Recall \eqref{lemma can connection crystal} that the morphism $U\to T$ induces an isomorphism $U\xrightarrow{\sim} \underline{T}$. Then we obtain an object $(U,\mathfrak{T})$ of $\CRIS$. The morphism $\varphi_{T/k}:T\to X'$ \eqref{varphi Tk} is the same as the composition $T\xrightarrow{\phi_{T/k}} U'\to X'$ \eqref{C inverse underlying}. 
	Moreover, the $\varphi_{T/k}^{*}(\rT_{X'/k})$-torsor $\mathscr{L}_{\XX'_{2},\varphi_{T/k}}$ of $T_{\zar}$ \eqref{torsor of liftings} is the same as the $\phi_{T/k}^{*}(\rT_{U'/k})$-torsor $\mathscr{R}_{\XX',(U',\mathfrak{T},\phi_{T/k})}$ of $T_{\zar}$ \eqref{Relevement X}. Recall that $\mathscr{F}_{\XX_{2}',(U,T)}$ is the sheaf of affine functions on $\mathscr{L}_{\XX_{2}',\varphi_{T/k}}$. By \eqref{underlineF F}, we deduce a canonical isomorphism of $\mathscr{O}_{T}$-modules
	\begin{equation} \label{iso sheaves of affine functions}
		\mathscr{F}_{\XX_{2}',(U,T)}\xrightarrow{\sim} \underline{\mathscr{F}}_{\XX,(U,\mathfrak{T})}.
	\end{equation}
	
	Let $g:(U_{1},T_{1},\delta_{1})\to (U_{2},T_{2},\delta_{2})$ be a morphism of $\Cris(X/k)$. Suppose that there exists an $\SS$-morphism $\mathfrak{g}:\mathfrak{T}_{1}\to \mathfrak{T}_{2}$ of flat formal $\SS$-schemes with special fiber $g:T_{1}\to T_{2}$. Then we obtain a morphism $\mathfrak{g}:(U_{1},\mathfrak{T}_{1})\to (U_{2},\mathfrak{T}_{2})$ of $\CRIS$. In view of \ref{torsor Frob liftings}(ii), the transition morphism of the crystal of $\mathscr{O}_{\CRIS,1}$-modules $\underline{\mathscr{F}}_{\XX}$ of $\widetilde{\CRIS}$ associated to $\mathfrak{g}$ (\ref{cg for F pCRIS}, \ref{C inverse underlying})
	\begin{equation}
		c_{\mathfrak{g}}:g^{*}(\underline{\mathscr{F}}_{\XX,(U_{2},\mathfrak{T}_{2})})\xrightarrow{\sim} \underline{\mathscr{F}}_{\XX,(U_{1},\mathfrak{T}_{1})}
	\end{equation}
	is compatible with the transition morphism of the crystal $\mathscr{F}_{\XX_{2}'}$ of $\mathscr{O}_{X/k}$-modules associated to $g$ 
	\begin{equation}
		c_{g}:g^{*}(\mathscr{F}_{\XX'_{2},(U_{2},T_{2})})\xrightarrow{\sim} \mathscr{F}_{\XX'_{2},(U_{1},T_{1})}.
	\end{equation}
	Applying above facts to projections from the PD-envelop of diagonal immersion $(X,\PP_{X})$ to $(X,X)$, we deduce the following lemma. 
\end{nothing}

\begin{lemma} \label{crystal of affine functions concides}
	The isomorphism $\mathcal{F}_{\XX_{2}'}\xrightarrow{\sim} \underline{\mathcal{F}}_{\XX}$ \eqref{iso sheaves of affine functions} is compatible with the actions of $\rD_{X/k}$.
\end{lemma}

\begin{coro} \label{coro two splitting}
	The isomorphism $\mathcal{F}_{\XX_{2}'}\xrightarrow{\sim} \underline{\mathcal{F}}_{\XX}$ induces a canonical isomorphism of $\mathscr{O}_{X}$-algebras $\mathcal{A}_{\XX_{2}'}\xrightarrow{\sim} \underline{\mathcal{B}}_{\XX}$ compatible with the actions of $\rD_{X/k}^{\gamma}$.
\end{coro}
\begin{proof}
By \ref{crystal of affine functions concides}, we obtain an isomorphism $\mathcal{A}_{\XX_{2}'}\xrightarrow{\sim} \underline{\mathcal{B}}_{\XX}$ compatible with actions of $\rD_{X/k}$. 
By \ref{Gamma action on A} and \eqref{base change Gamma T}, the $\mathscr{O}_{X}$-module $\underline{\mathcal{F}}_{\XX}$ (resp. $\mathscr{O}_{X}$-algebra $\underline{\mathcal{B}}_{\XX}$) is equipped with a $F_{X/k}^{*}(\widehat{\Gamma}(\rT_{X'/k}))$-module structure that we denote by $\vartheta_{\underline{\mathcal{F}}}$ (resp. $\vartheta_{\underline{\mathcal{B}}}$). 
By \ref{alg splitting mod} and the following lemma, the isomorphism $\mathcal{A}_{\XX_{2}'}\xrightarrow{\sim} \underline{\mathcal{B}}_{\XX}$ is compatible with actions of $F_{X/k}^{*}(\widehat{\Gamma}(\rT_{X'/k}))$. 
Then the assertion follows.
\end{proof}
\begin{lemma} \label{two Gamma modules actions}
	Two $F_{X/k}^{*}(\widehat{\Gamma}(\rT_{X'/k}))$-module structures $\psi_{\underline{\mathcal{F}}}$ and $\vartheta_{\underline{\mathcal{F}}}$ on $\underline{\mathcal{F}}_{\XX}$ \eqref{C inverse underlying} (resp. $\psi_{\underline{\mathcal{B}}}$ and $\vartheta_{\underline{\mathcal{F}}}$ on $\underline{\mathcal{B}}_{\XX}$) coincide.
\end{lemma}
\begin{proof} It suffices to show the assertion for $\underline{\mathcal{F}}_{\XX}$. We have $\underline{\mathcal{F}}_{\XX}=\mathscr{F}_{\XX',(X',\XX,F_{X/k})}$. 

The question being local, we can reduce to case where there exists an $\SS$-morphism $F:\XX\to \XX'$ which lifts the relative Frobenius morphism $F_{X/k}$ of $X$. Then we obtain a morphism $F:(X',\XX,F_{X/k})\to (X',\XX')$ of $\pCRIS'$ and an isomorphism \eqref{cg for F pCRIS}
\begin{equation}
	\eta_{F}:F_{X/k}^{*}(\mathcal{F}_{\XX'})\xrightarrow{\sim} \underline{\mathcal{F}_{\XX}}.
\end{equation}

Let $\psi_{\mathcal{F}'}$ be the $\widehat{\Gamma}(\rT_{X'/k})$-module structure on $\mathcal{F}_{\XX'}$ induced by the crystal $\mathscr{F}_{\XX'}$. By \ref{lifting F Cartier}, the isomorphism $\eta_{F}:F_{X/k}^{*}(\mathcal{F}_{\XX'})\xrightarrow{\sim} \underline{\mathcal{F}}_{\XX}$ underlies an isomorphism of $F_{X/k}^{*}(\widehat{\Gamma}(\rT_{X'/k}))$-modules
\begin{equation} \label{eta F lemma 2actions}
	\eta_{F}:F_{X/k}^{*}(\mathcal{F}_{\XX'},\iota^{*}(\psi_{\mathcal{F}'}))\xrightarrow{\sim} (\underline{\mathcal{F}}_{\XX},\psi_{\underline{\mathcal{F}}}).
\end{equation}

On the other hand, regarding $\mathcal{F}_{\XX'}$ as a sheaf of affine functions, the $\widehat{\Gamma}(\rT_{X'/k})$-action on $\mathcal{F}_{\XX'}$ defined in \ref{Gamma action on A} is equal to $\iota^{*}(\psi_{\mathcal{F}'})$ \eqref{GammaT FBXX}. 
By \ref{affine inverse image num}, $\eta_{F}$ induces an isomorphism of $F_{X/k}^{*}(\widehat{\Gamma}(\rT_{X'/k}))$-modules
\begin{equation} \label{FX canonical iso 2actions}
	\eta_{F}: F_{X/k}^{*}( \mathcal{F}_{\XX'},\iota^{*}(\psi_{\mathcal{F}'})) \xrightarrow{\sim} (\underline{\mathcal{F}}_{\XX}, \vartheta_{\underline{\mathcal{F}}}).
\end{equation}
The assertion follows. 
\end{proof}

\begin{nothing}
	\textit{Proof of \ref{comparison Cartier}}.
	Let $\mathscr{N}$ be a quasi-coherent crystal of $\mathscr{O}_{\underline{\mathscr{E}},1}$-modules of $\widetilde{\underline{\mathscr{E}}}$, $N$ the associated $\rD_{X/k}^{\gamma}$-module, $\mathscr{M}'=\rmC_{X/\SS*}(\mathscr{N})$ and $M'$ the associated $\widehat{\Gamma}(\rT_{X'/k})$-module. 
By \ref{thm pullback Cartier}, we have a canonical isomorphism $\rmC^{*}_{X/\SS}(\mathscr{M}')\xrightarrow{\sim} \mathscr{N}$. The $\rD_{X/k}^{\gamma}$-module structure $(\nabla_{\underline{\mathcal{B}}_{\XX}},\psi_{\underline{\mathcal{B}}_{\XX}})$ on $\underline{\mathcal{B}}_{\XX}$ induces a $\rD_{X/k}^{\gamma}$-module structure on its dual $\underline{\mathcal{B}}_{\XX}^{\vee}$ \eqref{FHom Dgamma}. 
By \eqref{Cartier OV compare} and \ref{coro two splitting}, we have a canonical isomorphism of $\widehat{\Gamma}(\rT_{X'/k})$-modules
\begin{eqnarray}\label{recall COV N calcul}
	\rmC_{\XX_{2}'}(N)\xrightarrow{\sim} \iota^{*}(\FHom_{\rD_{X/k}^{\gamma}}( \underline{\mathcal{B}}_{\XX}^{\vee},N))
\end{eqnarray}
where the $\widehat{\Gamma}(\rT_{X'/k})$-action on the right hand side is given by that of $\underline{\mathcal{B}}_{\XX}^{\vee}$. The canonical morphism
\begin{equation} \label{Hom Cartan}
	\FHom_{\mathscr{O}_{X}}( \mathscr{O}_{X},\underline{\mathcal{B}}_{\XX}\otimes_{\mathscr{O}_{X}} N) \to \FHom_{\mathscr{O}_{X}}( \underline{\mathcal{B}}_{\XX}^{\vee},N)
\end{equation}
sends a local section $\varphi$ of $\FHom_{\mathscr{O}_{X}}( \mathscr{O}_{X},\underline{\mathcal{B}}_{\XX}\otimes_{\mathscr{O}_{X}} N)$ to $\chi:\underline{\mathcal{B}}_{\XX}^{\vee}\to N$ defined for every local section $f$ of $\underline{\mathcal{B}}_{\XX}^{\vee}$ by
\begin{equation}
	\chi(f)=(f\otimes \id_{N})(\varphi(1)). \label{definition phi to chi}
\end{equation}
Since $\underline{\mathcal{B}}_{\XX}$ is locally a direct sum of free $\mathscr{O}_{X}$-modules of finite type \eqref{eq thm ad 2}, the morphism \eqref{Hom Cartan} is an isomorphism. We equip $\mathscr{O}_{X}$ the $\rD_{X/k}^{\gamma}$-module structure $(d,0)$. By \eqref{formula hom connection} and \eqref{Hom psi}, one verifies that a local section $\varphi: \mathscr{O}_{X}\to \underline{\mathcal{B}}_{\XX}\otimes_{\mathscr{O}_{X}} N$ is $\rD_{X/k}^{\gamma}$-equivariant if and only if $\chi$ \eqref{definition phi to chi} is $\rD_{X/k}^{\gamma}$-equivariant.
Then we deduce an isomorphism 
\begin{equation}
	\FHom_{\rD_{X/k}^{\gamma}}( (\mathscr{O}_{X},d,0),\underline{\mathcal{B}}_{\XX}\otimes_{\mathscr{O}_{X}} N) \xrightarrow{\sim}\FHom_{\rD_{X/k}^{\gamma}}( \underline{\mathcal{B}}_{\XX}^{\vee},N).
\end{equation}
In view of \eqref{Hom psi}, the above isomorphism induces an isomorphism of $\widehat{\Gamma}(\rT_{X'/k})$-modules \eqref{prep iso ad}
\begin{equation}
	(\FHom_{\rD_{X/k}^{\gamma}}( \mathscr{O}_{X},\underline{\mathcal{B}}_{\XX}\otimes_{\mathscr{O}_{X}} N), \psi_{\underline{\mathcal{B}}_{\XX}})\xrightarrow{\sim}\iota^{*}(\FHom_{\rD_{X/k}^{\gamma}}( \underline{\mathcal{B}}_{\XX}^{\vee},N)).
\end{equation}
Then the assertion follows from \ref{D-inv theta-inv}(ii) and \eqref{recall COV N calcul}. \hfill $\qed$
\end{nothing}

\section{Fontaine modules} \label{Fil modules}
	Inspired by the work of Fontaine--Laffaille \cite{FL}, Faltings introduced a notion of \textit{Fontaine module} on a smooth scheme over $\rW$ \cite{Fal89}. 
	In this section, we propose a new definition of Fontaine module using Cartier equivalence \eqref{definition MF}. 
	Compared to Faltings' original definition, our definition avoids the choice of (local) liftings of Frobenius to talk about ``divided Frobenius structures'' \eqref{def MF Tsuji} and encodes them in the Cartier equivalence.  
	We also show the compatibility between various definitions in (\ref{phiF familly phi}, \ref{coro MF1 MF2}, \ref{MF1 MF2 diagonal}). 

	Let $\XX$ denotes a smooth formal $\SS$-scheme and $X$ its special fiber. 

	\begin{definition} \label{filtered mods}
	Let $n$ be an integer $\ge 1$. We define the category of filtered modules with quasi-nilpotent integrable connection $\MIC_{\rF}(\XX_{n}/\SS_{n})$ as follows. \textit{A filtered module with quasi-nilpotent integrable connection} is a triple $(M,\nabla,M^{\bullet})$ consisting of an object $(M,\nabla)$ of $\MIC^{\qn}(\XX_{n}/\SS_{n})$ and a decreasing filtration $\{M^{i}\}_{i\in \mathbb{Z}}$
	\begin{equation}
		\cdots \subseteq M^{2} \subseteq M^{1} \subseteq M^{0}=M=M^{-1}\cdots
	\end{equation}
	satisfying Griffiths' transversality 
	\begin{equation} \label{Griffiths transversality}
		\nabla(M^{i})\subset M^{i-1}\otimes \Omega_{\XX_{n}/\SS_{n}}^{1} \quad \forall ~ i\ge 0.
	\end{equation}	
	
	Given two objects $(M_{1},\nabla_{1},M_{1}^{\bullet})$ and $(M_{2},\nabla_{2},M_{2}^{\bullet})$ of $\MIC_{\rF}(\XX_{n}/\SS_{n})$, a morphism from $(M_{1},\nabla_{1},M_{1}^{\bullet})$ to $(M_{2},\nabla_{2},M_{2}^{\bullet})$ is a horizontal $\mathscr{O}_{\XX_{n}}$-linear morphism $f:M_{1}\to M_{2}$ compatible with the filtrations.

	For any $\ell\ge 0$, we denote by $\MIC^{\ell}_{\rF}(\XX_{n}/\SS_{n})$ the full subcategory of $\MIC_{\rF}(\XX_{n}/\SS_{n})$ consisting of objects with length $\le \ell$ (i.e. the filtration satisfies $M^{\ell+1}=0$). 
\end{definition}

\begin{nothing} \label{data pr MF}
	Let $\ell$ be an integer $\ge 0$, $(M,\nabla,M^{\bullet})$ an object of $\MIC^{\ell}_{\rF}(\XX_{n}/\SS_{n})$. 
	We consider the $\mathscr{O}_{\XX_{n}}$-linear morphism
	\begin{equation} \label{morphism FL}
		g: \oplus_{i=1}^{\ell} M^{i}\to \oplus_{i=0}^{\ell}M^{i}
	\end{equation}
	defined for every local section $m_{i}$ of $M^{i}$ by $g(m_{i})=(m_{i},-pm_{i})$ in $M^{i-1}\oplus M^{i}$. We set  
	\begin{equation} \label{definition M tilde}
		\widetilde{M}=\Coker(g).
	\end{equation}
	For any $0\le j\le \ell$, we denote by $(-)_{j}$ the composition $M^{j}\to \oplus_{i=0}^{\ell}M^{i}\to \widetilde{M}$ and by $\widetilde{M}^{-j}$ the canonical image of $\oplus_{i=0}^{j}M^{i}$ in $\widetilde{M}$. We obtain a decreasing filtration
	\begin{equation}
		\widetilde{M}^{0}\subseteq \widetilde{M}^{-1}\subseteq \cdots \subseteq \widetilde{M}^{-\ell}=\widetilde{M}. \label{filtration N M tilde}
	\end{equation}
	
	We consider the $\rW_{n}$-linear morphism
	\begin{equation}
		h:\oplus_{i=0}^{\ell}M^{i}\to \bigl(\oplus_{i=0}^{\ell}M^{i}\bigr)\otimes_{\mathscr{O}_{\XX_{n}}}\Omega_{\XX_{n}/\SS_{n}}^{1}
	\end{equation}
	defined by 
	\begin{eqnarray}
		h|_{M^{i}}&=&\nabla:M^{i}\to M^{i-1}\otimes_{\mathscr{O}_{\XX_{n}}}\Omega_{\XX_{n}/\SS_{n}}^{1}, \quad \textnormal{for $1\le i\le \ell$}, \label{def tilde nabla}\\
		h|_{M^{0}}&=&p\nabla:M^{0}\to M^{0}\otimes_{\mathscr{O}_{\XX_{n}}}\Omega_{\XX_{n}/\SS_{n}}^{1}. \nonumber
	\end{eqnarray}
\end{nothing}

\begin{lemma} \label{lemma p connection}
	The $\rW_{n}$-linear morphism $h$ induces a quasi-nilpotent integrable $p$-connection $\widetilde{\nabla}$ on $\widetilde{M}$ such that for any $-\ell \le i\le -1$, we have
	\begin{equation} \label{transversality N fil}
		\widetilde{\nabla}(\widetilde{M}^{i})\subset \widetilde{M}^{i+1}\otimes_{\mathscr{O}_{\XX_{n}}}\Omega_{\XX_{n}/\SS_{n}}^{1}.
	\end{equation}
\end{lemma}
\begin{proof} It follows from the definition that the composition
\begin{equation}
	\oplus_{i=1}^{\ell} M^{i} \xrightarrow{g} \oplus_{i=0}^{\ell}M^{i}\xrightarrow{h} \bigl(\oplus_{i=0}^{\ell}M^{i}\bigr)\otimes_{\mathscr{O}_{\XX_{n}}}\Omega_{\XX_{n}/\SS_{n}}^{1} \to \widetilde{M}\otimes_{\mathscr{O}_{\XX_{n}}}\Omega_{\XX_{n}/\SS_{n}}^{1}
\end{equation}
is zero. Hence the morphism $h$ induces a $\rW_{n}$-linear morphism 
\begin{equation} 
	\widetilde{\nabla}:\widetilde{M}\to \widetilde{M}\otimes_{\mathscr{O}_{\XX_{n}}}\Omega_{\XX_{n}/\SS_{n}}^{1}.
\end{equation}
We show that $\widetilde{\nabla}$ is a $p$-connection. The restriction of $h$ to $M^{0}$ is the $p$-connection $p\nabla$. Hence the restriction of $\widetilde{\nabla}$ to $\widetilde{M}^{0}$ is a $p$-connection. 
Let $f$ be a local section of $\mathscr{O}_{\XX}$, $i$ an integer $\in [1,\ell]$ and $m$ a local section of $M^{i}$. 
The morphism $h$ sends $fm\in M^{i}$ to
\begin{equation}
	\nabla(fm)=f\nabla(m)+m\otimes df \in M^{i-1}\otimes \Omega_{\XX_{n}/\SS_{n}}^{1}.
\end{equation}
Note that we have $(m)_{i-1}=(pm)_{i}$. Then we deduce that
\begin{equation}
	\widetilde{\nabla}(f(m)_{i})=f\widetilde{\nabla}( (m)_{i})+p(m)_{i}\otimes df.
\end{equation}
The assertion follows. The integrability of $\nabla$ implies easily that of $\widetilde{\nabla}$. 

It is clear from the definition that $\widetilde{\nabla}$ satisfies the condition \eqref{transversality N fil}. Since $M$ is $p^{n}$-torsion, the restriction of $\widetilde{\nabla}$ to $\widetilde{M}^{0}$ is quasi-nilpotent. In view of \eqref{transversality N fil}, we deduce that $\widetilde{\nabla}$ is quasi-nilpotent.
\end{proof}
\begin{nothing} \label{local description Tstrat tildeM}
	Let $\varepsilon_{\widetilde{M},T}$ be the $\mathcal{T}_{\XX}$-stratification on $\widetilde{M}$ associated to $\widetilde{\nabla}$ \eqref{qn pconnection T}. We present its local description. 
	Suppose that there exists an \'etale $\SS$-morphism $\XX\to \widehat{\mathbb{A}}^{d}_{\SS}=\Spf(\rW\{T_{1},\cdots,T_{d}\})$. We take again the notation of \ref{quasi-nilpotent coords} and of \ref{Local description T X}. 
	For any $0\le i \le \ell$, any local section $m$ of $M^{i}$ and any $I\in \mathbb{N}$, in view of \eqref{def tilde nabla}, we have 
\begin{equation}
	\widetilde{\nabla}_{\partial^{I}}( (m)_{i})=
	\left\{ \begin{array}{ll}
		(\nabla_{\partial^{I}}(m))_{i-|I|} & \textnormal{if } 0\le i\le |I|; \\
		p^{|I|-i}(\nabla_{\partial^{I}}(m))_{0} & \textnormal{if } |I|> i.
		\end{array} \right.
\end{equation}
By \eqref{stratification given by pconnection}, we deduce that
\begin{equation}
	\varepsilon_{\widetilde{M},T}(1\otimes (m)_{i})=
	\sum_{|I|\le i} \big(\nabla_{\partial^{I}}(m)\big)_{i-|I|}\otimes \biggl(\frac{\xi}{p}\biggr)^{[I]} + \sum_{|I|> i} p^{|I|-i}\big(\nabla_{\partial^{I}}(m)\big)_{0}\otimes \biggl(\frac{\xi}{p}\biggr)^{[I]}.
	\label{stratification T tildeM}
\end{equation}
\end{nothing}

\begin{prop} \label{prop construction R stra tildeM}
	Suppose that $\ell\le p-1$. 
	There exists an $\mathcal{R}_{\XX}$-stratification $\varepsilon_{\widetilde{M}}$ on $\widetilde{M}$, which induces the above $\mathcal{T}_{\XX}$-stratification $\varepsilon_{\widetilde{M},T}$ via the functor \eqref{varpi pullback}.
\end{prop}
\begin{proof} Let $\varepsilon_{M}$ be the $\mathcal{P}_{\XX}$-stratification on $M$ and $\theta_{M}:M\to M\otimes_{\mathscr{O}_{\XX}}\mathcal{P}_{\XX}$ the morphism defined by $\theta_{M}(m)=\varepsilon_{M}(1\otimes m)$. We denote by $J_{P}$ the PD-ideal of $\mathcal{P}_{\XX}$.
	By flatness of $J^{[\bullet]}_{P}$ over $\mathscr{O}_\XX$ \eqref{flatness PX}, $M^{i}\otimes_{\mathscr{O}_{\XX_{n}}}J^{[j]}_{P}$ is a submodule of $M\otimes_{\mathscr{O}_{\XX_{n}}}\mathcal{P}_{\XX}$ for any $i,j \ge 0$. By (\cite{Ogus94} 3.1.3) and Griffiths' transversality, we have 
\begin{equation}
	\theta_{M}(M^{i})\subset \sum_{j=0}^{i} M^{j}\otimes_{\mathscr{O}_{\XX_{n}}} J_{P}^{[i-j]},\qquad \forall~ 0\le i\le \ell.
\end{equation}
We denote the target by $(M\otimes_{\mathscr{O}_{\XX_{n}}}\mathcal{P}_{\XX})^{i}$ and the induced morphism by $\theta_{M}^{i}:M^{i}\to (M\otimes_{\mathscr{O}_{\XX_{n}}}\mathcal{P}_{\XX})^{i}$. 
Let $s:\mathcal{P}_{\XX}\to \mathcal{R}_{\XX} $ be the homomorphism defined in \ref{lemma alg P to R}, and for all $j\le p-1$, $s^{j}:J^{[j]}_{P}\to \mathcal{R}_{\XX}$ the induced morphism \eqref{P to R divided}. 
For any $j\le i$, we have a canonical morphism 
\begin{equation} \label{()i-j sj}
	(-)_{j}\otimes s^{i-j}:M^{j}\otimes_{\mathscr{O}_{\XX_{n}}}J_{P}^{[i-j]}\to \widetilde{M}\otimes_{\mathscr{O}_{\XX_{n}}}\mathcal{R}_{\XX}.
\end{equation}
For any $0\le j'<j \le i$, by flatness of $J^{[\bullet]}_{P}$ over $\mathscr{O}_{\XX}$ \eqref{flatness PX}, we have the intersection in $M\otimes_{\mathscr{O}_{\XX_{n}}}\mathcal{P}_{\XX}$
\begin{displaymath}
	 M^{j}\otimes_{\mathscr{O}_{\XX_{n}}}J_{P}^{[i-j]}\cap M^{j'}\otimes_{\mathscr{O}_{\XX_{n}}}J_{P}^{[i-j']}=M^{j}\otimes_{\mathscr{O}_{\XX_{n}}}J_{P}^{[i-j']}.
\end{displaymath}
Since $s^{i-j}|_{J_{P}^{[i-j']}}=p^{j-j'}s^{i-j'}$ \eqref{P to R divided}, we deduce that the morphisms \eqref{()i-j sj} are compatible and they induce an $\mathscr{O}_{\XX}$-linear morphism
\begin{equation}
	u^{i}: (M\otimes_{\mathscr{O}_{\XX_{n}}}\mathcal{P}_{\XX})^{i} \to \widetilde{M}\otimes_{\mathscr{O}_{\XX_{n}}}\mathcal{R}_{\XX}.
\end{equation}
By construction, we have $u^{i}|_{(M\otimes_{\mathscr{O}_{\XX_{n}}}\mathcal{P}_{\XX})^{i+1}}=pu^{i+1}$. Then the morphism
\begin{displaymath}
	\oplus_{i=0}^{\ell}u^{i}\circ\theta_{M}^{i} :\oplus_{i=0}^{\ell} M^{i}\to \widetilde{M}\otimes_{\mathscr{O}_{\XX_{n}}}\mathcal{R}_{\XX}
\end{displaymath}
induces an $\mathscr{O}_{\XX}$-linear morphism
\begin{equation}
	\theta_{\widetilde{M}}:\widetilde{M}\to \widetilde{M}\otimes_{\mathscr{O}_{\XX_{n}}}\mathcal{R}_{\XX}.
\end{equation}

We will show that the above morphism satisfies the conditions of \ref{lemma stratification} and hence induces an $\mathcal{R}_{\XX}$-stratification.
To do this, we can reduce to the case where there exists an \'etale $\SS$-morphism $\XX\to \widehat{\mathbb{A}}^{d}_{\SS}=\Spf(\rW\{T_{1},\cdots,T_{d}\})$. We take again the notation of \ref{B-O stratification MIC}. 
For any $0\le i \le \ell$, any local section $m$ of $M^{i}$, we have 
\begin{equation}
	\theta^{i}_{M}(m)=\sum_{I\in \mathbb{N}^{d}} \nabla_{\partial^{I}}(m)\otimes \xi^{[I]} ~\in\sum_{j=0}^{i} M^{j}\otimes_{\mathscr{O}_{\XX_{n}}} J_{P}^{[i-j]}.
\end{equation}
The $p$-adic valuation of $I!$ is less than $\sum_{k\ge 1}\lfloor \frac{|I|}{p^{k}}\rfloor$. If $i\le p-1$ and $|I|\ge i$, $\frac{p^{|I|-i}}{I!}$ is an element of $\mathbb{Z}_{p}$.
By \eqref{calcul P to R alg}, we have $s^{i}(\xi^{[I]})=\frac{p^{|I|-i}}{I!} (\frac{\xi}{p})^{I}$. 
Then we deduce that 
\begin{equation} \label{stratification R tildeM}
	\theta_{\widetilde{M}}((m)_{i})=
	\sum_{|I|\le i}\frac{1}{I!} \big(\nabla_{\partial^{I}}(m)\big)_{i-|I|}\otimes \biggl(\frac{\xi}{p}\biggr)^{I} + \sum_{|I|> i}\frac{p^{|I|-i}}{I!} \big(\nabla_{\partial^{I}}(m)\big)_{0}\otimes \biggl(\frac{\xi}{p}\biggr)^{I}.
\end{equation}
It is clear that $\theta_{\widetilde{M}}$ verifies condition (i) of \ref{lemma stratification}. 
By \eqref{stratification R tildeM} and the local description \eqref{local description Hopf RQ} of $\delta:\mathcal{R}_{\XX,n}\to \mathcal{R}_{\XX,n}\otimes_{\mathscr{O}_{\XX_{n}}}\mathcal{R}_{\XX,n}$, we deduce that
\begin{eqnarray*}
	&&\theta_{\widetilde{M}}\otimes\id_{\mathcal{R}} (\theta_{\widetilde{M}}((m)_{i}))\\
	&=& \sum_{|I|+|J|\le i}\frac{1}{I!J!} \bigr(\nabla_{\partial^{I+J}}(m)\bigl)_{i-|I|-|J|}\otimes \biggl(\frac{\xi}{p}\biggr)^{I} \otimes\biggl(\frac{\xi}{p}\biggr)^{J} \\ 
	&& +\sum_{|I|+|J|> i}\frac{p^{|I|+|J|-i}}{I!J!} \bigr(\nabla_{\partial^{I+J}}(m)\bigl)_{0} \otimes \biggl(\frac{\xi}{p}\biggr)^{I} \otimes \biggl(\frac{\xi}{p}\biggr)^{J} \\
	&=& \id_{\widetilde{M}}\otimes \delta (\theta_{\widetilde{M}}( (m)_{i})).
\end{eqnarray*}
By \ref{lemma stratification}, we obtain an $\mathcal{R}_{\XX}$-stratification $\varepsilon_{\widetilde{M}}$ on $\widetilde{M}$. By comparing \eqref{stratification T tildeM} and \eqref{stratification R tildeM}, $\varepsilon_{\mathscr{M}}$ extends $\varepsilon_{\widetilde{\mathscr{M}},T}$.
\end{proof}

\begin{nothing}\label{M pullback pi}
	By \ref{equi crystals stratification}, we associate to $(\widetilde{M},\varepsilon_{\widetilde{M}})$ a crystal of $\mathscr{O}_{\pCRIS,n}$-modules of $\widetilde{\pCRIS}$ that we denote by $\widetilde{\mathscr{M}}$.

	We put $\XX'=\XX\times_{\SS,\sigma}\SS$ \eqref{notations} and we denote by $\pi:\XX'\to \XX$ the canonical morphism. In view of \ref{prop dilatation flat} and \ref{prop R Q Hopf alg}, the morphism $\pi$ induces a morphism of formal groupoids $\pi_{R}:\RR_{\XX'}\simeq \RR_{\XX}\times_{\SS,\sigma}\SS\to \RR_{\XX}$ above $\pi$ \eqref{morphism of groupoids}. We denote by $\widetilde{\mathscr{M}}'$ the crystal of $\mathscr{O}_{\pCRIS',n}$-modules of $\widetilde{\pCRIS}'$ associated to the $\mathscr{O}_{\XX'_{n}}$-module with $\mathcal{R}_{\XX'}$-stratification $(\pi^{*}(\widetilde{M}),\pi_{R}^{*}(\varepsilon_{\widetilde{M}}))$ \eqref{morphism Hopf functor}. The $\mathscr{O}_{\XX'_{n}}$-module with integrable $p$-connection associated to $\widetilde{\mathscr{M}}'$ \eqref{functor mu} is the inverse image of $(\widetilde{M},\widetilde{\nabla})$ by $\pi$ (\cite{Shiho} page 6).

The previous construction is clearly functorial and it defines a functor \eqref{def crystals}
	\begin{eqnarray} \label{functor MIC to crys}
		\MIC^{p-1}_{\rF}(\XX_{n}/\SS_{n})&\to&\CC(\mathscr{O}_{\pCRIS',n})\\
		(M,\nabla,M^{\bullet}) &\mapsto& \widetilde{\mathscr{M}}'. \nonumber
	\end{eqnarray}
\end{nothing}

\begin{definition} \label{definition MF}
	(i) A ($p^{n}$-torsion) \textit{Fontaine module over $\XX$} is a quadruple $(M,\nabla,M^{\bullet},\varphi)$ consisting of an object $(M,\nabla,M^{\bullet})$ of $\MIC^{p-1}_{\rF}(\XX_{n}/\SS_{n})$ with quasi-coherent $\mathscr{O}_{\XX_{n}}$-module $M^{i}$, and a morphism of $\MIC^{\qn}(\XX_{n}/\SS_{n})$
	\begin{equation}
		\varphi:\nu(\rmC^{*}(\widetilde{\mathscr{M}}'))\to (M,\nabla) \label{morphism varphi}
	\end{equation}
	where $\rmC^{*}$ is the Cartier equivalence \eqref{Cartier transform} and $\nu:\CC(\mathscr{O}_{\CRIS,n})\to \MIC^{\qn}(\XX_{n}/\SS_{n})$ is defined in \eqref{functor nu}.

(ii) We say that a Fontaine module is \textit{strongly divisible} if $\varphi$ \eqref{morphism varphi} is an isomorphism. 

(iii) Given two Fontaine modules $(M_{1},\nabla_{1},M_{1}^{\bullet},\varphi_{1})$ and $(M_{2},\nabla_{2},M_{2}^{\bullet},\varphi_{2})$, a morphism from $(M_{1},\nabla_{1},$ $M_{1}^{\bullet},\varphi_{1})$ to $(M_{2},\nabla_{2},M_{2}^{\bullet},\varphi_{2})$ is a morphism $f:(M_{1},\nabla_{1},M_{1}^{\bullet}) \to (M_{2},\nabla_{2},M_{2}^{\bullet})$ of $\MIC^{p-1}_{\rF}(\XX_{n}/\SS_{n})$ such that the following diagram commutes
\begin{equation}
	\xymatrix{
		\nu(\rmC^{*}(\widetilde{\mathscr{M}}'_{1}))\ar[r]^{\varphi_{1}} \ar[d]_{\nu(\rmC^{*}(\widetilde{f}'))}& (M_{1},\nabla_{1}) \ar[d]^{f} \\	
		\nu(\rmC^{*}(\widetilde{\mathscr{M}}'_{2}))\ar[r]^{\varphi_{2}}& (M_{2},\nabla_{2})
	}
\end{equation}
where $\widetilde{f}':\widetilde{\mathscr{M}}'_{1}\to \widetilde{\mathscr{M}}'_{2}$ is the $\mathscr{O}_{\pCRIS',n}$-linear morphism induced by $f$ \eqref{functor MIC to crys}.
\end{definition}

We denote by $\MFb(\XX)$ the category of Fontaine modules over $\XX$ and by $\MF(\XX)$ the full subcategory of $\MFb(\XX)$ of strongly divisible Fontaine modules $(M,\nabla,M^{\bullet},\varphi)$ such that $M$ is coherent.

\begin{rem} \label{def MF OV}
	In the $p$-torsion case, given an object $(M,\nabla,M^{\bullet})$ of $\MIC^{\ell}_{\rF}(X/k)$, $\widetilde{M}=\gr(M)$ and $\widetilde{\nabla}$ is the Higgs field on $\gr(M)$ induced by $\nabla$ and Griffiths' transversality which is of length $\le \ell$. In (\cite{OV07} 4.16), using their Cartier transform $\rmC_{\XX_{2}'}^{-1}$ \eqref{thm OV}, Ogus and Vologodsky define a $p$-torsion Fontaine module as an object $(M,\nabla,M^{\bullet})$ of $\MIC^{p-1}_{\rF}(X/k)$ together with a horizontal isomorphism
	\begin{equation} \label{def MF OV phi}
		\varphi: \rmC^{-1}_{\XX_{2}'}(\pi^{*}(\Gr(M),\theta))\xrightarrow{\sim} (M,\nabla).
	\end{equation}
	By \ref{comparison Cartier}, our definition \ref{definition MF} is compatible with theirs.
\end{rem}

\begin{nothing}
	In the remainder of this section, we compare definition \ref{definition MF} with Faltings' definition \cite{Fal89} and Tsuji's definition (in a broader context) \cite{Tsu96}. 
	We will formulate their definitions using the notion of \textit{(proto-)T-crystals}, the crystalline counterpart of filtered modules with quasi-nilpotent integrable connection \eqref{filtered mods}, introduced by Ogus \cite{Ogus94}. 
	
	Let $n$ be an integer and $\mathcal{X}$ a smooth scheme over $\SS_{n}$.
	We denote by $\Cris(\mathcal{X}/\SS_{n})$ (resp. $(\mathcal{X}/\SS_{n})_{\cris}$) the crystalline site (resp. topos) of $\mathcal{X}$ over $\SS_{n}$, by $\mathscr{O}_{\mathcal{X}/\SS_{n}}$ the structural sheaf and by $J_{\mathcal{X}/\SS_{n}}$ its PD-ideal. 
\end{nothing}
\begin{definition}[\cite{Ogus94} 2.1.2 and 3.1]
	\textnormal{(i)} Let $(U,T)$ be an object of $\Cris(\mathcal{X}/\SS_{n})$, $J_{T}$ the PD-ideal of $U$ in $T$ and $M$ an $\mathscr{O}_{T}$-module. We say that a decreasing filtration $\{M^{i}\}_{i\in \mathbb{Z}}$ of $M$ by $\mathscr{O}_{T}$-submodules is \textit{G-transversal to $J_{T}$} if for any $i\in \mathbb{Z}$, we have 
	\begin{equation}
		J_{T}M\cap M^{i}=J_{T}^{[1]}M^{i-1}+ J_{T}^{[2]}M^{i-2} +\cdots
	\end{equation}
	In particular, we see that such a filtration is $J_{T}$-saturated, i.e. $J_{T}^{[i]}M^{j}\subset M^{i+j}$ for all $i\ge 0, j$.

	\textnormal{(ii)} Let $\mathscr{M}$ be a crystal of $\mathscr{O}_{\mathcal{X}/\SS_{n}}$-modules. We say that a decreasing filtration $\{\mathscr{M}^{i}\}_{i\in \mathbb{Z}}$ of $\mathscr{M}$ by $\mathscr{O}_{\mathcal{X}/\SS_{n}}$-submodules of $\mathscr{M}$ is \textit{G-transversal to $J_{\mathcal{X}/\SS_{n}}$} if for every object $T$ of $\Cris(\mathcal{X}/\SS_{n})$, the filtration $\{\mathscr{M}^{i}_{T}\}_{i\in \mathbb{Z}}$ of $\mathscr{M}_{T}$ is G-transversal to $J_{\mathcal{X}/\SS_{n},T}$.
\end{definition}

\begin{lemma}[\cite{Ogus94} 3.1.1]\label{lemma Ogus G transversal}
	Let $\mathscr{M}$ be a crystal of $\mathscr{O}_{\mathcal{X}/\SS_{n}}$-modules endowed with a filtration $\{\mathscr{M}^{i}\}_{i\in \mathbb{Z}}$ G-transversal to $J_{\mathcal{X}/\SS_{n}}$. For any morphism $f:T_{2}\to T_{1}$ of $\Cris(\mathcal{X}/\SS_{n})$, via the transition isomorphism $f^{*}(\mathscr{M}_{T_{1}})\xrightarrow{\sim} \mathscr{M}_{T_{2}}$, $\mathscr{M}_{T_{2}}^{i}$ coincides with
	\begin{equation} \label{filtration pullback by f}
		\sum_{i_{1}+i_{2}=i}J^{[i_{1}]}_{T_{2}}\IM(f^{*}(\mathscr{M}_{T_{1}}^{i_{2}})\to f^{*}(\mathscr{M}_{T_{1}})).
	\end{equation}
\end{lemma}

\begin{theorem}[\cite{Ogus94} 3.1.2, 3.2.3] \label{thm Ogus Griffiths tran}
	Let $\mathcal{Y}$ be a smooth $\SS_{n}$-scheme, $\iota:\mathcal{X}\to \mathcal{Y}$ a closed $\SS_{n}$-immersion and $\mathcal{D}$ the PD-envelope of $\iota$ compatible with $\gamma$. 
	Let $\mathscr{M}$ be a quasi-coherent crystal of $\mathscr{O}_{\mathcal{X}/\SS_{n}}$-modules, $M=\mathscr{M}_{\mathcal{D}}$ and $\nabla:M\to M\otimes_{\mathscr{O}_{\mathcal{Y}}}\Omega_{\mathcal{Y}/\SS_{n}}^{1}$ the associated quasi-nilpotent integrable connection on $M$ (as an $\mathscr{O}_{\mathcal{Y}}$-module) (\cite{BO} 6.6). 
	Then the evaluation of sheaves of $(\mathcal{X}/\SS_{n})_{\cris}$ on $\mathcal{D}$ induces an equivalence of following sets of data:

	\textnormal{(i)} A decreasing filtration $\{\mathscr{M}^{i}\}_{i\in \mathbb{Z}}$ by quasi-coherent $\mathscr{O}_{\mathcal{X}/\SS_{n}}$-modules on $\mathscr{M}$ which is G-transversal to $J_{\mathcal{X}/\SS_{n}}$.

	\textnormal{(ii)} A decreasing filtration $\{M^{i}\}_{i\in \mathbb{Z}}$ by quasi-coherent $\mathscr{O}_{\mathcal{D}}$-modules on $M$ which is G-transversal to $J_{\mathcal{X}/\SS_{n},\mathcal{D}}$ and which satisfies Griffiths' transversality i.e. $\nabla(M^{i})\subset M^{i-1}\otimes_{\mathscr{O}_{\mathcal{Y}}}\Omega_{\mathcal{Y}/\SS_{n}}^{1}$ for all $i$. 
\end{theorem}

We briefly review the construction of the data (i) from the data (ii) and we refer to \cite{Ogus94} for more details. Let $\{M_{i}\}_{i\in \mathbb{Z}}$ be a filtration as in (ii). Let $\mathcal{D}(1)$ be the PD-envelope of the immersion $X\xrightarrow{\iota}Y\xrightarrow{\Delta} Y\times_{\SS_{n}}Y$ compatible with $\gamma$, $p_{1},p_{2}:\mathcal{D}(1)\to \mathcal{D}$ the canonical projections and $\varepsilon:p_{2}^{*}(M)\xrightarrow{\sim} p_{1}^{*}(M)$ the $\mathscr{O}_{\mathcal{D}(1)}$-stratification induced by $\nabla$. We define a filtration $\{A^{i}_{j}\}_{i\in \mathbb{Z}}$ on $p_{j}^{*}(M)$ by the formula \eqref{filtration pullback by f} for $j=1,2$. 
By Griffiths' transversality, one verifies that the isomorphism $\varepsilon:p_{2}^{*}(M)\xrightarrow{\sim}p_{1}^{*}(M)$ induces for any $i$, an isomorphism (\cite{Ogus94} 3.1.3)
\begin{equation}
	A^{i}_{2}\xrightarrow{\sim} A^{i}_{1}. \label{filtered iso}
\end{equation}

Given an object $T$ of $\Cris(\mathcal{X}/\SS_{n})$, there exists locally a morphism $r:T\to \mathcal{D}$ of $\Cris(\mathcal{X}/\SS_{n})$. Using the formula \eqref{filtration pullback by f}, we obtain a filtration $\{\mathscr{M}_{T}^{i}\}_{i\in \mathbb{Z}}$ on $r^{*}(M)\xrightarrow{\sim} \mathscr{M}_{T}$. Using the fact that $\varepsilon$ is a filtered isomorphism, one verifies that the filtration $\{\mathscr{M}_{T}^{i}\}_{i\in \mathbb{Z}}$ on $\mathscr{M}_{T}$ is independent of the choice of $r$ up to isomorphisms which come from the stratification and is well-defined. Then we obtain a filtration $\{\mathscr{M}^{i}\}_{i\in \mathbb{Z}}$ of $\mathscr{M}$ by quasi-coherent $\mathscr{O}_{\mathcal{X}/\SS_{n}}$-modules. By (\cite{Ogus94} 2.2.1.2, 2.3.1), one verifies that $\{\mathscr{M}^{i}\}_{i\in \mathbb{Z}}$ is G-transversal to $J_{\mathcal{X}/\SS_{n}}$.

\begin{definition}[\cite{Ogus94} 3.2.1, 3.2.3] \label{def T-crystal}
	(i) A \textit{proto-T-crystal} is a pair $(\mathscr{M},(\mathscr{M}^{i})_{i\in \mathbb{Z}})$ consisting of a quasi-coherent crystal $\mathscr{M}$ with a filtration $\{\mathscr{M}^{i}\}_{i\in \mathbb{Z}}$ G-transversal to $J_{\mathcal{X}/\SS_{n}}$. 
	
	(ii) A proto-T-crystal $(\mathscr{M},\mathscr{M}^{\bullet})$ is called \textit{T-crystal} if for $m>0$ and $i$, the canonical morphism 
	\begin{displaymath}
	\mathscr{M}^{i}_{\mathcal{X}}\otimes_{\mathscr{O}_{\mathcal{X}}}(\mathscr{O}_{\mathcal{X}}/p^m\mathscr{O}_{\mathcal{X}})\to 
	\mathscr{M}_{\mathcal{X}}\otimes_{\mathscr{O}_{\mathcal{X}}}(\mathscr{O}_{\mathcal{X}}/p^m\mathscr{O}_{\mathcal{X}})
	\end{displaymath}
	is injective. 
\end{definition}

We say a proto-T-crystal has length $\le\ell$ if $\mathscr{M}=\mathscr{M}^0$ and the evaluation of $(\mathscr{M}^i)_{i\in \mathbb{Z}}$ at $\mathcal{X}$ has length $\le \ell$ as in \eqref{filtered mods}.

\begin{nothing} \label{setting MF F}
	Let $\YY$ be a smooth formal $\SS$-scheme, $Y$ its special fiber and $\iota:\XX\hookrightarrow \YY$ a closed $\SS$-immersion.  
	For any $n\ge 1$, we denote by $\DD_{n}$ the PD-envelope of $\iota_{n}$ compatible with $\gamma$, by $J_{\DD_{n}}$ the PD-ideal of $\mathscr{O}_{\DD_{n}}$. 
	
	Recall that $\sigma$ denotes the Frobenius endomorphism of $\rW$. We suppose that there exists a $\sigma$-morphism $F_{\YY}:\YY\to \YY$ lifting the Frobenius morphism $F_{Y}:Y\to Y$. For any $n\ge 1$, the $\mathscr{O}_{\YY_{n}}$-linear morphism $\frac{dF_{\YY}}{p}$ \eqref{dFn+1 Omega} induces by adjunction a semi-linear morphism with respect to $F_{\YY}$ that we abusively denote by 
	\begin{equation} \label{dFY p}
		\frac{dF_{\YY}}{p}:\Omega_{\YY_{n}/\SS_{n}}^{1}\to F_{\YY*}(\Omega_{\YY_{n}/\SS_{n}}^{1})=\Omega_{\YY_{n}/\SS_{n}}^{1}.
	\end{equation}
	
	Since $\DD_{n}$ is equal to the PD-envelope of the immersion $X\to \YY_{n}$ compatible with $\gamma$, the morphism $F_{\YY}$ induces a $\sigma$-morphism $F_{\DD}:\DD\to \DD$ lifting the Frobenius morphism of $\DD_{1}$. 
	We denote by $\varphi_{\DD_{n}}$ the homomorphism $\mathscr{O}_{\DD_{n}}= F_{\DD}^{-1}(\mathscr{O}_{\DD_{n}})\to \mathscr{O}_{\DD_{n}}$ induced by $F_{\DD}$. Since $\varphi_{\DD_{1}}(J_{\DD_{1}})=0$, we deduce that for any $0\le r <p$, we have $\varphi_{\DD_{n}}(J^{[r]})\subset p^{r}\mathscr{O}_{\DD_{n}}$. Since $\DD_{n}$ is flat over $\SS_{n}$ (\cite{Ber} I 4.5.1), dividing $\varphi_{\DD_{n+r}}$ by $p^{r}$, we obtain a semi-linear morphism with respect to $F_{\DD}$
	\begin{equation}
		\varphi_{\DD_{n}}^{r}:J_{\DD_{n}}^{[r]}\to \mathscr{O}_{\DD_{n}}\qquad \forall~0\le r\le p-1.
	\end{equation}
\end{nothing}

\begin{definition}[\cite{Fal89} II c), \cite{Tsu96}, 2.1.7] \label{def MF Tsuji}
	(i) Let $(\mathscr{M},\mathscr{M}^{\bullet})$ be a proto-T-crystal over $\XX_{n}/\SS_n$ of length $<p$ \eqref{def T-crystal}, $(M=\mathscr{M}_{\XX_{n}},\nabla)$ the associated quasi-coherent $\mathscr{O}_{\DD_{n}}$-module with integrable connection and $M^{\bullet}$ the associated filtration \eqref{thm Ogus Griffiths tran}. 
	A family of \textit{divided Frobenius morphisms} on $(\mathscr{M},\mathscr{M}^{\bullet})$ with respect to $(\iota,F_{\YY})$ is a family of semi-linear morphisms $\{\varphi^{r}_{\mathfrak{M}}:M^{r}\to M\}_{r<p}$ with respect to $F_{\DD}$ satisfying the following conditions:
	\begin{itemize}
		\item[(a)] $\varphi^{r}_{\mathfrak{M}}|M^{r+1}=p\varphi^{r+1}_{\mathfrak{M}}, \forall~ r<p-1$ (in particular, $\varphi_{\mathfrak{M}}^{-i}=p^{i}\varphi_{\mathfrak{M}}^{0}$ for $i>0$).

		\item[(b)] For any integers $r,s\ge 0$ such that $r+s<p$ and any local sections $a$ of $J_{\DD}^{[r]}$, $x$ of $M^{s}$, we have
			\begin{displaymath}
				\varphi^{r+s}_{\mathfrak{M}}(ax)=\varphi_{\DD_{n}}^{r}(a)\varphi^{s}_{\mathfrak{M}}(x).
			\end{displaymath}

		\item[(c)] The following diagram commutes for all $r<p$ \eqref{dFY p}:
			\begin{equation} \label{Fal horizontal Frob div}
				\xymatrixcolsep{5pc}\xymatrix{
					M^{r} \ar[r]^-{\nabla|_{M^{r}}} \ar[d]_{\varphi^{r}_{\mathfrak{M}}}& M^{r-1}\otimes_{\mathscr{O}_{\YY_{n}}}\Omega_{\YY_{n}/\SS_{n}}^{1} \ar[d]^{\varphi_{\mathfrak{M}}^{r-1}\otimes \frac{d F_{\YY}}{p}} \\
					M \ar[r]^-{\nabla} & M\otimes_{\mathscr{O}_{\YY_{n}}}\Omega_{\YY_{n}/\SS_{n}}^{1} }
			\end{equation}
	\end{itemize}

	(ii) A \textit{Fontaine module over $\XX$ with respect to $(\iota,F_{\YY})$} is a pair of a proto-T-crystal $(\mathscr{M},\mathscr{M}^{\bullet})$ over $\XX_{n}/\SS_n$ of length $<p$ for some integer $n\ge 1$, together with a family of divided Frobenius morphisms with respect to $(\iota,F_{\YY})$.
	We can equivalently write such a data as a quadruple $\mathfrak{M}=(M,\nabla,M^{\bullet},\varphi_{\mathfrak{M}}^{\bullet})$ given by the evaluation of $(\mathscr{M},\mathscr{M}^{\bullet})$ at $\DD$ \eqref{thm Ogus Griffiths tran}. 

	(iii) A morphism between two Fontaine modules $\mathfrak{M}$ and $\mathfrak{N}$ is a morphism of crystals compatible with the filtrations and the divided Frobenius morphisms. 
\end{definition}

	We denote by $\MFb(\XX;\iota,F_{\YY})$ the category of Fontaine modules over $\XX$ with respect to $(\iota,F_{\YY})$. The quadruple $(\mathscr{O}_{\DD_{n}},\nabla_{\DD_{n}},J_{\DD_{n}}^{[\bullet]},\varphi_{\DD_{n}}^{\bullet})$ is an object of $\MFb(\XX;\iota,F_{\YY})$. 

	Let $\mathfrak{M}=(M,M^{\bullet},\nabla_{M},\varphi^{\bullet}_{\mathfrak{M}})$ be an object of $\MFb(\XX;\iota,F_{\YY})$. 
Since $M$ is an $\mathscr{O}_{\DD}$-module, the de Rham complexe $M\otimes_{\mathscr{O}_{\YY_{n}}}\Omega_{\YY_{n}/\SS_{n}}^{\bullet}$ is concentrated on $X$. 
By Griffith's transversality, for any $r\le p-1$, we have a subcomplex $(M^{r-q}\otimes_{\mathscr{O}_{\YY_{n}}}\Omega_{\YY_{n}/\SS_{n}}^{q})_{q\ge 0}$ of the de Rham complex $M\otimes_{\mathscr{O}_{\YY_{n}}}\Omega_{\YY_{n}/\SS_{n}}^{\bullet}$. 
By \eqref{Fal horizontal Frob div}, the divided Frobenius morphisms $\{\varphi_{\mathfrak{M}}^{\bullet}\}$ and $\frac{dF_{\YY}}{p}$ induce a $\rW$-linear morphism of complexes
	\begin{equation}\label{Frob divided morphisms complexes}
		(M^{r-\bullet}\otimes_{\mathscr{O}_{\YY_{n}}}\Omega_{\YY_{n}/\SS_{n}}^{\bullet})\otimes_{\sigma,\rW}\rW \to M\otimes_{\mathscr{O}_{\YY_{n}}}\Omega_{\YY_{n}/\SS_{n}}^{\bullet}.
	\end{equation}

\begin{nothing} \label{M tilde Tsuji}
	Let $\mathfrak{M}=(M,M^{\bullet},\nabla_{M},\varphi^{\bullet}_{\mathfrak{M}})$ be a $p^{n}$-torsion object of $\MFb(\XX;\iota,\YY)$. We define the $\mathscr{O}_{\DD}$-module $\widetilde{\mathfrak{M}}$ as the quotient of $\oplus_{r<p}F_{\DD}^{*}(M^{r})$ by the $\mathscr{O}_{\DD_{n}}$-submodule generated by local sections of the following forms: 
	\begin{eqnarray*}
		&\textnormal{(i)}& (1\otimes x)_{r-1}-(1\otimes px)_{r} \qquad \textnormal{for all } x\in M^{r}, r<p, \\
		&\textnormal{(ii)}& (\varphi_{\DD_{n}}^{r}(a)\otimes x)_{s}-(1\otimes ax)_{r+s} \qquad \textnormal{for all } a\in J_{\DD_{n}}^{[r]}, x\in M^{s}, r\ge 0, r+s<p,
	\end{eqnarray*}
	where $(-)_{r}$ denotes the canonical inclusion $F_{\DD}^{*}(M^{r})\to \oplus_{r<p}F_{\DD}^{*}(M^{r})$.
	In view of condition (d) of \ref{def MF Tsuji}, the morphisms $\{ \varphi_{\mathfrak{M}}^{r}\}_{r<p}$ induce an $\mathscr{O}_{\DD}$-linear morphism
	\begin{equation} \label{phi M tilde}
		\varphi_{\mathfrak{M}}:\widetilde{\mathfrak{M}}\to M.
	\end{equation}
\end{nothing}

\begin{definition} \label{strong divisible}
	We say that an object $\mathfrak{M}$ of $\MFb(\XX;\iota,F_{\YY})$ is \textit{strongly divisible} if $\varphi_{\mathfrak{M}}$ is an isomorphism. 
\end{definition}

\begin{nothing} \label{MF id F}
	In the case $\XX=\YY$ and $\iota=\id$, we have $\DD_{n}=\XX_{n}$. 
	We write simply $\MFb(\XX;F_{\XX})$ for $\MFb(\XX;\id,F_{\XX})$ and we denote by $\MF(\XX;F_{\XX})$ the full subcategory of $\MFb(\XX;F_{\XX})$ of strongly divisible objects whose underlying $\mathscr{O}_{\XX}$-modules are coherent. 

	Let $\mathfrak{M}=(M,\nabla,M^{\bullet},\varphi_{\mathfrak{M}}^{\bullet})$ be an object of $\MFb(\XX;F_{\XX})$. 
	The condition (i-b) of \ref{def MF Tsuji} and the relation (ii) of \ref{M tilde Tsuji} are empty and we have $\widetilde{\mathfrak{M}}=F_{\XX}^{*}(\widetilde{M})$ \eqref{definition M tilde}.
\end{nothing}

\begin{nothing} \label{phi and phi F}
	Suppose that there exists a $\sigma$-lifting $F_{\XX}:\XX\to \XX$ of the Frobenius morphism $F_{X}$. The morphism $F_{\XX}$ induces an $\SS$-morphism $F:\XX\to \XX'$.
	Let $(M,\nabla,M^{\bullet})$ be an object of $\MIC^{\qn}_{\rF}(\XX_{n}/\SS_{n})$. By \eqref{eta F MIC} and \ref{M pullback pi}, the morphism $F$ induces a functorial isomorphism of $\MIC^{\qn}(\XX_{n}/\SS_{n})$:
	\begin{equation} \label{eta F MF module}
		\eta_{F}:\Phi_{n}(\pi^{*}(\widetilde{M},\widetilde{\nabla})) \xrightarrow{\sim} \nu(\rmC^{*}(\widetilde{\mathscr{M}}')),
	\end{equation}
	where $\Phi_{n}$ is Shiho's functor \eqref{Shiho Phi} defined by $F$. The underlying $\mathscr{O}_{\XX_{n}}$-module of $\Phi_{n}(\pi^{*}(\widetilde{M},\widetilde{\nabla}))$ is $F_{\XX}^{*}(\widetilde{M})$ \eqref{Shiho Phi}. We denote the connection on $F_{\XX}^{*}(\widetilde{M})$ by $\nabla_{F}$. 
	
	Given a horizontal morphism $\varphi:\nu(\rmC^{*}(\widetilde{\mathscr{M}}'))\to (M,\nabla)$, we obtain a horizontal morphism $\varphi_{F}$ and a family of morphisms $\{\varphi_{F}^{i}:M^{i}\to M\}_{i=0}^{p-1}$:
	\begin{eqnarray}
		&\varphi_{F}&=\varphi\circ \eta_{F}: (F_{\XX}^{*}(\widetilde{M}),\nabla_{F})\to (M,\nabla), \label{phi F}\\
		&\varphi_{F}^{i}&:M^{i}\xrightarrow{(-)_{i}} \widetilde{M}\xrightarrow{\varphi_{F}} M. \label{Frobenius div morphisms F}
	\end{eqnarray}
	For $i>0$, we set $\varphi_{F}^{-i}=p^{i}\varphi_{F}^{0}$. Then we obtain a functor
	\begin{equation} \label{lambda F MF}
		\lambda_{F}: \MFb(\XX) \to \MFb(\XX;F_{\XX}) \qquad (M,\nabla,M^{\bullet},\varphi)\mapsto (M,\nabla,M^{\bullet},\varphi_{F}^{\bullet}).
	\end{equation}

	Conversely, a family of divided Frobenius morphisms $\{\varphi_{F}^{i}:M^{i}\to M\}_{i\le p-1}$ satisfying (i-a,c) of \ref{def MF Tsuji} induces an $\mathscr{O}_{\XX_{n}}$-linear morphism $\varphi_{F}:F_{\XX}^{*}(\widetilde{M})\to M$ \eqref{phi M tilde}. Then we obtain a morphism $\varphi:\rmC^{*}(\widetilde{\mathscr{M}}')_{(X,\XX)}\to M$ by composing with $\eta_{F}^{-1}$. We define a functor
	\begin{eqnarray} 
		\label{chi F MF} \chi_{F}: \MFb(\XX;F_{\XX}) \to \MFb(\XX) \qquad (M,\nabla,M^{\bullet},\varphi_{F}^{\bullet}) \mapsto (M,\nabla,M^{\bullet},\varphi). 
	\end{eqnarray}
\end{nothing}

\begin{prop}\label{phiF familly phi}
	The functors $\lambda_{F}$ \eqref{lambda F MF} and $\chi_{F}$ \eqref{chi F MF} are well-defined. They induce equivalences of categories quasi-inverse to each other and preserve the strong divisibility condition (\ref{definition MF}, \ref{strong divisible}). 
\end{prop}
\begin{proof} Let $(M,\nabla,M^{\bullet},\varphi)$ be an object of $\MFb(\XX)$.
It follows from the definition of $\widetilde{M}$ that the morphisms $\{\varphi_{F}^{\bullet}\}$ satisfy condition (i-b) of \ref{def MF Tsuji}. We show that they also satisfy condition (i-c). 
	Recall \ref{lemma verification shiho} that for any local sections $m$ of $\widetilde{M}$ and $f$ of $\mathscr{O}_{\XX_{n}}$, we have \eqref{shiho formula}
	\begin{equation} \label{def nabla F}
		\nabla_{F}(f F_{\XX}^{*}(m))=f\zeta\bigl(F_{\XX}^{*}(\widetilde{\nabla}(m))\bigr)+F_{\XX}^{*}(m)\otimes df
	\end{equation}
	where $\zeta$ denotes the composition 
	\begin{displaymath}
		F_{\XX}^{*}(\widetilde{M}\otimes_{\mathscr{O}_{\XX_{n}}}\Omega_{\XX_{n}/\SS_{n}}^{1})\xrightarrow{\sim} F_{\XX}^{*}(\widetilde{M})\otimes_{\mathscr{O}_{\XX_{n}}}F_{\XX}^{*}(\Omega_{\XX_{n}/\SS_{n}}^{1})\xrightarrow{\id\otimes dF_{\XX}/p} F_{\XX}^{*}(\widetilde{M})\otimes_{\mathscr{O}_{\XX_{n}}}\Omega_{\XX_{n}/\SS_{n}}^{1}.
	\end{displaymath}

	For any $0\le i\le p-1$ and any local section $m$ of $M^{i}$, we denote by $(\nabla(m))_{i-1}$ the image of $\nabla(m)$ via $M^{i-1}\otimes_{\mathscr{O}_{\XX_{n}}}\Omega_{\XX_{n}/\SS_{n}}^{1}\to \widetilde{M}\otimes_{\mathscr{O}_{\XX_{n}}}\Omega_{\XX_{n}/\SS_{n}}^{1}$ for $1\le i\le p-1$ and by $(\nabla(m))_{-1}$ the image of $p\nabla(m)$ in $\widetilde{M}\otimes_{\mathscr{O}_{\XX_{n}}}\Omega_{\XX_{n}/\SS_{n}}^{1}$ for $i=0$.
	In view of the definition of $\widetilde{\nabla}$ \eqref{lemma p connection}, we have $\widetilde{\nabla}( (m)_{i})=(\nabla(m))_{i-1}$ and
	\begin{eqnarray} \label{calcul commutative}
		\nabla(\varphi_{F}^{i}(m)) &=& \varphi_{F}\otimes\id (\nabla_{F}(F_{\XX}^{*}((m)_{i}))) \qquad \eqref{Frobenius div morphisms F} \\
		&=& \varphi_{F}\otimes\id (\zeta(F_{\XX}^{*}(\widetilde{\nabla}( (m)_{i})))) \quad ~\eqref{def nabla F} \nonumber \\
		&=& \varphi_{F}\otimes\frac{dF_{\XX}}{p}((\nabla(m))_{i-1})  \nonumber \\
		&=& \varphi_{F}^{i-1}\otimes \frac{dF_{\XX}}{p} (\nabla(m)).  \nonumber
	\end{eqnarray}
	The commutativity of \eqref{Fal horizontal Frob div} follows. The functor \eqref{lambda F MF} is well-defined and preserves the strong divisibility conditions. 

	Conversely, let $(M,\nabla,M^{\bullet},\varphi_{F}^{\bullet})$ be an object of $\MFb(\XX;F_{\XX})$. 
	In view of \eqref{calcul commutative}, the associated morphism $\varphi_{F}:F_{\XX}^{*}(\widetilde{M})\to M$ is compatible with the connections $\nabla_{F}$ and $\nabla$. 
	Via $\eta_{F}^{-1}$ \eqref{eta F MF module}, we obtain a morphism $\varphi:\nu(\rmC^{*}(\widetilde{\mathscr{M}}'))\to (M,\nabla)$ as \eqref{morphism varphi}. Hence the functor $\chi_{F}$ \eqref{chi F MF} is well-defined and is clearly quasi-inverse to $\lambda_{F}$.
\end{proof}
\begin{rem} \label{def Fontaine Laffaille}
	In particular, we see that the notion of Fontaine module over $\SS$ \eqref{definition MF} is compatible with the notion of Fontaine modules over $\rW$ introduced by Fontaine and Laffaille (cf. \cite{FL} 1.2 or \cite{Wach} 2.2.1). 
\end{rem}
\begin{nothing}\label{coro MF1 MF2}
	Let $F_{1},F_{2}:\XX\to \XX'$ be two liftings of the relative Frobenius morphism $F_{X/k}$ of $X$ and let $F_{i,\XX}=\pi\circ F_{i}:\XX\to \XX$. 
	In (\cite{Fal89} proof of Thm. 2.3), Faltings proposed a Taylor formula to construct an equivalence of categories between $\MFb(\XX;F_{2,\XX})$ and $\MFb(\XX;F_{1,\XX})$. 
	We will show that this Taylor formula \eqref{Taylor formula of Faltings} is naturally encoded in the Cartier equivalence. 
	
	We present an explicit description of the equivalence of categories
	\begin{equation} \label{equivalence beta FF}
		\lambda_{F_{2}}\circ \chi_{F_{1}}:\MFb(\XX;F_{1,\XX}) \xrightarrow{\sim} \MFb(\XX;F_{2,\XX}).
	\end{equation}
	In particular, we will see that Faltings' construction coincides with the above equivalence. 

	The morphisms $F_{1}$ and $F_{2}$ induce a morphism of $\pCRIS'$ \eqref{morphism alpha}
	\begin{equation} \label{morphism alpha second}
		\alpha:\rho(X,\XX)\to (X',\RR_{\XX'}).
	\end{equation}
	Let $(M,\nabla,M^{\bullet})$ be an object of $\MIC^{p-1}_{\rF}(\XX_{n}/\SS_{n})$. Recall that the morphism $\alpha$ induces a functorial isomorphism of $\MIC^{\qn}(\XX_{n}/\SS_{n})$ \eqref{alpha iso MIC}
	\begin{equation} \label{isomorphisms alpha MF}
		\alpha^{*}(\pi^{*}_{R}(\varepsilon_{\widetilde{M}})):(F_{2,\XX}^{*}(\widetilde{M}),\nabla_{F_{2}})\xrightarrow{\sim} (F_{1,\XX}^{*}(\widetilde{M}),\nabla_{F_{1}})
	\end{equation}
	such that $\eta_{F_{2}}=\eta_{F_{1}}\circ \alpha^{*}(\pi^{*}_{R}(\varepsilon_{\widetilde{M}}))$ \eqref{alpha iso stra}. In view of the proof of \ref{phiF familly phi}, a family of divided Frobenius morphisms $\{\varphi_{F_{j}}^{i}\}_{i\le p-1}$ is equivalent to a horizontal morphism $\varphi_{F_{j}}$ \eqref{phi F} for $j=1,2$. Then the functor \eqref{equivalence beta FF} is given by 
	\begin{eqnarray}
		\MFb(\XX;F_{1,\XX}) &\to& \MFb(\XX;F_{2,\XX}) \\
		(M,\nabla,M^{\bullet},\varphi_{F_{1}})&\mapsto& (M,\nabla,M^{\bullet},\varphi_{F_{1}}\circ \alpha^{*}(\pi^{*}_{R}(\varepsilon_{\widetilde{M}}))). \nonumber
	\end{eqnarray}
	
	Let us describe the isomorphism \eqref{isomorphisms alpha MF} in terms of a system of local coordinates. Assume that there exists an \'etale $\SS$-morphism $f:\XX\to \widehat{\mathbb{A}}_{\SS}^{d}=\Spf(\rW\{T_{1},\cdots,T_{d}\})$ and put $t_{i}$ the image of $T_{i}$ in $\mathscr{O}_{\XX}$ and $\xi_{i}=1\otimes t_{i}- t_{i}\otimes 1\in \mathscr{O}_{\XX^{2}}$. 
	Let $i$ be an integer $\in [0,p-1]$, $m$ an element of $M^{i}$ and $(m)_{i}$ its image in $\widetilde{M}$. We have \eqref{stratification R tildeM} 
	\begin{equation} \label{stratification Mtilde}
		\varepsilon_{\widetilde{M}}(1\otimes (m)_{i})= \sum_{|I|\le i} \frac{1}{|I|!}\big(\nabla_{\partial^{I}}(m)\big)_{i-|I|}\otimes \biggl(\frac{\xi}{p}\biggr)^{I} + \sum_{|I|> i} \frac{p^{|I|-i}}{|I|!}\big(\nabla_{\partial^{I}}(m)\big)_{0}\otimes \biggl(\frac{\xi}{p}\biggr)^{I}.
	\end{equation}

	Recall that the morphism $\alpha:\XX\to \mathcal{R}_{\XX'}$ \eqref{morphism alpha second} induces a homomorphism $a:\mathcal{R}_{\XX'}\to \mathscr{O}_{\XX}$ \eqref{homomorphism a} which sends $\frac{\xi'_{i}}{p}$ to $\frac{F_{2,\XX}^{*}(t_{i})-F_{1,\XX}^{*}(t_{i})}{p}$ \eqref{homomorphism a formula}. We deduce that 
	\begin{eqnarray}
		\alpha^{*}(\pi^{*}_{R}(\varepsilon_{\widetilde{M}}))(1\otimes_{F_{2}} (m)_{i})=& \sum_{|I|\le i} \biggl(\frac{F_{2,\XX}^{*}(\underline{t})-F_{1,\XX}^{*}(\underline{t})}{p}\biggr)^{I} \bigotimes_{F_{1}} \frac{1}{|I|!}\big(\nabla_{\partial^{I}}(m)\big)_{i-|I|} \label{Taylor formula of Faltings}& \\
		&+ \sum_{|I|> i} \biggl(\frac{F_{2,\XX}^{*}(\underline{t})-F_{1,\XX}^{*}(\underline{t})}{p}\biggr)^{I} \bigotimes_{F_{1}} \frac{p^{|I|-i}}{|I|!}\big(\nabla_{\partial^{I}}(m)\big)_{0}& \nonumber
	\end{eqnarray}
	where 
	\begin{displaymath}
		\biggl(\frac{F_{2,\XX}^{*}(\underline{t})-F_{1,\XX}^{*}(\underline{t})}{p}\biggr)^{I}=\prod_{j=1}^{d}\biggl(\frac{F_{2,\XX}^{*}(t_{j})-F_{1,\XX}^{*}(t_{j})}{p}\biggr)^{i_{j}}\qquad \forall~ I=(i_{1},\cdots,i_{d})\in \mathbb{N}^{d}. \\
	\end{displaymath}
\end{nothing} 

We review some basic properties about Fontaine modules. 
The following result is fisrt showed in \cite{Fal89} and we refer to \cite{Ogus94} for another approach. 

\begin{prop}[\cite{Fal89} 2.1;\cite{Ogus94} 5.3.3] \label{MF ab Ogus}
	Suppose that there exists a $\sigma$-lifting $F_{\XX}:\XX\to \XX$ of the Frobenius morphism. 
Let $(M,\nabla,M^{\bullet},\varphi)$ be an object of $\MF(\XX;F_{\XX})$ \eqref{MF id F}.

	\textnormal{(i)} Then each $M^{i}$ is locally a direct sum of sheaves of the form $\mathscr{O}_{\XX_{n}}$; each morphism $M^i\to M^{i-1}$ locally split. 
	In particular, $(M,\nabla,M^{\bullet})$ forms a T-crystal \eqref{def T-crystal}. 

	\textnormal{(ii)} Any morphism of $\MF(\XX;F_{\XX})$ is strictly compatible with the filtrations (\cite{HodgeII} 1.1.5). 
	
	\textnormal{(iii)} The category $\MF(\XX;F_{\XX})$ is abelian.
\end{prop}

Then we deduce the corresponding statement for (global) Fontaine modules. 

\begin{coro}\label{abelian cat MF}
	\textnormal{(i)} For every object $(M,\nabla,M^{\bullet},\varphi)$ of $\MF(\XX)$, $(M,\nabla,M^{\bullet})$ forms a T-crystal. 
	
	\textnormal{(ii)} Any morphism of $\MF(\XX)$ is strictly compatible with the filtrations. 
	
	\textnormal{(iii)} The category $\MF(\XX)$ is abelian. 
\end{coro}
\begin{proof} Assertions (i) and (ii) being local, they follow from \ref{phiF familly phi} and \ref{MF ab Ogus}. 

For any morphism $f:(M_{1},\nabla_{1},M_{1}^{\bullet},\varphi_{1})\to (M_{2},\nabla_{2},M_{2}^{\bullet},\varphi_{2})$ of $p^{n}$-torsion objects of $\MF(\XX)$, we denote by $(L,\nabla_{L})$ and $(N,\nabla_{N})$ the kernel and the cokernel of $f$ in $\MIC^{\qn}(\XX_{n}/\SS_{n})$. We denote by $L^{\bullet}$ (resp. $N^{\bullet}$) the filtration on $L$ (resp. $N$) induced by $M_{1}^{\bullet}$ (resp. $M_{2}^{\bullet}$) (\cite{HodgeII} 1.1.8). Since $f$ is strictly compabitlbe, for any $i<p$, we have $f(M_{1}^{i})=f(M_{1})\cap M_{2}^{i}$ (\cite{HodgeII} 1.1.11) and an exact sequence $0\to L^{i}\to M_{1}^{i}\to M_{2}^{i}\to N^{i}\to 0$. By the snake lemma, we deduce an exact sequence 
\begin{displaymath}
	0\to \widetilde{L}\to \widetilde{M_{1}}\to \widetilde{M_{2}} \to \widetilde{N}\to 0.
\end{displaymath}
Then we deduce a commutative diagram
\begin{displaymath}
	\xymatrix{
		0\ar[r] & \nu(\rmC^{*}(\widetilde{\mathscr{L}}'))\ar[r] & \nu(\rmC^{*}(\widetilde{\mathscr{M}}'_{1})) \ar[d]^{\wr}_{\varphi_{1}} \ar[r] & \nu(\rmC^{*}(\widetilde{\mathscr{M}}'_{2})) \ar[r] \ar[d]_{\wr}^{\varphi_{2}} & \nu(\rmC^{*}(\widetilde{\mathscr{N}}')) \ar[r] & 0 \\
		0\ar[r] & (L,\nabla_{L}) \ar[r] & (M_{1},\nabla_{1}) \ar[r] & (M_{2},\nabla_{2}) \ar[r] & (N,\nabla_{N}) \ar[r] & 0 }
\end{displaymath}
Hence we can define $\Ker(f)$ and $\Coker(f)$ in $\MF(\XX)$. We deduce that the category $\MF(\XX)$ is abelian. 
\end{proof}

\begin{nothing} \label{morphisms setting Tsuji}
	In the following, we compare categories $\MF(\XX;\iota, F_{\YY})$ with respect to different choice of data $(\iota:\XX\to \YY,F_{\YY})$ following Tsuji \cite{Tsu96}. 

	Let $(\iota_{1}:\XX\to \YY_{1},F_{\YY_{1}})$ and $(\iota_{2}:\XX\to \YY_{2},F_{\YY_{2}})$ be two data as in \ref{setting MF F} and suppose that there exists a \textit{smooth} $\SS$-morphism $g:\YY_{2}\to \YY_{1}$ compatible with $\iota_{1},\iota_{2}$ and the Frobenius morphisms $F_{\YY_{1}}$ and $F_{\YY_{2}}$. Then $g$ induces a PD-morphism $g_{D}:\DD_{2}\to\DD_{1}$ compatible with $F_{\DD_{2}}$ and $F_{\DD_{1}}$. Note that $g_{D}$ induces an isomorphism on the underlying topological spaces.
\end{nothing}
\begin{lemma}[\cite{Tsu96} 2.2.2]\label{lemma local isomorphism PDenvelop}
	Let $x$ be a point of $\XX_{n}$ and let $t_{1},\cdots,t_{d}$ be a family of local sections of $\mathscr{O}_{\YY_{2,n}}$ in a neighborhood of $\iota_{2}(x)$ such that $\{dt_{1},\ldots,dt_{d}\}$ form a basis of $\Omega_{\YY_{2,n}/\YY_{1,n},x}^{1}$ and that $\iota_{2}^{*}(t_{i})=0$ (the existence follows from the fact that $\iota_{1}$ is a closed immersion). Then there exists an $\mathscr{O}_{\DD_{1,n},x}$-PD-isomorphism
	\begin{equation} \label{isomorphism local D1 D2}
		\mathscr{O}_{\DD_{1,n},x}\langle T_{1},\cdots,T_{d}\rangle \xrightarrow{\sim} \mathscr{O}_{\DD_{2,n},x}
	\end{equation}
	which sends $T_{i}$ to $t_{i}$.
\end{lemma}

\begin{prop}[\cite{Tsu96} Proof of 2.2.1]
	Keep the assumption of \ref{morphisms setting Tsuji}. The morphism $g:\YY_{2}\to \YY_{1}$ induces equivalences of categories quasi-inverse to each other:
\begin{equation} \label{pullback MF}
	g^{*}:\MFb(\XX;\iota_{1},F_{\YY_{1}})\to\MFb(\XX;\iota_{2},F_{\YY_{2}}),\qquad g_{*}:\MFb(\XX;\iota_{2},F_{\YY_{2}})\to\MFb(\XX;\iota_{1},F_{\YY_{1}}).
\end{equation}
\end{prop}

We present the construction of pull-back functor and we refer to \cite{Tsu96} for the construction of the push-forward functor. 
Let $\mathfrak{M}=(M_{1},\nabla_{1},M^{\bullet}_{1},\varphi^{\bullet}_{1})$ be a $p^n$-torsion object of $\MFb(\XX;\iota_{1},F_{\YY_{1}})$ and $(\mathscr{M},\mathscr{M}^{\bullet})$ the associated proto-T-crystal. 

The evaluation $(M_2,\nabla_2,M_2^{\bullet})$ of $(\mathscr{M},\mathscr{M}^{\bullet})$ at $\DD_{2}$ is given by $M_{2}=\mathscr{O}_{\DD_{2,n}}\otimes_{\mathscr{O}_{\DD_{1,n}}}M_{1}$ and for $r\in \mathbb{Z}$, $M_2^r=\sum_{r_{1}\ge 0, r_{1}+r_{2}=r}\IM(J_{\DD_{2,n}}^{[r_{1}]}\otimes_{\mathscr{O}_{\DD_{1,n}}}M^{r_{2}}_{1}\to M_{2})$ \eqref{lemma Ogus G transversal}.

	The connection $\nabla_{2}:M_{2}\to M_{2}\otimes_{\mathscr{O}_{\YY_{2,n}}}\Omega_{\YY_{2,n}/\SS_{n}}^{1}$ is defined, for any local sections $a$ of $\mathscr{O}_{\DD_{2,n}}$ and $m$ of $M$, by
\begin{equation} \label{pullback connection}
	\nabla_{2}(a\otimes m)=a\cdot g^{*}(\nabla_{1}(m))+m\otimes \nabla_{\DD_{2,n}}(a).
\end{equation}

Let $x$ be a point of $\XX$. With the notation and assumption of \ref{lemma local isomorphism PDenvelop}, for any $I=(i_{1},\cdots,i_{d})\in \mathbb{N}^{d}$, we set $\underline{T}^{[I]}=\prod_{j=1}^{d}T_{j}^{[i_{j}]}$. In view of \ref{def MF Tsuji}(i-c) and \eqref{isomorphism local D1 D2}, $M_{2,x}^{r}$ can be written as a direct sum of $\mathscr{O}_{\DD_{1,n},x}$-modules $M_{2,x}^{r}=\bigoplus_{|I|=s} M_{1,x}^{r-s}\cdot \underline{T}^{[I]}$.
For any $r< p$, we deduce that the semi-linear morphisms 
\begin{displaymath}
	\varphi_{J_{\DD_{2,n}}}^{r_{1}}\otimes \varphi_{1}^{r_{2}}:J_{\DD_{2,n}}^{[r_{1}]}\otimes_{\mathscr{O}_{\DD_{1,n}}}M_{1}^{r_{2}}\to M_{2},\qquad r_{1}+ r_{2}=r,~ 0\le r_{1}\le p-1
\end{displaymath}
are compatible and induce a family of divided Frobenius morphisms with respect to $F_{\DD_{2}}$: 
\begin{equation}
	\varphi_{2}^{r}:M_{2}^{r}\to M_{2}.
\end{equation}
Conditions \ref{def MF Tsuji}(i a-c) follow from those of $\varphi_{1}^{\bullet}$ and of $\varphi_{\DD_{2,n}}^{\bullet}$. 

\begin{rem}
By \eqref{pullback connection}, the morphism $g$ induces a morphism of de Rham complexes
\begin{equation}
	M_{1}\otimes_{\mathscr{O}_{\YY_{1,n}}}\Omega_{\YY_{1,n}/\SS_{n}}^{\bullet} \to M_{2}\otimes_{\mathscr{O}_{\YY_{2,n}}}\Omega_{\YY_{2,n}/\SS_{n}}^{\bullet}.
\end{equation}
For any $r\le p-1$, we have a commutative diagram \eqref{Frob divided morphisms complexes}
\begin{equation} \label{Frob divided morphisms complexes compatible}
	\xymatrix{
		(M_{1}^{r-\bullet}\otimes_{\mathscr{O}_{\YY_{1,n}}}\Omega_{\YY_{1,n}/\SS_{n}}^{\bullet})\otimes_{\sigma,\rW}\rW \ar[r] \ar[d]& M_{1}\otimes_{\mathscr{O}_{\YY_{1,n}}}\Omega_{\YY_{1,n}/\SS_{n}}^{\bullet} \ar[d] \\
		(M_{2}^{r-\bullet}\otimes_{\mathscr{O}_{\YY_{2,n}}}\Omega_{\YY_{2,n}/\SS_{n}}^{\bullet})\otimes_{\sigma,\rW}\rW \ar[r] & M_{2}\otimes_{\mathscr{O}_{\YY_{2,n}}}\Omega_{\YY_{2,n}/\SS_{n}}^{\bullet} }
\end{equation}
\end{rem}
\begin{lemma}[\cite{Tsu96} 2.3.2]
	The functor $g^{*}$ \eqref{pullback MF} sends strongly divisible objects to strongly divisible objects \eqref{strong divisible}.
\end{lemma}
\begin{proof} The canonical morphisms $\mathscr{O}_{\DD_{2,n}}\otimes_{\mathscr{O}_{\DD_{1,n}}}M_{1}^{r}\to M_{2}^{r}$ induce an $\mathscr{O}_{\DD_{n}}$-linear morphism \eqref{M tilde Tsuji}
	\begin{equation} \label{ug compare Mtilde}
		u_{g}:g_{D}^{*}(\widetilde{\mathfrak{M}})\to \widetilde{g^{*}(\mathfrak{M})}.
	\end{equation}
	In view of condition \ref{M tilde Tsuji}(ii), the above morphism is surjective. 
	We have a commutative diagram:
	\begin{equation} \label{pullback Frob divided}
		\xymatrix{
			g_{D}^{*}(\widetilde{\mathfrak{M}})\ar@{->>}[r]^{u_{g}} \ar[d]_{g_{D}^{*}(\varphi_{\mathfrak{M}})}& \widetilde{g^{*}(\mathfrak{M})} \ar[d]^{\varphi_{g^{*}(\mathfrak{M})}} \\
			g_{D}^{*}(M_{1}) \ar@{=}[r] & M_{2} }
	\end{equation}
	If $\varphi_{\mathfrak{M}}$ is an isomorphism, then $u_{g}$ is an isomorphism and so is $\varphi_{g^{*}(\mathfrak{M})}$. The lemma follows. 
\end{proof}

In the end, we construct a natural functor from $\MFb(\XX)$ to the category of Fontaine modules with respect to the diagonal immersion and two liftings of Frobenius morphism on $\XX$. 

\begin{prop} \label{MF1 MF2 diagonal}
	Suppose that $\XX$ is separated over $\SS$. 
	Take again the notation of \ref{coro MF1 MF2} and let $\Delta:\XX\to \XX^{2}$ be the diagonal immersion, $q_{1},q_{2}:\XX^{2}\to \XX$ the canonical projections and $F_{\XX^{2}}=(F_{1,\XX},F_{2,\XX}):\XX^{2}\to \XX^{2}$. Then the diagram
	\begin{equation}
		\xymatrix{
			\MFb(\XX) \ar[r]^-{\lambda_{F_{1}}} \ar[d]_{\lambda_{F_{2}}} & \MFb(\XX;F_{1,\XX}) \ar[d]^{q_{1}^{*}} \\
			\MFb(\XX;F_{2,\XX}) \ar[r]^-{q_{2}^{*}} & \MFb(\XX;\Delta,F_{\XX^{2}})
		}
	\end{equation}
	is commutative up to a canonical isomorphism.
\end{prop}
\begin{proof} 
Let $\mathfrak{M}=(M,\nabla,M^{\bullet},\varphi)$ be a Fontaine module over $\XX$ and $\varphi_{F_{i}}=\varphi\circ \eta_{F_{i}}:F_{i,\XX}^{*}(\widetilde{M})\to M$ \eqref{phi F} for $i=1,2$. 
Then we have $\lambda_{F_{i}}(\mathfrak{M})=(M,\nabla,M^{\bullet},\varphi_{F_{i}})\in \MFb(\XX;F_{i,\XX})$. 
We denote abusively by $q_{i}$ the composition $\PP_{\XX}\to \XX^{2}\xrightarrow{q_{i}} \XX$ and we set $\mathfrak{M}_{i}=q_{i}^{*}(\lambda_{F_{i}}(\mathfrak{M}))=(q_{i}^{*}(M),\nabla_{i},A_{i}^{\bullet},\varphi_{\mathfrak{M}_{i}}^{\bullet}) \in \MFb(\XX;\Delta,F_{\XX^{2}})$.

	Let $\varepsilon:q_{2}^{*}(M)\xrightarrow{\sim} q_{1}^{*}(M)$ be the $\mathcal{P}_{\XX}$-stratification on $M$. By \ref{thm Ogus Griffiths tran}, it is a filtered isomorphism with respect to the filtrations $A_{2}^{\bullet}$ and $A_{1}^{\bullet}$. 
	It remains to show that $\varepsilon$ is compatible with the divided Frobenius morphisms $\varphi_{\mathfrak{M}_{i}}^{\bullet}$ on both sides.

Since $\varphi:\rmC^{*}(\widetilde{\mathscr{M}}')_{(X,\XX)}\to M$ is a horizontal morphism, the following diagram commutes
\begin{equation}\label{outer diagram 1}
	\xymatrix{
		q_{2}^{*}(\rmC^{*}(\widetilde{\mathscr{M}}')_{(X,\XX)}) \ar[r]^{\sim} \ar[d]_{q_{2}^{*}(\varphi)} &
		\rmC^{*}(\widetilde{\mathscr{M}}')_{(X,\PP_{\XX})} &
		q_{1}^{*}(\rmC^{*}(\widetilde{\mathscr{M}}')_{(X,\XX)}) \ar[d]^{q_{1}^{*}(\varphi)} \ar[l]_{\sim}& \\
		q_{2}^{*}(M) \ar[rr]^{\varepsilon}_{\sim}&& q_{1}^{*}(M)  }
\end{equation}
We have $\varphi_{F_{i}}=\eta_{F_{i}}\circ \varphi$ and the following commutative diagrams:
\begin{equation} \label{outer diagram 2}
	\xymatrix{
		q_{2}^{*}(F_{2,\XX}^{*}(\widetilde{M})) \ar[r]^-{q_{2}^{*}(\eta_{F_{2}})}_-{\sim} \ar[rd]_{q_{2}^{*}(\varphi_{F_{2}})} & q_{2}^{*}(\rmC^{*}(\widetilde{\mathscr{M}}')_{(X,\XX)}) \ar[d]^{q_{2}^{*}(\varphi)}\\
		&q_{2}^{*}(M)  } \qquad
	\xymatrix{
		q_{1}^{*}(F_{1,\XX}^{*}(\widetilde{M})) \ar[r]^-{q_{1}^{*}(\eta_{F_{1}})}_-{\sim} \ar[rd]_{q_{1}^{*}(\varphi_{F_{1}})} & q_{1}^{*}(\rmC^{*}(\widetilde{\mathscr{M}}')_{(X,\XX)}) \ar[d]^{q_{1}^{*}(\varphi)}\\
		&q_{1}^{*}(M)  }
\end{equation}
The filtered isomorphism $\varepsilon:A_{2}^{\bullet}\xrightarrow{\sim} A_{1}^{\bullet}$ induces an isomorphism $\widetilde{\varepsilon}:\widetilde{\mathfrak{M}_{2}}\xrightarrow{\sim} \widetilde{\mathfrak{M}_{1}}$ \eqref{M tilde Tsuji}. 
In the diagrams
\begin{equation} \label{diagram Cartier Tsuji}
	\xymatrix{
		q_{2}^{*}(F_{2,\XX}^{*}(\widetilde{M})) \ar[r]^-{\sim} \ar@{->>}[d]^{u_{q_{2}}} \ar@/_3pc/[dd]_{q_{2}^{*}(\varphi_{F_{2}})}
		& \rmC^{*}(\widetilde{\mathscr{M}}')_{(X,P_{\XX})} \ar@{}[d]|{(1)} &
		q_{1}^{*}(F_{1,\XX}^{*}(\widetilde{M}))  \ar[l]_-{\sim} \ar@{->>}[d]_{u_{q_{1}}} \ar@/^3pc/[dd]^{q_{1}^{*}(\varphi_{F_{1}})} \\
		\widetilde{\mathfrak{M}_{2}} \ar[rr]^{\widetilde{\varepsilon}} \ar[d]^{\varphi_{\mathfrak{M}_{2}}} & \ar@{}[d]|{(2)} &\widetilde{\mathfrak{M}_{1}} \ar[d]_{\varphi_{\mathfrak{M}_{1}}} \\
		q_{2}^{*}(M) \ar[rr]^{\varepsilon} & & q_{1}^{*}(M) }
\end{equation}
the left-hand side diagram and the right-hand side diagram are commutative \eqref{pullback Frob divided}. 

To prove the assertion, it suffices to show that diagram (2) is commutative. By \eqref{outer diagram 1} and \eqref{outer diagram 2}, the outer diagram of \eqref{diagram Cartier Tsuji} commutes. Since $u_{q_{i}}$ is surjective \eqref{ug compare Mtilde}, it suffices to prove the following lemma.
\end{proof}
\begin{lemma}
	Diagram \textnormal{(1)} of \eqref{diagram Cartier Tsuji} is commutative.
\end{lemma}
\begin{proof} Recall \eqref{construction psi12 alpha} that $F_{1},F_{2}$ induce a $\SS$-morphism $\QQ_{\XX}\to \RR_{\XX'}$. We denote the composition $\rho(X,\PP_{\XX})\to \rho(X,\QQ_{\XX})\to (X',\RR_{\XX'})$ of morphisms of $\mathscr{E}'$ \eqref{functor rho} by $f$. It fits into the following commutative diagram:
\begin{displaymath}
	\xymatrix{
		\rho(X,\XX) \ar[d]_{F_{2}} & \rho(X,\PP_{\XX}) \ar[d]_{f} \ar[l]_{\rho(q_{2})} \ar[r]^{\rho(q_{1})} & \rho(X,\XX) \ar[d]^{F_{1}} \\
		(X',\XX') & (X',\RR_{\XX'}) \ar[l]_{q_{2}'} \ar[r]^{q_{1}'} & (X',\XX') }
\end{displaymath}
Hence the composition $q_{2}^{*}(F_{2,\XX}^{*}(\widetilde{M}))\xrightarrow{\sim}q_{1}^{*}(F_{1,\XX}^{*}(\widetilde{M}))$ of the upper arrows of diagram (1) coincides with the pull-back of the $\mathcal{R}_{\XX}$-stratification $\varepsilon_{\widetilde{M}}$ on $\widetilde{M}$ \eqref{prop construction R stra tildeM} via the composition $\PP_{\XX}\xrightarrow{f}\RR_{\XX'}\xrightarrow{\pi_{R}}\RR_{\XX}$ \eqref{M pullback pi}:
\begin{equation} \label{top arrow 1}
	f^{*}(\pi_{R}^{*}(\varepsilon_{\widetilde{M}})): \mathcal{P}_{\XX}\otimes_{\mathcal{R}_{\XX}}(\mathcal{R}_{\XX}\otimes_{\mathscr{O_{\XX_{n}}}}\widetilde{M})\xrightarrow{\sim}\mathcal{P}_{\XX}\otimes_{\mathcal{R}_{\XX}}(\widetilde{M}\otimes_{\mathscr{O}_{\XX_{n}}}\mathcal{R}_{\XX}).
\end{equation}

To show the lemma, we may suppose that there exists an \'etale morphism $\XX\to \widehat{\mathbb{A}}_{\SS}^{d}$ and we take again the notation of \ref{coro MF1 MF2}. For $1\le k\le d$, we set $F_{1}^{*}(t_{k}')=t_{k}^{p}+pa_{k}$ and $F_{2}^{*}(t_{k}')=t_{k}^{p}+pb_{k}$. By \eqref{calcul F1F2 pullback}, the homomorphism $\mathcal{R}_{\XX'}\to \mathcal{P}_{\XX}$ induced by $f$ sends
\begin{equation} \label{description of morphism f}
	\biggl(\frac{\xi_{k}'}{p}\biggr) \mapsto z_{k}=(p-1)!\xi_{k}^{[p]}+\sum_{j=1}^{p-1}\frac{(p-1)!}{j!(p-j)!}\xi_{k}^{j}(t_{k}\otimes1)^{p-j}+(1\otimes b_{k}-a_{k}\otimes 1).
\end{equation}
For any multi-index $I=(i_{1},\cdots,i_{d})$, we set $z^{I}=\prod_{k=1}^{d}z_{k}^{i_{k}}$.
Let $i$ be an integer $\in [0,p-1]$, $m$ a local section of $M^{i}$ and $(m)_{i}$ its image in $\widetilde{M}$. By \eqref{stratification Mtilde} and \eqref{description of morphism f}, the isomorphism \eqref{top arrow 1} sends 
\begin{equation} \label{cal top arrow 1}
	1\otimes(1\otimes_{q_{2}}(m)_{i})\mapsto \sum_{|I|\le i} \frac{z^{I}}{I!}\otimes \bigl((\nabla_{\partial^{I}}(m))_{i-|I|}\otimes_{q_{1}}1\bigr) +\sum_{|I|> i} \frac{p^{|I|-i} z^{I}}{I!}\otimes \bigl((\nabla_{\partial^{I}}(m))_{0}\otimes_{q_{1}}1\bigr)
\end{equation}

The $\mathcal{P}_{\XX}$-linear morphism $u_{q_{2}}$ (resp. $u_{q_{1}}$) \eqref{ug compare Mtilde} sends 
\begin{equation} \label{descrpition uqj}	
	1\otimes(1\otimes_{q_{2}}(m)_{i})\mapsto \bigl(1\otimes_{F_{\PP_{\XX}}}(1\otimes_{q_{2}}m)\bigr)_{i} \qquad \textnormal{(resp. }  1\otimes( (m)_{i}\otimes_{q_{1}}1)\mapsto \bigl(1\otimes_{F_{\PP_{\XX}}}(m\otimes_{q_{1}}1)\bigr)_{i}
\end{equation}
where $(-)_{i}:F_{\PP_{\XX}}^{*}(A_{j}^{i})\to \widetilde{\mathfrak{M}_{j}}$ denotes the canonical morphism for $j=1,2$.

On the other hand, the isomorphism $\varepsilon:A_{2}^{i}\xrightarrow{\sim} A_{1}^{i}$ sends $1\otimes m$ to $\sum_{I}\nabla_{\partial^{I}}(m)\otimes \xi^{[I]}$. 
Hence the isomorphism $\widetilde{\varepsilon}:\widetilde{\mathfrak{M}_{2}}\xrightarrow{\sim} \widetilde{\mathfrak{M}_{1}}$ in diagram (1) sends 
\begin{equation} \label{cal lower arrow 1}
	\bigl(1\otimes_{F_{\PP_{\XX}}}(1\otimes_{q_{2}}m)\bigr)_{i}\mapsto \sum_{I} \bigl(1\otimes_{F_{\PP_{\XX}}} (\nabla_{\partial^{I}}(m)\otimes_{q_{1}} \xi^{[I]})\bigr)_{i}
\end{equation}
With the notation of \ref{setting MF F}, the divided Frobenius morphisms on $(\mathscr{O}_{\PP_{\XX_{n}}},J^{[\bullet]}_{\PP_{\XX_{n}}})$ satisfies
\begin{equation}
	\varphi_{\PP_{\XX_{n}}}(\xi_{k})=pz_{k},~ \forall~ 1\le k\le d,\qquad \varphi_{\PP_{\XX_{n}}}^{i}(\xi^{[I]})=\frac{p^{|I|-i} z^{I}}{I!}, \quad \forall ~ i<p, |I|\ge i.
\end{equation}
By conditions (i) and (ii) of \ref{M tilde Tsuji}, we deduce that
\begin{equation} \label{condition 2 lower arrow}
	(1\otimes_{F_{\PP_{\XX}}} (\nabla_{\partial^{I}}(m)\otimes \xi^{[I]}))_{i}=\left\{ \begin{array}{ll}
		(\frac{z^{I}}{I!}\otimes_{F_{\PP_{\XX}}}(\nabla_{\partial^{I}}(m)\otimes 1))_{i-|I|} & \textnormal{if } |I|\le i\\
		(\frac{p^{|I|-i}z^{I}}{I!}\otimes_{F_{\PP_{\XX}}}(\nabla_{\partial^{I}}(m)\otimes 1))_{0} & \textnormal{if } |I|> i		
	\end{array} \right.
\end{equation}
By comparing \eqref{cal top arrow 1}, \eqref{descrpition uqj}, \eqref{cal lower arrow 1} and \eqref{condition 2 lower arrow}, the assertion follows.
\end{proof}
\section{The Fontaine module structure on the crystalline cohomology of a Fontaine module} \label{coh fil mod}

To illustrates our definition of Fontaine modules, we reprove the following result of Faltings on the crystalline cohomology of a Fontaine module. 

\begin{theorem}[\cite{Fal89} IV 4.1] \label{fil mod coh fil mod}
	Let $\XX$ be a smooth proper formal $\SS$-scheme of relative dimension $d$ and $(M,\nabla,M^{\bullet},\varphi)$ a $p^{n}$-torsion object of $\MF(\XX)$ \eqref{definition MF} of length $\le \ell\le p-1$ (i.e. $M^{\ell+1}=0$). We denote by $\rF^{i}$ the subcomplex $M^{i-\bullet}\otimes_{\mathscr{O}_{\XX_{n}}}\Omega_{\XX_{n}/\SS_{n}}^{\bullet}$ of the de Rham complex $M\otimes_{\mathscr{O}_{\XX_{n}}}\Omega_{\XX_{n}/\SS_{n}}^{\bullet}$.	

	\textnormal{(i)} Let $m$ be an integer such that $\min\{m,d\}+\ell \le p-1$ and $i$ an integer $\le p-1$. Then the canonical morphism $\mathbb{H}^{m}(\rF^{i})\to \mathbb{H}^{m}(M\otimes_{\mathscr{O}_{\XX_{n}}}\Omega_{\XX_{n}/\SS_{n}}^{\bullet})$ is injective. 

	\textnormal{(ii)} Let $m$ be an integer such that $\min\{m,d-1\}+\ell \le p-2$ and $i$ an integer $\le p-1$. The isomorphism $\varphi$ induces a family of semi-linear morphisms $\phi_{\rH}^{m,i}:\mathbb{H}^{m}(\rF^{i})\to \mathbb{H}^{m}(M\otimes_{\mathscr{O}_{\XX_{n}}}\Omega_{\XX_{n}/\SS_{n}}^{\bullet})$ (with respect to $\sigma$). 
	Then the data 
	\begin{displaymath}
		(\mathbb{H}^{m}(M\otimes_{\mathscr{O}_{\XX_{n}}}\Omega_{\XX_{n}/\SS_{n}}^{\bullet}), (\mathbb{H}^{m}(\rF^{i}))_{i=0}^{p-1},(\phi^{m,i}_{\rH})_{i=0}^{p-1})
	\end{displaymath}
	form an object of $\MF(\SS)$ \eqref{def Fontaine Laffaille}. 

	\textnormal{(iii)} In the spectral sequence of the filtered de Rham complex $(M\otimes_{\mathscr{O}_{\XX_{n}}}\Omega_{\XX_{n}/\SS_{n}}^{\bullet},\rF^{i})$ \textnormal{(\cite{HodgeII} 1.4.5)}
	\begin{equation}
		\rE^{r,s}_{1}=\mathbb{H}^{r+s}(\gr^{r}_{F}(M\otimes_{\mathscr{O}_{\XX_{n}}}\Omega_{\XX_{n}/\SS_{n}}^{\bullet}))\Rightarrow \mathbb{H}^{r+s}(M\otimes_{\mathscr{O}_{\XX_{n}}}\Omega_{\XX_{n}/\SS_{n}}^{\bullet}),
	\end{equation}
	the differential morphism $d_{1}^{r,s}$ vanishes for $\min\{r+s,d-1\}+\ell\le p-2$. 
\end{theorem}
 
The data $(\mathscr{O}_{\XX_{n}},d)$ defines a Fontaine module of length $0$ over $\XX$. By \ref{fil mod coh fil mod}(iii), we deduce that:

\begin{coro}[\cite{FM} 2.8; \cite{Fal89} IV 4.1]
	If $d\le p-1$, the Hodge to de Rham spectral sequence of $\XX_{n}/\SS_{n}$ degenerates at $\rE_{1}$.
\end{coro}

We defer the proof of the theorem to \ref{proof of coh theo} and we begin with some preparations on crystalline cohomology of T-cyrstals \eqref{def T-crystal} following Ogus \cite{Ogus94}. 

\begin{nothing}
	Let $\iota:\mathcal{X}\to \mathcal{Y}$ a closed immersion of smooth $\SS_{n}$-schemes and $\mathcal{D}$ the PD-envelope of $\iota$ compatible with $\gamma$. 
	Let $(\mathscr{M},\mathscr{M}^{\bullet})$ a T-crystal of $(\mathcal{X}/\SS_{n})_{\cris}$ \eqref{def T-crystal} and $(M,\nabla)$ (resp. $M^{\bullet}$) the associated $\mathscr{O}_{\mathcal{D}}$-module with integrable connection (resp. filtration) \eqref{thm Ogus Griffiths tran}.
	Since $M^{\bullet}$ satisfy Griffiths' transversality, it induces a filtration on the de Rham complex $M\otimes_{\mathscr{O}_{\mathcal{Y}}}\Omega_{\mathcal{Y}/\SS_{n}}^{\bullet}$ defined for every $i\in \mathbb{Z}$ and $q\in \mathbb{Z}_{\ge 0}$ by
	\begin{equation} \label{filtration dR complex}
		\rF^{i}(M\otimes_{\mathscr{O}_{\mathcal{Y}}}\Omega_{\mathcal{Y}/\SS_{n}}^{q})=M^{i-q}\otimes_{\mathscr{O}_{\mathcal{Y}}}\Omega_{\mathcal{Y}/\SS_{n}}^{q}. 
	\end{equation}
	Since $M$ is an $\mathscr{O}_{\mathcal{D}}$-module, the de Rham complex $M\otimes_{\mathscr{O}_{\mathcal{Y}}}\Omega_{\mathcal{Y}/\SS_{n}}^{\bullet}$ is concentrated on $\mathcal{X}$. 

	Ogus showed that the above subcomplex computes the crystalline cohomology of $\mathscr{M}^i$. 
\end{nothing}

\begin{theorem}[\cite{Ogus94} 6.1.1] \label{thm Poincare lemma}
	We keep the above notation and we denote by $u_{\mathcal{X}/\SS_{n}}$ the canonical morphism of topoi $(\mathcal{X}/\SS_{n})_{\cris}\to \mathcal{X}_{\zar}$. 
	
	\textnormal{(i)} There exists an isomorphism in the derived category $\rD(\mathcal{X}_{\zar},\rW_{n})$
	\begin{equation}
		\rR u_{\mathcal{X}/\SS_{n}*}(\mathscr{M})\xrightarrow{\sim} \mathscr{M}_{\mathcal{D}}\otimes_{\mathscr{O}_{\mathcal{Y}}}\Omega_{\mathcal{Y}/\SS_{n}}^{\bullet}. \label{Poincaré lemma}
	\end{equation}

	\textnormal{(ii)} For every $i$, there exists an isomorphism in $\rD(\mathcal{X}_{\zar},\rW_{n})$ compatible with \eqref{Poincaré lemma}
	\begin{equation}
		\rR u_{\mathcal{X}/\SS_{n}*}(\mathscr{M}^{i})\xrightarrow{\sim} \rF^{i}(\mathscr{M}_{\mathcal{D}}\otimes_{\mathscr{O}_{\mathcal{Y}}}\Omega_{\mathcal{Y}/\SS_{n}}^{\bullet}). \label{filtered Poincaré lemma}
	\end{equation}
\end{theorem}

\begin{coro}[\cite{Ogus94} 6.1.7] \label{compatibility poincare lemma}
	Let $\iota_{1}:\mathcal{X}\to \mathcal{Y}_{1}$ and $\iota_{2}:\mathcal{X}\to \mathcal{Y}_{2}$ be two closed $\SS_{n}$-immersions of $\mathcal{X}$ into smooth $\SS_{n}$-schemes, $\mathcal{D}_{1},\mathcal{D}_{2}$ the PD-envelopes of $\iota_{1},\iota_{2}$ compatible with $\gamma$ and $f:\mathcal{Y}_{2}\to \mathcal{Y}_{1}$ an $\SS_{n}$-morphism such that $\iota_{2}=f\circ \iota_{1}$.
	Then the morphisms of complexes induced by $f$
	\begin{eqnarray*}
		\mathscr{M}_{\mathcal{D}_{1}}\otimes_{\mathscr{O}_{\mathcal{Y}_{1}}}\Omega_{\mathcal{Y}_{1}/\SS_{n}}^{\bullet} &\to& \mathscr{M}_{\mathcal{D}_{2}}\otimes_{\mathscr{O}_{\mathcal{Y}_{2}}}\Omega_{\mathcal{Y}_{2}/\SS_{n}}^{\bullet}, \\
		\rF^{i}(\mathscr{M}_{\mathcal{D}_{1}}\otimes_{\mathscr{O}_{\mathcal{Y}_{1}}}\Omega_{\mathcal{Y}_{1}/\SS_{n}}^{\bullet}) &\to& \rF^{i}(\mathscr{M}_{\mathcal{D}_{2}}\otimes_{\mathscr{O}_{\mathcal{Y}_{2}}}\Omega_{\mathcal{Y}_{2}/\SS_{n}}^{\bullet}),~ \forall~ i\in \mathbb{Z},
	\end{eqnarray*}
	are quasi-isomorphisms and compatible with \eqref{Poincaré lemma} and \eqref{filtered Poincaré lemma}.
\end{coro}

\begin{nothing} \label{construction bicomplex}
	We suppose that $\mathcal{X}$ is separated over $\SS_{n}$ and we take a Zariski covering $\mathscr{U}=\{U_{i}\}_{i\in I}$ of $\mathcal{X}$ consisting of affine schemes. 
	For any integer $r\ge 0$ and any element $J=(j_{0},j_{1},\cdots,j_{r})$ of $I^{r+1}$, we denote by $U^{J}$ the intersection $\cap_{i=0}^{r}U_{j_{i}}$,  by $U_{J}$ the product $\{U_{j_{i}}\}_{i=0}^{r}$ over $\rW_{n}$ and by $P_{J}$ the PD-envelope of the diagonal immersion $U^{J} \to U_{J}$ compatible with $\gamma$. Note that $U^{J}$ and $P_{J}$ are also affine. 
	Then we obtain two compatible simplicial objects:
	\begin{eqnarray} \label{simplicial object U}
		\xymatrix{
			\bigsqcup_{i\in I}U_{i} \ar@<-.8ex>[r]_{s} & \bigsqcup_{J\in I^{2}}U^{J} \ar@<-.8ex>[l]_{d} \ar[l] \ar@<-.8ex>[r] \ar@<-1.6ex>[r]_{s}& \bigsqcup_{J\in I^{3}}U^{J} \cdots  \ar@<-1.6ex>[l]_{d} \ar[l] \ar@<-.8ex>[l]} \\
			\label{simplicial object U prod} \xymatrix{
			\bigsqcup_{i\in I}U_{i} \ar@<-.8ex>[r]_{s} & \bigsqcup_{J\in I^{2}}U_{J} \ar@<-.8ex>[l]_{d} \ar[l] \ar@<-.8ex>[r] \ar@<-1.6ex>[r]_{s}& \bigsqcup_{J\in I^{3}}U_{J} \cdots  \ar@<-1.6ex>[l]_{d} \ar[l] \ar@<-.8ex>[l]}
	\end{eqnarray}
	where $d,s$ denote the faces and degeneracy morphisms.
	Then the faces and degeneracy morphisms of \eqref{simplicial object U prod} induce a simplicial object compatible with \eqref{simplicial object U}:
	\begin{equation}\label{simplicial morphisms}
		\xymatrix{
			\bigsqcup_{i\in I}U_{i} \ar@<-.8ex>[r]_{s} & \bigsqcup_{J\in I^{2}}P_{J} \ar@<-.8ex>[l]_{d} \ar[l] \ar@<-.8ex>[r] \ar@<-1.6ex>[r]_{s}& \bigsqcup_{J\in I^{3}}P_{J} \cdots  \ar@<-1.6ex>[l]_{d} \ar[l] \ar@<-.8ex>[l]}
	\end{equation}

	Let $(\mathscr{M},\mathscr{M}^{\bullet})$ be a T-crystal of $(\mathcal{X}/\SS_{n})_{\cris}$. We associate to $(\mathscr{M},\mathscr{M}^{\bullet})$ a bicomplex $C_{\mathscr{U}}^{\bullet,\bullet}(\mathscr{M})$ and for any $i\in \mathbb{Z}$, a bicomplex $\rF^{i}\bigl(C_{\mathscr{U}}^{\bullet,\bullet}(\mathscr{M})\bigr)$) by setting 
	\begin{eqnarray}
		C_{\mathscr{U}}^{r,s}(\mathscr{M})&=& \bigoplus_{J\in I^{r+1}} \Gamma(U^{J},\mathscr{M}_{P_{J}}\otimes_{\mathscr{O}_{U_{J}}}\Omega^{s}_{U_{J}/\SS_{n}}), \quad r,s\ge 0\\
		\rF^{i}\bigr(C_{\mathscr{U}}^{r,s}(\mathscr{M})\bigl)&=& \bigoplus_{J\in I^{r+1}} \Gamma(U^{J},\mathscr{M}_{P_{J}}^{i-s}\otimes_{\mathscr{O}_{U_{J}}}\Omega^{s}_{U_{J}/\SS_{n}}),\quad r,s\ge 0,
	\end{eqnarray}
	the horizontal differential morphism $\partial_{1}^{r,s}$ is the alternating sum of the resctriction morphisms induced by the faces morphisms \eqref{simplicial morphisms} and the vertical differential morphism $\partial_{2}^{r,s}$ is given by the connection on $\mathscr{M}_{P_{J}}$. 
\end{nothing}

\begin{prop}[\cite{Fal89} IV a)] \label{cohomology bicomplex}
	There exist canonical isomorphisms of cohomology groups:
	\begin{eqnarray}
		\rH^{\bullet}( (\mathcal{X}/\SS_{n})_{\cris},\mathscr{M}) &\xrightarrow{\sim}& \rH^{\bullet}(\Tot(C^{\bullet,\bullet}_{\mathscr{U}}(\mathscr{M}))) \\
		\rH^{\bullet}( (\mathcal{X}/\SS_{n})_{\cris},\mathscr{M}^{i}) &\xrightarrow{\sim}& \rH^{\bullet}(\Tot(\rF^{i}\bigl(C_{\mathscr{U}}^{\bullet,\bullet}(\mathscr{M})\bigr))),~ \forall i\in \mathbb{Z}. \nonumber
	\end{eqnarray}
\end{prop}

\begin{proof} Let $e$ be the final object of $(\mathcal{X}/\SS_{n})_{\cris}$. An open subscheme $U$ of $\mathcal{X}$ defines a subobject $\widetilde{U}$ of $e$ by
\begin{displaymath}
		\widetilde{U}(V,T)= \left\{ \begin{array}{ll}
			e(V,T) & \textnormal{if } V\subset U,\\
			\emptyset & \textnormal{otherwise}.
		\end{array} \right.			
\end{displaymath}
Moreover, there exists a canonical equivalence of topoi (\cite{Ber} IV 3.1.2)
\begin{equation}
	(\mathcal{X}/\SS_{n})_{\cris/\widetilde{U}}\xrightarrow{\sim} (U/\SS_{n})_{\cris}
\end{equation}
which identifies the localisation morphism with respect to $\widetilde{U}$ and the functoriality morphism induced by $U\to X$. 
Then, the morphisms $\{\widetilde{U}_{i}\to e\}_{i\in I}$ form a covering in $(\mathcal{X}/\SS_{n})_{\cris}$. With the notation of \ref{construction bicomplex}, we have $\prod_{j\in J} \widetilde{U}_{j}=\widetilde{U^{J}}$. By cohomological descent, we have a spectral sequence (\cite{HodgeIII} 5.3.3.2, \cite{SGAIV} Vbis 2.5.5)
\begin{equation} \label{ss coh descent}
	\rE^{r,s}_{1}=\bigoplus_{J\in I^{r+1}}\rH^{s}( (U^{J}/\SS_{n})_{\cris},\mathscr{M}|\widetilde{U^{J}})\Rightarrow \rH^{r+s}( (\mathcal{X}/\SS_{n})_{\cris}, \mathscr{M}),
\end{equation}
whose differential morphism $d^{r,s}_{1}:\rE^{r,s}_{1}\to \rE^{r+1,s}_{1}$ is the alternating sum of the morphisms induced by the faces morphisms of \eqref{simplicial object U}. 

On the other hand, we calculate $\rH^{\bullet}(\Tot(C_{\mathscr{U}}^{\bullet,\bullet}(\mathscr{M})))$ by filtering this bicomplex by rows (\cite{EGAIII} 0.11.3.2). By \ref{thm Poincare lemma}(i), the vertical cohomology of $C_{\mathscr{U}}^{\bullet,\bullet}(\mathscr{M})$ is isomorphic to a direct sum of crystalline cohomology groups
\begin{eqnarray}
	\rH^{r,s}_{\rI}&=& \Ker \partial^{r,s}_{2}/\IM(\partial^{r,s-1}_{2}) \\
	&\simeq & \bigoplus_{J\in I^{r+1}}\rH^{s}( (U^{J}/\SS_{n})_{\cris},\mathscr{M}|\widetilde{U^{J}}).
\end{eqnarray}
Recall that we have a spectral sequence (\cite{Weib} 5.6.1)
\begin{equation} \label{ss bicomplex rows filtration}
	\rE'^{r,s}_{1}=\rH^{r,s}_{\rI} \Rightarrow \rH^{r+s}(\Tot(C^{\bullet,\bullet}_{\mathscr{U}}(\mathscr{M}))),
\end{equation}
whose differential morphism $d'^{r,s}_{1}:\rE'^{r,s}_{1}\to \rE'^{r+1,s}_{1}$ is induced by the morphism of complexes
\begin{equation}
	\partial_{1}^{r,\bullet}: C^{r,\bullet}_{\mathscr{U}}(\mathscr{M})\to C^{r+1,\bullet}_{\mathscr{U}}(\mathscr{M}).
\end{equation}
In view of \ref{compatibility poincare lemma} and \ref{construction bicomplex}, the morphism $d'^{r,s}_{1}$ coincides with $d^{r,s}_{1}$ \eqref{ss coh descent}. Then the assertion for $\mathscr{M}$ follows.

Using \ref{thm Poincare lemma}(ii) and \ref{compatibility poincare lemma}, one verifies the assertion for $\mathscr{M}^{i}$ in the same way.
\end{proof}

%

\begin{nothing}\label{lambda M}
	We consider the injective $\mathscr{O}_{\mathcal{X}/\SS_{n}}$-linear morphism
	\begin{equation}
		g:\oplus_{i=1}^{p-1}\mathscr{M}^{i}\to \oplus_{i=0}^{p-1}\mathscr{M}^{i}
	\end{equation}
	defined for every local section $m$ of $\mathscr{M}^{i}$ by $g(m)=(m,-pm)$ in $\mathscr{M}^{i-1}\oplus \mathscr{M}^{i}$. We set $\Lambda_{\mathscr{M}}=\Coker(g)$. 
	For any $r,s\in \mathbb{Z}$, the morphism $g$ induces an injective morphism
	\begin{equation}
		\bigoplus_{i=1}^{p-1} \bigoplus_{J\in I^{r+1}} \Gamma(U^{J},\mathscr{M}_{P_{J}}^{i-s}\otimes_{\mathscr{O}_{U_{J}}}\Omega^{s}_{U_{J}/\SS_{n}}) \to \bigoplus_{i=0}^{p-1} \bigoplus_{J\in I^{r+1}} \Gamma(U^{J},\mathscr{M}_{P_{J}}^{i-s}\otimes_{\mathscr{O}_{U_{J}}}\Omega^{s}_{U_{J}/\SS_{n}})
	\end{equation}
	compatible with the differential morphisms and hence an injective morphism of bicomplexes:
	\begin{equation} \label{morphism gC}
		g_{C}:\bigoplus_{i=1}^{p-1} \rF^{i}\bigl(C_{\mathscr{U}}^{\bullet,\bullet}(\mathscr{M})\bigr)\to \bigoplus_{i=0}^{p-1} \rF^{i}\bigl(C_{\mathscr{U}}^{\bullet,\bullet}(\mathscr{M})\bigr).
	\end{equation}
	We denote its quotient by $C_{\mathscr{U}}^{\bullet,\bullet}(\Lambda_{\mathscr{M}})$. By \ref{cohomology bicomplex}, the crystalline cohomology groups of $\Lambda_{\mathscr{M}}$ are canonically isomorphic to cohomology groups of $\Tot(C_{\mathscr{U}}^{\bullet,\bullet}(\Lambda_{\mathscr{M}}))$.
\end{nothing}

\begin{nothing}
	Let $\XX$ be a smooth and separated formal $\SS$-scheme. 
	In the following, we will use \ref{cohomology bicomplex} to construct a Fontaine module structure on the crystalline cohomology of an object of $\MF(\XX)$ \eqref{definition MF}. For this purpose, we first show how to associate an object of $\MFb(\XX)$ to a Fontaine module with respect to a family of Frobenius liftings \eqref{def MF Tsuji}.

	Suppose that there exist $m+1$ liftings $F_{1},\cdots,F_{m+1}$ of the relative Frobenius morphism $F_{X/k}$ of $X$. We set $\pi:\XX'\to \XX$ the canonical morphism, $F_{\XX^{m+1}}=\pi^{m+1}\circ(F_{1},\cdots,F_{m+1}):\XX^{m+1}\to \XX^{m+1}$ and $\Delta:\XX\to \XX^{m+1}$ the diagonal closed immersion. For $1\le i\le m+1$, the projection $q_{i}:\XX^{m+1}\to \XX$ on the $i$-th component induces a functor \eqref{pullback MF}
	\begin{displaymath}
		q_{i}^{*}:\MFb(\XX;F_{i,\XX})\to \MFb(\XX;\Delta,F_{\XX^{m+1}}). 
	\end{displaymath}
	By composing with the functor $\lambda_{F_{i}}$ \eqref{lambda F MF}, we obtain a functor
	\begin{equation} \label{MF to MFI}
		q_{i}^{*}\circ \lambda_{F_{i}}: \MFb(\XX)\to \MFb(\XX;\Delta,F_{\XX^{m+1}}).
	\end{equation}

	We denote by $\PP_{\XX}(m)$ the PD-envelope of the diagonal immersion $\Delta$ compatible with $\gamma$ \eqref{PD envelop P} and by $F_{\PP_{\XX}(m)}:\PP_{\XX}(m)\to \PP_{\XX}(m)$ the lifting of Frobenius morphism induced by $F_{\XX^{m+1}}:\XX^{m+1}\to \XX^{m+1}$. We denote abusively the composition $\PP_{\XX}(m)\to \XX^{m+1}\xrightarrow{q_{i}}\XX$ by $q_{i}$. 
\end{nothing}

\begin{prop} \label{MF to MF diagonal}
	Let $n$ be an integer $\ge 1$, $\mathfrak{M}=(M,\nabla,M^{\bullet},\varphi)$ a $p^{n}$-torsion Fontaine module over $\XX$, $\mathscr{M}$ the crystal of $\mathscr{O}_{\XX_{n}/\SS_{n}}$-modules associated to $(M,\nabla)$ and $\{\mathscr{M}^{i}\}$ the filtration on $\mathscr{M}$ associated to $\{M^{i}\}$ \eqref{thm Ogus Griffiths tran}. 
	Then there exists a functor 
\begin{eqnarray}
	\MFb(\XX)&\to& \MFb(\XX;\Delta,F_{\XX^{m+1}}) \\
	(M,\nabla,M^{\bullet},\varphi) &\mapsto& (\mathscr{M}_{\PP_{\XX}(m)},\nabla_{\PP_{\XX}(m)},\mathscr{M}^{\bullet}_{\PP_{\XX}(m)},\varphi_{\PP_{\XX}(m)}^{\bullet}) \nonumber
\end{eqnarray}
which is isomorphic to the functor $q_{i}^{*}\circ \lambda_{F_{i}}$ \eqref{MF to MFI} via the transition morphism
\begin{equation} \label{transition iso qi}
	c_{q_{i}}:q_{i}^{*}(M)\xrightarrow{\sim} \mathscr{M}_{\PP_{\XX}(m)}.
\end{equation}
\end{prop}
\begin{proof} It suffices to show that the divided Frobenius morphisms $\varphi_{\mathfrak{M}_{i}}^{\bullet}$ constructed by $q_{i}^{*}\circ \lambda_{F_{i}}$ \eqref{MF to MFI} are compatible via $q_i$ \eqref{transition iso qi}. 
We can reduce to the case $m=1$. In this case, the proposition follows from \ref{MF1 MF2 diagonal}.
\end{proof}
\begin{nothing} \label{XX smooth proper}
	In the rest of this section, we suppose that $\XX$ is a smooth proper formal $\SS$-scheme of relative dimension $d$. 
	Let $n$ an integer $\ge 1$, $(M,\nabla,M^{\bullet},\varphi)$ a $p^{n}$-torsion object of $\MF(\XX)$ of length $\ell<p$ (i.e. $M^{\ell}\neq 0$, $M^{\ell+1}=0$) and $(\mathscr{M},\mathscr{M}^{\bullet})$ the associated T-crystal \eqref{abelian cat MF}.
	We write simply $\rH^{\bullet}_{\cris}(-)$ for the crystalline cohomology groups $\rH^{\bullet}( (\XX_{n}/\SS_{n})_{\cris},-)$. 
\end{nothing}
\begin{lemma}
	Keep the notation of \ref{lambda M} and of \ref{XX smooth proper}. 
	The morphism $\varphi$ induces $\rW$-linear morphisms of bicomplexes
	\begin{eqnarray} 
		&\phi_{C}^{i}:\rF^{i}\big(C^{\bullet,\bullet}_{\mathscr{U}}(\mathscr{M})\big)\otimes_{\sigma,\rW}\rW \to C^{\bullet,\bullet}_{\mathscr{U}}(\mathscr{M}),&\quad \forall~0\le i\le p-1, \label{morphism of bicomplexes phi i C}\\
		\label{morphism of bicomplexes psi C}
		&\psi_{C}:C^{\bullet,\bullet}_{\mathscr{U}}(\Lambda_{\mathscr{M}})\otimes_{\sigma,\rW}\rW \to C^{\bullet,\bullet}_{\mathscr{U}}(\mathscr{M}),&
	\end{eqnarray}	
	which are functorial in $(M,\nabla,M^{\bullet},\varphi)\in \MF(\XX)$.
	In particular, we obtain for $m\ge 0$, $\rW$-linear morphisms 
	\begin{eqnarray} 
		&\phi_{\rH}^{m,i}: \rH^{m}_{\cris}(\mathscr{M}^{i})\otimes_{\sigma,\rW}\rW\to \rH^{m}_{\cris}(\mathscr{M}),& \quad \forall~0\le i\le p-1, \label{morphisms phi mi coh}\\ 
		\label{morphism psi m coh}
		&\psi^{m}:\rH^{m}_{\cris}(\Lambda_{\mathscr{M}})\otimes_{\sigma,\rW}\rW\to \rH^{m}_{\cris}(\mathscr{M}).&
	\end{eqnarray}
\end{lemma}
\begin{proof} We take a Zariski covering $\mathscr{U}=\{\mathfrak{U}_{i}\}_{i\in I}$ of $\XX$ consisting of affine formal schemes and for each $i\in I$ a lifting $F_{i}:\mathfrak{U}_{i}\to \mathfrak{U}'_{i}$ of the relative Frobenius morphism of $\mathfrak{U}_{i,1}$. 
For any integer $r\ge 0$ and $J=(j_{0},\cdots,j_{r})\in I^{r+1}$, we denote by $\UU^{J}$ the intersection $\cap_{i=0}^{r}\UU_{j_{i}}$, by $\UU_{J}$ the product of $(r+1)$-copies of $\UU^{J}$ and by $\PP_{J}$ the PD-envelope of the diagonal closed immersion $\UU^{J}\to \UU_{J}$. 
Note that $\PP_{J}$ is equal to the PD-envelope of the immersion of $\UU^{J}$ in the product of $\{\UU_{j_{i}}\}_{i=0}^{r}$ over $\SS$. 
We denote by $F_{\UU_{J}}:\UU_{J}\to \UU_{J}$ the morphism induced by $\{F_{j_{i}}\}_{i=0}^{n}$ and by $F_{P_{J}}:P_{J}\to P_{J}$ the lifting of the Frobenius morphism induced by $F_{\UU_{J}}$.
By \ref{MF to MF diagonal}, we associate to $(M,\nabla,M^{\bullet},\varphi)$ a family of divided Frobenius morphisms with respect to $(\mathfrak{U}^{J}\to \mathfrak{U}_{J}, F_{\UU_{J}})$:
\begin{equation}
	\varphi_{P_{J}}^{i}:\mathscr{M}^{i}_{P_{J}}\to \mathscr{M}_{P_{J}} \qquad \forall~i\le p-1.
\end{equation}

For any $r,s\ge 0$, in view of \eqref{Frob divided morphisms complexes} and \eqref{Frob divided morphisms complexes compatible}, the $\rW$-linear morphism
\begin{eqnarray*}\label{definition psi C}
	\bigoplus_{J\in I^{r+1}} \varphi_{P_{J}}^{i-s}\otimes \wedge^{s}\biggl(\frac{dF_{U_{J}}}{p}\biggr) :
	\bigoplus_{J\in I^{r+1}}\Gamma(\UU^{J},\mathscr{M}_{P_{J}}^{i-s}\otimes_{\mathscr{O}_{\UU_{J,n}}}\Omega^{s}_{\UU_{J,n}/\SS_{n}})\otimes_{\sigma,\rW}\rW  \\
	\to 
	\bigoplus_{J\in I^{r+1}}\Gamma(\UU^{J},\mathscr{M}_{P_{J}}\otimes_{\mathscr{O}_{\UU_{J,n}}}\Omega^{s}_{\UU_{J,n}/\SS_{n}}) \nonumber
\end{eqnarray*}
is compatible with differential morphisms $\partial_{1}^{\bullet,\bullet},\partial_{2}^{\bullet,\bullet}$ of $\rF^{i}\bigr(C_{\mathscr{U}}^{\bullet,\bullet}(\mathscr{M})\bigl)$, $C^{\bullet,\bullet}_{\mathscr{U}}(\mathscr{M})$ respectively \eqref{construction bicomplex}. 
Then we obtain a $\rW$-linear morphism of bicomplexes
\begin{equation}
	\phi_{C}^{i}:\rF^{i}\bigr(C_{\mathscr{U}}^{\bullet,\bullet}(\mathscr{M})\bigl)\otimes_{\sigma,\rW}\rW\to C^{\bullet,\bullet}_{\mathscr{U}}(\mathscr{M}).
\end{equation}
By condition (i-a) of \ref{def MF Tsuji}, we see that the composition \eqref{morphism gC}
\begin{equation}
	\xymatrix{
		\bigoplus_{i=1}^{p-1}\rF^{i}\bigr(C_{\mathscr{U}}^{\bullet,\bullet}(\mathscr{M})\bigl)\otimes_{\sigma,\rW}\rW \ar[r]^-{g_{C}} & 
		\bigoplus_{i=0}^{p-1} \rF^{i}\bigr(C_{\mathscr{U}}^{\bullet,\bullet}(\mathscr{M})\bigl)\otimes_{\sigma,\rW}\rW \ar[r]^-{\oplus\phi_{C}^{i}} & C^{\bullet,\bullet}_{\mathscr{U}}(\mathscr{M}) }
\end{equation}
vanishes. Then we obtain a $\rW$-linear morphism of bicomplexes
\begin{equation}
	\psi_{C}:C^{\bullet,\bullet}_{\mathscr{U}}(\Lambda_{\mathscr{M}})\otimes_{\sigma,\rW}\rW \to C^{\bullet,\bullet}_{\mathscr{U}}(\mathscr{M}).
\end{equation}
It is clear that the above constructions are functorial.

By \ref{cohomology bicomplex} and \ref{lambda M}, we obtain morphisms of cohomology groups \eqref{morphisms phi mi coh} and \eqref{morphism psi m coh}. 
\end{proof}
\begin{prop} \label{psi m iso ptorsion}
	If $pM=0$ and $\min\{m,d-1\}+\ell \le p-2$, the morphism $\psi^{m}$ \eqref{morphism psi m coh} is an isomorphism.
\end{prop}
\begin{proof} We use $\rH^{r,s}_{\rI}(-)$ to denote the vertical cohomology of a bicomplex. 
Recall \eqref{ss bicomplex rows filtration} that we have two spectral sequences:
\begin{eqnarray}
	\rE'^{r,s}_{1}=\rH^{r,s}_{\rI}(C^{\bullet,\bullet}_{\mathscr{U}}(\Lambda_{\mathscr{M}})) &\Rightarrow& \rH^{r+s}(\Tot(C^{\bullet,\bullet}_{\mathscr{U}}(\Lambda_{\mathscr{M}}))) \\
	\rE^{r,s}_{1}=\rH^{r,s}_{\rI}(C^{\bullet,\bullet}_{\mathscr{U}}(\mathscr{M})) &\Rightarrow& \rH^{r+s}(\Tot(C^{\bullet,\bullet}_{\mathscr{U}}(\mathscr{M})))
\end{eqnarray}
The morphism of bicomplexes $\psi_{C}$ induces a $k$-linear morphism of spectral sequences $\rE'\otimes_{\sigma,k}k\to \rE$. Then it is enough to prove that for every $r\ge 0$ and $s$ satisfying $\min\{s,d-1\}+\ell \le p-2$, the induced morphism
\begin{equation}
	\rH^{r,s}_{\rI}(C^{\bullet,\bullet}_{\mathscr{U}}(\Lambda_{\mathscr{M}}))\otimes_{\sigma,k}k \to \rH^{r,s}_{\rI}(C^{\bullet,\bullet}_{\mathscr{U}}(\mathscr{M}))	
\end{equation}
is an isomorphism. 
Since $\mathscr{M}$ is $p$-torsion, we have $\Lambda_{\mathscr{M}}=\oplus_{i=0}^{p-2}\gr^{i}(\mathscr{M})\oplus \mathscr{M}^{p-1}$. For any $s\ge 0$, we set $\Lambda_{\mathscr{M}}^{-s}=\oplus_{i=0}^{p-2}\gr^{i-s}(\mathscr{M})\oplus \mathscr{M}^{p-1-s}$. Since $\mathfrak{U}^{J}$ is affine, $C^{r,s}_{\mathscr{U}}(\Lambda_{\mathscr{M}})$ can be written as a direct sum
\begin{displaymath}
	\bigoplus_{J\in I^{r+1}} \Gamma(\mathfrak{U}^{J},\Lambda_{\mathscr{M},\PP_{J}}^{-s}\otimes_{\mathscr{O}_{\UU_{J,1}}}\Omega^{s}_{\UU_{J,1}/k}).
\end{displaymath}
Recall \eqref{definition psi C} that the divided Frobenius morphisms $\varphi_{\PP_{J}}^{i}:\mathscr{M}^{i}_{\PP_{J}}\to \mathscr{M}_{\PP_{J}}$ and the semi-linear morphism $\wedge^{s}(\frac{dF_{\UU}}{p}):\Omega^{s}_{\UU_{J,1}/k} \to \Omega^{s}_{\UU_{J,1}/k}$ induce a $k$-linear morphism
\begin{eqnarray*}
	\psi_{J,C}^{r,s}: 
	\Gamma(\UU^{J},\Lambda^{-s}_{\mathscr{M},P_{J}}\otimes_{\mathscr{O}_{\UU_{J,1}}}\Omega^{s}_{\UU_{J,1}/k})\otimes_{\sigma,k}k 
	\to \Gamma(\UU^{J},\mathscr{M}_{P_{J}}\otimes_{\mathscr{O}_{\UU_{J,1}}}\Omega^{s}_{\UU_{J,1}/k})
\end{eqnarray*}
and that the morphism $\psi_{C}^{r,s}$ \eqref{morphism of bicomplexes psi C} is defined by a direct sum of morphisms $\oplus_{J\in I^{r+1}} \psi_{J,C}^{r,s}$.

Then the assertion follows from the following lemma.
\end{proof}
\begin{lemma} \label{lemma Lamda M M complex}
	For any $r\ge 0$, $J\in I^{r+1}$, the morphism of complexes
	\begin{equation} \label{psi JC complex}
		\psi_{J,C}^{r,\bullet}: \Gamma(\UU^{J},\Lambda^{-\bullet}_{\mathscr{M},P_{J}}\otimes_{\mathscr{O}_{\UU_{J,1}}}\Omega^{\bullet}_{\UU_{J,1}/k})\otimes_{\sigma,k}k \to \Gamma(\UU^{J},\mathscr{M}_{P_{J}}\otimes_{\mathscr{O}_{\UU_{J,1}}}\Omega^{\bullet}_{\UU_{J,1}/k})
	\end{equation}
	induces an isomorphism on the $m$-th cohomology group for any integer $m$ satisfying $\min\{m,d-1\}+\ell \le p-2$.
\end{lemma}
\begin{proof} In view of \eqref{Frob divided morphisms complexes compatible} and \eqref{compatibility poincare lemma}, we can reduce to the case where $r=0,~J\in I$. To simplify the notation, we write $\UU$ for $\UU_{J}$, $F$ for the lifting of Frobenius $F_{J}:\UU\to \UU'$, $U$ the special fiber of $\UU$ and $M$ (resp. $M^{i}$) for $\mathscr{M}_{\UU}=M|U$ (resp. $\mathscr{M}^{i}_{\UU}=M^{i}|U$). 

We set $\gr(M)=\oplus_{i=0}^{\ell} M^{i}/M^{i+1}$. By Griffiths' transversality, $\nabla$ induces a Higgs field on $\gr(M)$:
\begin{displaymath}
	\theta:\gr(M)\to \gr(M)\otimes_{\mathscr{O}_{X}}\Omega_{X/k}^{1}.
\end{displaymath}
The source of \eqref{psi JC complex} can be written as
\begin{equation}
	\Gamma(U,(\oplus_{i=0}^{p-2-s}\gr^{i}(M)\oplus M^{p-1-s})\otimes_{\mathscr{O}_{U}}\Omega_{U/k}^{s})\otimes_{\sigma,k}k
\end{equation}
which is equal to 
\begin{equation}
	\Gamma(U,\gr(M)\otimes_{\mathscr{O}_{U}}\Omega_{U/k}^{s})\otimes_{\sigma,k}k \qquad \textnormal{if } s\le p-1-\ell.
\end{equation}
The differential morphism of the source is induced by $\theta$ for $s\le p-1-\ell$. 

Since $pM=0$, we have $(\widetilde{M},\widetilde{\nabla})=(\gr(M),\theta)$ \eqref{def MF OV}.
The isomorphism $\varphi$ and the lifting $F$ induce an isomorphism of $\MIC(\XX_{n}/\SS_{n})$ \eqref{phi F}
\begin{equation}
	\varphi_{F}:\Phi_{1}\bigr( (\gr(M),\theta)\otimes_{\sigma,k}k\bigl)\xrightarrow{\sim} (M,\nabla).
\end{equation}
Recall \eqref{phiF familly phi} that $\varphi_{F}$ induces a family of divided Frobenius morphisms $\varphi_{F}^{\bullet}$. 
The morphism $\psi_{J,C}^{r,\bullet}$ \eqref{psi JC complex}, which is induced by $\varphi_{F}^{\bullet}$, coincides with the composition of morphism of complexes \eqref{morphism Cartier iso nilpotent} induced by $\Phi_{1}$ and the isomorphism of de Rham complexes induced by $\varphi_{F}$
\begin{equation}
	(\gr(M)\otimes_{\mathscr{O}_{U}}\Omega_{U/k}^{\bullet})\otimes_{\sigma,k}k\to M\otimes_{\mathscr{O}_{U}}\Omega_{U/k}^{\bullet}
\end{equation}
in degrees $\le p-1-\ell$. Then the lemma follows from \ref{Cartier iso nilpotent}.
\end{proof}
\begin{rem}
	Suppose that $pM=0$. Let $\gr(M)=\oplus_{i=0}^{\ell}M^{i}/M^{i+1}$ and $\theta$ the Higgs field on $\gr(M)$ induced by $\nabla$ and Griffiths' transversality. By \ref{thm Poincare lemma} and a similar argument of \ref{lemma Lamda M M complex}, we deduce for $m\le p-2-\ell$, an isomorphism
	\begin{equation} \label{iso grM complex Lambda M}
		\rH^{m}_{\cris}(\Lambda_{\mathscr{M}})\xrightarrow{\sim}\mathbb{H}^{m}( \gr(M)\otimes\Omega_{X/k}^{\bullet}).
	\end{equation}
\end{rem}

\begin{prop} \label{psi m iso}
	\textnormal{(i)} If $d\le p-1-\ell$, the morphism $\psi^{m}$ \eqref{morphism psi m coh} is an isomorphism for all $m$.  
	
	\textnormal{(ii)} If $d> p-1-\ell$, the morphism $\psi^{m}$ \eqref{morphism psi m coh} is an isomorphism for $m+\ell \le p-3$, and is a monomorphism for $m+\ell=p-2$.
\end{prop}

\begin{proof} We prove it by induction on $n$. In the case $n=1$, i.e. $M$ is $p$-torsion, it follows from \ref{psi m iso ptorsion}. Suppose that the proposition is true for $n-1$. By \ref{abelian cat MF}, the quadruple $(pM,\nabla|_{pM},pM^{\bullet},\varphi|_{p\rmC^{*}(\widetilde{\mathscr{M}}')})$ is a subobject of $(M,\nabla,M^{\bullet},\varphi)$ in $\MF(\XX)$ and we denote its quotient by $(\overline{M},\overline{\nabla},\overline{M}^{\bullet},\overline{\varphi})$. 
If $\overline{\mathscr{M}}$ and $\overline{\mathscr{M}}^{\bullet}$ denote the crystal of $\mathscr{O}_{\XX_{n}/\SS_{n}}$-modules and the filtration associated to $(\overline{M},\overline{\nabla},\overline{M}^{\bullet})$ \eqref{thm Ogus Griffiths tran}, we have an exact sequence:
\begin{equation}
	0\to p\mathscr{M}^{i}\to \mathscr{M}^{i} \to \overline{\mathscr{M}}^{i}\to 0 \qquad \forall ~ i\le p-1.
\end{equation}
By the snake lemma, we have a commutative diagram:
\begin{equation}
	\xymatrix{
		&0\ar[d] & 0\ar[d] & 0\ar[d] & \\
		0\ar[r] & \oplus_{i=1}^{p-1} p\mathscr{M}^{i} \ar[r] \ar[d] & \oplus_{i=0}^{p-1} p\mathscr{M}^{i} \ar[r] \ar[d] & \Lambda_{p\mathscr{M}} \ar[r] \ar[d]& 0\\
		0\ar[r] & \oplus_{i=1}^{p-1} \mathscr{M}^{i} \ar[r] \ar[d] & \oplus_{i=0}^{p-1} \mathscr{M}^{i} \ar[r] \ar[d] & \Lambda_{\mathscr{M}} \ar[r] \ar[d]& 0 \\
		0\ar[r] & \oplus_{i=1}^{p-1} \overline{\mathscr{M}}^{i} \ar[r] \ar[d] & \oplus_{i=0}^{p-1} \overline{\mathscr{M}}^{i} \ar[r] \ar[d] & \Lambda_{\overline{\mathscr{M}}} \ar[r] \ar[d]& 0 \\
		&0&0&0& }
\end{equation}
Since $\psi_m$ is functorial in $\MF(\XX)$, the assertion follows by dévissage from the induction hypothesis. 
\end{proof}
\begin{prop} \label{prop exact seq Lambda M}
	\textnormal{(i)} Let $m$ be an integer such that $\min\{m,d-1\}\le p-2-\ell$, the exact sequence \eqref{lambda M}
	\begin{displaymath}
		0\to \oplus_{i=1}^{p-1}\mathscr{M}^{i}\xrightarrow{g} \oplus_{i=0}^{p-1}\mathscr{M}^{i}\to \Lambda_{\mathscr{M}}\to 0
	\end{displaymath}
	induces an exact sequence of cohomology groups
	\begin{equation}
		0\to \oplus_{i=1}^{p-1}\rH^{m}_{\cris}(\mathscr{M}^{i})\to \oplus_{i=0}^{p-1}\rH^{m}_{\cris}(\mathscr{M}^{i})\to \rH^{m}_{\cris}(\Lambda_{\mathscr{M}})\to 0.
	\end{equation}

	\textnormal{(ii)} If $m+\ell=p-2$, morphism $\psi^{m}$ \eqref{morphism psi m coh} is an isomorphism. 
\end{prop}

\begin{proof} We prove the proposition by induction on $m$. The case $m=-1$ is trivial. Suppose that the assertion is true for $m-1$ and we will prove it for $m$. By hypothesis of induction, we have an exact sequence (which automatically exists for $m=0$)
\begin{equation} \label{surjectivity to prove}
		0\to \oplus_{i=1}^{p-1}\rH^{m}_{\cris}(\mathscr{M}^{i})\to \oplus_{i=0}^{p-1}\rH^{m}_{\cris}(\mathscr{M}^{i})\to \rH^{m}_{\cris}(\Lambda_{\mathscr{M}}).	
\end{equation}
Since $\mathscr{M},\mathscr{M}^{i}$ are coherent, by \ref{thm Poincare lemma}, we see that the cohomology groups in the above sequence are finite type $\rW_{n}$-modules.
By \ref{psi m iso}, we have the following inequalities on the length of $\rW_{n}$-modules 
\begin{eqnarray}
	\length_{\rW_{n}} \rH^{m}_{\cris}(\Lambda_{\mathscr{M}})&\le& \length_{\rW_{n}} \rH^{m}_{\cris}(\mathscr{M}) \label{equality FM}\\
	&=& \sum_{i=0}^{p-1} \length_{\rW_{n}} \rH^{m}_{\cris}(\mathscr{M}^{i}) -\sum_{i=1}^{p-1} \length_{\rW_{n}} \rH^{m}_{\cris}(\mathscr{M}^{i}). \nonumber
\end{eqnarray}
Then we deduce the surjectivity of the last arrow of \eqref{surjectivity to prove}, and that the above inequality is an equality. The assertion (i) follows.
The assertion (ii) follows from \ref{psi m iso} and the equality \eqref{equality FM}. 
\end{proof}
\begin{nothing} \label{proof of coh theo}
	\textit{Proof of \ref{fil mod coh fil mod}}. Assertion (i) follows from \ref{thm Poincare lemma} and \ref{prop exact seq Lambda M}.

We take for $\phi_{\rH}^{m,i}$ the morphism \eqref{morphisms phi mi coh}. Then assertion (ii) follows from (i), \ref{psi m iso} and \ref{prop exact seq Lambda M}. 

Note that the complex $\rF^{i}=0$ if $i>d+\ell$. For any $r,s$ satisfying $\min\{r+s,d-1\}+\ell\le p-2$, we deduce from (i) that 
\begin{equation}
	\mathbb{H}^{r+s}(\gr^{r}_{\rF}(M\otimes_{\mathscr{O}_{\XX_{n}}}\Omega_{\XX_{n}/\SS_{n}}^{\bullet}))\xrightarrow{\sim} \mathbb{H}^{r+s}(\rF^{r})/\mathbb{H}^{r+s}(\rF^{r+1}). 
\end{equation}
Then assertion (iii) follows by comparing the $\rW_{n}$-length of $\rE^{r,s}_{1}$ and of $\rE^{r+s}$.
\end{nothing}
\begin{rem}\label{final remark on coh}
	Using the comparison isomorphism between the de Rham and the Dolbeault complexes \eqref{iso complexes OV}, Ogus and Vologodsky proved \ref{fil mod coh fil mod} for $p$-torsion Fontaine modules \eqref{def MF OV} (\cite{OV07} 4.17). 
	More precisely, let $(M,\nabla,M^{\bullet},\varphi)$ be a $p$-torsion object of $\MF(\XX)$ of length $\ell$ and $\theta$ the Higgs field on $\gr(M)$ induced by $\nabla$ and Griffiths' transversality. By \ref{iso complexes OV}, the isomorphism \eqref{iso complexes OV eq} induces via \eqref{def MF OV phi}
	\begin{displaymath}
		\varphi:\rmC_{\XX_{2}'}^{-1}(\pi^{*}(\Gr(M),\theta))\xrightarrow{\sim} (M,\nabla),
	\end{displaymath}
	for $m\le p-1-\ell$, an isomorphism:
	\begin{equation} \label{iso cohomology hig dr}
		\mathbb{H}^{m}( \gr(M)\otimes\Omega_{X/k}^{\bullet}) \otimes_{\sigma,k}k\xrightarrow{\sim} \mathbb{H}^{m}(M\otimes\Omega_{X/k}^{\bullet}). 
	\end{equation}
	By \eqref{iso grM complex Lambda M}, we obtain for $m\le p-2-\ell$, an isomorphism
	\begin{equation} 
		\rH_{\cris}^{m}(\Lambda_{\mathscr{M}})\otimes_{\sigma,k}k\xrightarrow{\sim} \rH_{\cris}^{m}(\mathscr{M}).
	\end{equation}
	The above isomorphism is an analogue of the isomorphism $\psi^{m}$ (\ref{morphism psi m coh}, \ref{psi m iso ptorsion}) and allows us to deduce \ref{fil mod coh fil mod} for $p$-torsion objects. We don't know whether these two isomorphisms coincide or not.

	Theorem \ref{fil mod coh fil mod} provides a generalisation of Ogus--Vologodsky's result \eqref{iso complexes OV} for $p^{n}$-torsion Fontaine modules. However, we don't know how to generalise \ref{iso complexes OV} for the Cartier equivalence modulo $p^{n}$.
\end{rem}

\end{document}